\documentclass[10pt]{article}
\usepackage{amsmath,amssymb}
\usepackage[dvips]{graphicx}
\usepackage{theorem}
\usepackage{esint}
\usepackage{color}
\usepackage{hyperref}
\hypersetup{hidelinks}
\usepackage{xcolor}
\usepackage{booktabs}
\usepackage{tabularx}
\usepackage{multirow}
\usepackage{lineno}
\usepackage{bm}
\numberwithin{equation}{section}
{\theorembodyfont{\upshape}}
\newtheorem{remark}{\textbf{\emph{Remark}}}[section]

\theoremstyle{plain}
\theoremstyle{plain} 
\theoremstyle{plain}
\numberwithin{equation}{section}
{\theorembodyfont{\upshape}\theoremstyle{plain}}

\usepackage[left=3cm,right=2.5cm]{geometry}
\usepackage[ruled,linesnumbered]{algorithm2e}
\usepackage{lineno}
\usepackage{float}
\begin{document}

\title{WINO: A Weak-Form Physics Informed Neural Operator for Hyperelasticity on Variable Domains}

\author{Bokai Zhu \thanks{School of Science, Harbin Institute of Technology, Shenzhen, P. R. China. Email address: 200810224@stu.hit.edu.cn.}
\and Yizheng Wang \thanks{Department of Engineering Mechanics, Tsinghua University, Beijing, P. R. China.}
\and Qinghui Zhang \thanks{School of Science, Harbin Institute of Technology, Shenzhen, Guangdong. Email address: zhangqh@hit.edu.cn. This research was partially supported by the Natural Science Foundation of China under grant 12471370.
}
\and Timon Rabczuk \thanks{Corresponding author. Institute of Structural Mechanics, Bauhaus-Universität Weimar, Weimar, 99423, Germany. Email address: timon.rabczuk@uni-weimar.de.}
}

\date{}

\maketitle

\begin{abstract}

	We propose a Weak-form Physics-Informed Neural Operator (WINO), a data-free framework that combines the efficiency of neural operators with the geometric flexibility of the $\varphi$-finite element method ($\varphi$-FEM). $\varphi$-FEM is an unfitted method that accommodates geometric variations without body-fitted meshes, where the domain geometry is represented by the level-set function $\varphi$. To impose the boundary conditions, Dirichlet problems adopt the $\varphi$-FEM lifting so only the homogeneous displacement contribution is learned, whereas traction-driven Neumann problems additionally predict the auxiliary fields necessary for the unfitted weak formulation. Parameters are trained by minimizing squared weak-form residuals aligned with $\varphi$-FEM together with squared penalties on the cut-cell auxiliary equations, which removes the need for large paired datasets of converged reference solutions. When labeled reference data are available, an optional data-augmented variant (WINO+data) can further combine this physics-informed loss with a supervised term. After training, WINO outputs can seed the nonlinear $\varphi$-FEM solvers as neural operator warm starts (NOWS), which reduce iteration counts relative to traditional cold-started solvers. Numerical benchmarks show substantial accuracy of WINO together with total training times of about $15\%$--$70\%$ of those of supervised $\varphi$-FEM-FNO across all cases, without requiring reference-solution generation.

\end{abstract}

{\bf Key words}. Hyperelasticity, Physics informed neural operator, Level set method, Partial differential equation, Weak form

\section{Introduction}

Hyperelasticity represents a typical class of nonlinear elastic problems, in which the stress-strain relationship is both nonlinear and material-dependent \cite{holzapfel2002nonlinear}. Due to their distinctive constitutive behavior, hyperelastic materials (e.g., rubbers and biological tissues) are widely used in various engineering applications \cite{beatty1987topics}. The numerical solution of hyperelastic problems has long been a challenging task in both scientific research and engineering practice. Conventional numerical methods, such as the finite element method (FEM) \cite{zienkiewicz1977finite, reddy2026introduction, weiss1996finite} and the finite difference method (FDM) \cite{thomas2013numerical}, often suffer from high computational complexity and low efficiency when applied to such nonlinear problems. Furthermore, these mesh-based methods require domain discretization; for physical domains with irregular geometries, mesh regeneration or remeshing becomes necessary, which significantly increases the computational burden. Therefore, solving hyperelastic problems on irregularly shaped domains remains a challenging issue that demands more effective and efficient numerical approaches.

In recent years, operator learning has emerged as a prominent neural-network paradigm for PDEs \cite{lu2019deeponet, lu2021learning, li2020fourier}. Discretized samples of the inputs are passed through an operator network such as deep operator network (DeepONet) \cite{lu2019deeponet, lu2021learning, kumar2025synergistic} or fourier neural operator (FNO) \cite{li2020fourier}, which returns a discrete approximation of the solution, thereby amortizing training over families of parametric PDE problems. At inference time, only the prescribed input fields need to be supplied, and the forward evaluation is typically orders of magnitude faster than a repeated conventional FEM solver or full network retraining, often on the order of milliseconds per query in reported regimes. 

Although neural operators can cover broad families of parametric PDEs, most pipelines remain supervised and thus rely on training sets whose labels may be expensive to obtain from classical solvers. Physics-informed neural operators \cite{li2024physics,wang2021learning, ahmadi2026physics} mitigate this cost by building the governing residuals into the loss, reducing or removing the need for precomputed solution pairs and enabling label-free or weakly supervised training. Nonetheless, many such designs penalize the strong form of PDEs, which introduces high-order derivatives and collocation-only enforcement. This issue is especially acute for hyperelasticity, where stresses depend nonlinearly on deformation gradients. Weak-form neural operators are better aligned with finite-element statements: they avoid unnecessary differentiation orders and are typically more stable for large-deformation problems. Representative examples include the variational physics-informed neural operator (VINO) \cite{eshaghi2025variational} and the finite element operator network \cite{lee2023finite}, which minimize an energy-based objective and a weak form residual squared loss function, respectively. Furthermore, physics-informed neural operators can generate preliminary approximate solutions that serve as an effective initial guess for traditional iterative solvers. By providing this near-optimal starting point with a very small inference time, neural operator warm start (NOWS) \cite{eshaghi2025nows, wang2026pretrain} effectively initializes the conventional solver, significantly reducing the number of iterations required to converge to a tight tolerance.

Many of these approaches, however, still assume a fixed domain or a fixed spatial layout of material parameters, so markedly different geometries can require extensive new data and retraining. A compact geometric parametrization is therefore needed to learn operators across families of domains. For complex domains, immersed boundary methods (IBM) provide an alternative to body-fitted mesh generation by embedding the physical boundary into a fixed background mesh and treating interface conditions through weak enforcement or implicit geometric descriptors \cite{peskin2002immersed}. Representative IBM families include fictitious-domain, cut finite element (Cut-FEM), and level-set based unfitted discretizations \cite{burman2025cut}. The level-set idea \cite{osher2001level, osher2004level} represents geometry implicitly through a scalar function $\varphi$ whose sign distinguishes interior from exterior (e.g., $\varphi<0$ inside and $\varphi>0$ outside). $\varphi$-FEM \cite{duprez2020phi,duprez2020varphi, duprez2025varphi} is an efficient level-set based unfitted method on a fixed Cartesian background mesh. Unlike typical Cut-FEM schemes \cite{burman2025cut}, it defines the domain through $\varphi$ and avoids heavily trimmed quadrature on cut cells. This setting is compatible with Cartesian-grid neural operators. However, the standard $\varphi$-FEM weak formulation is generally nonsymmetric and does not coincide with a pure energy minimization principle, so energy-based variational operator frameworks are not directly applicable without reformulation.

To address these issues, we propose the Weak-form Physics-Informed Neural Operator (WINO) for large-deformation hyperelasticity on families of domains. The learnable map follows a FNO architecture. In Dirichlet problems we adopt the $\varphi$-FEM lifting $\mathbf{u}_h=\varphi_h\mathbf{w}_h+\mathbf{g}_h$, where $\mathbf{w}_h$ is the model prediction and $\mathbf{g}_h$ is a suitable extension of the boundary condition. In traction-driven Neumann problems the same backbone maps the level-set description and loads to the displacement together with the auxiliary variables needed on cut cells to close the $\varphi$-FEM weak formulation. We minimize squared weak-form momentum residuals consistent with $\varphi$-FEM, augmented by squared residuals that enforce the auxiliary conditions in intersected elements, which does not rely on precomputed displacement labels. We evaluate the integrals on the unchanged rectangular mesh and compute via Gaussian quadrature. When the training progress is finished, WINO outputs can seed the nonlinear $\varphi$-FEM Newton iteration as neural operator warm starts (NOWS), reducing iteration counts relative to cold starts or classical load continuation. We report numerical results on five benchmark families, including pure Dirichlet problems and mixed Dirichlet and Neumann problems. On the pure Dirichlet benchmarks, WINO is competitive with $\varphi$-FEM-FNO under the same metrics. On the mixed-boundary benchmarks, label-free WINO achieves substantial accuracy while requiring about $25\%$--$70\%$ of the total training time of supervised $\varphi$-FEM-FNO, without requiring converged reference-solution generation. When labeled data are available, the WINO+data variant can further improve accuracy, and trained WINO outputs also serve as warm starts for the high-fidelity $\varphi$-FEM solver.
%%%%%%%%%%%%%%%%%%%%%%%%%%%%%%%%%%
\section{Methods}\label{section_Method}

In this section, we introduce the proposed framework for solving parametric hyperelasticity problems over complex and varying geometries. We first define the governing equations and the underlying physical model in Section 2.1. To seamlessly handle varying domains without the computational overhead of generating body-fitted meshes, we then formulate the boundary value problem using the $\varphi$-FEM approach in Section 2.2. Subsequently, Section 2.3 details the Fourier neural operator architecture. Section 2.4 presents the loss functions used to train the operator together with neural operator warm starts (NOWS) for accelerating the $\varphi$-FEM Newton/Krylov solver.

\subsection{Problem setting}

\begin{figure}[t]
	\centering
	\includegraphics[width=0.4\textwidth]{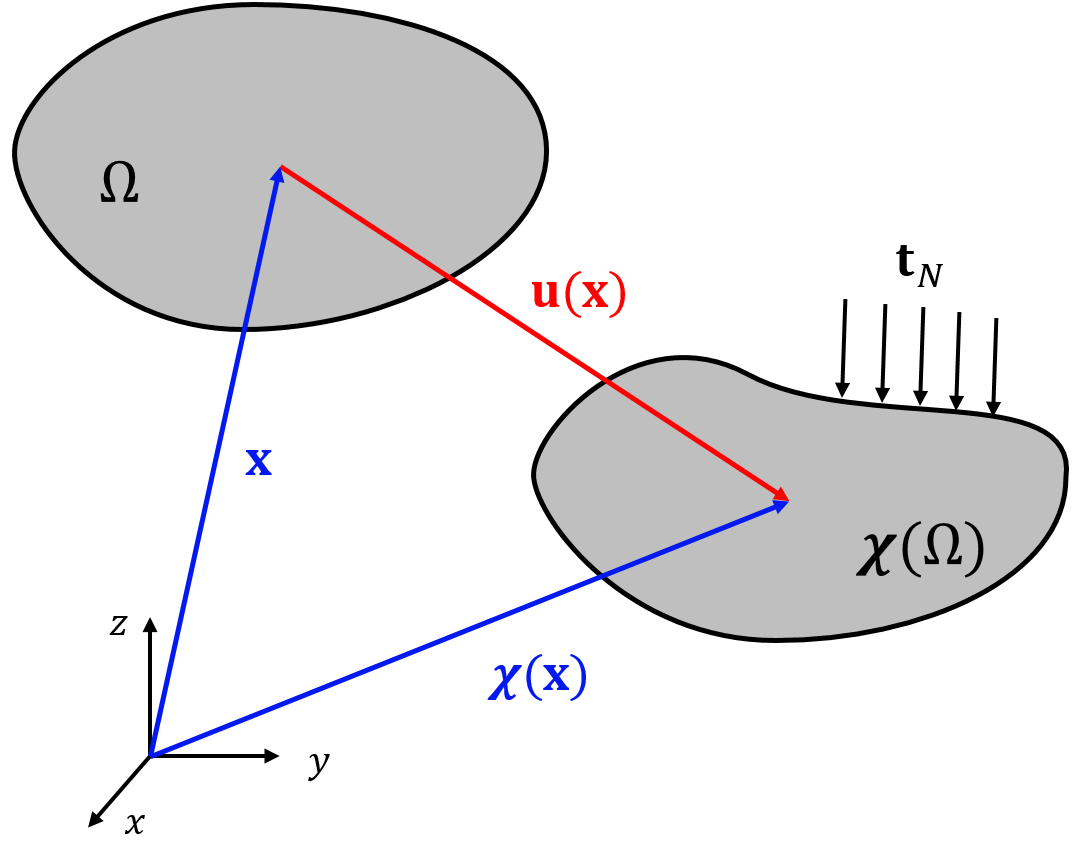}
	\caption{The deformation mapping $\boldsymbol{\chi}$ under the load $\mathbf{t}_N$.}
	\label{fig:deformation_mapping}
\end{figure}

In this subsection, we introduce the governing equations for hyperelastic materials undergoing large deformations. We identify the body with a domain $\Omega \subset \mathbb{R}^d$ in its undeformed configuration (typically $d=2$ or $3$), whose material points are labeled by $\mathbf{x}\in\Omega$. The deformation mapping of each material point from the undeformed configuration to the deformed configuration is described by a smooth mapping $\boldsymbol{\chi}:\Omega\to\Omega_{\mathbf{t}}$, where $\Omega_{\mathbf{t}}$ is the deformed configuration under the load $\mathbf{t}_N$, see Fig. \ref{fig:deformation_mapping}. In the displacement-based formulation used throughout this work, we write
\begin{equation}
	\boldsymbol{\chi}(\mathbf{x})=\mathbf{x}+\mathbf{u}(\mathbf{x}),
\end{equation}
where $\mathbf{u}:\Omega\to\mathbb{R}^d$ is the unknown displacement field. The deformation gradient tensor is
\begin{equation}
	\mathbf{F}(\mathbf{x})=\nabla\boldsymbol{\chi}(\mathbf{x})=\mathbf{I}+\nabla\mathbf{u}(\mathbf{x}),
\end{equation}
which characterizes the local deformation induced by $\boldsymbol{\chi}$. Therefore, we define the boundary value problem for hyperelastic materials as
\begin{alignat}{2}
	-\nabla \cdot \mathbf{P}(\mathbf{F}) &= \mathbf{f}, &\quad &\text{in } \Omega, \nonumber \\
	\mathbf{u} &= \mathbf{g}, &\quad &\text{on } \Gamma_D, \label{eq:Hyperelasticity} \\
	\mathbf{P}(\mathbf{F}) \cdot \mathbf{n} &= \mathbf{t}_N, &\quad &\text{on } \Gamma_N, \nonumber
\end{alignat}
where $\mathbf{P}$ is the first Piola--Kirchhoff stress tensor and $\mathbf{f}$ is the body force density on $\Omega$. The boundary $\partial\Omega$ is partitioned into a Dirichlet part $\Gamma_D$ and a Neumann part $\Gamma_N$ with $\Gamma_D\cap\Gamma_N=\emptyset$ and ${\Gamma_D\cup\Gamma_N}=\partial\Omega$. The field $\mathbf{g}$ prescribes the displacement on $\Gamma_D$, $\mathbf{t}_N$ is the nominal traction on $\Gamma_N$, and $\mathbf{n}$ is the outward unit normal on $\partial\Omega$ in the reference configuration. In this paper we focus on planar problems, $\Omega\subset\mathbb{R}^2$.

For hyperelastic materials, their nonlinear mechanical behavior is fully characterized by a strain energy density function $W(\mathbf{F}(\mathbf{u}))$. The first Piola-Kirchhoff stress tensor is directly derived from the partial derivative of this energy density function with respect to the deformation gradient:
\begin{equation}
	\mathbf{P}(\mathbf{F}(\mathbf{u})) = \frac{\partial W(\mathbf{F}(\mathbf{u}))}{\partial \mathbf{F}}. \label{eq:constitutive_law}
\end{equation}

In this study, we consider the compressible Neo-Hookean hyperelastic constitutive model as a representative example. The strain energy density function for this material is defined as follows (see \cite{bonet1997nonlinear}):
\begin{equation}
	W(\mathbf{F}(\mathbf{u})) = \frac{\mu}{2}(\text{tr}(\mathbf{C}) - 3 - 2\ln(J)) + \frac{\lambda}{2}(\ln(J))^2, \label{eq:strain_energy}
\end{equation}
where $\mathbf{C} = \mathbf{F}^T \mathbf{F}$ is the right Cauchy-Green deformation tensor, and $J = \det(\mathbf{F})$ is the Jacobian determinant of the deformation gradient, representing the local volume ratio of the material before and after deformation. $\mu$ and $\lambda$ are the Lamé parameters that characterize the intrinsic mechanical properties of the material. In this paper, we only formulate the hyperelasticity problem under the plane strain conditions.

Through the strong form of the nonlinear partial differential equations described above, we define the mechanical response of the physical system under specific boundary conditions and constitutive relations. In the subsequent subsections, we will explore how to transform this boundary value problem into a weak form suitable for neural network solvers in the context of varying geometric domains.

%%%%%%%%%%%%%%%%%%%%%%%%%%%%%%%%%%
\subsection{$\varphi$-FEM formulation}\label{subsec:phi_fem}

We first briefly describe the $\varphi$-FEM method introduced in \cite{duprez2025varphi, duprez2020varphi} to solve \eqref{eq:Hyperelasticity}. As the method is specifically designed for Neumann problems. We suppose the physical domain $\Omega$ is embedded within a simpler background Cartesian grid denoted by $\mathcal{O}$. The varying geometry and its boundaries are implicitly described using a level-set function $\varphi$. The active physical domain and its boundary are defined as:
\begin{equation}
	\begin{aligned}
		\Omega &:= \{\mathbf{x} \in \mathcal{O} : \varphi(\mathbf{x}) < 0\}, \\
		\Gamma &:= \{\mathbf{x} \in \mathcal{O} : \varphi(\mathbf{x}) = 0\}.
	\end{aligned}
\end{equation}
To ensure the mathematical rigor and optimal convergence of the $\varphi$-FEM formulation, we state the regularity assumptions for the level-set function $\varphi$ following standard practices \cite{duprez2020varphi, duprez2025varphi}. We assume that $\varphi$ is sufficiently smooth in a tubular neighborhood $\mathcal{O}_{\Gamma}$ around the boundary $\Gamma$. Specifically, $\varphi \in C^{k+1}(\mathcal{O}_{\Gamma})$, where $k \ge 1$ is the polynomial degree of the finite element space used for the spatial discretization. Furthermore, to guarantee that the boundary is well-defined and to prevent degenerate cut-cells, $\varphi$ behaves locally as a signed distance function such that its gradient is bounded away from zero, i.e., $|\nabla \varphi| \ge m > 0$ on $\mathcal{O}_{\Gamma}$ for some constant $m$.

\begin{figure}[t]
	\centering
	\includegraphics[width=0.4\textwidth]{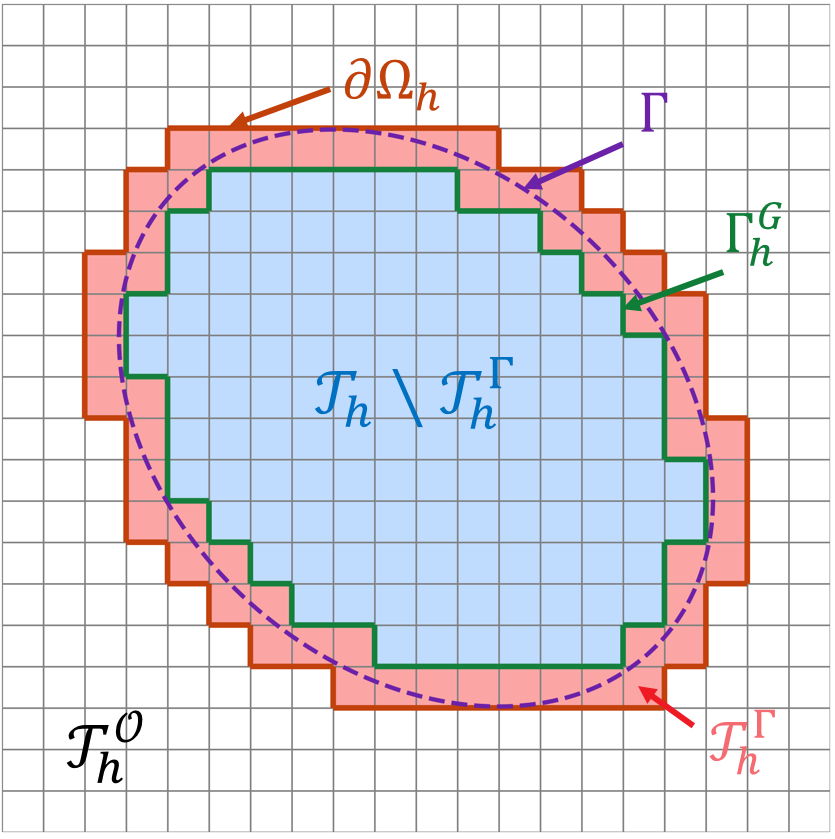}
	\caption{$\varphi$-FEM meshes and facet sets for an elliptic domain: the exact boundary $\Gamma$ (purple); facets belonging to $\Gamma_h^G$ (green); facets on $\partial\Omega_h$ (red); active interior cells $\mathcal{T}_h\setminus\mathcal{T}_h^\Gamma$ (blue); cut cells $\mathcal{T}_h^\Gamma$ (red fill); Cartesian background mesh $\mathcal{T}_h^{\mathcal{O}}$ (overall grid).}
	\label{fig:phi_fem_mesh_tags}
\end{figure}

Let $\mathcal{T}_h^{\mathcal{O}}$ be a Cartesian background mesh of the box $\mathcal{O}$ with a characteristic element size $h$. We define the computational mesh $\mathcal{T}_h$ as the collection of cells in $\mathcal{T}_h^{\mathcal{O}}$ that intersect the physical domain:
\begin{equation}
	\mathcal{T}_h := \{ T \in \mathcal{T}_h^{\mathcal{O}} : T \cap \{\varphi_h < 0\} \neq \emptyset \},
\end{equation}
where $\varphi_h$ is obtained by standard finite-element interpolation of $\varphi$ on $\mathcal{T}_h^{\mathcal{O}}$ through the mesh shape functions: nodal values are set to $(\varphi_h)_i=\varphi(\mathbf{x}_i)$, and on each cell $\varphi_h$ is reconstructed as a linear combination of the local shape functions with these nodal coefficients \cite{duprez2020varphi,duprez2025varphi}. Furthermore, to properly treat the boundary conditions, we introduce a boundary sub-mesh $\mathcal{T}_h^\Gamma \subset \mathcal{T}_h$ containing the cells cut by the approximate boundary:
\begin{equation}
	\mathcal{T}_h^\Gamma := \{ T \in \mathcal{T}_h : T \cap \{\varphi_h = 0\} \neq \emptyset \}.
\end{equation}
We denote the continuous domains covered by these meshes as
\begin{equation}
	\Omega_h := \cup_{T \in \mathcal{T}_h} T, \quad \Omega_h^\Gamma := \cup_{T \in \mathcal{T}_h^\Gamma} T.
	\label{eq:omega_h}
\end{equation}
Additionally, we define $\Gamma_h^G := \partial\Omega_h^\Gamma \setminus \partial\Omega_h$ as the facets between the $\mathcal{T}_h\setminus\mathcal{T}_h^\Gamma$ and the $\mathcal{T}_h^\Gamma$. To better understand $\varphi$-FEM, we illustrate the meshes and facets for an elliptic domain in Fig. \ref{fig:phi_fem_mesh_tags}.

For the spatial discretization, let $k \ge 1$ be an integer. The planar displacement $\mathbf{u}_h:\Omega_h\to\mathbb{R}^2$ is sought in the vector-valued finite element space on the active mesh,
\begin{equation}\label{eq:phi_fem_space}
	V_h^{(k)} := \{ \mathbf{v}_h \in H^1(\Omega_h;\mathbb{R}^2) : \mathbf{v}_h|_T \in (\mathbb{P}_k(T))^2 \,\, \forall T \in \mathcal{T}_h \}.
\end{equation}

\subsubsection*{Pure Dirichlet problems}
For pure Dirichlet problems, where $\Gamma=\Gamma_D$ in \eqref{eq:Hyperelasticity}, we follow the standard $\varphi$-FEM strategy \cite{duprez2020phi} and write the discrete displacement as
\begin{equation}
	\mathbf{u}_h=\varphi_h\mathbf{w}_h+\mathbf{g}_h,
	\label{eq:phifem_dirichlet_compose}
\end{equation}
where $\mathbf{g}_h$ extends the boundary data $\mathbf{g}$ to $\Omega_h$ by Lagrange interpolation on $\mathcal{T}_h^{\mathcal{O}}$ from nodal values of a prescribed extension of $\mathbf{g}$, and $\mathbf{w}_h$ is the new unknown. Since $\varphi_h$ vanishes on $\Gamma$, the condition $\mathbf{u}_h=\mathbf{g}$ holds on the implicit boundary while all computations use the fixed background mesh. Homogeneous trial and test components are taken from $V_h^{(k)}$ in \eqref{eq:phi_fem_space}. With $\mathbf{u}_h=\varphi_h\mathbf{w}_h+\mathbf{g}_h$ and test functions $\mathbf{v}_h=\varphi_h\mathbf{s}_h$, where $\mathbf{w}_h,\mathbf{s}_h\in V_h^{(k)}$, the stabilized $\varphi$-FEM weak formulation of the planar hyperelastic problem reads: find $\mathbf{w}_h\in V_h^{(k)}$ such that, for all $\mathbf{s}_h\in V_h^{(k)}$,
\begin{equation}
	\int_{\Omega_h}\mathbf{P}(\mathbf{F}(\mathbf{u}_h)) : \nabla \mathbf{v}_h
	-\int_{\partial\Omega_h}(\mathbf{P}(\mathbf{F}(\mathbf{u}_h))\mathbf{n})\cdot\mathbf{v}_h  {}+ G_h(\mathbf{u}_h,\mathbf{v}_h)
	= \int_{\Omega_h} \mathbf{f}_h\cdot\mathbf{v}_h + G_h^{\mathrm{rhs}}(\mathbf{v}_h),
\label{eq:phifem_dirichlet_weak}
\end{equation}
with
\begin{equation}
\begin{aligned}
	G_h(\mathbf{u}_h,\mathbf{v}_h) := {}& \sigma_D h\sum_{E\in\mathcal{F}_h^\Gamma}\int_E
	\bigl[\!\bigl[\mathbf{P}(\mathbf{F}(\mathbf{u}_h))\mathbf{n}_E\bigr]\!\bigr] \cdot
	\bigl[\!\bigl[D_u(\mathbf{P} \circ \mathbf{F})(\mathbf{u}_h)[\mathbf{v}_h]\mathbf{n}_E\bigr]\!\bigr] \\
	&{}+ \sigma_D h^2\sum_{T\in\mathcal{T}_h^\Gamma}\int_T (\nabla\cdot\mathbf{P}(\mathbf{F}(\mathbf{u}_h))) \cdot (\nabla\cdot D_u(\mathbf{P} \circ \mathbf{F})(\mathbf{u}_h)[\mathbf{v}_h]),
\end{aligned}
\label{eq:phifem_dirichlet_Gh}
\end{equation}
Here, $\mathcal{F}_h^\Gamma$ denotes the set of internal facets that lie in the red cut-cell region of Fig.~\ref{fig:phi_fem_mesh_tags} (the band where $\mathcal{T}_h^\Gamma$ is shown with red fill), and
\begin{equation}
	G_h^{\mathrm{rhs}}(\mathbf{v}_h):=-\sigma_D h^2\sum_{T\in\mathcal{T}_h^\Gamma}\int_T \mathbf{f}_h\cdot\bigl(\nabla\cdot D_u(\mathbf{P} \circ \mathbf{F})(\mathbf{u}_h)[\mathbf{v}_h]\bigr),
	\label{eq:phifem_dirichlet_Grhs}
\end{equation}
where $\sigma_D>0$ is a stabilization parameter and $\bigl[\!\bigl[\cdot\bigr]\!\bigr]$ denotes the jump across each facet $E\in\mathcal{F}_h^\Gamma$ in the first line of $G_h$. 

\subsubsection*{Pure Neumann problems}
For pure Neumann problems, imposing boundary conditions accurately on an implicitly defined boundary $\Gamma$ without conforming elements is a primary challenge. To address this, the $\varphi$-FEM approach introduces two auxiliary variables defined exclusively on the boundary sub-mesh domain $\Omega_h^\Gamma$ \cite{duprez2020varphi}. The core boundary constraints are weakly enforced through the following relations:

\begin{equation}\label{eq:strong_yp}
	\begin{aligned} 
		\mathbf{y}_h + \mathbf{P}(\mathbf{F}(\mathbf{u}_h)) &= \mathbf{0}, & \text{on } \Omega_h^\Gamma, \\ 
		\mathbf{y}_h \nabla \varphi_h + \frac{1}{h} \mathbf{p}_h \varphi_h + \mathbf{t}_h |\nabla \varphi_h| &= \mathbf{0}, & \text{on } \Omega_h^\Gamma.
	\end{aligned}
\end{equation}
In this formulation, the tensor variable $\mathbf{y}_h$ approximates the negative stress tensor $-\mathbf{P}(\mathbf{F}(\mathbf{u}_h))$, while the vector variable $\mathbf{p}_h$ acts as a Lagrange multiplier to stabilize the boundary constraints. $\mathbf{t}_h$ denotes an extension of the prescribed Neumann traction $\mathbf{t}_N$ from the boundary $\Gamma$ to the cut-cell region $\Omega_h^\Gamma$. The corresponding finite element spaces for these auxiliary variables, $Z_h^{(k)}$ and $Q_h^{(k-1)}$, are defined as:
\begin{align}
	Z_h^{(k)} &:= \{ \mathbf{z}_h \in L^2(\Omega_h^\Gamma) : \mathbf{z}_h|_T \in \mathbb{P}_k(T) \,\, \forall T \in \mathcal{T}_h^\Gamma \}, \\
	Q_h^{(k-1)} &:= \{ \mathbf{q}_h \in L^2(\Omega_h^\Gamma) : \mathbf{q}_h|_T \in \mathbb{P}_{k-1}(T) \,\, \forall T \in \mathcal{T}_h^\Gamma \}.
\end{align}

Consequently, the strong form of the hyperelastic boundary value problem is recast into the following fully coupled weak form formulation: Find $(\mathbf{u}_h, \mathbf{y}_h, \mathbf{p}_h) \in V_h^{(k)} \times Z_h^{(k)} \times Q_h^{(k-1)}$ such that for all test functions $(\mathbf{v}_h, \mathbf{z}_h, \mathbf{q}_h) \in V_h^{(k)} \times Z_h^{(k)} \times Q_h^{(k-1)}$,

\begin{equation}
	\begin{aligned}
		& \int_{\Omega_h} \mathbf{P}(\mathbf{F}(\mathbf{u}_h)) : \nabla \mathbf{v}_h 
		+ \int_{\partial\Omega_h} (\mathbf{y}_h \mathbf{n}) \cdot \mathbf{v}_h - \int_{\Omega_h} \mathbf{f}_h\cdot\mathbf{v}_h
		+ G_h(\mathbf{u}_h, \mathbf{v}_h) \\
		& \quad + \frac{\gamma_p}{h^2} \int_{\Omega_h^{\Gamma}} \left(\mathbf{y}_h \nabla \varphi_{h} + \frac{1}{h} \mathbf{p}_h \varphi_{h} + \mathbf{t}_h |\nabla \varphi_h| \right) \cdot \left(\mathbf{z}_h \nabla \varphi_{h} + \frac{1}{h} \mathbf{q}_h \varphi_{h} \right) \\
		& \quad + \gamma_u \int_{\Omega_h^{\Gamma}} (\mathbf{y}_h + \mathbf{P}(\mathbf{F}(\mathbf{u}_h))) : (\mathbf{z}_h + D_u(\mathbf{P} \circ \mathbf{F})(\mathbf{u}_h)[\mathbf{v}_h]) \\
		& \quad + \gamma_{div} \int_{\Omega_h^{\Gamma}} (\nabla \cdot \mathbf{y}_h) \cdot (\nabla \cdot \mathbf{z}_h) = 0,
	\end{aligned}
	\label{eq:phifem_weak_form}
\end{equation}
where $\gamma_u$, $\gamma_p$, and $\gamma_{div}$ are penalty parameters that weakly enforce the auxiliary Neumann constraints \eqref{eq:strong_yp} on the cut-cell region $\Omega_h^\Gamma$, where the implicit boundary is not aligned with mesh facets. The operator $D_u(\mathbf{P} \circ \mathbf{F})(\mathbf{u}_h)[\mathbf{v}_h]$ computes the directional derivative of the stress tensor with respect to the displacement field $\mathbf{u}_h$ in the direction of $\mathbf{v}_h$. To guarantee numerical stability and coercivity across the irregular cut cells, a ghost penalty stabilization term $G_h(\mathbf{u}_h, \mathbf{v}_h)$ is incorporated over the internal facets $\Gamma_h^G$ of the boundary sub-mesh. This term is defined as:
\begin{equation}
	\label{Gh_definition}
	G_h(\mathbf{u}_h, \mathbf{v}_h) := \sigma_N h  \int_{\Gamma_h^G} \bigl[\!\bigl[\mathbf{P}(\mathbf{F}(\mathbf{u}_h)) \mathbf{n}\bigr]\!\bigr] \cdot \bigl[\!\bigl[D_u(\mathbf{P} \circ \mathbf{F})(\mathbf{u}_h)[\mathbf{v}_h] \mathbf{n}\bigr]\!\bigr],
\end{equation}
where $\sigma_N > 0$ is a stabilization parameter, $\mathbf{n}$ is the unit normal vector to the facet $\Gamma_h^G$, and the brackets $\bigl[\!\bigl[\cdot\bigr]\!\bigr]$ denote the jump of the enclosed quantity across the facet.

\vspace{1em}
\begin{remark}
	From a computational perspective, $\varphi$-FEM provides several practical advantages over both classical body-fitted finite element methods and standard cut finite element schemes. Unlike conforming discretizations, it does not require generating or updating a boundary-fitted mesh whenever the physical domain changes. The geometry is specified solely by a level-set field on a fixed Cartesian background, which makes it particularly attractive for parametric neural operator learning. Compared with typical Cut-FEM pipelines that often rely on exact clipping of elements and quadrature rules adapted to arbitrary cut-cell intersections, the $\varphi$-FEM formulation avoids those geometrically complex volume integrals on trimmed elements. Volume and interface terms are instead assembled from collections of background cells together with stabilization on the internal facets ($\mathcal{F}_h^\Gamma$ in the Dirichlet branch and $\Gamma_h^G$ in the Neumann branch), yielding a simpler implementation while retaining robustness for implicit boundaries.
\end{remark}

\begin{remark}
	The jump notation $[\![\cdot]\!]$ in \eqref{eq:phifem_dirichlet_Gh} and \eqref{Gh_definition} denotes the discontinuity of the enclosed quantity across an internal facet. In the pure Dirichlet formulation, the ghost-penalty sum in \eqref{eq:phifem_dirichlet_Gh} is taken over $E\in\mathcal{F}_h^\Gamma$, whereas in the Neumann formulation it is taken over $E\in\Gamma_h^G$ as defined in Fig.~\ref{fig:phi_fem_mesh_tags}. These terms form a ghost-penalty stabilization on the cut-cell band: because the implicit boundary is not aligned with the background mesh, mesh conformity is lost in $\mathcal{T}_h^\Gamma$, and penalizing flux jumps on these facet sets suppresses spurious oscillations and restores stability.
\end{remark}

% \vspace{1em}
% \begin{remark}
% For mixed boundary value problems, recent literature \cite{cotin2023phi} has introduced a ''lazy approach'', in which no explicit boundary condition is imposed on an element if the Dirichlet/Neumann transition falls inside that element. Stability is then maintained through stabilization terms (e.g., ghost penalty). In the present work, we adopt a direct hybrid strategy: the level-set function $\varphi=0$ is used to define the complex Neumann part, while Dirichlet boundaries are aligned with exterior edges of the Cartesian background mesh and imposed directly by standard conforming finite elements.
% \end{remark} 

\subsection{The structure of the neural operator}

We employ the Fourier Neural Operator (FNO) \cite{li2020fourier} as a surrogate mapping between grid-valued fields on a Cartesian background of size $N_x\times N_y$. In the pure Dirichlet lift framework of Section~2.2, the network predicts only the homogeneous coefficient $\mathbf{w}_h$ from the level-set field $\varphi_h$, the discrete extension $\mathbf{g}_h$ of the boundary data, and the field channel $\mathbf{f}_h$:
\begin{equation}
	\begin{split}
		\mathcal{G}_\theta : \mathbb{R}^{N_x \times N_y \times d_{in}} &\rightarrow \mathbb{R}^{N_x \times N_y \times d_{out}},\\
		(\varphi_h, \mathbf{g}_h, \mathbf{f}_h) &\mapsto \mathbf{w}_h,
	\end{split}
	\label{eq:fno_mapping_dirichlet}
\end{equation}
where $d_{in}$ stacks $\varphi_h$, the discrete Dirichlet boundary data $\mathbf{g}_h$, and the body force $\mathbf{f}_h$, while $d_{out}=2$ in the plane (two components of $\mathbf{w}_h$). The discrete displacement is then recovered pointwise on the mesh by the $\varphi$-FEM lift $\mathbf{u}_h=\varphi_h\mathbf{w}_h+\mathbf{g}_h$. For the Neumann-oriented $\varphi$-FEM weak form \eqref{eq:phifem_weak_form}, the same FNO skeleton instead maps $(\varphi_h, \mathbf{t}_h, \mathbf{f}_h)$ to the displacement and the two auxiliary fields on the background grid:
\begin{equation}
	\begin{split}
		\mathcal{G}_\theta : \mathbb{R}^{N_x \times N_y \times d_{in}} &\rightarrow \mathbb{R}^{N_x \times N_y \times d_{out}},\\
		(\varphi_h, \mathbf{t}_h, \mathbf{f}_h) &\mapsto (\mathbf{u}_h, \mathbf{y}_h, \mathbf{p}_h),
	\end{split}
	\label{eq:fno_mapping}
\end{equation}
where  $d_{out}=8$ in two dimensions (2 channels for $\mathbf{u}_h$, 4 for $\mathbf{y}_h$, and 2 for $\mathbf{p}_h$).

To enhance training stability and accelerate convergence, the input tensor is initially scaled by a standardization operator $\mathcal{S}$, and the corresponding output is rescaled via its inverse $\mathcal{S}^{-1}$. The core neural architecture is built upon an encoding layer $\mathcal{E}_\theta$, a stack of $L$ iterative integral kernel layers $\mathcal{K}_\theta^{(l)}$ (where $L$ is typically set to 4), and a decoding layer $\mathcal{D}_\theta$. The global feed-forward pipeline is expressed as:
\begin{equation}
	\mathcal{G}_\theta = \mathcal{S}^{-1} \circ \mathcal{D}_\theta \circ \mathcal{K}_\theta^{(L)} \circ \cdots \circ \mathcal{K}_\theta^{(2)} \circ \mathcal{K}_\theta^{(1)} \circ \mathcal{E}_\theta \circ \mathcal{S}. \label{eq:fno_pipeline}
\end{equation}
The encoder $\mathcal{E}_\theta$ operates as a point-wise fully connected network across the spatial grid. Its primary function is to lift the low-dimensional physical features into a higher-dimensional hidden channel space of width $d_h$. Consequently, the initial hidden state is formulated as $\mathbf{h}^{(0)} = \mathcal{E}_\theta(\mathcal{S}(\varphi_h, \mathbf{t}_h, \mathbf{f}_h)) \in \mathbb{R}^{N_x \times N_y \times d_h}$.

The dominant mathematical transformations are executed within the iterative kernel layers $\mathcal{K}_\theta^{(l)}$. At each step $l \in \{1, \dots, L\}$, the hidden representation $\mathbf{h}^{(l-1)}$ is processed through a dual-path mechanism: a global spectral convolution $\mathcal{G}_{conv}^{(l)}$ and a local linear transformation $\mathcal{T}_{loc}^{(l)}$. These are subsequently passed through a non-linear activation function $\sigma$ (e.g., GELU), yielding the updated state
\begin{equation}
	\mathbf{h}^{(l)} = \mathcal{K}_\theta^{(l)}(\mathbf{h}^{(l-1)}) = \sigma \left( \mathcal{G}_{conv}^{(l)}(\mathbf{h}^{(l-1)}) + \mathcal{T}_{loc}^{(l)}(\mathbf{h}^{(l-1)}) \right). \label{eq:kernel_layer}
\end{equation}

The local operation $\mathcal{T}_{loc}^{(l)}$ is defined by a standard point-wise affine transformation parameterized by a weight matrix $\mathbf{W}^{(l)} \in \mathbb{R}^{d_h \times d_h}$ and a bias vector $\mathbf{b}^{(l)} \in \mathbb{R}^{d_h}$
\begin{equation}
	\mathcal{T}_{loc}^{(l)}(\mathbf{h}^{(l-1)}) = \mathbf{W}^{(l)} \mathbf{h}^{(l-1)} + \mathbf{b}^{(l)}. \label{eq:local_linear}
\end{equation}
Concurrently, the global spectral operator $\mathcal{G}_{conv}^{(l)}$ extracts long-range spatial correlations by learning a parameterized mapping in the frequency domain. By employing the discrete Fast Fourier Transform ($\mathcal{F}$) alongside its inverse ($\mathcal{F}^{-1}$), this non-local convolution is evaluated as
\begin{equation}
	\mathcal{G}_{conv}^{(l)}(\mathbf{h}^{(l-1)}) = \mathcal{F}^{-1} \left( \mathbf{R}^{(l)} \mathcal{F}(\mathbf{h}^{(l-1)}) \right), \label{eq:global_conv}
\end{equation}
where $\mathbf{R}^{(l)}$ acts as a complex-valued parameter tensor in the Fourier space.

An intrinsic advantage of this architecture is its integrated low-pass filtering property. Because the mechanical responses of hyperelastic PDEs typically exhibit high spatial smoothness, high-frequency spectral modes contain marginal physical relevance. Thus, following the forward transformation $\mathcal{F}(\mathbf{h}^{(l-1)})$, the frequency spectrum is truncated, preserving only the lowest $k_{max} \times k_{max}$ Fourier modes. The parameterized complex tensor $\mathbf{R}^{(l)}$ is exclusively applied to these truncated components. This strategy not only diminishes the total parameter count substantially but also provides an implicit regularization effect that prevents overfitting to high-frequency discretization noise.

\begin{figure}[t]
	\centering
	\includegraphics[width=0.9\textwidth]{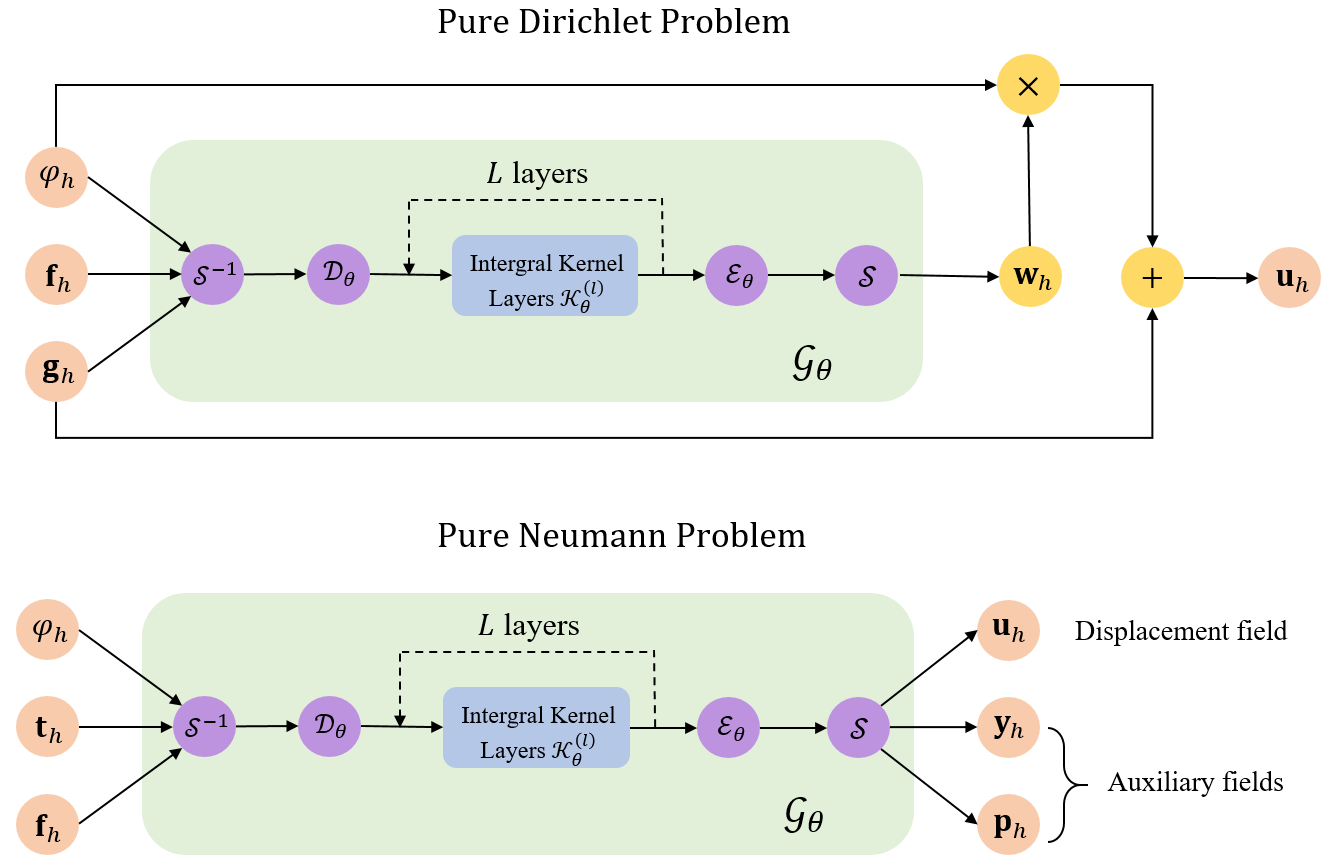}
	\caption{The architecture of the Fourier neural operator for Dirichlet boundary conditions and Neumann boundary conditions.}
	\label{fig:fno_architecture}
\end{figure}

Ultimately, after $L$ iterations of coupled global-local processing, the decoder network $\mathcal{D}_\theta$—which generally comprises a two-layer multi-layer perceptron (MLP)—projects the terminal hidden state $\mathbf{h}^{(L)}$ from the high-dimensional latent space $d_h$ back into the physical output space $d_{out}$. This final projection reconstructs either $\mathbf{w}_h$ or the tuple $(\mathbf{u}_h, \mathbf{y}_h, \mathbf{p}_h)$ over the background Cartesian grid, depending on the stacked inputs. The overall architecture of the Fourier neural operator is illustrated in Fig.~\ref{fig:fno_architecture}.

\begin{remark}
	In many FNO implementations, the spatial coordinates are concatenated with the physical fields as additional input channels (often written as $\mathbf{X}$ in the literature). This positional information helps the network distinguish points that may have similar local field values but different geometric locations, which is especially important on irregular domains represented by level sets. In practice, adding coordinate channels typically improves both optimization stability and final accuracy, particularly for extrapolation to unseen geometries or load distributions.
\end{remark}

\subsection{Loss functions and solver warm starts}

\subsubsection{Data-driven loss function}

When paired training data are available, a purely data-driven objective can be defined in the spirit of standard operator learning on parametric domains described by level-set fields on a fixed background mesh. We assume a set of paired samples is given by
\begin{equation}
	\mathcal{D}_{data}=\{(\mathcal{D}_{\mathrm{input}}^{(m)};\mathbf{u}_h^{ref,(m)})\}_{m=1}^{N_s}.
	\label{eq:paired_data}
\end{equation}
In this supervised setting, we denote the corresponding data-driven neural operator by $\mathcal{G}^{data}_\theta$. The input--output map follows the same $\varphi$-FEM branch structure as WINO in Section~2.3. For pure Dirichlet problems, we write $\mathcal{D}_{\mathrm{input}}^{(m)}=(\varphi_h^{(m)},\mathbf{g}_h^{(m)},\mathbf{f}_{h}^{(m)})$ and the operator predicts the homogeneous coefficient,
\begin{equation}
	\mathcal{G}^{data}_\theta : \mathcal{D}_{\mathrm{input}}^{(m)} \mapsto \mathbf{w}_h^{(m)},
	\label{eq:supervised_mapping_dirichlet}
\end{equation}
after which the discrete displacement is recovered by the $\varphi$-FEM lift \eqref{eq:phifem_dirichlet_compose},
\begin{equation}
	\mathbf{u}_h^{(m)}=\varphi_h^{(m)}\mathbf{w}_h^{(m)}+\mathbf{g}_h^{(m)}.
	\label{eq:supervised_lift}
\end{equation}
For the Neumann-oriented $\varphi$-FEM branch, we instead set $\mathcal{D}_{\mathrm{input}}^{(m)}=(\varphi_h^{(m)},\mathbf{t}_h^{(m)},\mathbf{f}_{h}^{(m)})$ and the operator predicts the displacement field directly,
\begin{equation}
	\mathcal{G}^{data}_\theta : \mathcal{D}_{\mathrm{input}}^{(m)} \mapsto \mathbf{u}_h^{(m)}.
	\label{eq:supervised_mapping}
\end{equation}
Denoting by $\mathbf{u}_h^{(m)}$ the supervised prediction obtained either from the Dirichlet lift \eqref{eq:supervised_lift} or from the Neumann map \eqref{eq:supervised_mapping}, and following the $\varphi$-FEM-FNO framework \cite{duprez2025varphi}, we define the supervised loss
\begin{equation}
	\mathcal{L}_{\text{FNO}}=\frac{1}{N_s}\sum_{m=1}^{N_s}
	\left(
	\left\|
	\mathbf{u}_h^{(m)}-\mathbf{u}_h^{ref,(m)}
	\right\|_{L^2(\Omega_h^{(m)})}^2
	+
	\left\|
	\nabla\mathbf{u}_h^{(m)}-\nabla\mathbf{u}_h^{ref,(m)}
	\right\|_{L^2(\Omega_h^{(m)})}^2
	\right),
	\label{eq:fno_loss}
\end{equation}
where $\Omega_h^{(m)}$ is the active computational domain induced by $\varphi_h^{(m)}$ according to \eqref{eq:omega_h}, and the discrete gradients in the second term are evaluated by the centered finite difference approximation. However, generating high-fidelity labeled data (e.g., by $\varphi$-FEM or other numerical solvers) may be prohibitively time-consuming for large-scale parametric studies. Therefore, in the next subsection, we introduce a label-free model and its weak-form physics-informed loss function.

\subsubsection{Weak form physics-informed loss function}

In this subsection, we firstly introduce a weak-form physics-informed loss that enables training without generating labeled solution data for pure Neumann boundary hyperealsticity problems. We collect input tuples in the same indexed form as in Section~2.4.1, but without reference solutions:
\begin{equation}
	\mathcal{D}_{phys}=\{\mathcal{D}_{\mathrm{input}}^{(m)}\}_{m=1}^{N_s},
	\label{eq:phys_data}
\end{equation}
Building upon this idea, we propose the Weak-form Physics-Informed Neural Operator (WINO), whose parameters are learned by minimizing a loss consistent with the $\varphi$-FEM formulation. For each sample $m$, the map
\begin{equation}
	\mathcal{G}_\theta : \mathcal{D}_{\mathrm{input}}^{(m)} \mapsto (\mathbf{u}_h^{(m)},\mathbf{y}_h^{(m)},\mathbf{p}_h^{(m)})
	\label{eq:wino_sample_map}
\end{equation}
produces the discrete fields on which the residuals are evaluated. Specifically, the total WINO loss combines the weak-form residual of the momentum equation \eqref{eq:phifem_weak_form} with the strong-form residuals of the auxiliary constraints \eqref{eq:strong_yp} introduced by $\varphi$-FEM, averaged over $m=1,\ldots,N_s$.

Let $T^e = [x_i,x_{i+1}] \times [y_j,y_{j+1}]$ be a rectangular cell of the background mesh, and let $(\xi,\eta) \in [-1,1]^2$ denote the local coordinates on the reference square. The mapping between the physical coordinates $(x,y)$ and the reference coordinates $(\xi,\eta)$ is given by
\begin{equation*}
	x=\frac{1-\xi}{2}x_i+\frac{1+\xi}{2}x_{i+1}, \qquad
	y=\frac{1-\eta}{2}y_j+\frac{1+\eta}{2}y_{j+1}.
\end{equation*}
The four bilinear shape functions on the reference square are given by
\begin{equation}
	\begin{aligned}
		N_1(\xi,\eta) &= \frac{1}{4}(1-\xi)(1-\eta), \qquad
		N_2(\xi,\eta) = \frac{1}{4}(1+\xi)(1-\eta), \\
		N_3(\xi,\eta) &= \frac{1}{4}(1+\xi)(1+\eta), \qquad
		N_4(\xi,\eta) = \frac{1}{4}(1-\xi)(1+\eta).
	\end{aligned}
	\label{eq:rect_shape_function}
\end{equation}
Accordingly, for any sample $m$ and any field $\psi \in \{\mathbf{u}_h^{(m)},\mathbf{y}_h^{(m)},\mathbf{p}_h^{(m)},\varphi_h^{(m)},\mathbf{t}_{h}^{(m)}, \mathbf{f}_{h}^{(m)}\}$, we denote its restriction to the element $T^e$ by
\begin{equation*}
	\psi^{e} := \psi|_{T^e}.
\end{equation*}
The corresponding bilinear interpolation can be written as
\begin{equation}
	\psi^{e}(\xi,\eta) \approx \mathbf{N}^{e}(\xi,\eta)^T\,\bar{\boldsymbol{\psi}}^{e},
	\label{eq:element_interpolation}
\end{equation}
where $\mathbf{N}^{e}(\xi,\eta)=\bigl[N_1,\ldots,N_4\bigr]^T$ lists the bilinear shape functions on the reference element and $\bar{\boldsymbol{\psi}}^{e}=\bigl[\bar{\psi}_1^{e},\ldots,\bar{\psi}_4^{e}\bigr]^T$ collects the nodal values of $\psi$ at the vertices of $T^e$. A second-order $Q9$ element uses the same matrix form with $\mathbf{N}^{Q9}\in\mathbb{R}^9$ and $\bar{\boldsymbol{\psi}}^{e}\in\mathbb{R}^9$, see \cite{eshaghi2025variational} for more details. On each element, the physical gradient is assembled in matrix form as
\begin{equation}
	\nabla\psi^{e}(\xi,\eta)=\mathbf{B}^{e}(\xi,\eta)\,\bar{\boldsymbol{\psi}}^{e},
	\label{eq:grad_psi_element}
\end{equation}
where $\mathbf{B}^{e}(\xi,\eta)=\bigl[\nabla N_1\ \cdots\ \nabla N_4\bigr]\in\mathbb{R}^{2\times 4}$ has the $a$-th column $\nabla N_a$. Each $\nabla N_a$ is computed from $\partial_\xi N_a$ and $\partial_\eta N_a$ together with the Jacobian of the mapping $(\xi,\eta)\mapsto(x,y)$. In particular, the notation \eqref{eq:element_interpolation} applies to the neural-operator outputs $\mathbf{u}_h^{(m)}$, $\mathbf{y}_h^{(m)}$, and $\mathbf{p}_h^{(m)}$, as well as to the discrete level-set field $\varphi_h^{(m)}$, the discrete Neumann loading $\mathbf{t}_{h}^{(m)}$ and the body force $\mathbf{f}_h^{(m)}$. 

To evaluate the weak formulation \eqref{eq:phifem_weak_form} in the discrete setting, we select test functions from a nodal basis on the background mesh. Since the displacement $\mathbf{u}_h^{(m)}$ is vector-valued, we use the corresponding vector-valued basis $\{\mathbf{v}_{i,j}\}$, with $\mathbf{v}_{i,j}:=v_i\,\mathbf{e}_j$ and $j\in\{1,2\}$, where $v_i$ denotes the scalar $Q4$ shape function associated with node $i$ and $\{\mathbf{e}_j\}_{j=1}^2$ is the canonical basis of $\mathbb{R}^2$. On each element $T^e$, the scalar part satisfies
\begin{equation}
	v_i|_{T^e}(\xi,\eta)=
	\begin{cases}
		N_a(\xi,\eta), & \text{if } \mathbf{x}_i \text{ is the } a\text{-th vertex of } T^e,\\
		0, & \text{if } \mathbf{x}_i \notin T^e,
	\end{cases}
	\qquad a\in\{1,2,3,4\},
\end{equation}
and the resulting global basis satisfies $v_i(\mathbf{x}_k)=\delta_{ik}$. Consequently, the vector basis satisfies the component-wise Kronecker-delta property
\begin{equation*}
	\mathbf{v}_{i,j}(\mathbf{x}_k)=v_i\mathbf{e}_j=\delta_{ik}\mathbf{e}_j.
\end{equation*}
Then the weak-form residual of the momentum equation \eqref{eq:phifem_weak_form} for sample $m$ is written as
\begin{equation}
	\mathcal{R}_u^{(m)}(\mathbf{v}_{i,j})=\int_{\Omega_h^{(m)}} \mathbf{P}(\mathbf{F}(\mathbf{u}_h^{(m)})):\nabla \mathbf{v}_{i,j}\, d\mathbf{x}
	+\int_{\partial\Omega_h^{(m)}} (\mathbf{y}_h^{(m)} \mathbf{n})\cdot \mathbf{v}_{i,j}\, ds - \int_{\Omega_h^{(m)}} \mathbf{f}_{h}^{(m)}\cdot\mathbf{v}_{i,j}\, d\mathbf{x}.
	\label{eq:weak_residual_u}
\end{equation}
Accordingly, the weak loss is defined by the mean over samples of the squared sum of these nodal residuals,
\begin{equation}
	\mathcal{L}_{weak}=\frac{1}{N_s}\sum_{m=1}^{N_s}\frac{1}{2N_{nodes}}\sum_{i=1}^{N_{nodes}}\sum_{j=1}^{2} \bigl|\mathcal{R}_u^{(m)}(\mathbf{v}_{i,j})\bigr|^2.
	\label{eq:weak_loss}
\end{equation}
For the auxiliary variables $\mathbf{y}_h^{(m)}$ and $\mathbf{p}_h^{(m)}$ defined on the cut-cell domain $\Omega_h^{\Gamma,(m)}$, we enforce the $\varphi$-FEM auxiliary constraints and the divergence condition through strong-form residuals, in accordance with \eqref{eq:strong_yp}:
\begin{align}\label{eq:strong_residuals}
	\mathcal{R}_y^{(m)} &= \mathbf{y}_h^{(m)} + \mathbf{P}(\mathbf{F}(\mathbf{u}_h^{(m)})),\nonumber \\
	\mathcal{R}_p^{(m)} &= \mathbf{y}_h^{(m)} \nabla \varphi_h^{(m)} + \frac{1}{h}\mathbf{p}_h^{(m)} \varphi_h^{(m)} + \mathbf{t}_h^{(m)} |\nabla \varphi_h^{(m)}|, \\
	\mathcal{R}_{d}^{(m)} &= \nabla \cdot \mathbf{y}_h^{(m)}.\nonumber
\end{align}
The corresponding squared residual loss is given by
\begin{equation}
\begin{aligned}
	\mathcal{L}_{sq} &= \lambda_1 \mathcal{L}_y + \lambda_2 \mathcal{L}_p + \lambda_3 \mathcal{L}_d \\
	&= \frac{1}{N_s}\sum_{m=1}^{N_s}\frac{1}{N_E^{(m)}}\int_{\Omega_h^{\Gamma,(m)}}
	\left(
	\lambda_1 \bigl\|\mathcal{R}_y^{(m)}\bigr\|^2
	+\lambda_2 \bigl\|\mathcal{R}_p^{(m)}\bigr\|^2
	+\lambda_3 \bigl\|\mathcal{R}_{d}^{(m)}\bigr\|^2
	\right)\, d\mathbf{x} \\
	&= \frac{1}{N_s}\sum_{m=1}^{N_s}\frac{1}{N_E^{(m)}}\sum_{E \subset \Omega_h^{\Gamma,(m)}}\int_{E}
	\left(
	\lambda_1 \bigl\|\mathcal{R}_y^{(m)}\bigr\|^2
	+\lambda_2 \bigl\|\mathcal{R}_p^{(m)}\bigr\|^2
	+\lambda_3 \bigl\|\mathcal{R}_{d}^{(m)}\bigr\|^2
	\right)\, d\mathbf{x},
\end{aligned}
\label{eq:strong_loss}
\end{equation}
where $\mathcal{L}_y$, $\mathcal{L}_p$, and $\mathcal{L}_d$ are defined by
\begin{equation*}
	\mathcal{L}_y=\frac{1}{N_s}\sum_{m=1}^{N_s}\frac{1}{N_E^{(m)}}\int_{\Omega_h^{\Gamma,(m)}}\bigl\|\mathcal{R}_y^{(m)}\bigr\|^2\, d\mathbf{x},
\end{equation*}
\begin{equation*}
	\mathcal{L}_p=\frac{1}{N_s}\sum_{m=1}^{N_s}\frac{1}{N_E^{(m)}}\int_{\Omega_h^{\Gamma,(m)}}\bigl\|\mathcal{R}_p^{(m)}\bigr\|^2\, d\mathbf{x},
\end{equation*}
\begin{equation*}
	\mathcal{L}_d=\frac{1}{N_s}\sum_{m=1}^{N_s}\frac{1}{N_E^{(m)}}\int_{\Omega_h^{\Gamma,(m)}}\bigl\|\mathcal{R}_{d}^{(m)}\bigr\|^2\, d\mathbf{x},
\end{equation*}
respectively; $\lambda_1$, $\lambda_2$, $\lambda_3$ are positive penalty parameters, $E$ are the elements of the cut cells, and $N_E^{(m)}$ is the number of such elements for sample $m$. The total WINO training loss function is given by
\begin{equation}
	\mathcal{L}_{total}=\mathcal{L}_{weak}+\mathcal{L}_{sq}.
	\label{eq:total_loss}
\end{equation}

\begin{figure}[t]
	\centering
	\includegraphics[width=0.88\textwidth]{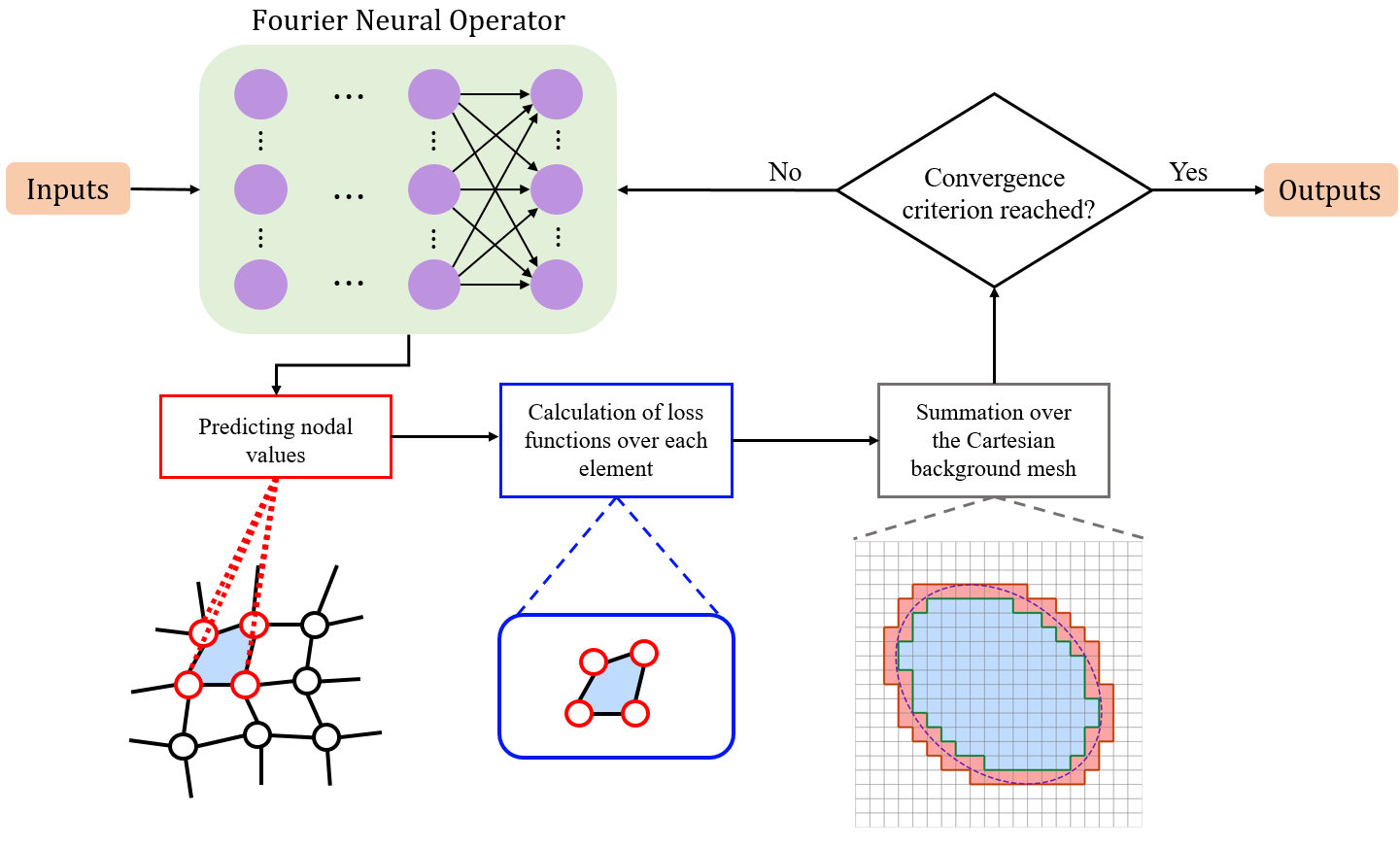}
	\caption{Schematic flowchart of WINO training: inputs without labeled data, FNO prediction on nodes, calculation of loss (weak form based residuals), and summation of loss over the Cartesian background mesh.}
	\label{fig:wino_schematic}
\end{figure}

Therefore, to train a solution operator for hyperelastic problems on varying domains, WINO only requires minimizing \eqref{eq:total_loss}. The spatial derivatives entering the finite-element assembly are evaluated at the Gaussian quadrature points of each element from the shape-function gradients and the isoparametric map \eqref{eq:grad_psi_element}, so the weak and squared residuals do not rely on additional numerical differentiation of the network outputs in physical space. Compared with a purely strong-form residual, the weak formulation requires lower regularity of the discrete fields, which reduces both the cost of loss evaluation and the sensitivity to errors introduced when approximating higher-order derivatives. In the same spirit as $\varphi$-FEM, WINO can learn solution maps for families of domains on a single fixed Cartesian background mesh, without generating body-fitted meshes or performing exact quadrature on arbitrarily cut cells. Many physics-informed neural schemes instead use an energy-based loss. However, some unfitted discretizations, including $\varphi$-FEM, are not tied to a symmetric variational principle or an equivalent scalar energy. Since WINO builds its objective directly from the $\varphi$-FEM weak residual and the auxiliary squared residuals, it provides an efficient training strategy in such nonsymmetric, non-energy settings. Domain and boundary terms are accumulated element-wise by quadrature over $\Omega_h^{(m)}$ and $\Omega_h^{\Gamma,(m)}$. Finally, Fig.~\ref{fig:wino_schematic} gives a schematic flowchart of the WINO training pipeline, and Algorithm~\ref{alg:WINO} lists the same procedure in pseudo-code.

\begin{remark}
	For pure Dirichlet problems, for each sample $m$ we use
	\begin{equation*}
		\mathcal{G}_\theta : \mathcal{D}_{\mathrm{input}}^{(m)} \mapsto \mathbf{w}_h^{(m)},
	\end{equation*}
	where $\mathcal{D}_{\mathrm{input}}^{(m)}:=(\varphi_h^{(m)},\mathbf{g}_h^{(m)},\mathbf{f}_{h}^{(m)})$ in this branch,
	together with the lift representation $\mathbf{u}_h^{(m)}=\varphi_h^{(m)}\mathbf{w}_h^{(m)}+\mathbf{g}_h^{(m)}$ according to \eqref{eq:phifem_dirichlet_compose}.
	Accordingly, we discretize these variables on each mesh using the same shape functions, and the weak residual keeps the same form as \eqref{eq:weak_residual_u}, except that the boundary flux variable $\mathbf{y}_h$ is replaced by the physical stress $-\mathbf{P}(\mathbf{F}(\mathbf{u}_h^{(m)}))$:
	\begin{equation}
		\begin{aligned}
			\mathcal{R}_{u,D}^{(m)}(\mathbf{v}_{i,j})={}&\int_{\Omega_h^{(m)}} \mathbf{P}(\mathbf{F}(\mathbf{u}_h^{(m)})):\nabla (\varphi_h\mathbf{v}_{i,j})\, d\mathbf{x}
			\\&-\int_{\partial\Omega_h^{(m)}} \bigl(\mathbf{P}(\mathbf{F}(\mathbf{u}_h^{(m)}))\mathbf{n}\bigr)\cdot (\varphi_h\mathbf{v}_{i,j})\, ds 
			-\int_{\Omega_h^{(m)}} \mathbf{f}_{h}^{(m)}\cdot(\varphi_h\mathbf{v}_{i,j})\, d\mathbf{x}.
		\end{aligned}
	\end{equation}
	Here, the test function is replaced by $\varphi_h\mathbf{v}_{i,j}$, following the definition in Section~\ref{subsec:phi_fem}. Hence the Dirichlet weak loss can be defined as
	\begin{equation}
		\mathcal{L}_{total}^D=\mathcal{L}_{weak}^{D}=\frac{1}{N_s}\sum_{m=1}^{N_s}\frac{1}{2N_{nodes}}\sum_{i=1}^{N_{nodes}}\sum_{j=1}^{2}
		\bigl|\mathcal{R}_{u,D}^{(m)}(\mathbf{v}_{i,j})\bigr|^2.
		\label{eq:total_loss_Dirichlet}
	\end{equation}
\end{remark}

\begin{remark}
	The WINO loss functions \eqref{eq:total_loss} and \eqref{eq:total_loss_Dirichlet} omit the ghost-penalty stabilization terms present in the reference $\varphi$-FEM discretization. Since the operator network produces globally smooth fields on the fixed background mesh, these cut-cell jump penalties are not essential to close the WINO training formulation.  We did experiment with adding analogous ghost-penalty residuals to the loss in Appendix~\ref{app:WINO_Gh}, but found that they did not improve accuracy on the reported benchmark, whereas the total training time increased substantially. Therefore, we retained the simpler loss formulation without ghost penalties.
\end{remark}

% \begin{remark}
% 	When labeled training data are available, one may combine the supervised loss \eqref{eq:data_loss} with the physics-informed loss \eqref{eq:total_loss} into the hybrid objective
% 	\begin{equation}
% 		\mathcal{L}=\mathcal{L}_{data}+\omega_{phys}\mathcal{L}_{total},
% 		\label{eq:hybrid_loss}
% 	\end{equation}
% 	where $\omega_{phys}\ge 0$ is a scalar weight.
% 	\end{remark}

\begin{algorithm}[t]
\caption{Training procedure of WINO (Dirichlet lift or Neumann $\varphi$-FEM branch)}
\label{alg:WINO}
\KwIn{Background mesh $\mathcal{T}_h$, initial parameters $\theta$, and a training set $\mathcal{D}=\{\mathcal{D}_{\mathrm{input}}^{(m)}\}_{m=1}^{N_s}$, where $\mathcal{D}_{\mathrm{input}}^{(m)}=(\varphi_h^{(m)},\mathbf{g}_h^{(m)},\mathbf{f}_{h}^{(m)})$ for the Dirichlet branch and $\mathcal{D}_{\mathrm{input}}^{(m)}=(\varphi_h^{(m)},\mathbf{t}_h^{(m)},\mathbf{f}_{h}^{(m)})$ for the Neumann branch.}
\KwOut{Trained parameters $\theta^\ast$.}
Initialize the neural-operator parameters $\theta$\;
\While{the stopping criterion is not satisfied}{
	Sample a mini-batch from $\mathcal{D}$ of size $N_{bs}$\;
	\For{each sample $m$ in the mini-batch}{
		\eIf{Dirichlet lift branch}{
			Evaluate $\mathcal{G}_\theta(\mathcal{D}_{\mathrm{input}}^{(m)})$ to obtain $\mathbf{w}_h^{(m)}$, set $\mathbf{u}_h^{(m)}=\varphi_h^{(m)}\mathbf{w}_h^{(m)}+\mathbf{g}_h^{(m)}$, and build element-wise interpolants using \eqref{eq:element_interpolation}\;
			For each element of $\mathcal{T}_h^{\mathcal{O}}$, assemble $\mathcal{R}_{u,D}^{(m)}(\mathbf{v}_{i,j})$ for all test functions as in the Dirichlet weak residual stated before \eqref{eq:total_loss_Dirichlet}\;
			Sum the contributions over all elements to update the mini-batch estimate of $\mathcal{L}_{\mathrm{total}}^{D}$ in \eqref{eq:total_loss_Dirichlet}\;
		}{
			Evaluate $\mathcal{G}_\theta(\mathcal{D}_{\mathrm{input}}^{(m)})$ to obtain $(\mathbf{u}_h^{(m)},\mathbf{y}_h^{(m)},\mathbf{p}_h^{(m)})$, and build element-wise interpolants using \eqref{eq:element_interpolation}\;
			For each element of $\mathcal{T}_h^{\mathcal{O}}$, assemble $\mathcal{R}_u^{(m)}(\mathbf{v}_{i,j})$ from \eqref{eq:weak_residual_u} for all test functions, and evaluate $\mathcal{R}_y^{(m)}$, $\mathcal{R}_p^{(m)}$, and $\mathcal{R}_{d}^{(m)}$ from \eqref{eq:strong_residuals} on cut cells\;
			Sum the residuals over all elements to update the mini-batch estimates of $\mathcal{L}_{\mathrm{weak}}$, $\mathcal{L}_{\mathrm{sq}}$, and $\mathcal{L}_{\mathrm{total}}$ according to \eqref{eq:weak_loss}, \eqref{eq:strong_loss}, and \eqref{eq:total_loss}\;
		}
	}
	Update the neural-operator parameters $\theta$ by backpropagation and an optimizer\;
}
\Return{the trained neural-operator parameters $\theta^\ast$}\;
\end{algorithm}

\subsubsection{Neural operator warm starts}
\label{sec:nows}

Discretizing the Dirichlet weak form \eqref{eq:phifem_dirichlet_weak} or the Neumann weak form \eqref{eq:phifem_weak_form} yields a finite-dimensional nonlinear system $\mathbf{R}(\mathbf{U})=\mathbf{0}$, where $\mathbf{U}\in\mathbb{R}^n$ collects the nodal degrees of freedom of the primary unknowns in the assembled $\varphi$-FEM system: for Neumann problems, $\mathbf{U}$ stacks the nodal values of $(\mathbf{u}_h,\mathbf{y}_h,\mathbf{p}_h)$; for Dirichlet problems with the lift \eqref{eq:phifem_dirichlet_compose}, $\mathbf{U}$ stacks those of $\mathbf{w}_h$. We solve this system by Newton iteration \cite{kelley2003solving, pernice1998nitsol}; each linearization is solved by Krylov methods (CG when the tangent matrix is SPD, otherwise restarted GMRES \cite{saad1986gmres,morgan2002gmres}), optionally with GAMG preconditioning in PETSc \cite{balay2019petsc}. Further details are given in Appendix~\ref{app:nows}. Neural operator warm starts (NOWS) \cite{eshaghi2025nows} evaluate the trained WINO mapping on the problem inputs to obtain an initial Newton iterate $\mathbf{U}^{(0)}$: in the Neumann branch, the predicted nodal fields $(\mathbf{u}_h,\mathbf{y}_h,\mathbf{p}_h)$ are stacked directly into $\mathbf{U}^{(0)}$ and may also supply the inner Krylov initial guess; in the Dirichlet branch, $\mathbf{w}_h$ is taken from WINO and $\mathbf{u}_h=\varphi_h\mathbf{w}_h+\mathbf{g}_h$ is recovered via \eqref{eq:phifem_dirichlet_compose}. A better $\mathbf{U}^{(0)}$ typically reduces $\|\mathbf{r}_0\|$, thereby lowering the number of outer Newton steps and inner Krylov iterations.

A critical advantage of embedding NOWS into our $\varphi$-FEM framework is that it preserves the stability, interpretability, and rigorous convergence guarantees native to traditional finite element solvers: the hybrid initialization accelerates convergence while strictly preserving the exactness of the final solution. Because NOWS intervenes only at initialization, it requires no modification to the $\varphi$-FEM discretization, the Newton/Krylov algorithms, or any preconditioning strategy, and thus provides a flexible and trustworthy bridge between rapid operator learning and high-fidelity computational mechanics. Benchmark-specific tolerances and solver settings are reported in Section~\ref{section_Results}.

\section{Results}\label{section_Results}

\begin{table}[t]
	\raggedright
	\footnotesize
	\setlength{\tabcolsep}{4pt}
	\caption{Training setup across the numerical benchmarks.}
	\label{tab:hyperparameters}
	\begin{tabularx}{\textwidth}{@{}l >{\raggedright\arraybackslash}X >{\raggedright\arraybackslash}X >{\raggedright\arraybackslash}X >{\raggedright\arraybackslash}X >{\raggedright\arraybackslash}X@{}}
		\toprule
		Hyperparameter & Elliptical shape & Random shape & Plate with a hole & Cook's membrane & Pressure vessel \\
		\midrule
		Training set size & 500 & $500$ & $300$ & $1000$ & $1000$ \\
		Test set size & 100 & $100$ & $100$ & $100$ & $100$ \\
		Epochs & 500 & $500$ & $1000$ & $5000$ & $2000$ \\
		Batch size & 50 & 50 & 50 & 50 & 50 \\
		Resolution size & $64\times64$ & $64\times64$ & $64\times64$ & $51\times51$ & $51\times51$ \\
		Learning rate & 0.005 & 0.005 & 0.005 & 0.01 & 0.002 \\
		Patience & 250 & 250 & 250 & 250 & 250 \\
		Decay factor & 0.7 & 0.7 & 0.5 & 0.8 & 0.7 \\
		Fourier modes & 16 & 16 & 16 & 16 & 16 \\
		Number of layers & 4 & 4 & 4 & 4 & 4 \\
		Width & 32 & 32 & 32 & 32 & 32 \\
		Parameter count & $4211714$ & $4211714$ & $4212424$ & $4212424$ & $4212552$ \\
		\bottomrule
	\end{tabularx}
\end{table}

\begin{table}[t]
	\raggedright
	\scriptsize
	\setlength{\tabcolsep}{4pt}
	\caption{Overall comparison of relative $L^2$ errors and computational costs (data generation + training time) across the numerical benchmarks.}
	\label{tab:performance}
	\begin{tabularx}{\textwidth}{@{}l >{\raggedright\arraybackslash}X >{\raggedright\arraybackslash}X >{\raggedright\arraybackslash}X >{\raggedright\arraybackslash}X >{\raggedright\arraybackslash}X@{}}
		\toprule
		\multirow{2}{*}{Method} 
		& \multicolumn{2}{c}{{Pure Dirichlet}} 
		& \multicolumn{3}{c}{{Mixed Dirichlet-Neumann}} \\
		\cmidrule(lr){2-3} \cmidrule(lr){4-6}
		& Elliptical shape & Random shape & Plate with a hole & Cook's membrane & Pressure vessel \\
		\midrule
		WINO & $0.79 \pm 0.43\%$ & $0.65 \pm 0.33\%$ & $2.26 \pm 1.10\%$ & $3.83 \pm 2.14\%$ & $3.71 \pm 0.78\%$ \\
		& $0.3 + 1248.8$ s & $0.3 + 1366.5$ s & $0.3 + 2912.8$ s & $0.6 + 19242.1$ s & $0.6 + 17868.4$ s \\
		\midrule
		$\varphi$-FEM-FNO & $0.86 \pm 0.53\%$ & $0.79 \pm 0.32\%$ & $0.78 \pm 0.40\%$ & $0.05 \pm 0.02\%$ & $0.03 \pm 0.01\%$ \\
		& $1392.5 + 644.4$ s & $8603.1 + 783.8$ s & $11635.8 + 1973.6$ s & $73096.6 + 12757.2$ s & $13498.2 + 14308.4$ s \\
		\midrule
		WINO+data & $0.82 \pm 0.54\%$ & $0.55 \pm 0.26\%$ & $0.54 \pm 0.20\%$ & $0.07 \pm 0.01\%$ & $0.01 \pm 0.01\%$ \\
		& $1392.5 + 1258.7$ s & $8603.1 + 1391.3$ s & $11635.8 + 2976.8$ s & $73096.6 + 19679.1$ s & $13498.2 + 18755.4$ s \\
		\bottomrule
	\end{tabularx}
\end{table}

In this section, we report numerical results for the proposed framework on parametric hyperelasticity problems with complex, varying geometries. The experiments are organized into two groups: pure Dirichlet problems and mixed Dirichlet and Neumann problems. In the second group, the level-set function $\varphi$ encodes only the Neumann (implicitly defined) part of the boundary, whereas the Dirichlet segments lie on edges of the Cartesian background mesh $\mathcal{T}_h^{\mathcal{O}}$. Accordingly, the high-fidelity $\varphi$-FEM discretization and the WINO model use the Neumann-oriented $\varphi$-FEM branch on cut cells for the level-set portion, while prescribed displacements on the mesh-aligned boundary are imposed with a standard conforming finite-element treatment. We compare three approaches: a purely data-driven $\varphi$-FEM-FNO, the physics-informed WINO, and WINO augmented with labeled data. All experiments are run on a single NVIDIA RTX~4090 GPU.

Let $\mathbf{u}^{\mathrm{ref}}$ denote the high-fidelity reference displacement and $\mathbf{u}^{\mathrm{pred}}$ the prediction, both restricted to the active physical domain $\Omega_h$ used for integration. We report the relative $L^2$ error
\begin{equation}
	\|e\|_{L^2}=\frac{\|\mathbf{u}^{\mathrm{pred}}-\mathbf{u}^{\mathrm{ref}}\|_{L^2(\Omega_h)}}{\|\mathbf{u}^{\mathrm{ref}}\|_{L^2(\Omega_h)}},
	\label{eq:rel_err_L2}
\end{equation}
the relative $H^1$ \emph{seminorm} error
\begin{equation}
	\|e\|_{H^1}=\frac{|\mathbf{u}^{\mathrm{pred}}-\mathbf{u}^{\mathrm{ref}}|_{H^1(\Omega_h)}}{|\mathbf{u}^{\mathrm{ref}}|_{H^1(\Omega_h)}},\qquad
	|\mathbf{v}|_{H^1(\Omega_h)}^2=\int_{\Omega_h}\nabla\mathbf{v}:\nabla\mathbf{v}\,\mathrm{d}\mathbf{x},
	\label{eq:rel_err_H1}
\end{equation}
To further highlight the physical relevance of the error, we also evaluate an energy-norm error measure. Let $\mathbf{U}^{\mathrm{ref}}, \mathbf{U}^{\mathrm{pred}} \in \mathbb{R}^n$ denote the vectors of finite-element displacement degrees of freedom associated with $\mathbf{u}^{\mathrm{ref}}$ and $\mathbf{u}^{\mathrm{pred}}$, respectively. 
Let $\mathbf{K} := \mathbf{J}(\mathbf{U}^{\mathrm{ref}})$ denote the tangent stiffness matrix, i.e., the Jacobian of the internal force vector in \eqref{eq:newton_linear}, evaluated at the converged reference configuration. We define the discrete energy seminorm
\begin{equation*}
\|\mathbf{U}^{\mathrm{pred}} - \mathbf{U}^{\mathrm{ref}}\|_{\mathbf{K}} 
= \sqrt{(\mathbf{U}^{\mathrm{pred}} - \mathbf{U}^{\mathrm{ref}})^{\mathrm{T}} 
\mathbf{K} 
(\mathbf{U}^{\mathrm{pred}} - \mathbf{U}^{\mathrm{ref}})},
\end{equation*}
which measures the error in terms of the elastic energy induced by the displacement difference. Note that this defines a seminorm, since $\mathbf{K}$ may not be strictly positive definite in general. The relative energy error is then given by
\begin{equation}
\|e\|_{E} = 
\frac{\|\mathbf{U}^{\mathrm{pred}} - \mathbf{U}^{\mathrm{ref}}\|_{\mathbf{K}}}
{\|\mathbf{U}^{\mathrm{ref}}\|_{\mathbf{K}}}, 
\qquad
\|\mathbf{U}^{\mathrm{ref}}\|_{\mathbf{K}} 
= \sqrt{(\mathbf{U}^{\mathrm{ref}})^{\mathrm{T}} 
\mathbf{K} 
\mathbf{U}^{\mathrm{ref}}},
\label{eq:rel_err_energy}
\end{equation}
which corresponds (up to a constant factor) to the relative error in stored elastic energy. Here, $\mathbf{K}$ is assembled using the same $\varphi$-FEM tangent operator as employed in the nonlinear solver. Neural-operator hyperparameters are listed in Table~\ref{tab:hyperparameters} unless stated otherwise, where \emph{patience} denotes the number of optimization steps between two learning-rate decays. In all benchmarks, the parameters of integral kernel layers are initialized with Xavier initialization, whereas the spectral weights in the global spectral convolutions are initialized from a normal distribution. Table~\ref{tab:performance} summarizes the performance comparison on these test cases. Overall, WINO achieves satisfactory accuracy across all benchmarks. Although $\varphi$-FEM-FNO attains higher accuracy in some cases, WINO has a substantially lower computational cost because it does not require generating $\varphi$-FEM reference solutions. The values quoted as $\mu\pm\sigma$ in the tables denote the sample mean and sample standard deviation of the relative errors over the held-out test dataset. Since the weak-form loss yields an ill-conditioned optimization landscape, we adopt SOAP optimizer (Shampoo with Adam in the preconditioner's eigenbasis) \cite{vyas2024soap} for faster and more stable training. SOAP is motivated by a formal connection between Shampoo (with $1/2$-power preconditioning) and Adafactor, and can be interpreted as performing Adam-like second-moment updates in the rotated coordinate system induced by Shampoo's preconditioner. Compared with standard Shampoo, this design improves computational efficiency while retaining strong preconditioning behavior, and introduces only one additional key hyperparameter over Adam, namely the preconditioning frequency.

\begin{remark}
	Table~\ref{tab:performance} indicates that WINO achieves smaller mean relative $L^2$ errors than $\varphi$-FEM-FNO on the pure Dirichlet benchmarks, whereas $\varphi$-FEM-FNO attains higher accuracy on the mixed Dirichlet-Neumann cases under the same reporting protocol.
	We attribute this to structural alignment of the training objective with the discretization. For pure Dirichlet problems, WINO uses the $\varphi$-FEM lift $\mathbf{u}_h=\varphi_h\mathbf{w}_h+\mathbf{g}_h$, so the network primarily learns the homogeneous component while the weak-form residual matches the primal $\varphi$-FEM equation without Neumann auxiliary fields.
	In contrast, mixed boundary settings couple the displacement with cut-cell auxiliary variables and penalty parameters, leading to a more intricate optimization landscape for WINO. Under these settings, supervised regression of converged $\varphi$-FEM solutions employed by $\varphi$-FEM-FNO appears to provide a learning advantage.
\end{remark}

\subsection{Pure Dirichlet problems}

We first consider pure Dirichlet problems, where the domain geometry is fixed and the boundary conditions are prescribed. We formulate the hyperelastic problem with pure Dirichlet boundary conditions on the elliptic domain. The boundary value problem reads
\begin{equation}
	\begin{split}
		-\nabla \cdot \mathbf{P}(\mathbf{F}(\mathbf{u})) &= \mathbf{f}, \quad \text{in } \Omega, \\
		\mathbf{u} &= \mathbf{g}, \quad \text{on } \Gamma,
	\end{split}
	\label{eq:case_dirichlet}
\end{equation}
where $\Omega$ is the physical domain represented by the level-set function $\varphi$, and $\Gamma=\partial\Omega$. The $\varphi$-FEM is computed by the equation \eqref{eq:phifem_dirichlet_weak} and the loss function of WINO is given by \eqref{eq:total_loss_Dirichlet}. In this subsection, we consider the elliptical shape and the random shape as the physical domain.

\subsubsection{Elliptical shape}\label{subsec:elliptical_shape}

\begin{figure}[t]
	\raggedright
	\begin{minipage}[t]{0.47\textwidth}
		\textbf{(a)}\par\vspace{0.3em}
		\includegraphics[width=\linewidth]{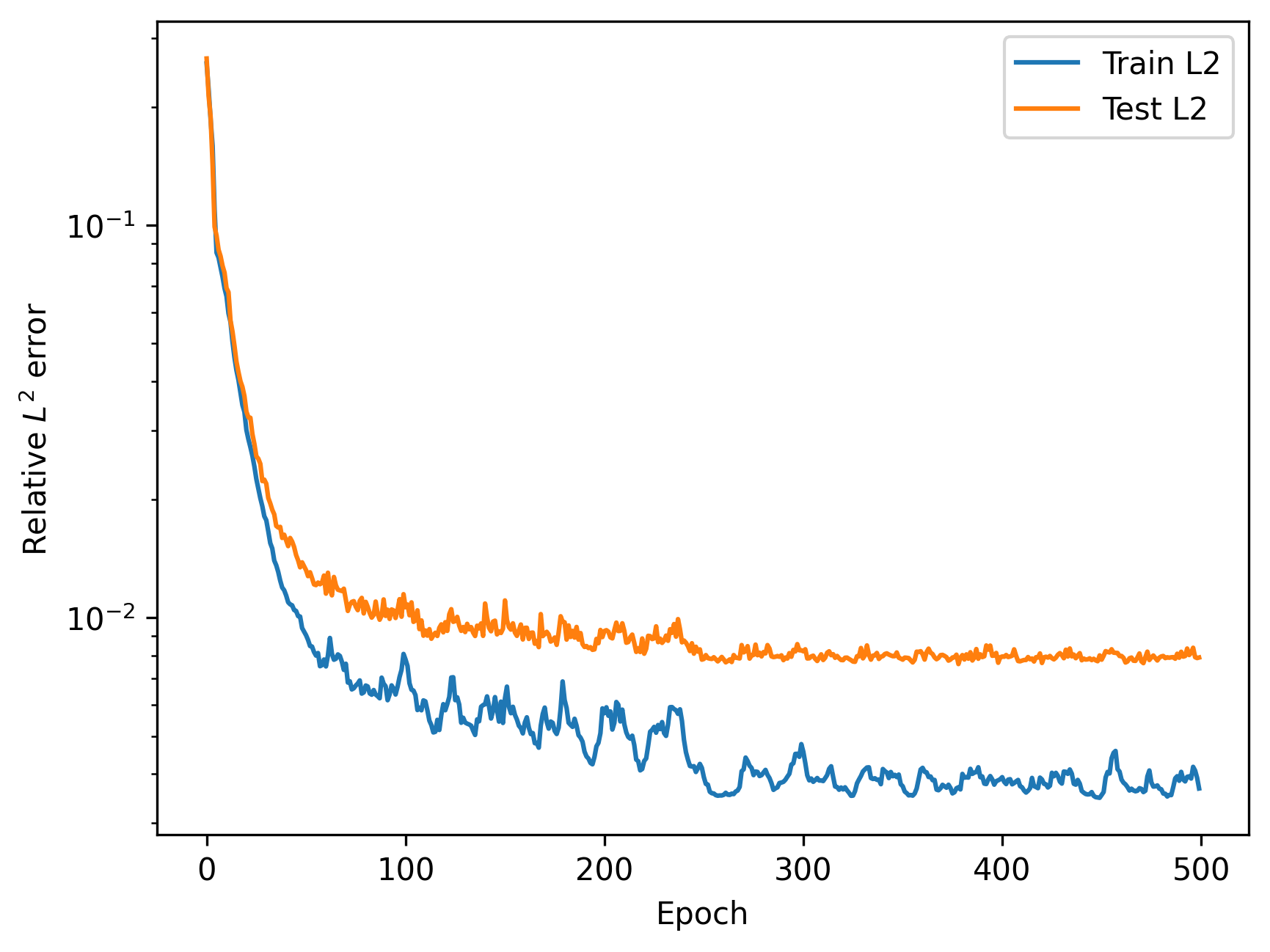}
	\end{minipage}
	\hspace{0.8em}
	\begin{minipage}[t]{0.46\textwidth}
		\textbf{(b)}\par\vspace{0.3em}
		\includegraphics[width=\linewidth]{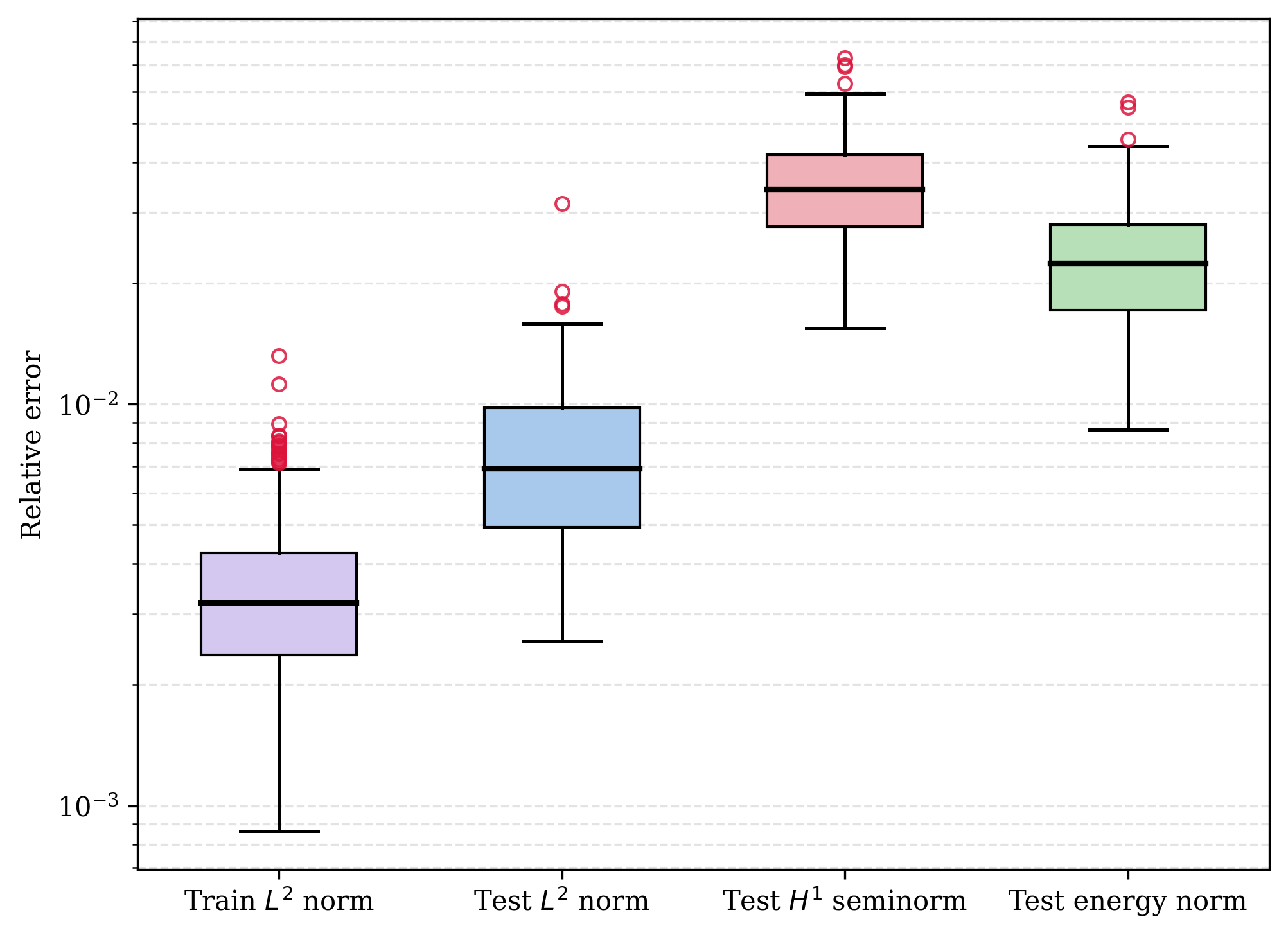}
	\end{minipage}
	\vspace{0.8em}
	\begin{minipage}[t]{\textwidth}
		\textbf{(c)}\par\vspace{0.5em}
		\includegraphics[width=0.18\textwidth]{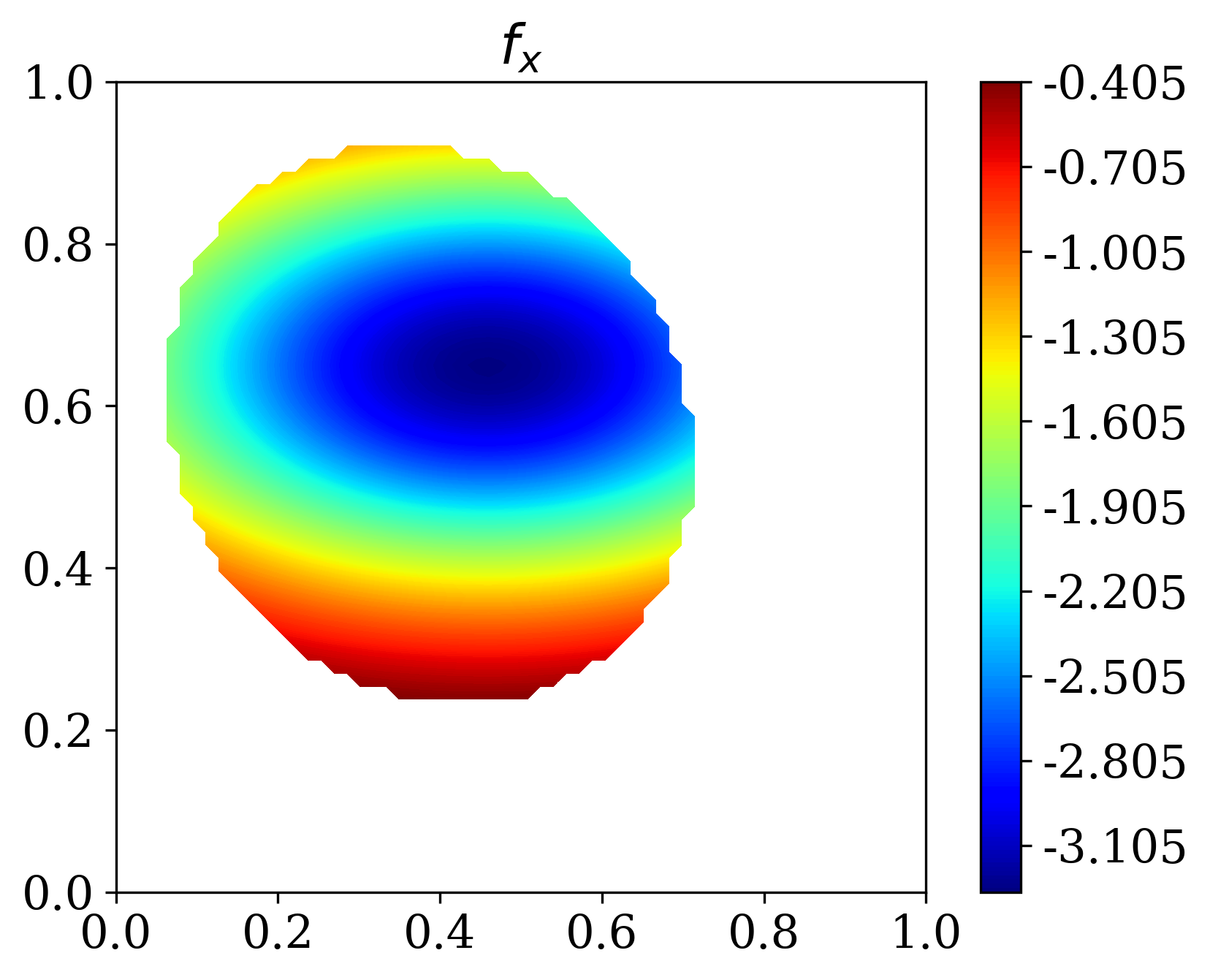}
		\hfill
		\includegraphics[width=0.18\textwidth]{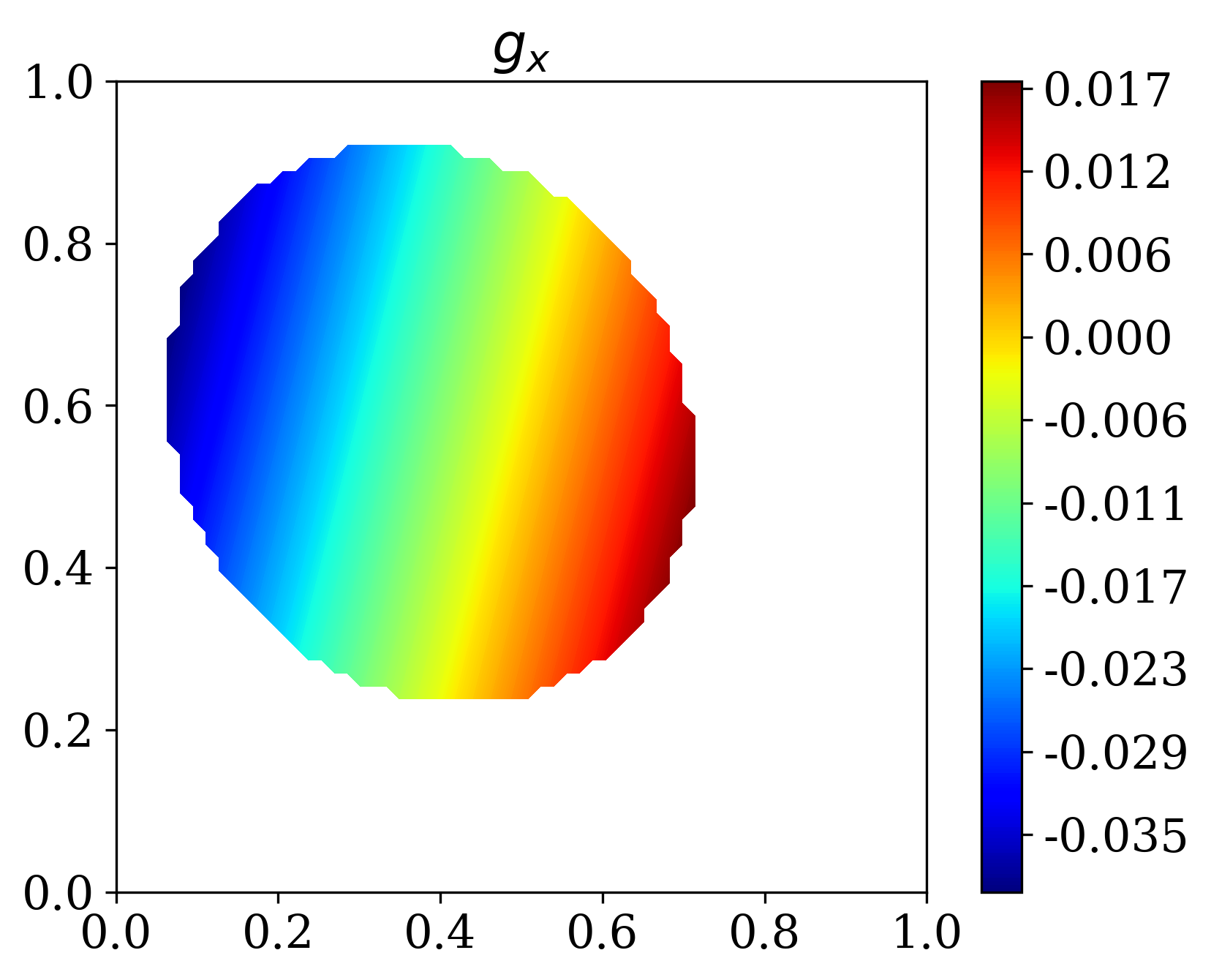}
		\hfill
		\includegraphics[width=0.18\textwidth]{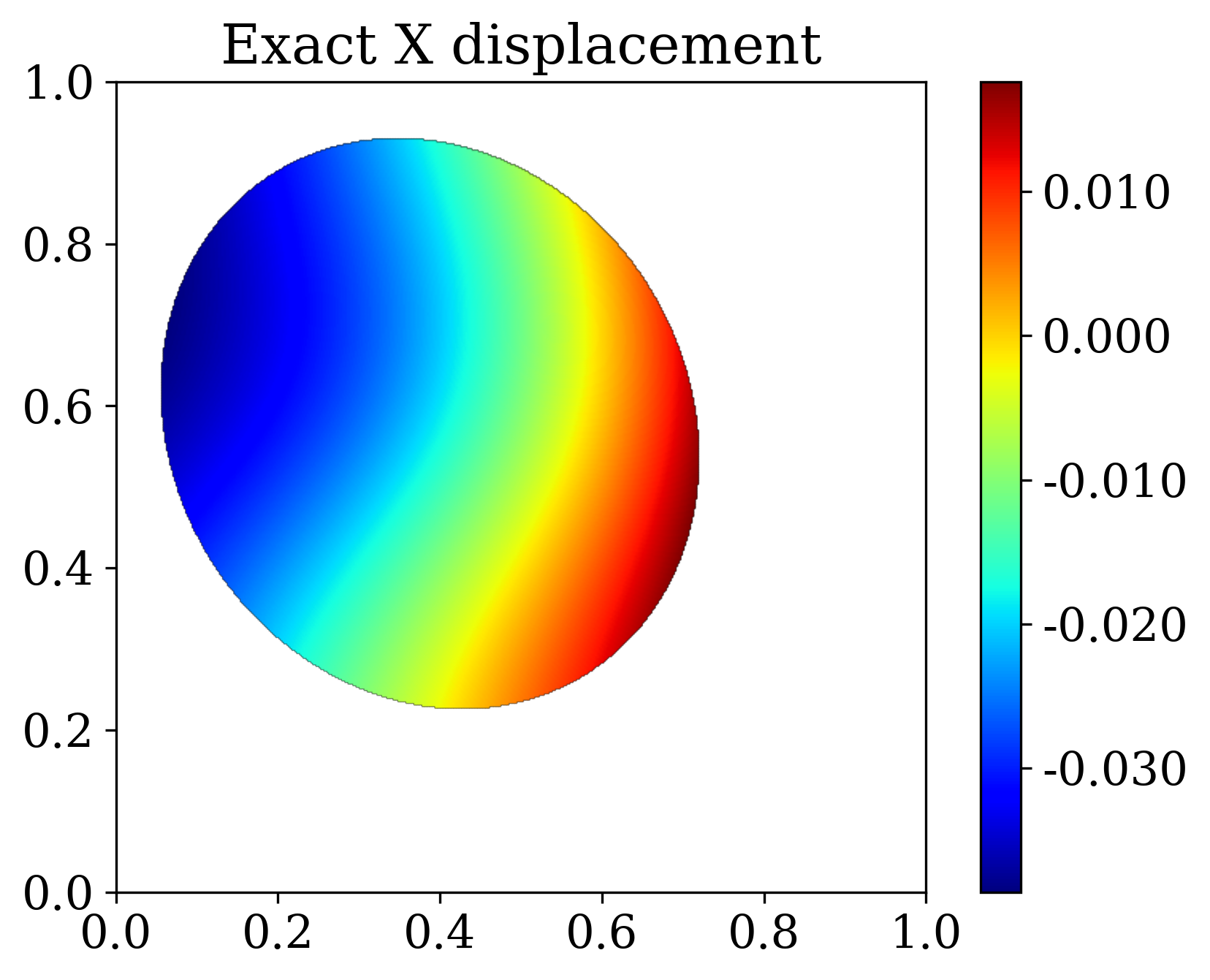}
		\hfill
		\includegraphics[width=0.18\textwidth]{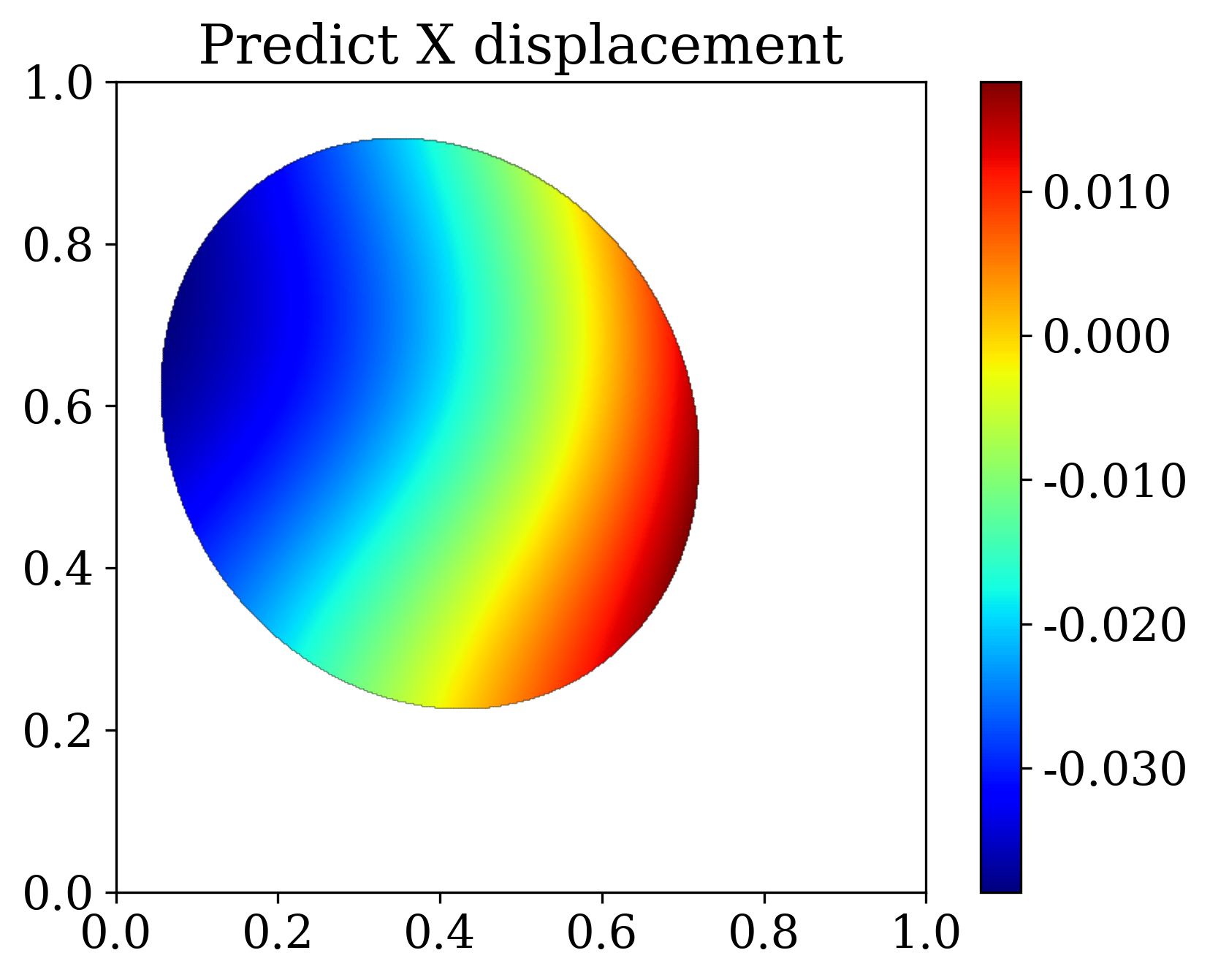}
		\hfill
		\includegraphics[width=0.19\textwidth]{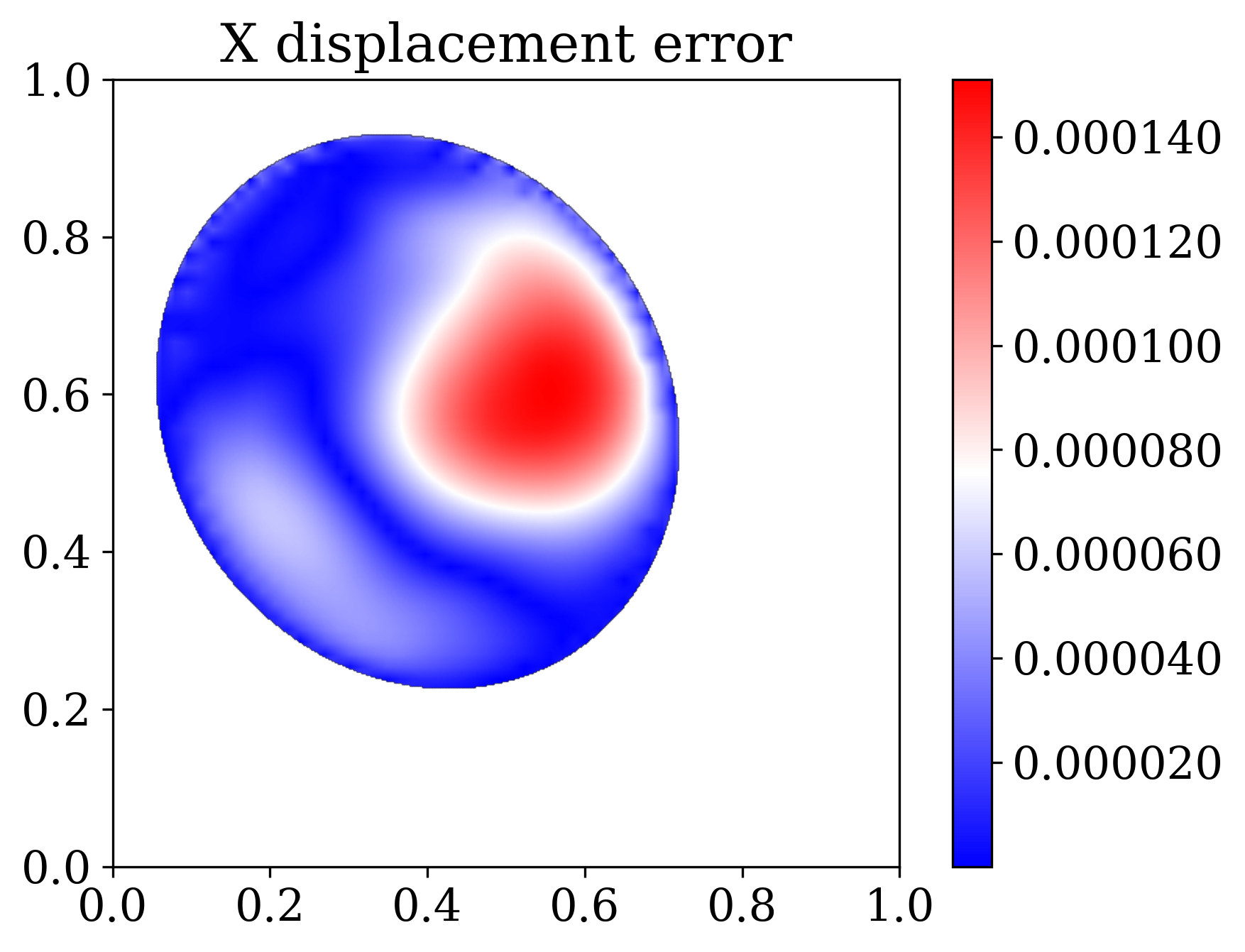}
		\vspace{0.5em}
		\includegraphics[width=0.18\textwidth]{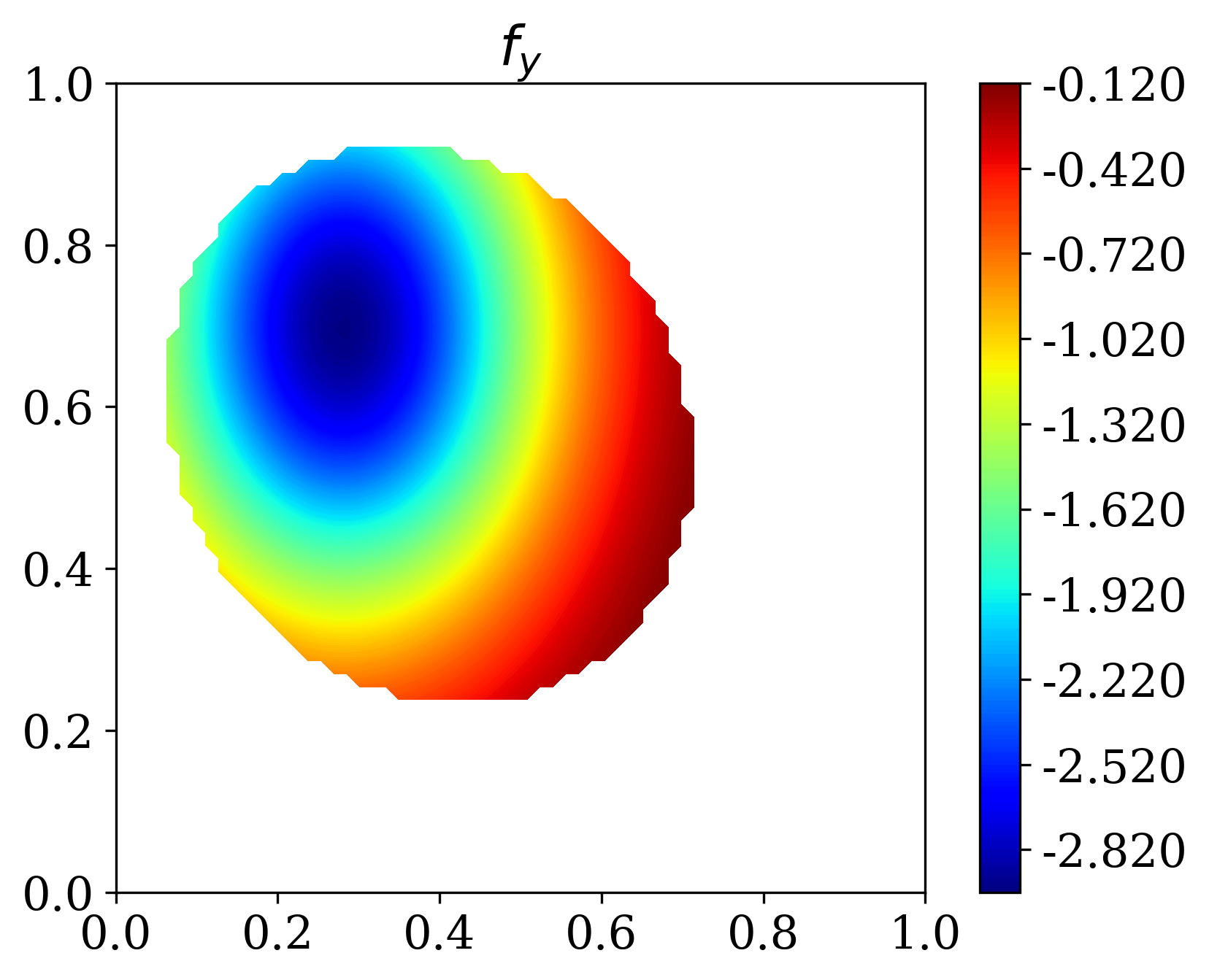}
		\hfill
		\includegraphics[width=0.18\textwidth]{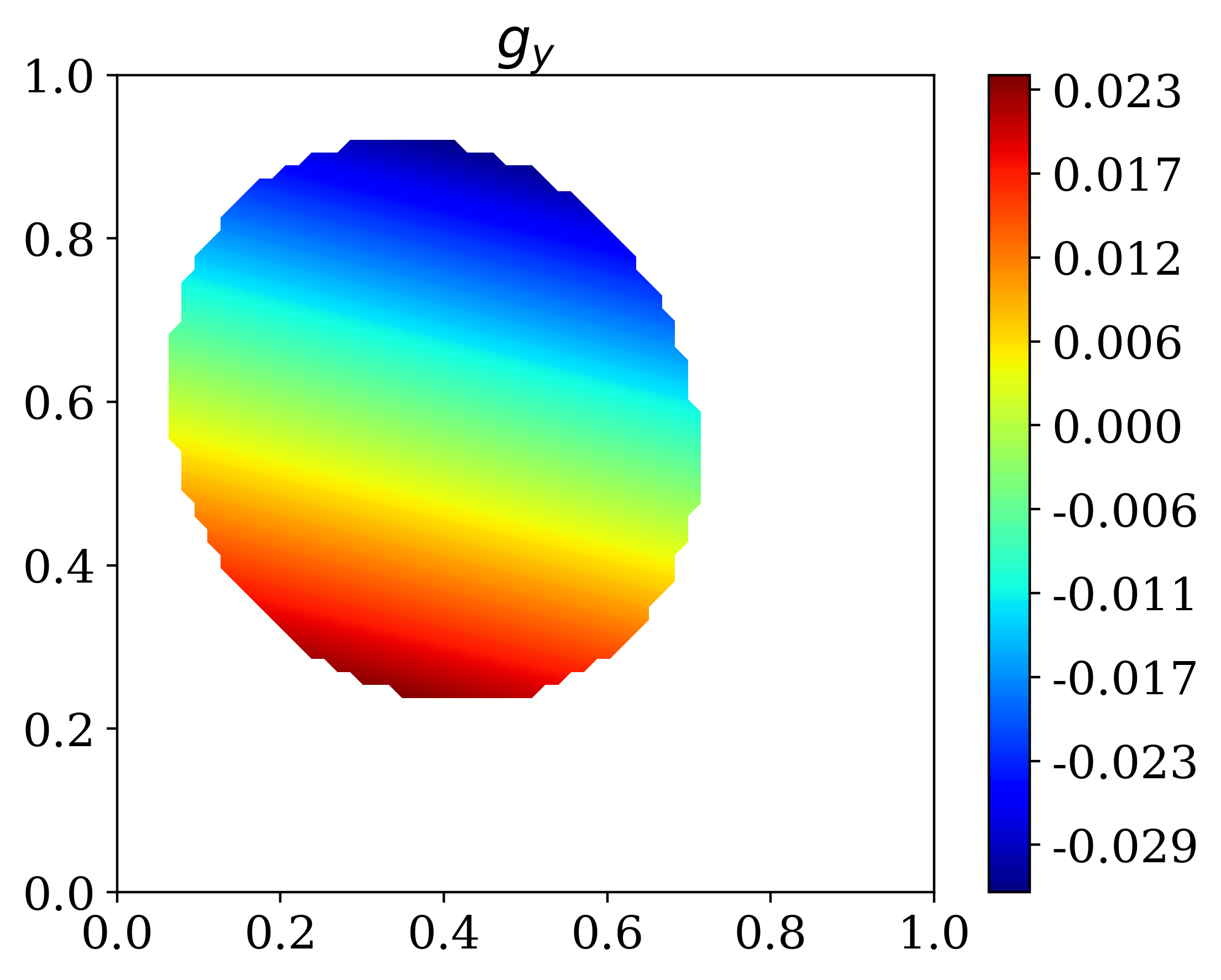}
		\hfill
		\includegraphics[width=0.18\textwidth]{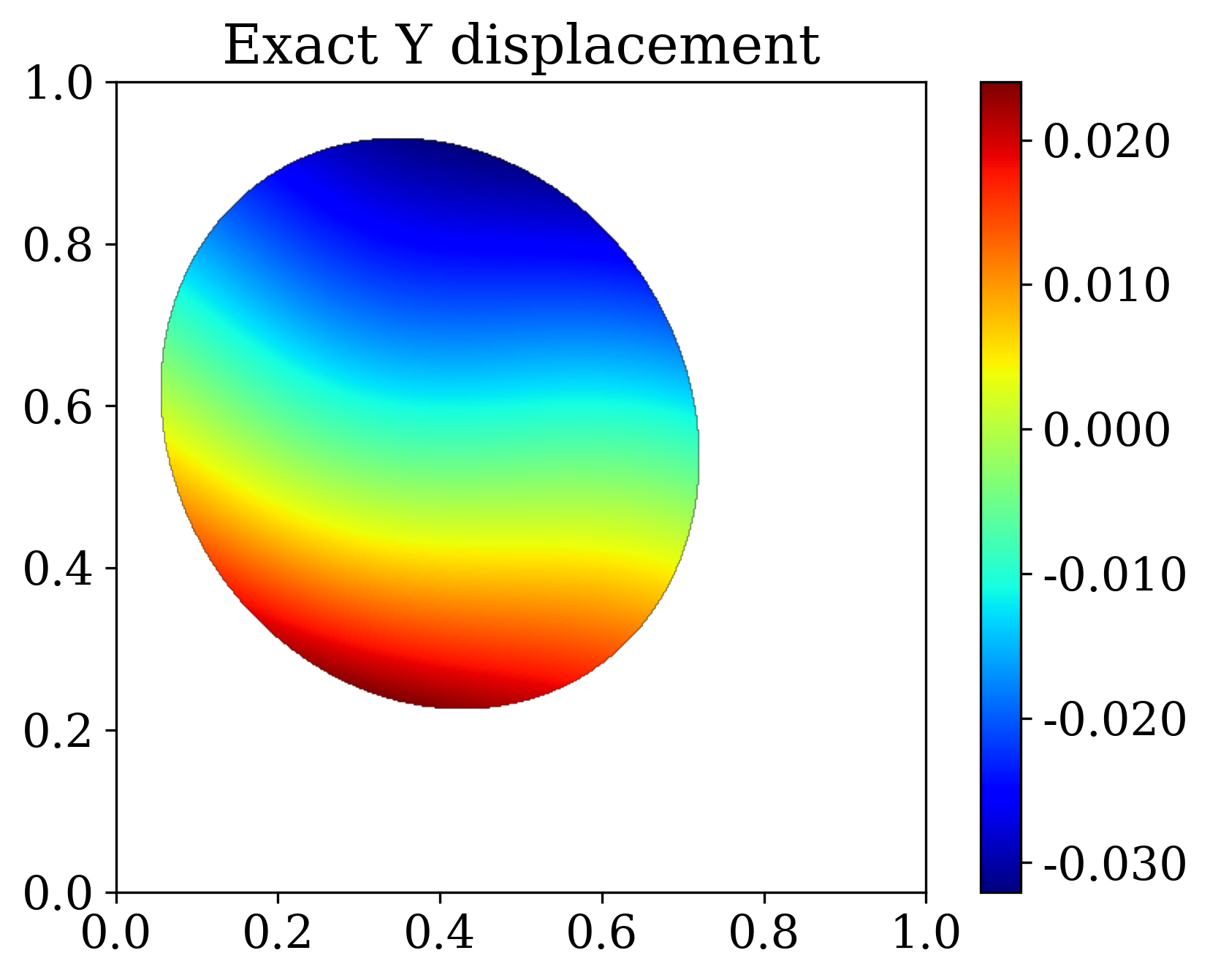}
		\hfill
		\includegraphics[width=0.18\textwidth]{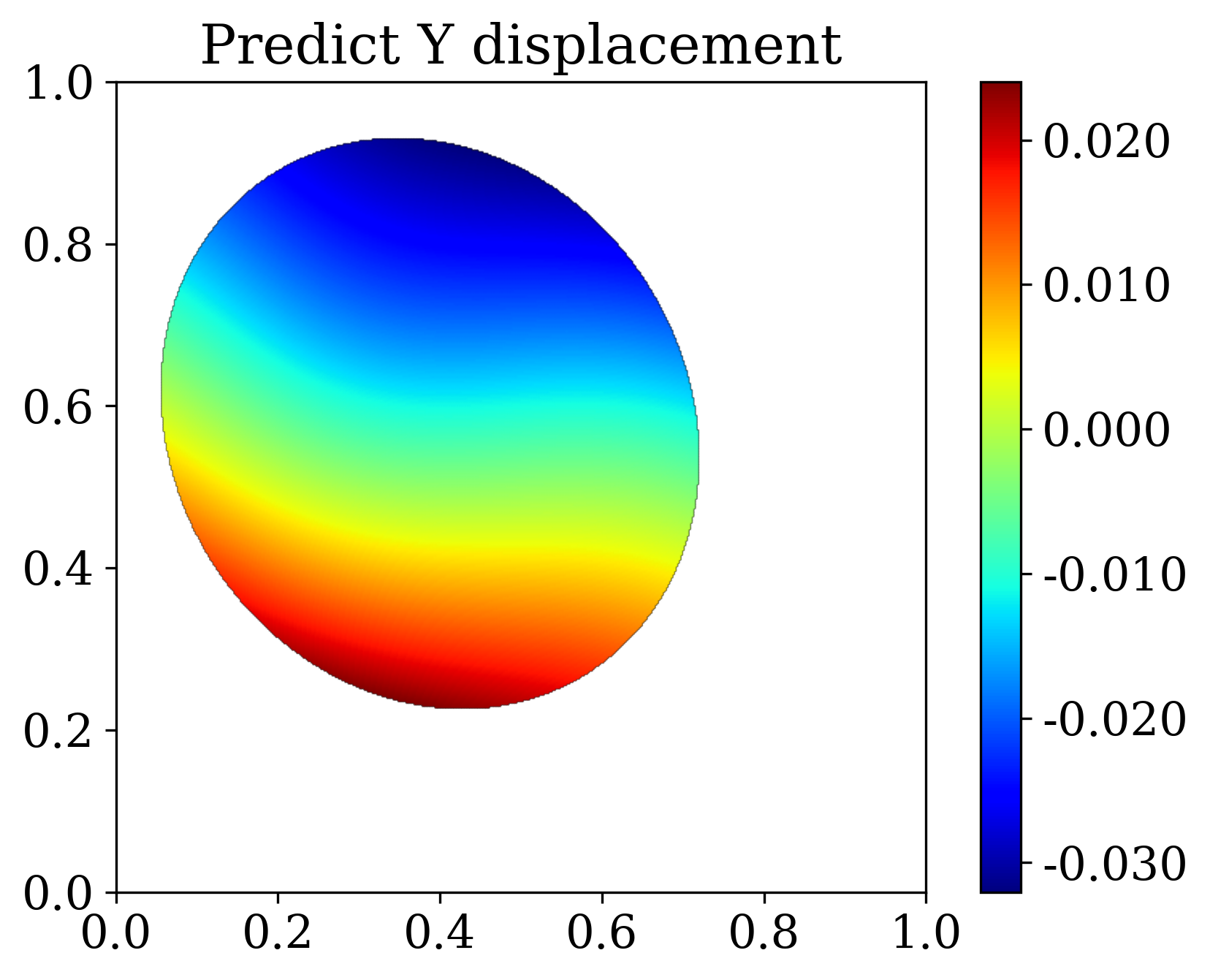}
		\hfill
		\includegraphics[width=0.19\textwidth]{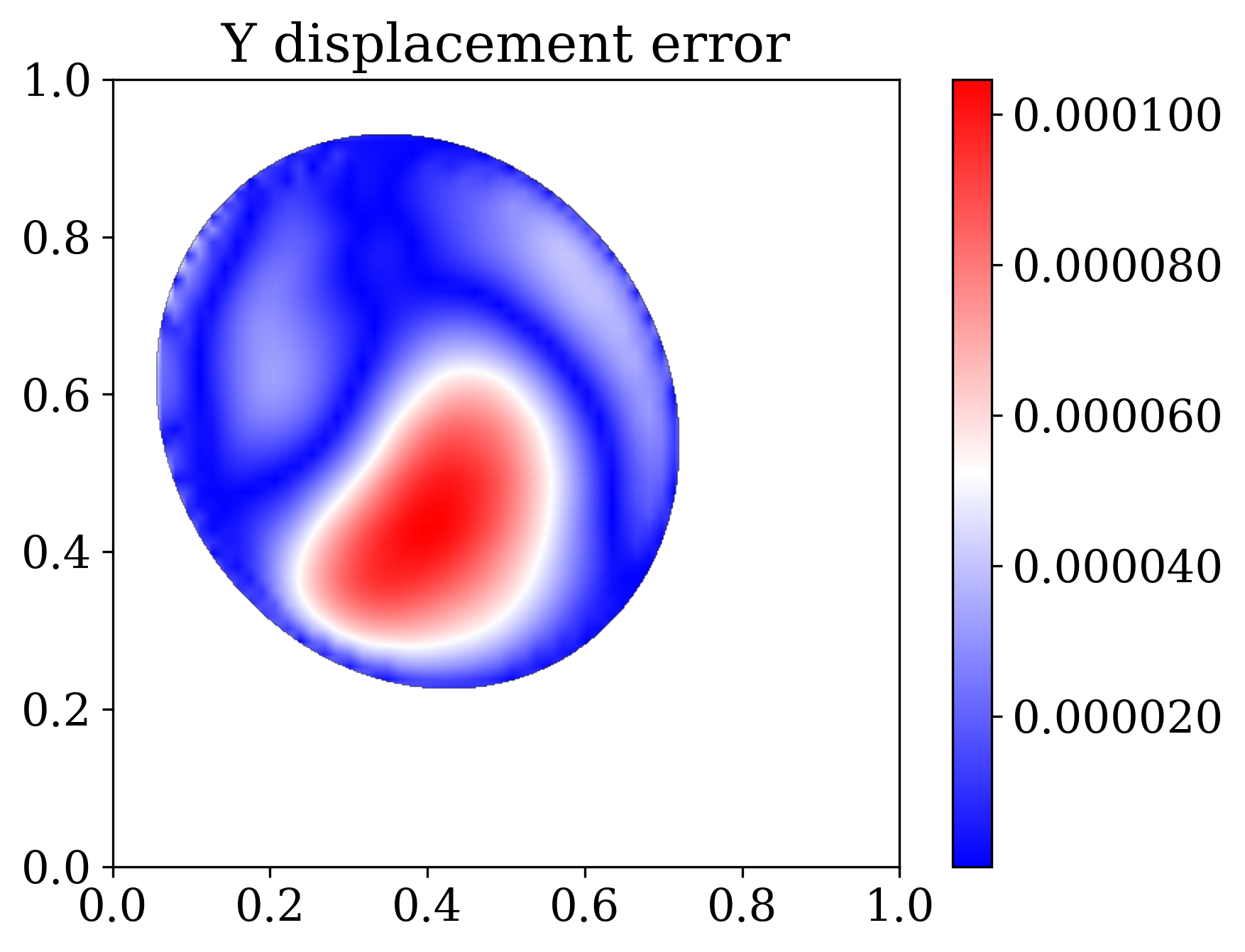}

	\end{minipage}

	\caption{Elliptical shape benchmark. (a) Evolution of the mean relative $L^2$ errors during training; (b) Relative errors after training; (c) Representative test sample: the top and bottom rows show the $x$- and $y$-components; the columns from left to right display $\mathbf{f}_h$, $\mathbf{g}_h$, exact displacement, WINO prediction, and absolute error contour, respectively.}
	\label{fig:case_ellip_error}
\end{figure}

We formulate the hyperelastic problem with pure Dirichlet boundary conditions on the elliptic domain. The boundary value problem is given by \eqref{eq:case_dirichlet}, where $\Omega$ is the physical domain represented by the level-set function $\varphi$, and $\Gamma=\partial\Omega$. For an elliptical geometry, we describe $\Omega$ as the interior of an ellipse via
\begin{equation}
	\varphi(x,y)=\left(\frac{x-c_x}{l_x}\right)^2+\left(\frac{y-c_y}{l_y}\right)^2-1,
	\label{eq:elliptical_phi}
\end{equation}
with center $(c_x,c_y)$ and semi-axis lengths $l_x,l_y$ sampled for each realization on the background domain $\mathcal{O}$, where $c_x,c_y\sim\mathcal{U}([0.3,0.7])$ and $l_x,l_y\sim\mathcal{U}([0.2,0.4])$. The body force and Dirichlet boundary data are given by
\begin{equation}
	\mathbf{f}(x,y)=
	\begin{bmatrix}
		A_1 \exp \left( -\frac{(x-x_{f1})^2}{2\sigma_{f1}} -\frac{(y-y_{f1})^2}{2\gamma_{f1}} \right) \\
		A_2 \exp \left( -\frac{(x-x_{f2})^2}{2\sigma_{f2}} -\frac{(y-y_{f2})^2}{2\gamma_{f2}} \right)
	\end{bmatrix},
	\label{eq:case_ellip_f}
\end{equation}
and
\begin{equation}
	\mathbf{g}(x,y)=
	\begin{bmatrix}
		0.1\left(\alpha x+0.25\,\beta y\right)\\
		0.1\left(\beta y-0.25\,\alpha x\right)
	\end{bmatrix},
	\label{eq:case_ellip_g}
\end{equation}
where $A_1,A_2 \sim\mathcal{U}([-4,-2]\cup[2,4])$, $x_{f1},y_{f1}, x_{f2},y_{f2}\sim\mathcal{U}([0.2, 0.8])$, $\sigma_{f1},\sigma_{f2},\gamma_{f1},\gamma_{f2}\sim\mathcal{U}([0.15,0.45])$, and $\alpha,\beta\sim\mathcal{U}([-0.8, -0.2]\cup[0.2, 0.8])$. The Young's modulus and the Poisson ratio are $E=10$, $\nu=0.3$.

The dataset contains 600 samples of $(\varphi_h,\mathbf{f}_h,\mathbf{g}_h)$ represented on a $64\times64$ Cartesian grid, partitioned into 500 training and 100 test cases. For each input triple, a reference displacement field is obtained with $\varphi$-FEM on the identical mesh and used for error evaluation. In the pure Dirichlet regime, the operator learns
\begin{equation}
\mathcal{G}_{\theta}:(\varphi_h,\mathbf{f}_h,\mathbf{g}_h)\mapsto\mathbf{w}_h,
\end{equation}
and the displacement is then calculated by $\mathbf{u}_h=\varphi_h\mathbf{w}_h+\mathbf{g}_h$. With this discretization, WINO is optimized for 500 epochs using SOAP under the loss function \eqref{eq:total_loss_Dirichlet}.

\begin{table}[t]
	\centering
	\caption{Computational costs (data generation + training time) and relative errors of three methods for the elliptical shape case. Cost ratio is total time relative to WINO.}
	\label{tab:performance_ellip}
	\footnotesize
	\begin{tabularx}{\linewidth}{@{}l >{\raggedright\arraybackslash}X >{\raggedright\arraybackslash}X >{\raggedright\arraybackslash}X >{\raggedright\arraybackslash}X >{\raggedright\arraybackslash}X@{}}
		\toprule
		Method & Total time & Cost ratio & $\|e\|_{L^2}$ & $\|e\|_{H^1}$ & $\|e\|_{E}$ \\
		\midrule
		WINO & $0.3 + 1248.8$s & $1.00\times$ & $0.79 \pm 0.43\%$ & $3.53 \pm 1.22\%$ & $2.38 \pm 0.96\%$ \\
		$\varphi$-FEM-FNO & $1392.5 + 644.4$s & $1.63\times$ & $0.86 \pm 0.53\%$ & $2.62 \pm 0.89\%$ & $2.13 \pm 0.87\%$ \\
		WINO+data & $1392.5 + 1258.7$s & $2.12\times$ & $0.82 \pm 0.54\%$ & $3.33 \pm 1.17\%$ & $2.31 \pm 0.97\%$ \\
		\bottomrule
	\end{tabularx}
\end{table}

Fig.~\ref{fig:case_ellip_error}a shows the training histories of the mean relative $L^2$ error, and Fig.~\ref{fig:case_ellip_error}b summarizes the final error distributions on the test set. In each box plot, the center line denotes the median and the box bounds indicate the interquartile range (IQR). These results confirm that WINO converges very fast in the first 100 epochs and attains accurate predictions. The representative example in Fig.~\ref{fig:case_ellip_error}c shows that the predicted displacement field is in close agreement with the reference solution, with absolute errors of $O(10^{-4})$. For a quantitative comparison, Table~\ref{tab:performance_ellip} reports both accuracy and computational cost for $\varphi$-FEM-FNO, WINO, and WINO with labeled data. While all three methods reach comparable error levels, $\varphi$-FEM-FNO's total wall-clock time is about $1.63\times$ that of WINO, mainly because supervised training requires generating large sets of $\varphi$-FEM reference solutions. This demonstrates the practical effectiveness of the proposed framework for parametric hyperelasticity problems with complex geometries. To assess training reproducibility, we additionally trained WINO with ten independent random seeds and evaluated each model on the same training and test datasets. The resulting mean test relative $L^2$ errors have a seed-averaged mean of $8.19\times10^{-3}$, with a $5$--$95\%$ percentile range of $[7.81\times10^{-3},\,8.67\times10^{-3}]$, indicating stable and repeatable performance across random initializations.

\begin{figure}[t]
	\centering
	\begin{minipage}[t]{0.58\textwidth}
		\centering
		\makebox[\linewidth][l]{\textbf{(a)}}\par\vspace{0.3em}
		\includegraphics[width=\linewidth]{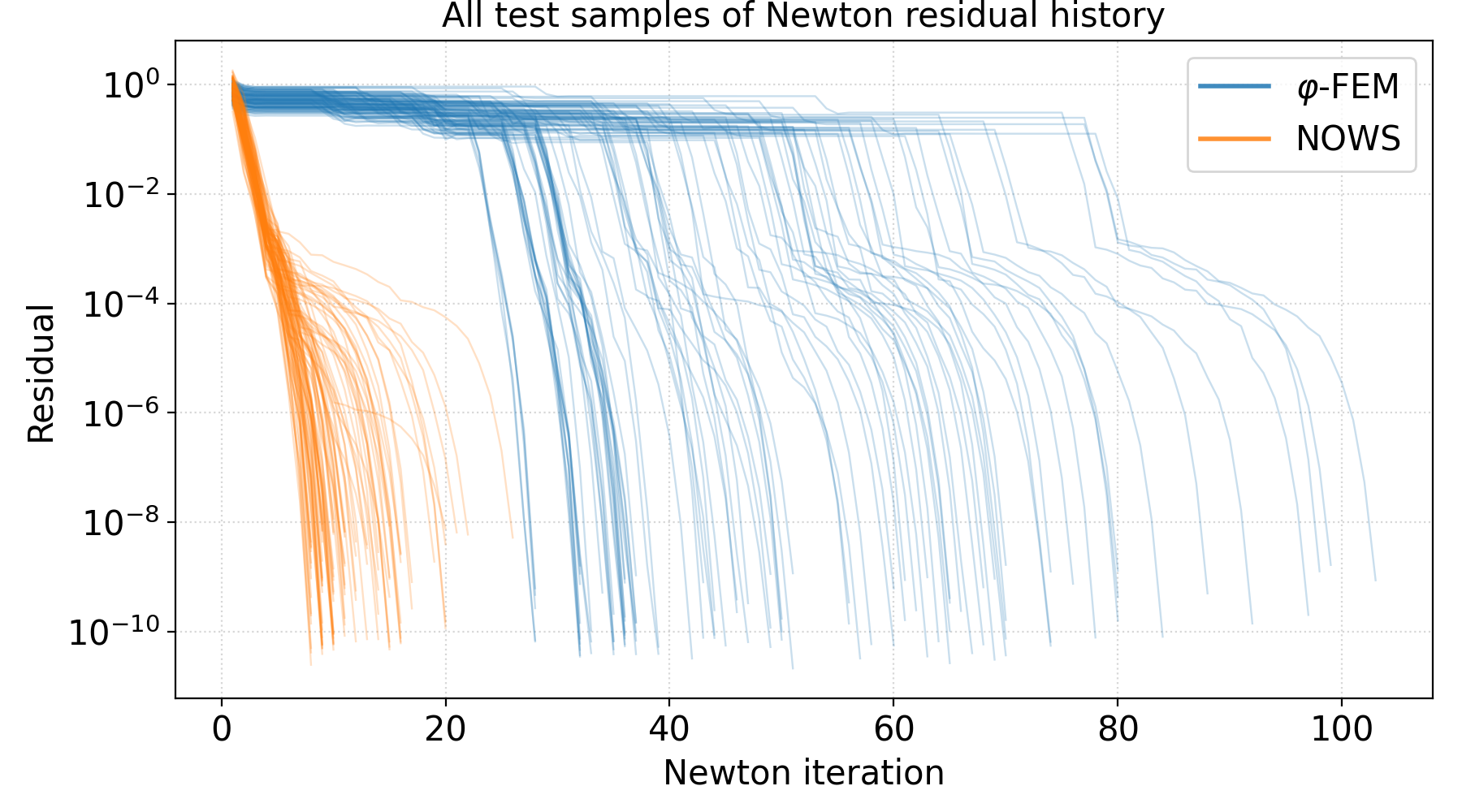}
	\end{minipage}
	\hfill
	\begin{minipage}[t]{0.39\textwidth}
		\centering
		\makebox[\linewidth][l]{\textbf{(b)}}\par\vspace{0.3em}
		\includegraphics[width=\linewidth]{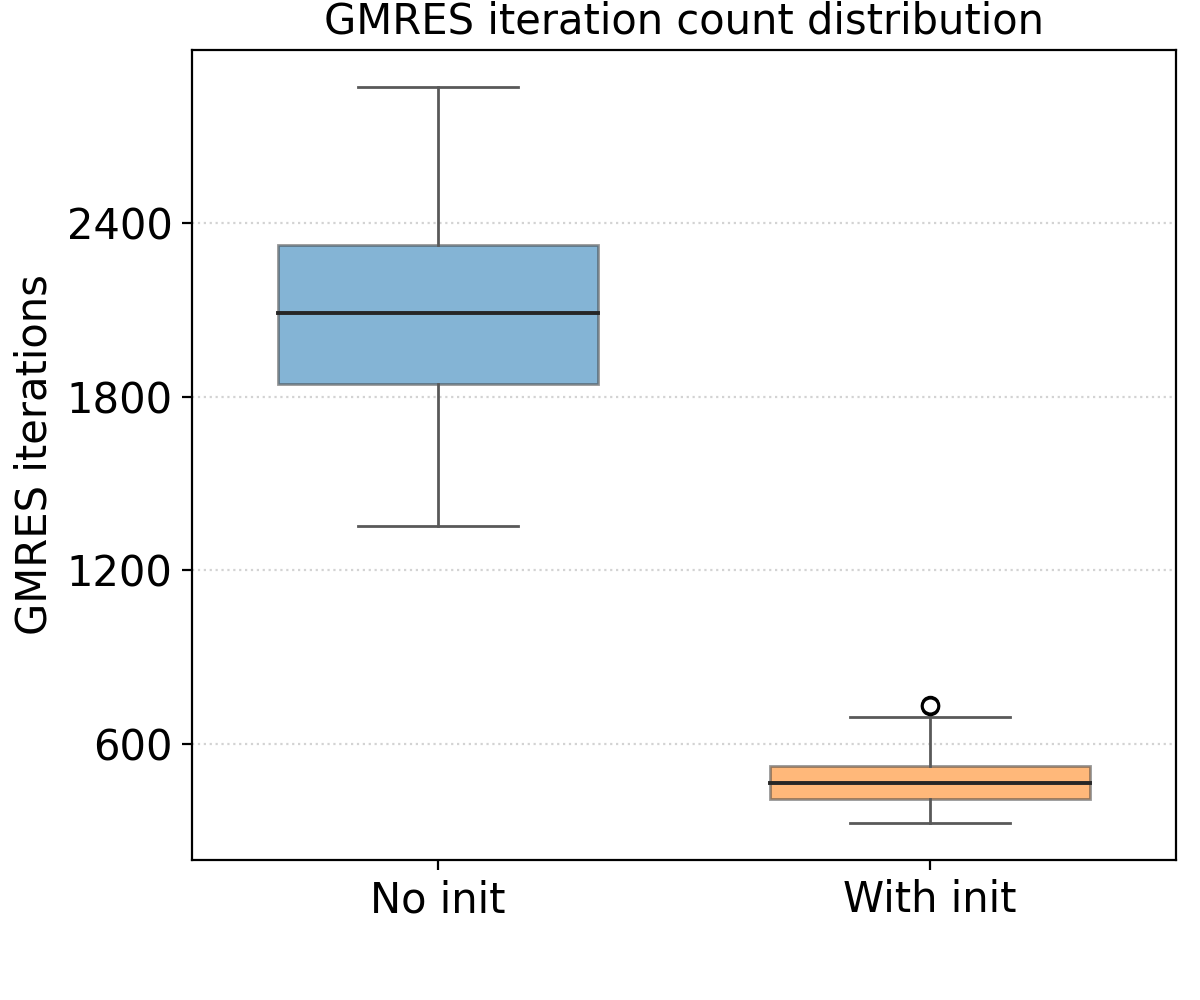}
	\end{minipage}
	\caption{NOWS comparison for the elliptical shape case. (a) Outer Newton residual and total Newton iterations for classical $\varphi$-FEM (4 load increments) and NOWS. (b) Box plots of total GMRES iteration counts for all load increments.}
	\label{fig:case_ellip_nows}
\end{figure}

In addition, WINO predictions can be used as an initial guess for the $\varphi$-FEM nonlinear solve, thereby accelerating convergence; this follows the neural-operator warm-start (NOWS) strategy introduced in Section~\ref{sec:nows}. For the hyperelastic Dirichlet formulation \eqref{eq:phifem_dirichlet_weak}, the discrete problem is a nonlinear residual equation. Newton's method linearizes it at each outer step and thus generates a sequence of sparse linear systems. For strongly nonlinear instances such as the present benchmark, plain Newton iteration from a cold start may converge slowly or even fail. A standard remedy is incremental load stepping, which solves a sequence of problems with scaled body forces $\mathbf{f}_h^{(k)}=\alpha_k\mathbf{f}_h$, where $0<\alpha_0<\cdots<\alpha_K=1$, and uses the solution at $\alpha_k$ as the initial guess for the problem at $\alpha_{k+1}$. WINO offers a simpler alternative, since its output can be used directly as the initial Newton iterate without continuation. We compare classical $\varphi$-FEM with four loading steps against NOWS without load stepping. For each Newton iteration, we use restarted GMRES as the inner solver with restart $m=120$ and relative tolerance $10^{-5}$. Fig.~\ref{fig:case_ellip_nows} reports the outer Newton residual histories, the total number of Newton iterations across all incremental loading steps, and box plots of the total GMRES iteration counts over all Newton steps for runs with and without neural initialization. The residual is evaluated with respect to the final load $\mathbf{f}_h$, so for $\varphi$-FEM with load stepping the curve exhibits a staircase pattern, where each plateau corresponds to convergence at an intermediate load level $\mathbf{f}_h^{(k)}$. NOWS substantially reduces the Newton iteration count and reaches fast convergence in a single solve pass. Averaged over the test instances, the wall-clock time for the nonlinear $\varphi$-FEM solve without warm start is $2.44\pm 1.09\,\mathrm{s}$, compared with $0.47\pm 0.16\,\mathrm{s}$ when WINO supplies the initial guess, which amounts to a substantial speedup. Moreover, the total GMRES iterations are significantly reduced with NOWS, to approximately one-third of the baseline level, which directly lowers the overall computational cost. In our setup, a WINO forward inference takes only about $5\,\mathrm{ms}$ to provide a high-quality initial guess, eliminating the need for incremental loading and significantly improving the overall computational efficiency of $\varphi$-FEM.

\subsubsection{Random shape}\label{subsec:random_shape}

\begin{figure}[t]
	\raggedright
	\begin{minipage}[t]{\textwidth}
		\textbf{(a)}\par\vspace{0.4em}

		\hspace*{0.01\textwidth}
		\includegraphics[width=0.19\textwidth]{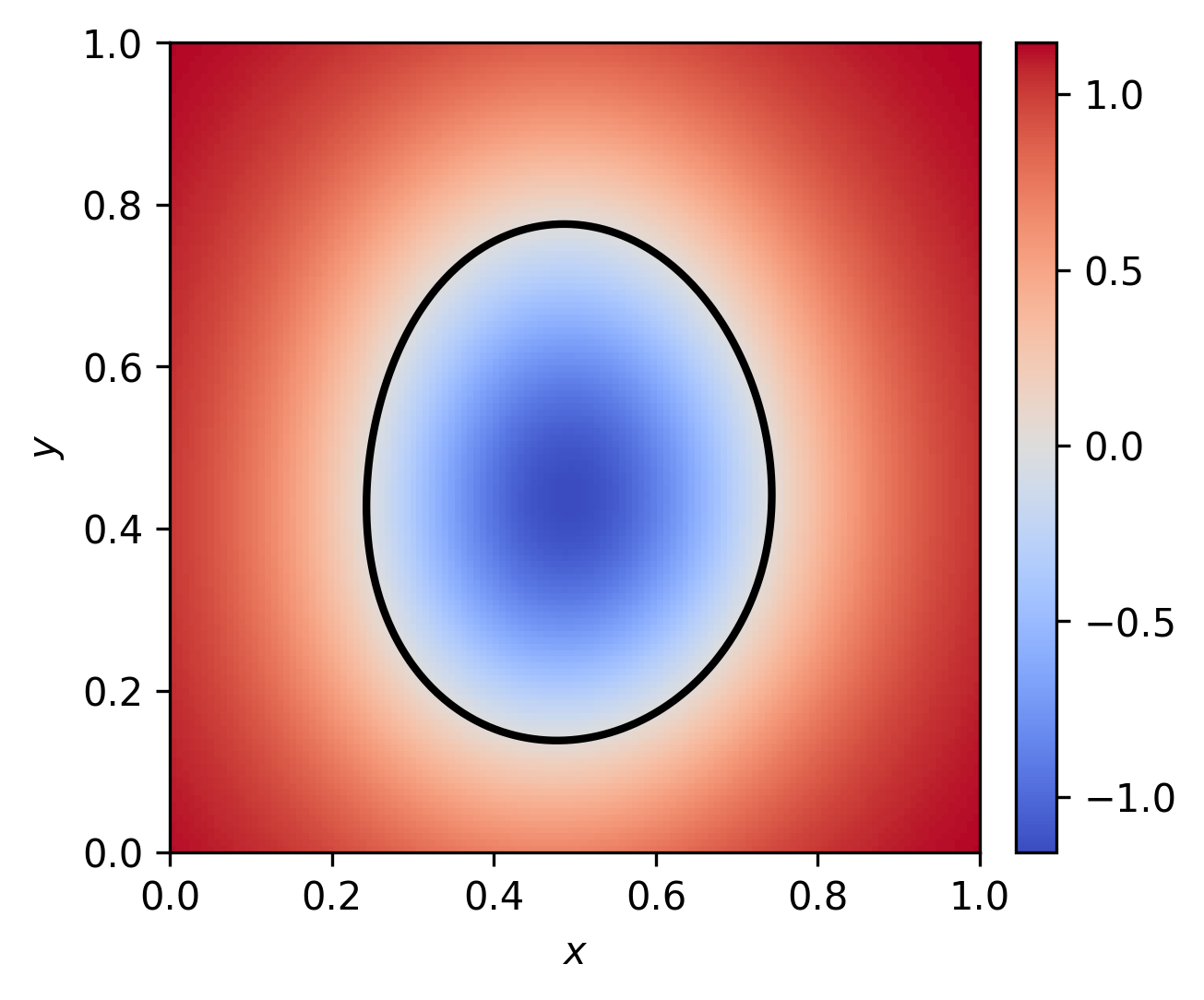}
		\hfill
		\includegraphics[width=0.19\textwidth]{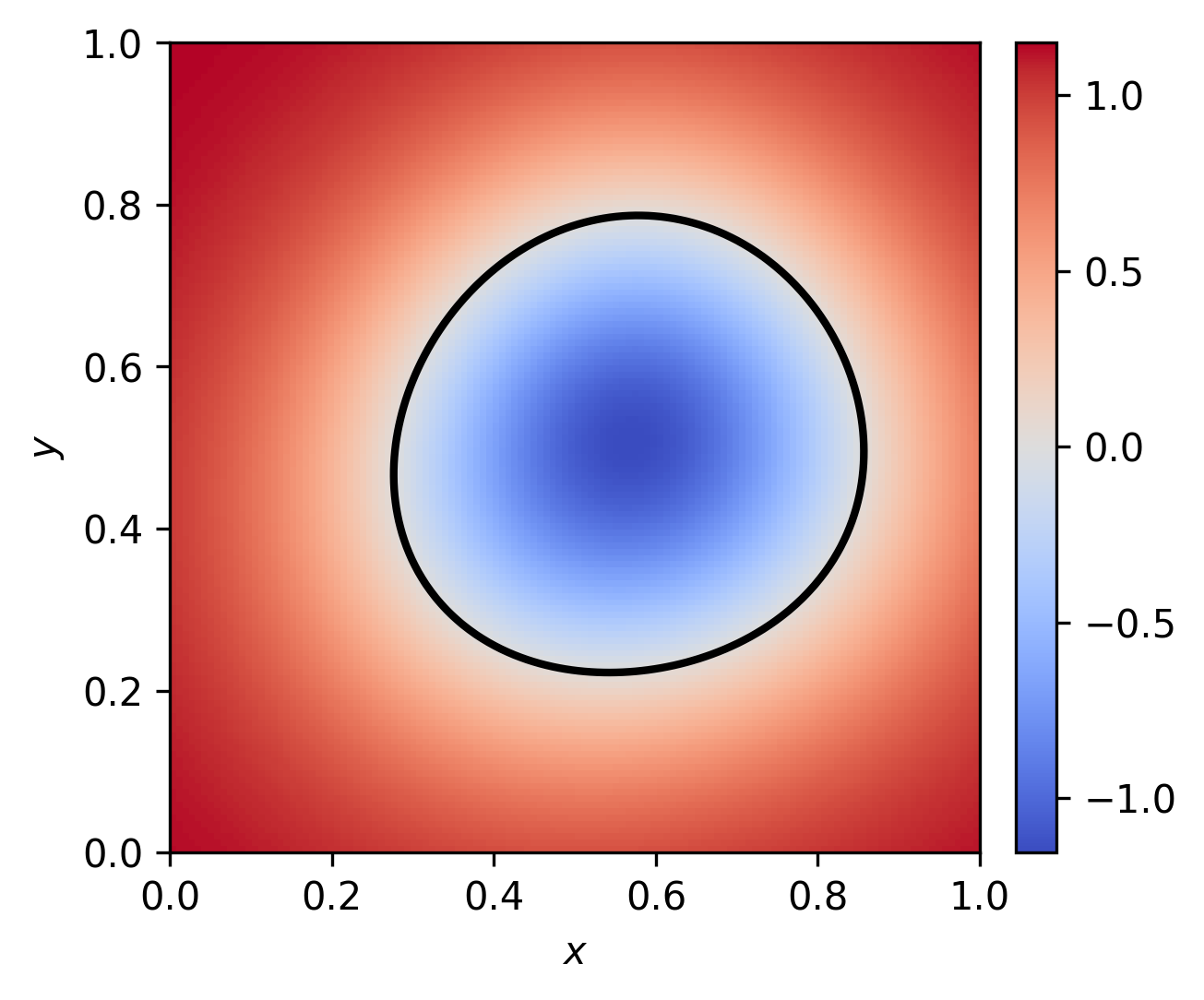}
		\hfill
		\includegraphics[width=0.19\textwidth]{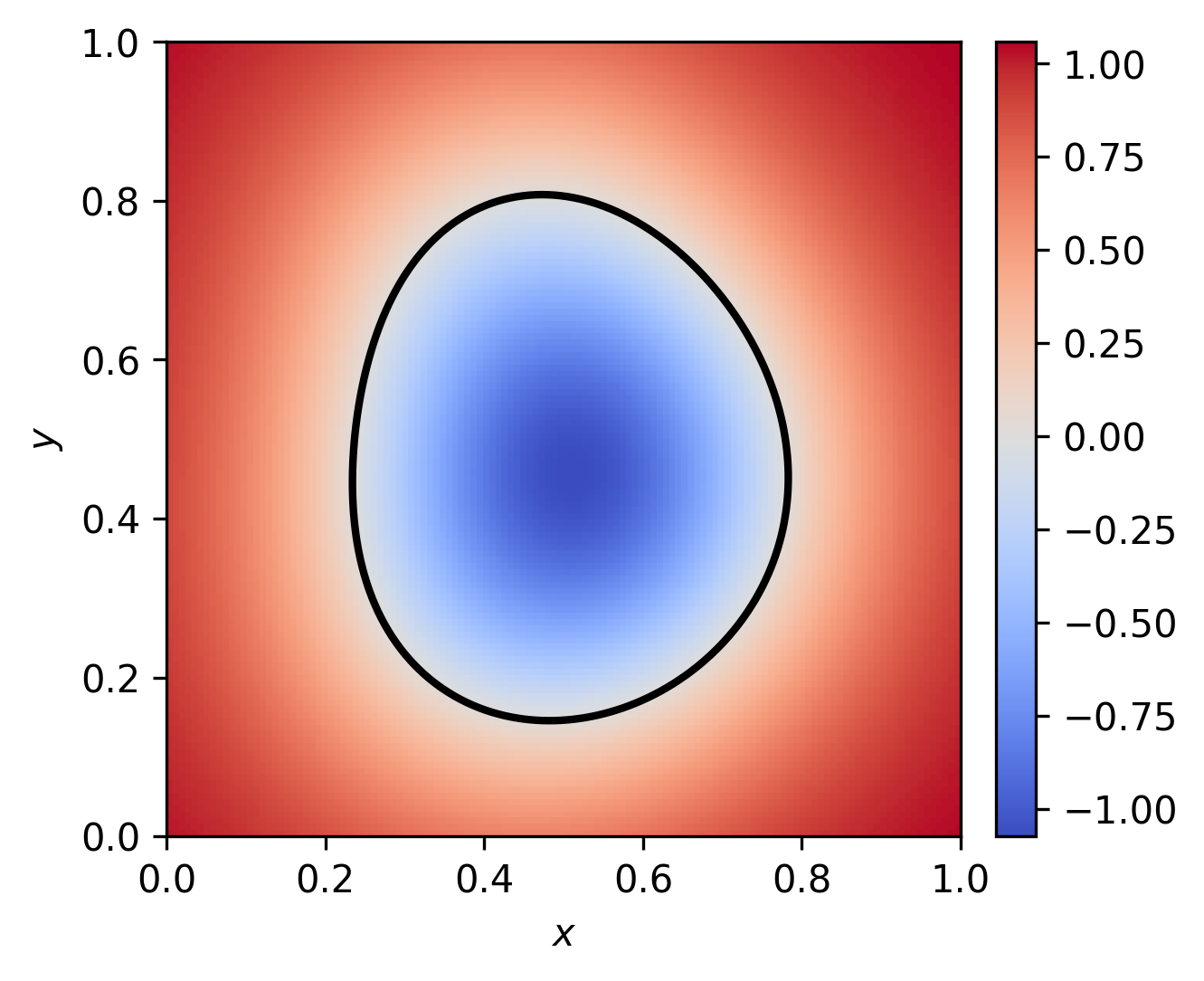}
		\hfill
		\includegraphics[width=0.19\textwidth]{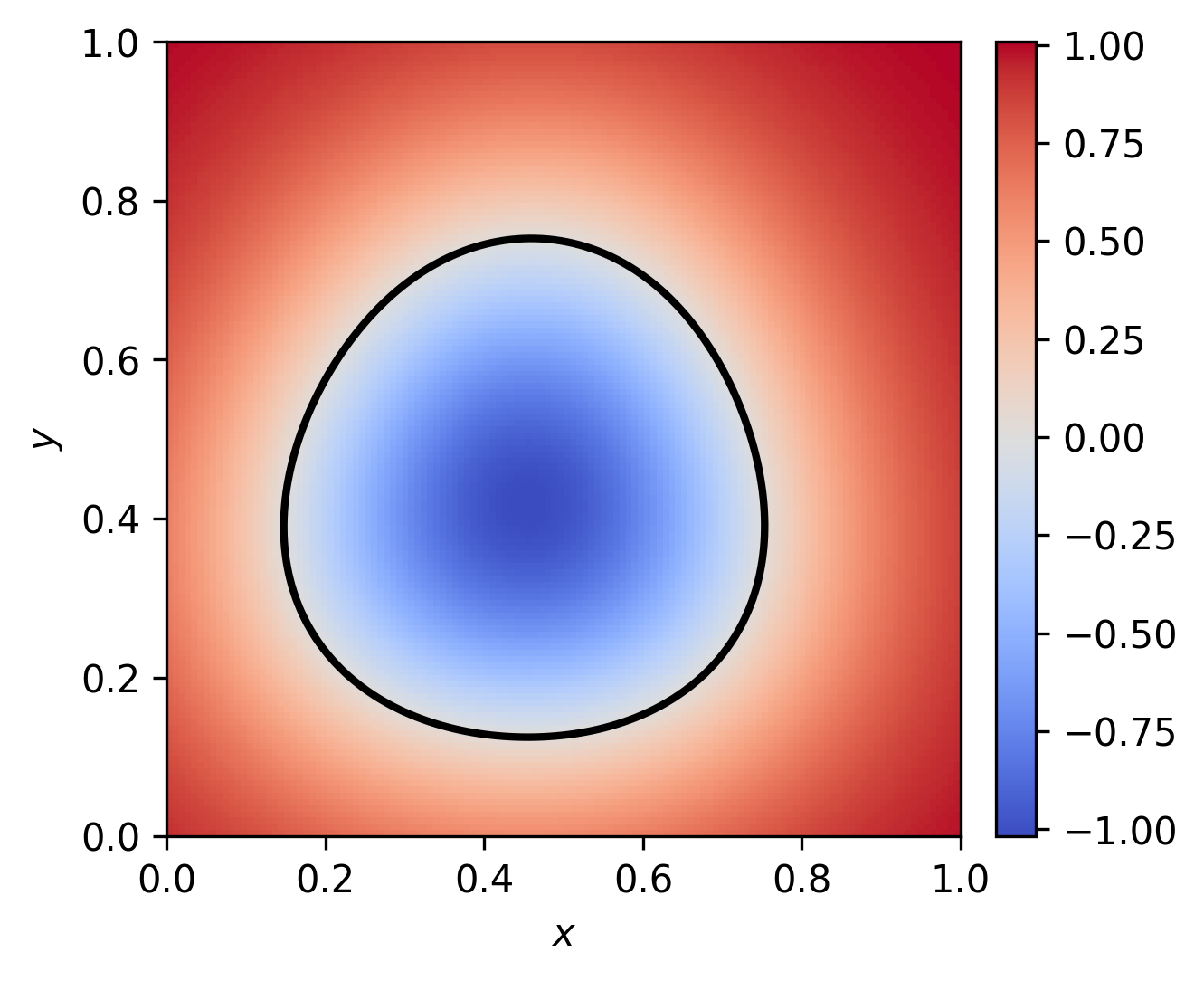}
		\hfill
		\includegraphics[width=0.19\textwidth]{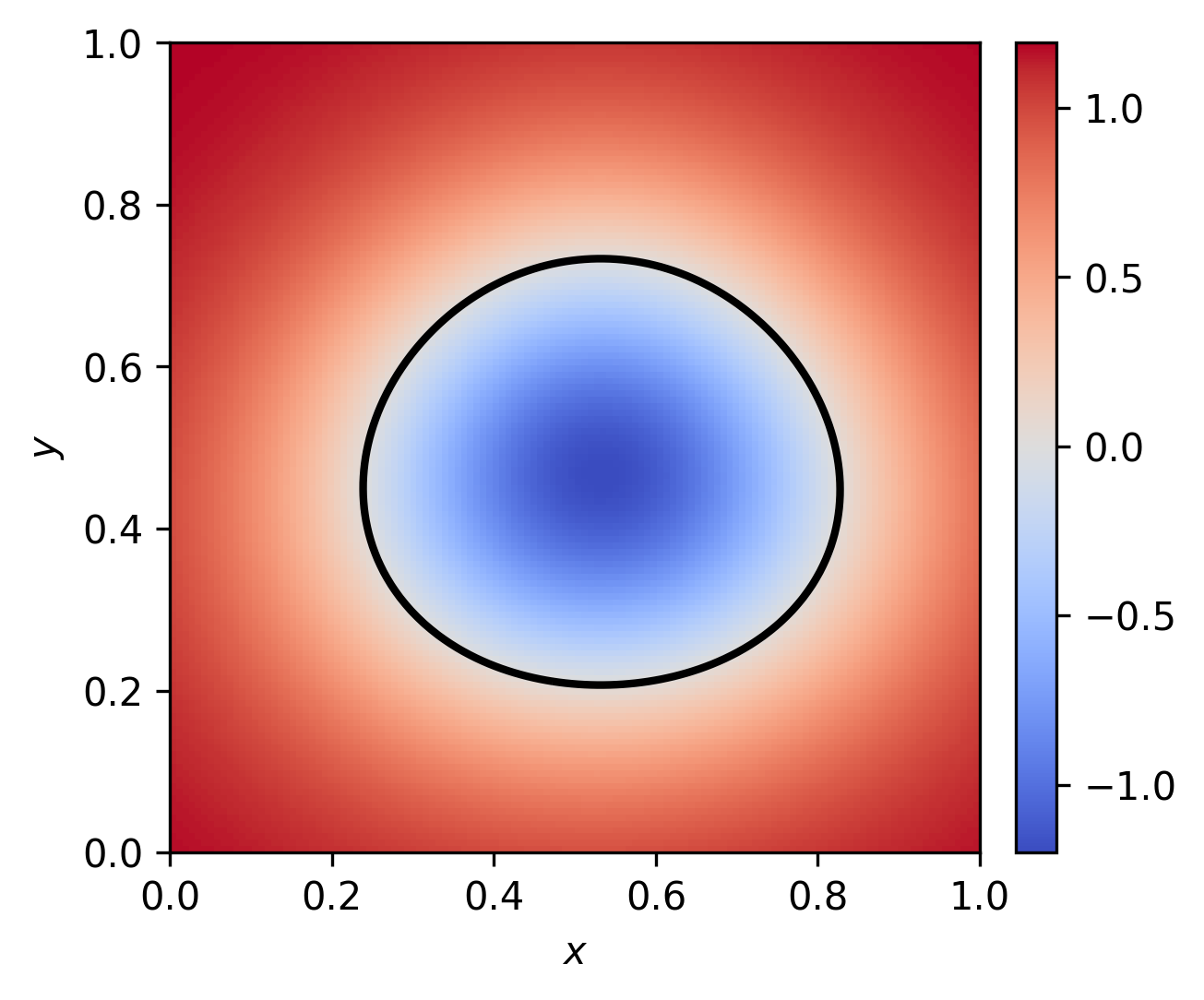}
	\end{minipage}

	\vspace{0.8em}
	\begin{minipage}[t]{0.31\textwidth}
		\textbf{(b)}\par\vspace{0.4em}
		\includegraphics[width=\linewidth]{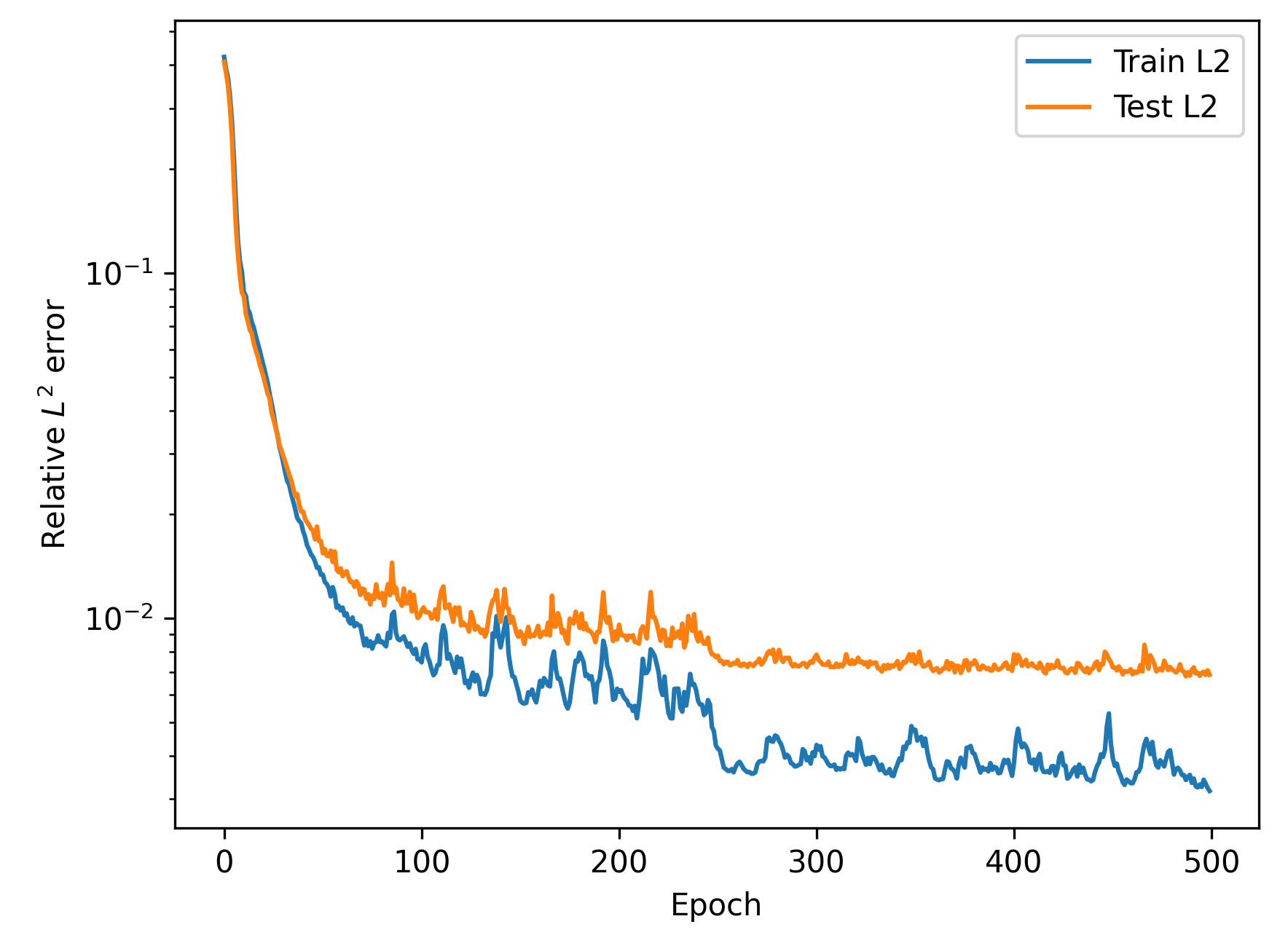}
	\end{minipage}
	\hspace{0.8em}
	\begin{minipage}[t]{0.3\textwidth}
		\textbf{(c)}\par\vspace{0.4em}
		\includegraphics[width=\linewidth]{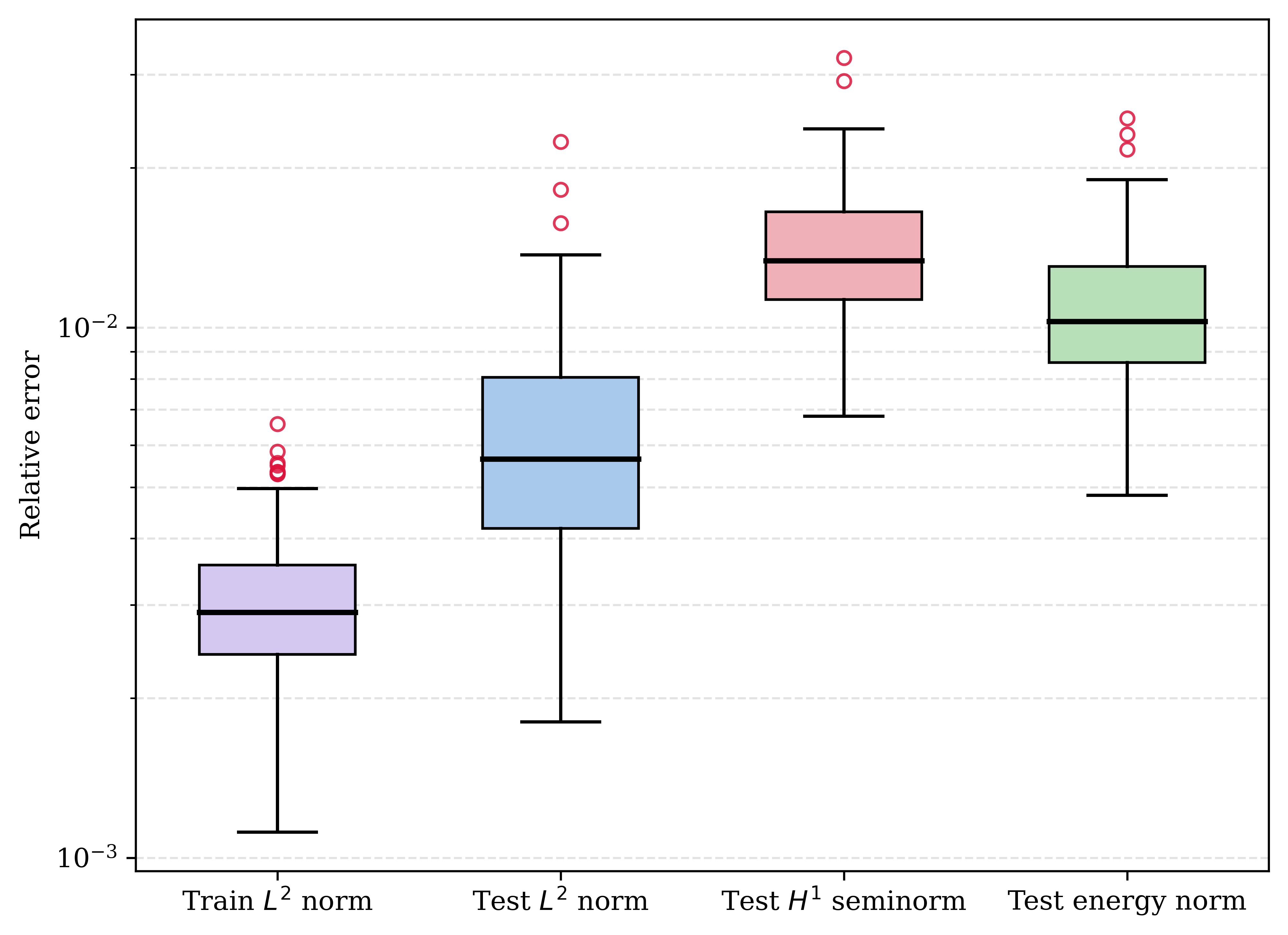}
	\end{minipage}
	\hspace{0.8em}
	\begin{minipage}[t]{0.31\textwidth}
		\textbf{(d)}\par\vspace{0.4em}
		\includegraphics[width=\linewidth]{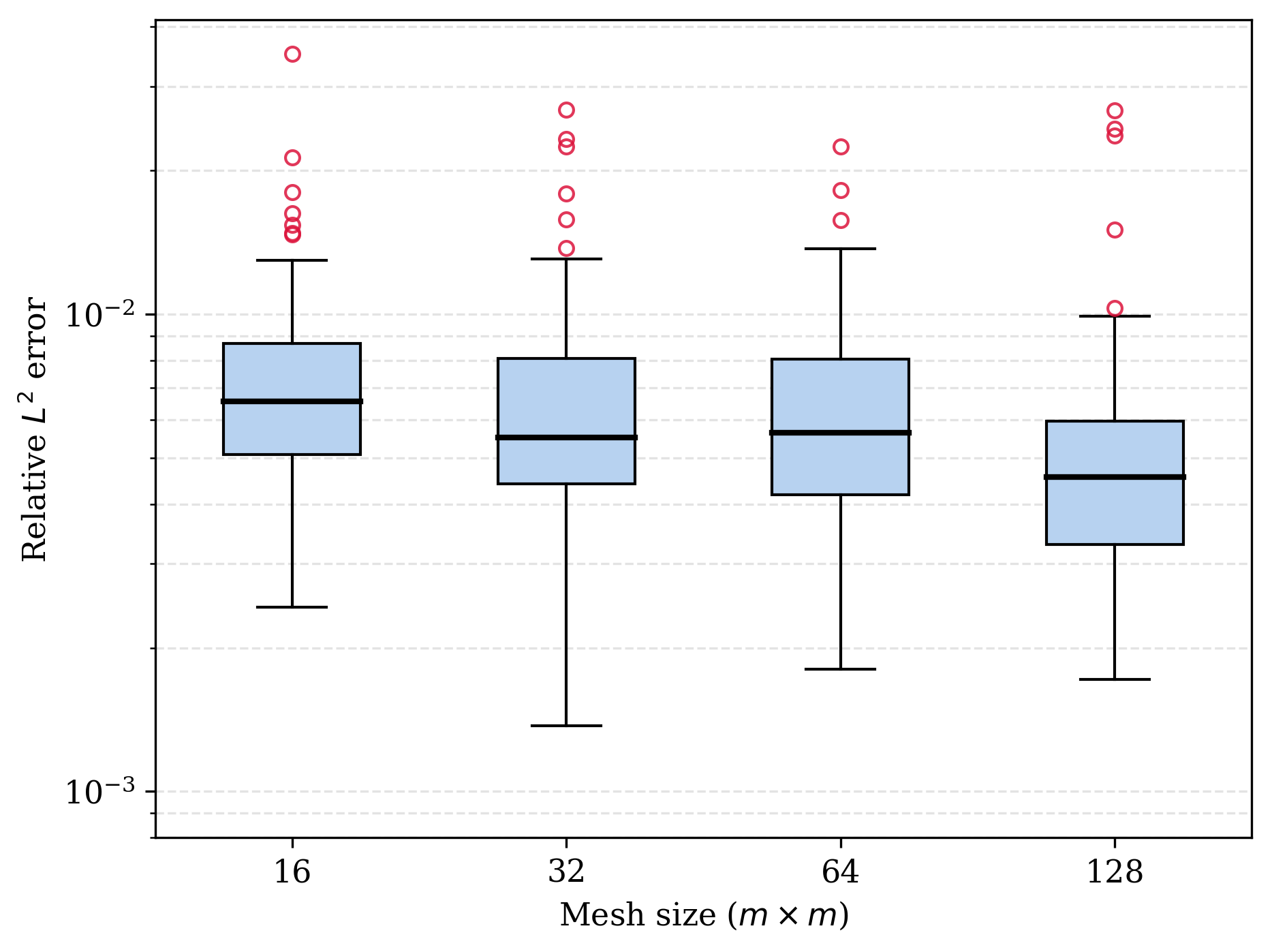}
	\end{minipage}

	\caption{Random shape benchmark. (a) Five test samples of the level-set functions and corresponding geometries; (b) Evolution of the mean relative $L^2$ errors during training; (c) Relative errors after training; (d) Box plots of the relative $L^2$ errors for WINO with different mesh sizes (from $16\times16$ to $128\times128$).}
	\label{fig:case_arbit_phi_error}
\end{figure}

In this part, we consider the hyperelastic problem \eqref{eq:case_dirichlet} on the random shape domain. The level-set function $\varphi$ is constructed from a sum of three Gaussian functions,
\begin{equation}
	\varphi(x, y)= -\varphi_g(x, y) + 0.5 \max_{(x,y)\in\Omega} \varphi_g(x, y),
\end{equation}
with
\begin{equation}
	\varphi_g(x, y) = \sum_{i=1}^{3} \exp \left( -\frac{(x-x_i)^2}{2\sigma_i} -\frac{(y-y_i)^2}{2\gamma_i} \right),
\end{equation}
where $(x_i, y_i)$ are the Gaussian centers and $(\sigma_i,\gamma_i)$ control the spread in the $x$- and $y$-directions. For $i=1,2,3$, we sample $\sigma_i,\gamma_i\sim\mathcal{U}([0.04,0.05])$, and independently
\[
\begin{aligned}
x_1&\sim\mathcal{U}([0.4,0.6]), & y_1&\sim\mathcal{U}([0.55,0.7]),\\
x_2&\sim\mathcal{U}([0.3,0.45]), & y_2&\sim\mathcal{U}([0.3,0.45]),\\
x_3&\sim\mathcal{U}([0.55,0.7]), & y_3&\sim\mathcal{U}([0.3,0.45]).\\
\end{aligned}
\]
We present five test samples of the level-set functions and their corresponding domains in Fig.~\ref{fig:case_arbit_phi_error}a. The definitions of $\mathbf{f}$ and $\mathbf{g}$ follow Section~\ref{subsec:elliptical_shape}; in this subsection, we set $A_1, A_2 \sim \mathcal{U}([-6,-3]\cup[3,6])$.

\begin{table}[t]
	\centering
	\footnotesize
	\caption{Computational costs (data generation + training time) and relative errors of three methods for the random shape case. Cost ratio is total time relative to WINO.}
	\label{tab:performance_arbit}
	\begin{tabularx}{\linewidth}{@{}l >{\raggedright\arraybackslash}X >{\raggedright\arraybackslash}X >{\raggedright\arraybackslash}X >{\raggedright\arraybackslash}X >{\raggedright\arraybackslash}X@{}}
		\toprule
		Method & Total time & Cost ratio & $\|e\|_{L^2}$ & $\|e\|_{H^1}$ & $\|e\|_{E}$ \\
		\midrule
		WINO & $0.3 + 1366.5$s & $1.00\times$ & $0.65 \pm 0.33\%$ & $1.41 \pm 0.42\%$ & $1.09 \pm 0.37\%$ \\
		$\varphi$-FEM-FNO & $8603.1 + 783.8$s & $6.87\times$ & $0.79 \pm 0.32\%$ & $2.05 \pm 0.58\%$ & $1.69 \pm 0.53\%$ \\
		WINO+data & $8603.1 + 1391.3$s & $7.32\times$ & $0.55 \pm 0.26\%$ & $1.23 \pm 0.38\%$ & $0.96 \pm 0.35\%$ \\
		\bottomrule
	\end{tabularx}
\end{table}

Using the sampling rules above, we randomly draw $600$ input realizations of $(\varphi_h,\mathbf{f}_h,\mathbf{g}_h)$, including $500$ training samples and $100$ test samples, represented on a $64\times64$ Cartesian grid. Reference solutions (ground truth) are computed by $\varphi$-FEM on a finer $128\times128$ grid. To evaluate errors consistently across discretizations, the coarse $64\times64$ fields are interpolated onto the $128\times128$ mesh, and all norms are evaluated on the active domain $\Omega_h$ implied by the fine-grid level-set $\varphi_h$. For this pure Dirichlet setting, we impose the lift $\mathbf{u}_h=\varphi_h\mathbf{w}_h+\mathbf{g}_h$, where $\mathbf{w}_h$ is the only field predicted by WINO. After shape-function discretization, WINO is trained with the loss \eqref{eq:total_loss_Dirichlet}. We use the SOAP optimizer for $500$ epochs with a learning rate of $0.05$. The corresponding operator mapping is
\begin{equation}
	\mathcal{G}_{\theta}:(\varphi_h,\mathbf{g}_h,\mathbf{f}_h)\mapsto\mathbf{w}_h.
\end{equation}

\begin{figure}[t]
	\noindent\textbf{(a)}\par\vspace{0.4em}
	{\centering
		\includegraphics[width=0.18\textwidth]{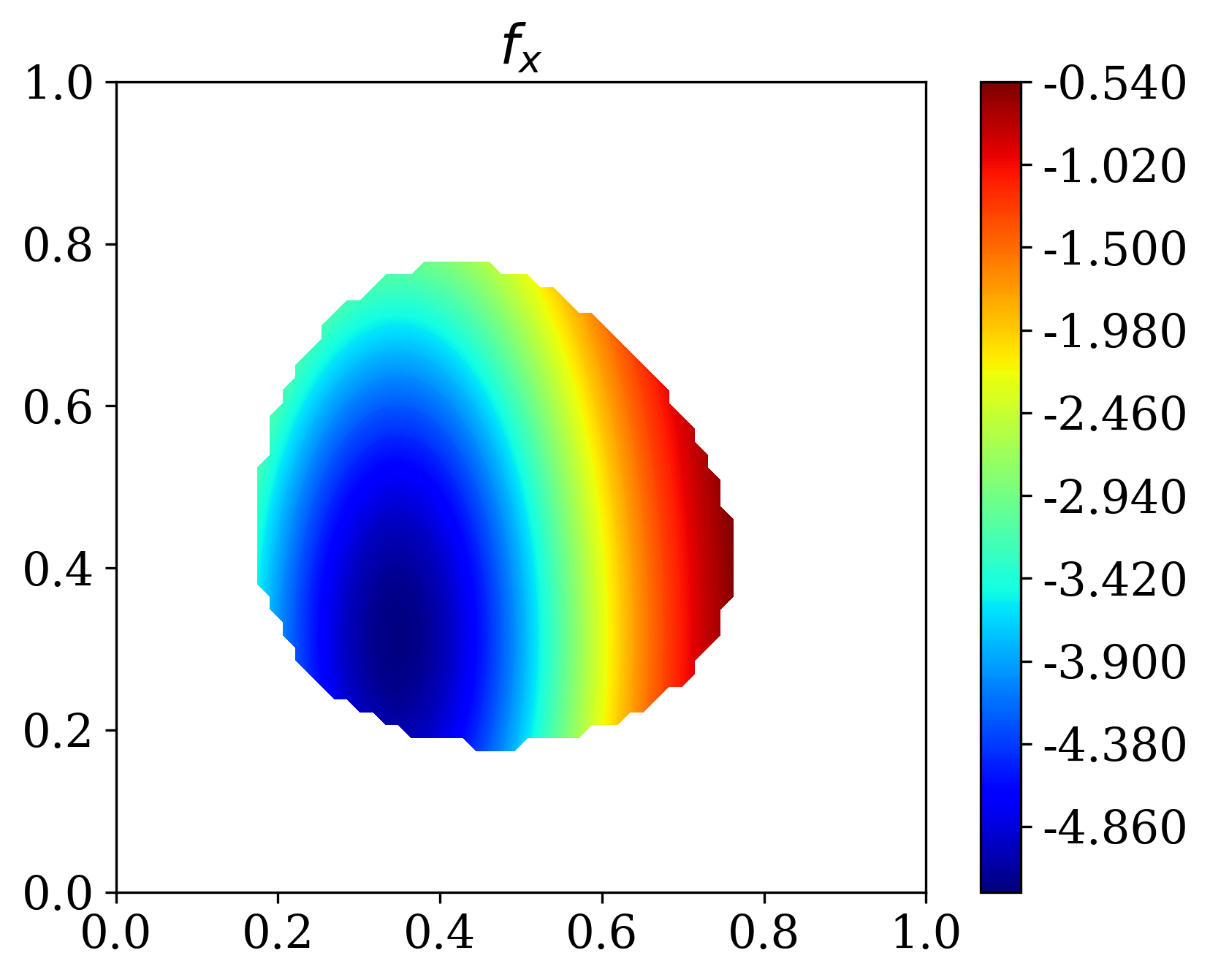}
		\includegraphics[width=0.18\textwidth]{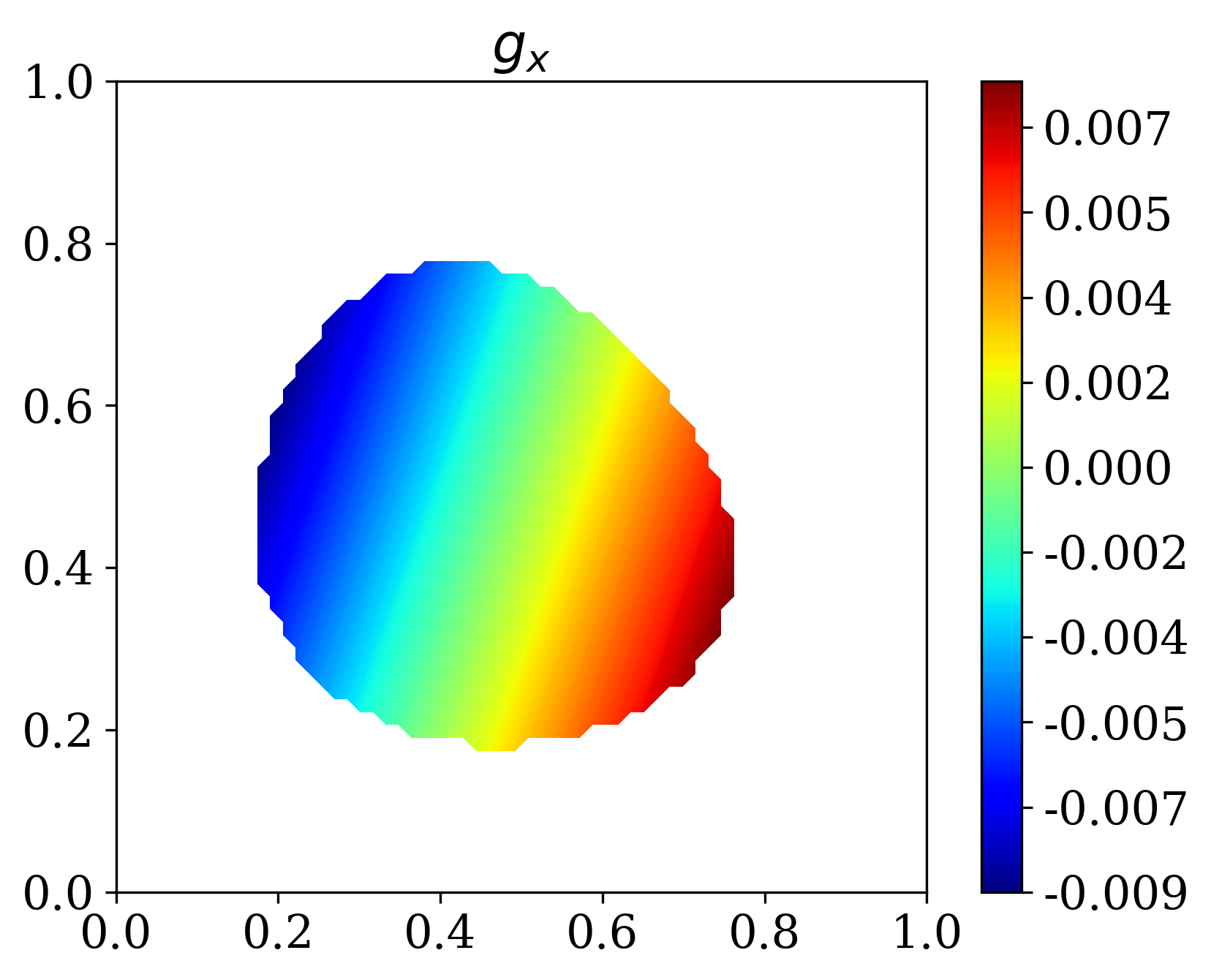}
		\includegraphics[width=0.18\textwidth]{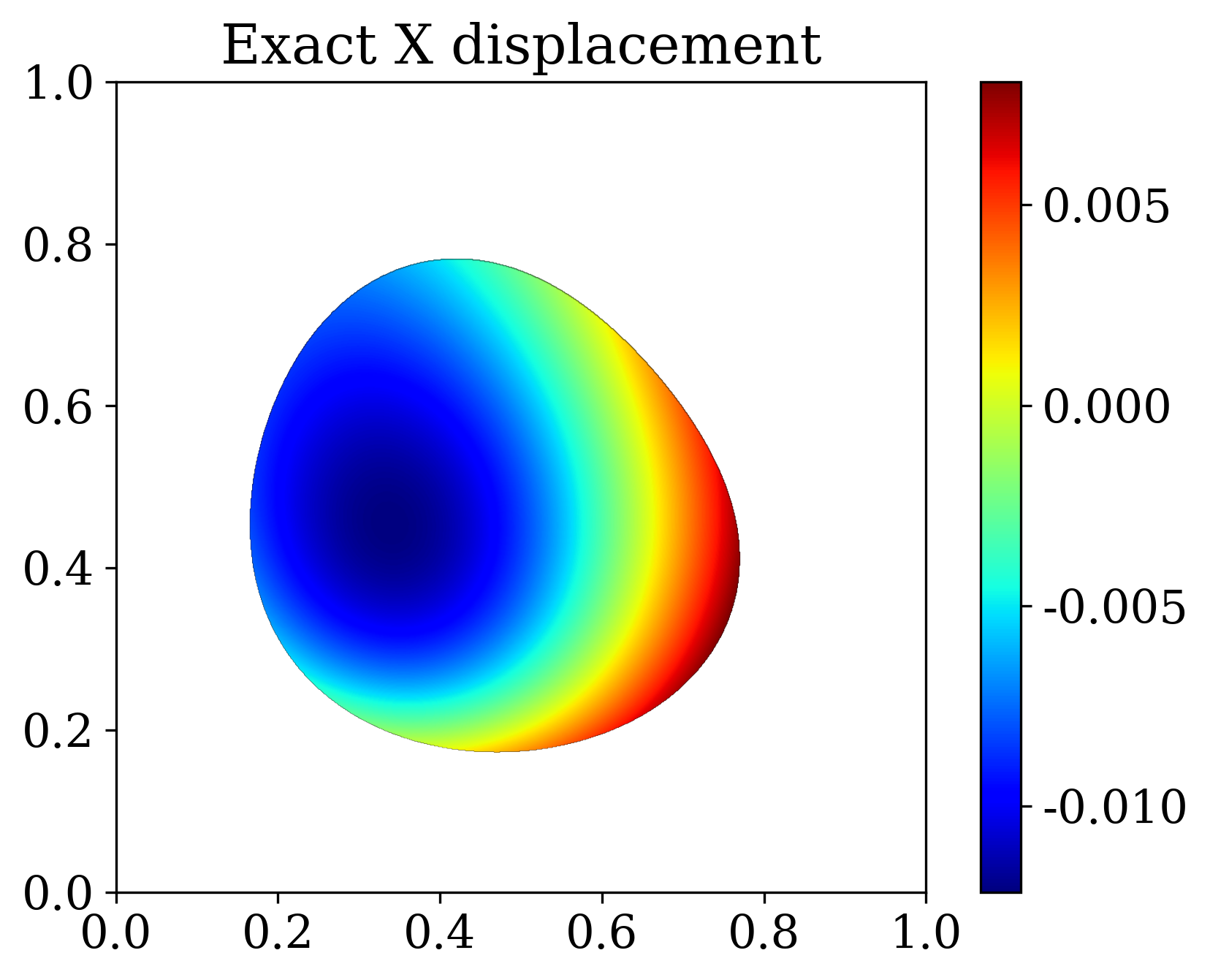}
		\includegraphics[width=0.18\textwidth]{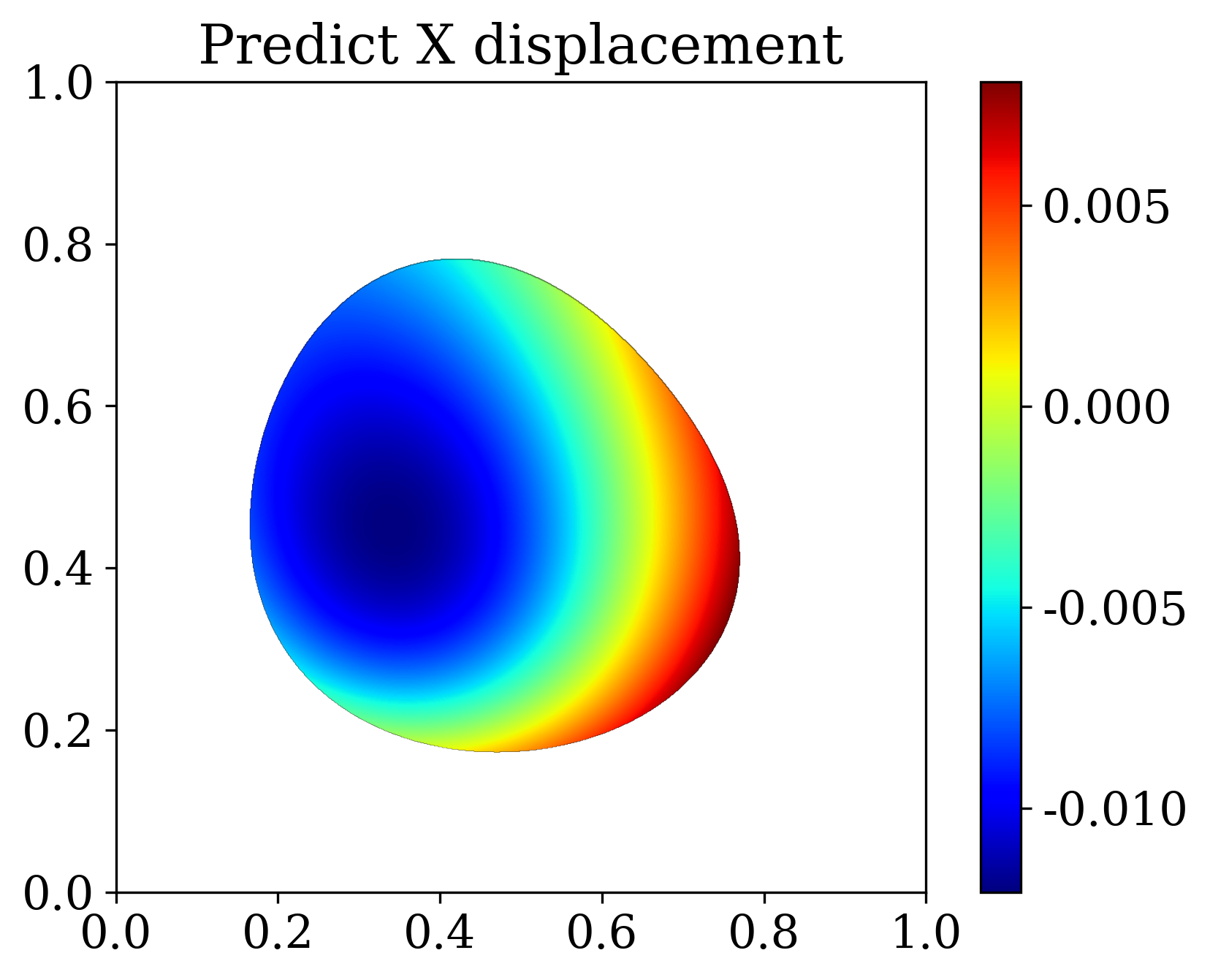}
		\includegraphics[width=0.19\textwidth]{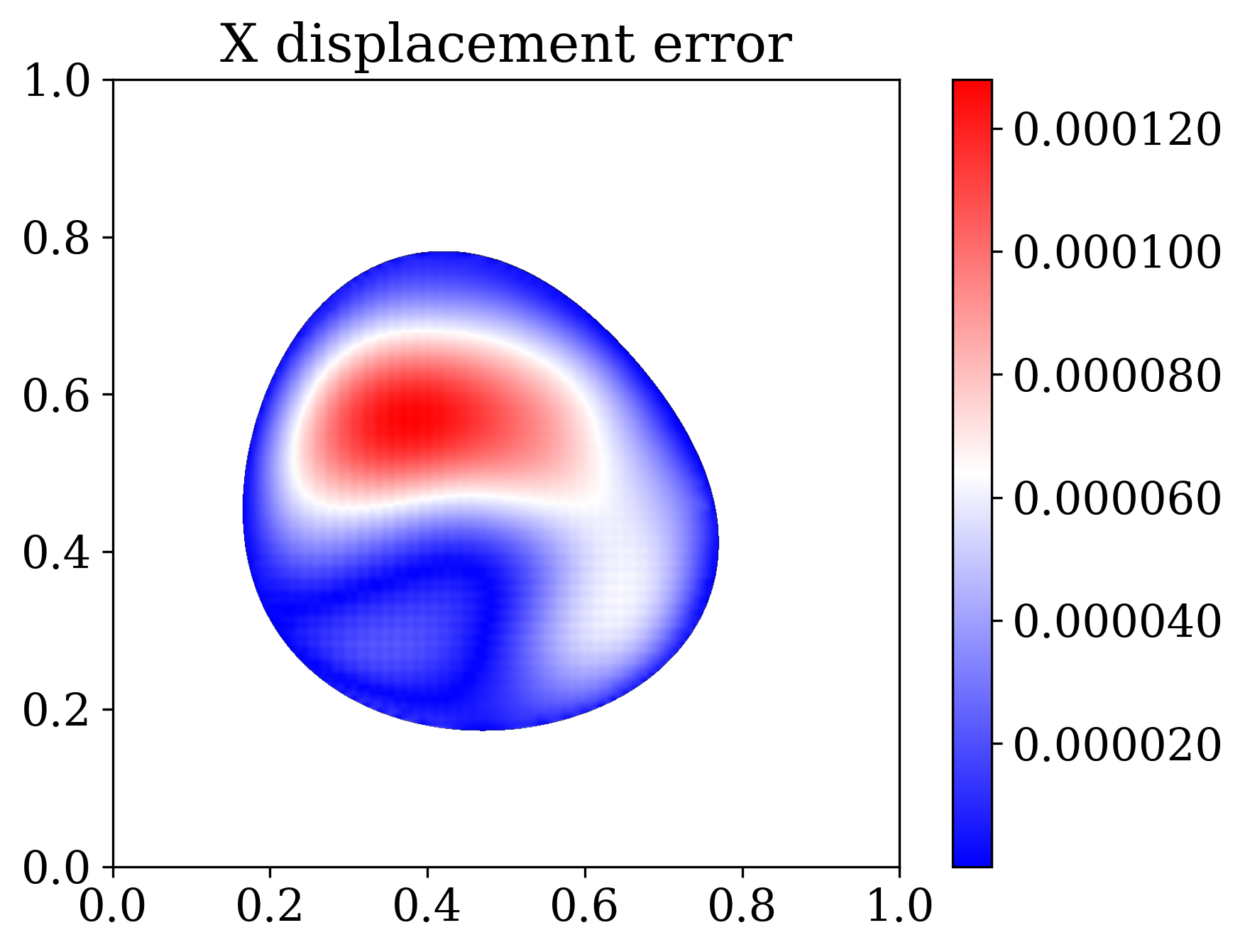}\\
		\vspace{0.1cm}
		\includegraphics[width=0.18\textwidth]{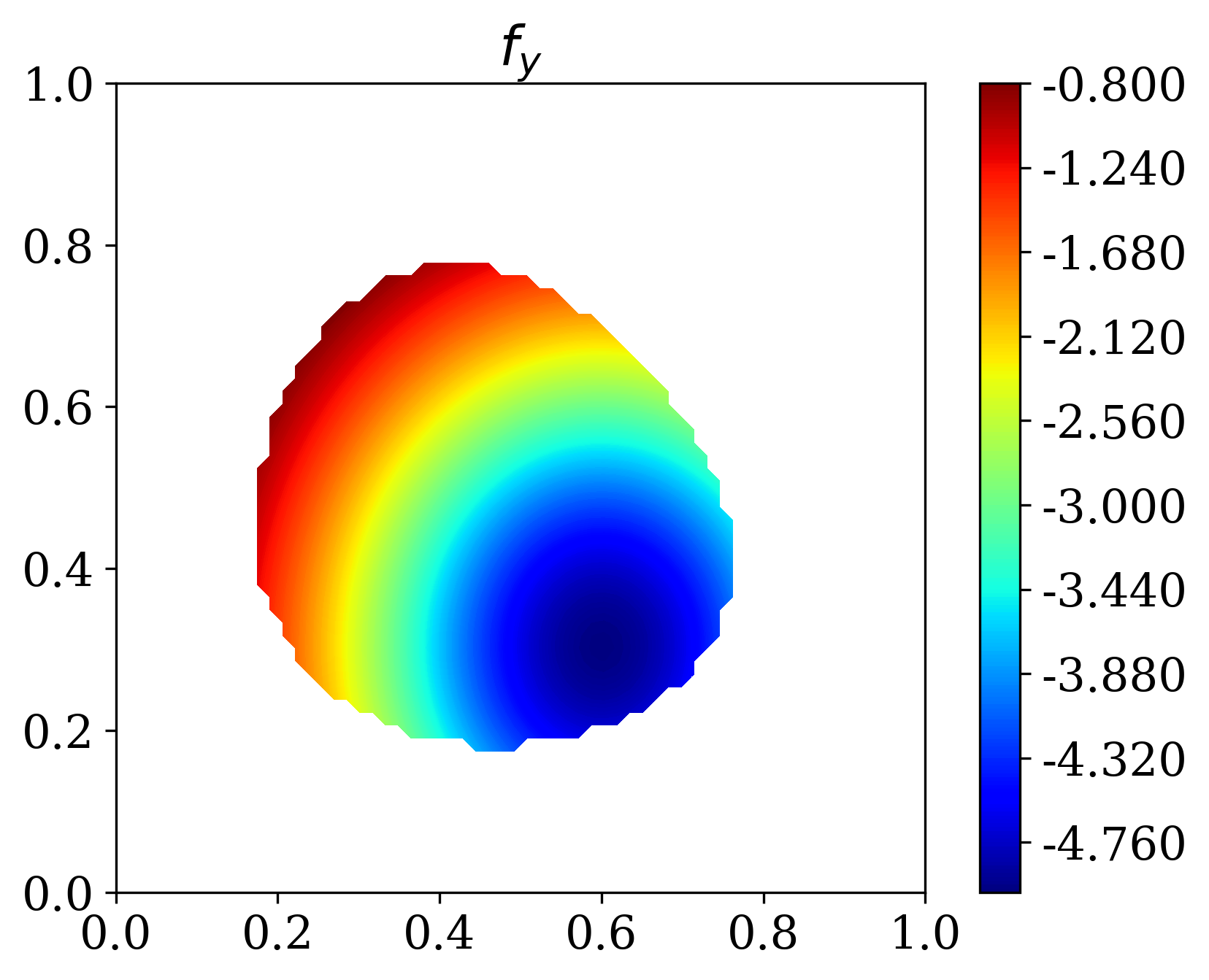}
		\includegraphics[width=0.18\textwidth]{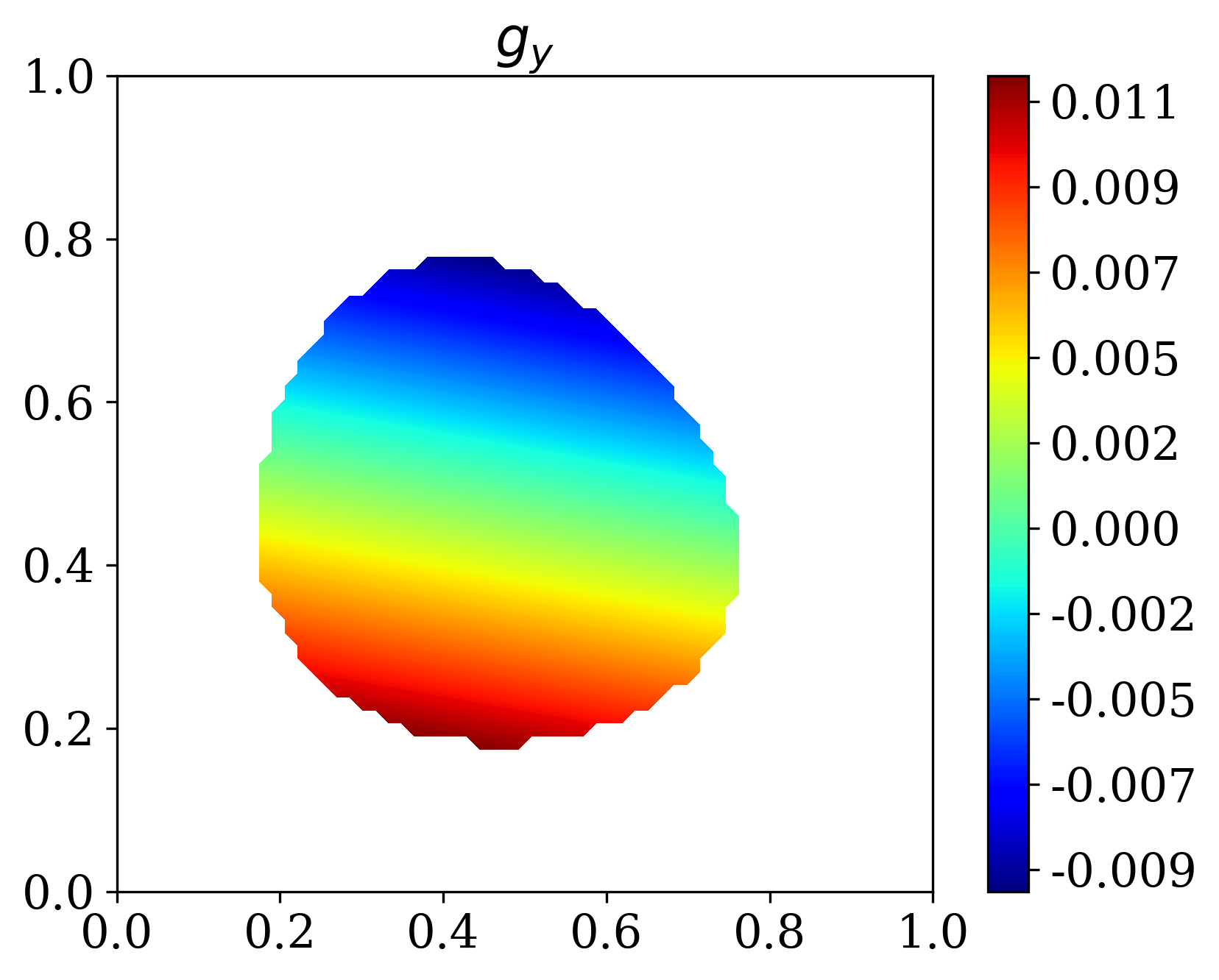}
		\includegraphics[width=0.18\textwidth]{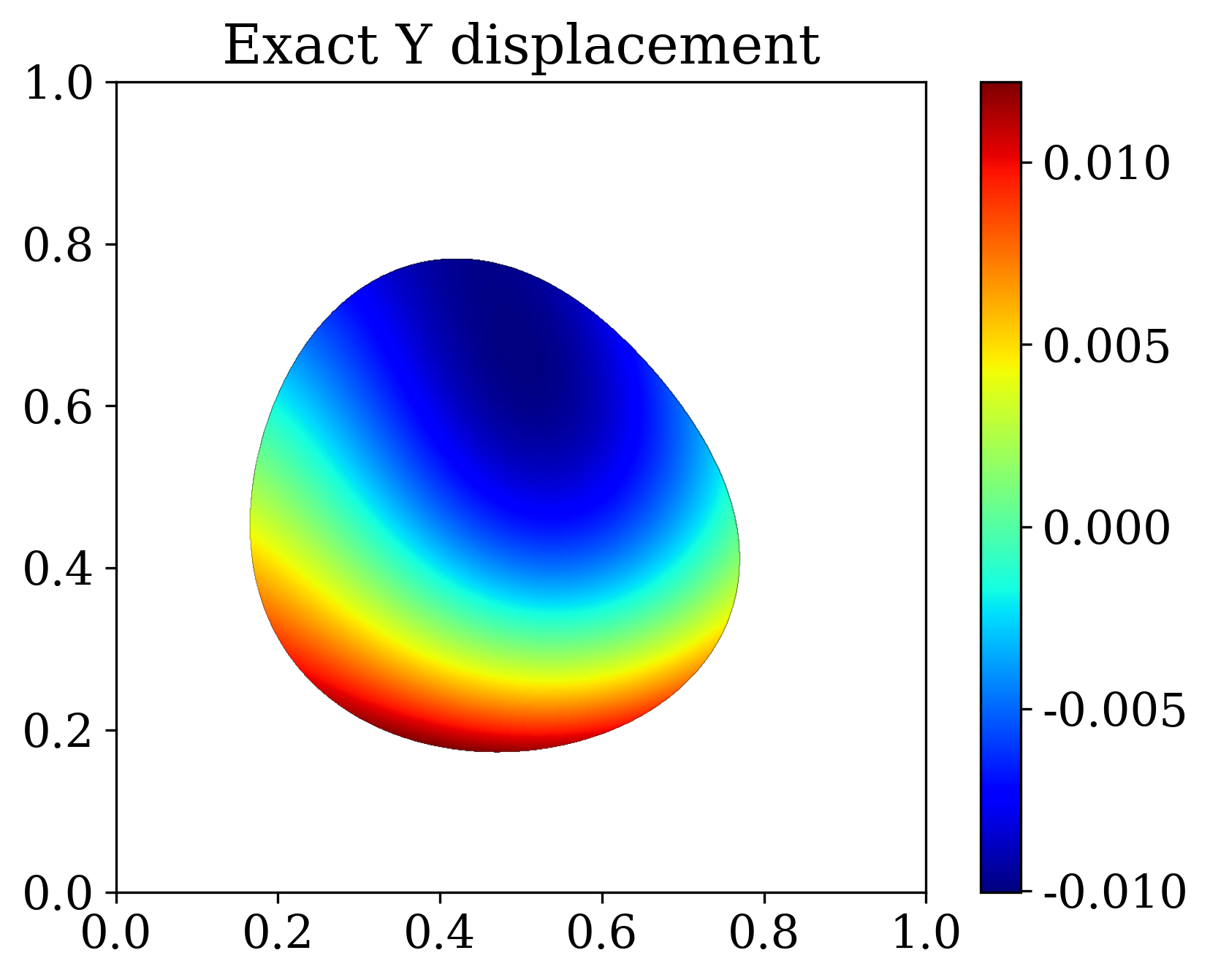}
		\includegraphics[width=0.18\textwidth]{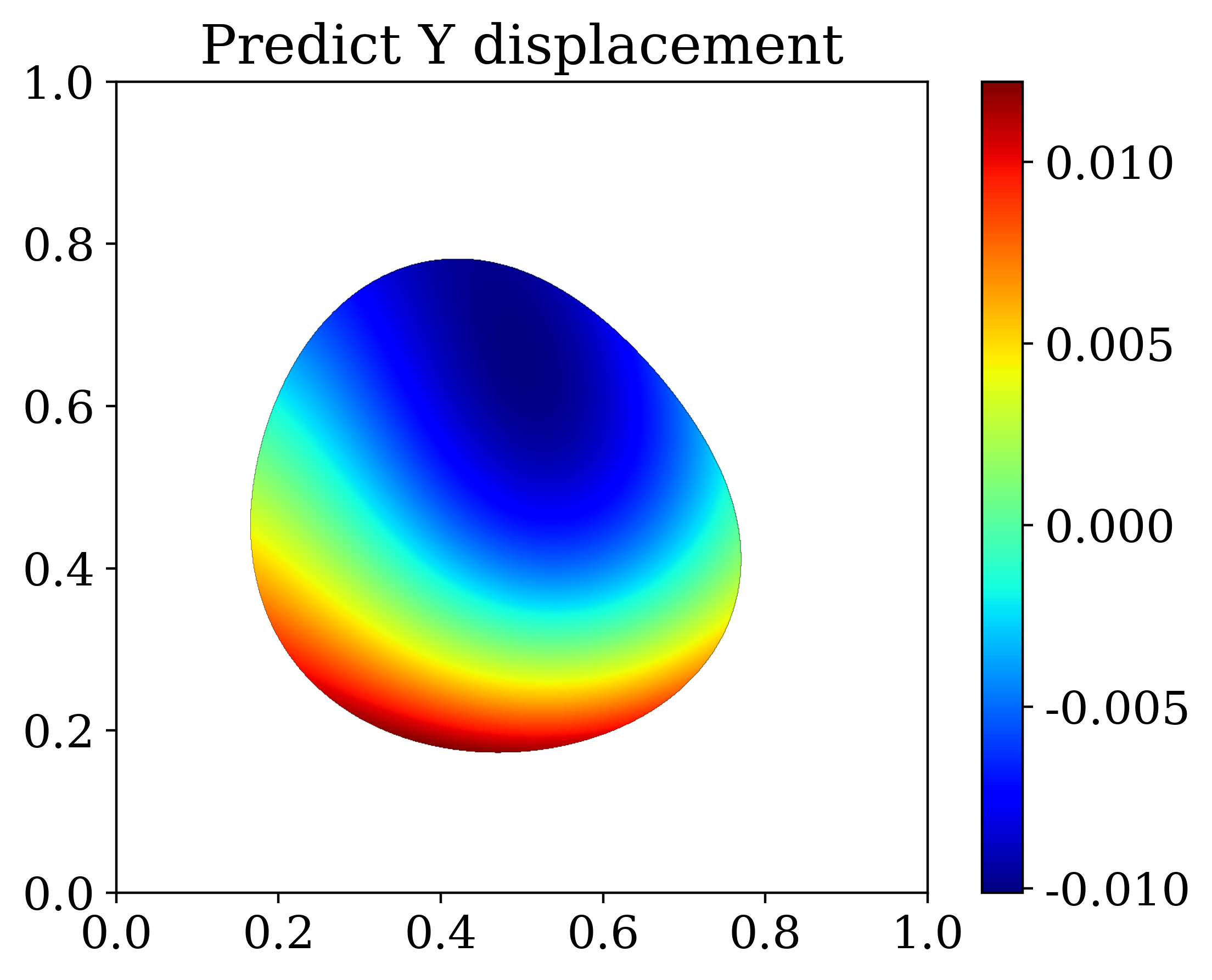}
		\includegraphics[width=0.19\textwidth]{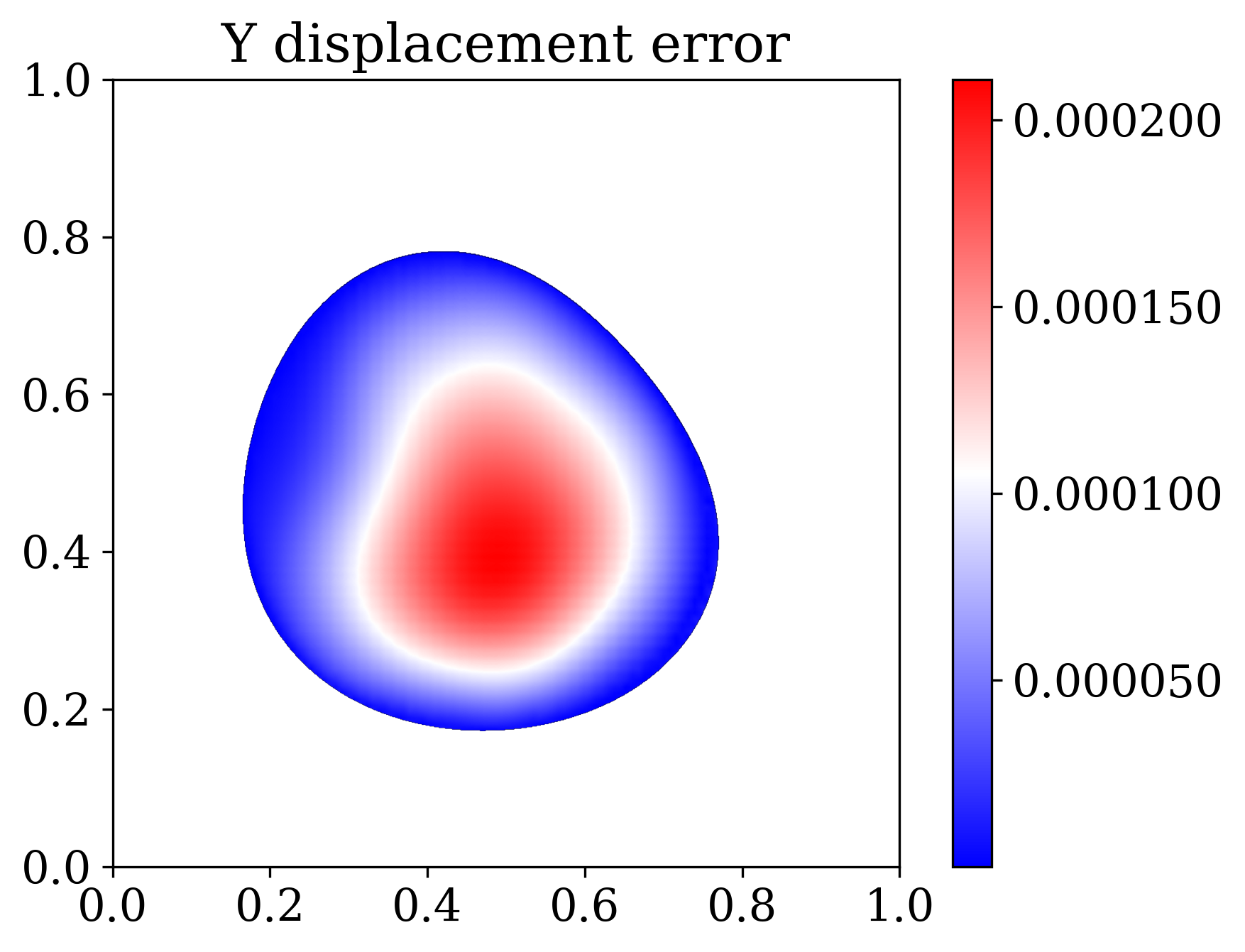}\par}
	\noindent\textbf{(b)}\par\vspace{0.4em}
	{\centering
		\includegraphics[width=0.18\textwidth]{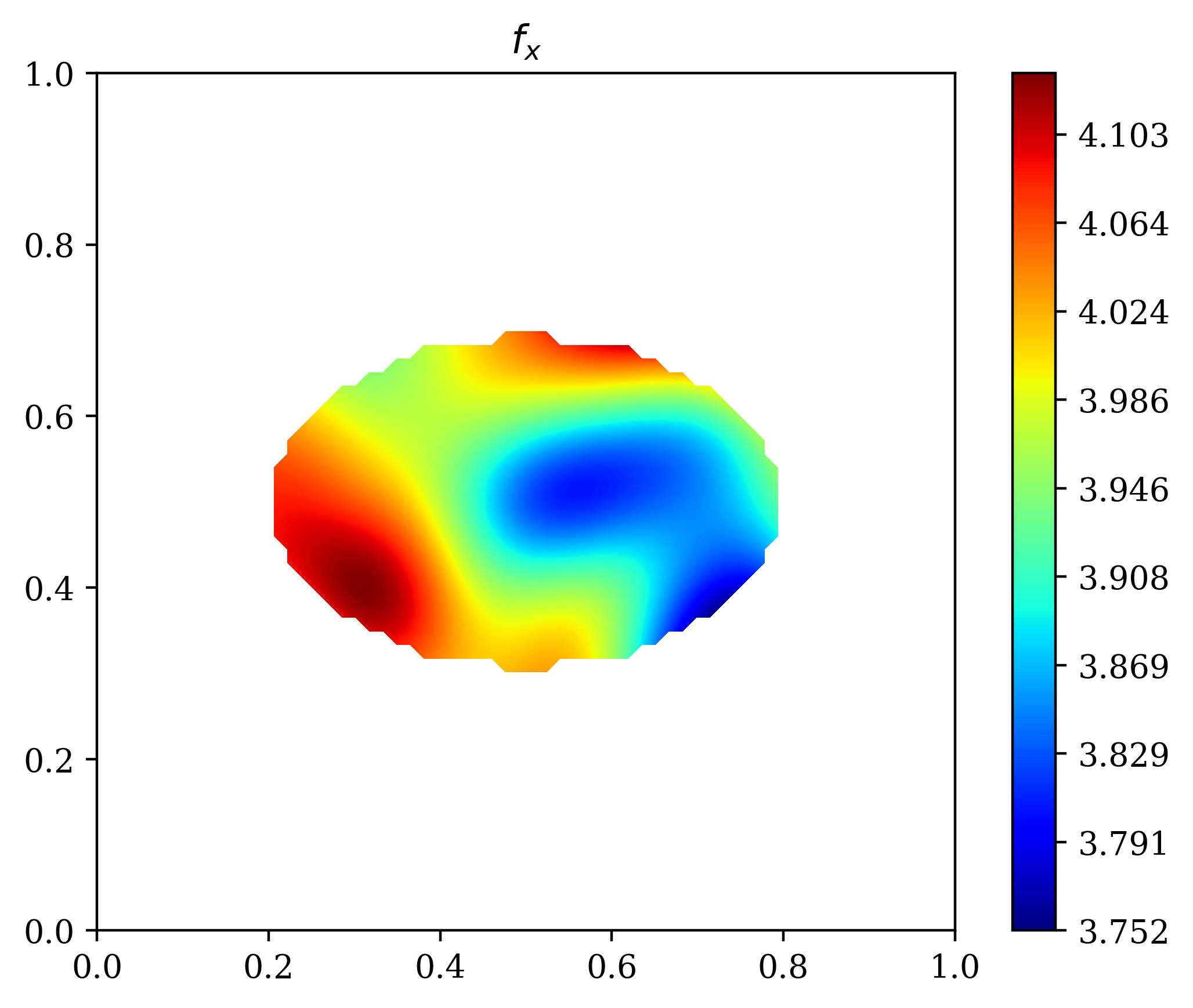}
		\includegraphics[width=0.18\textwidth]{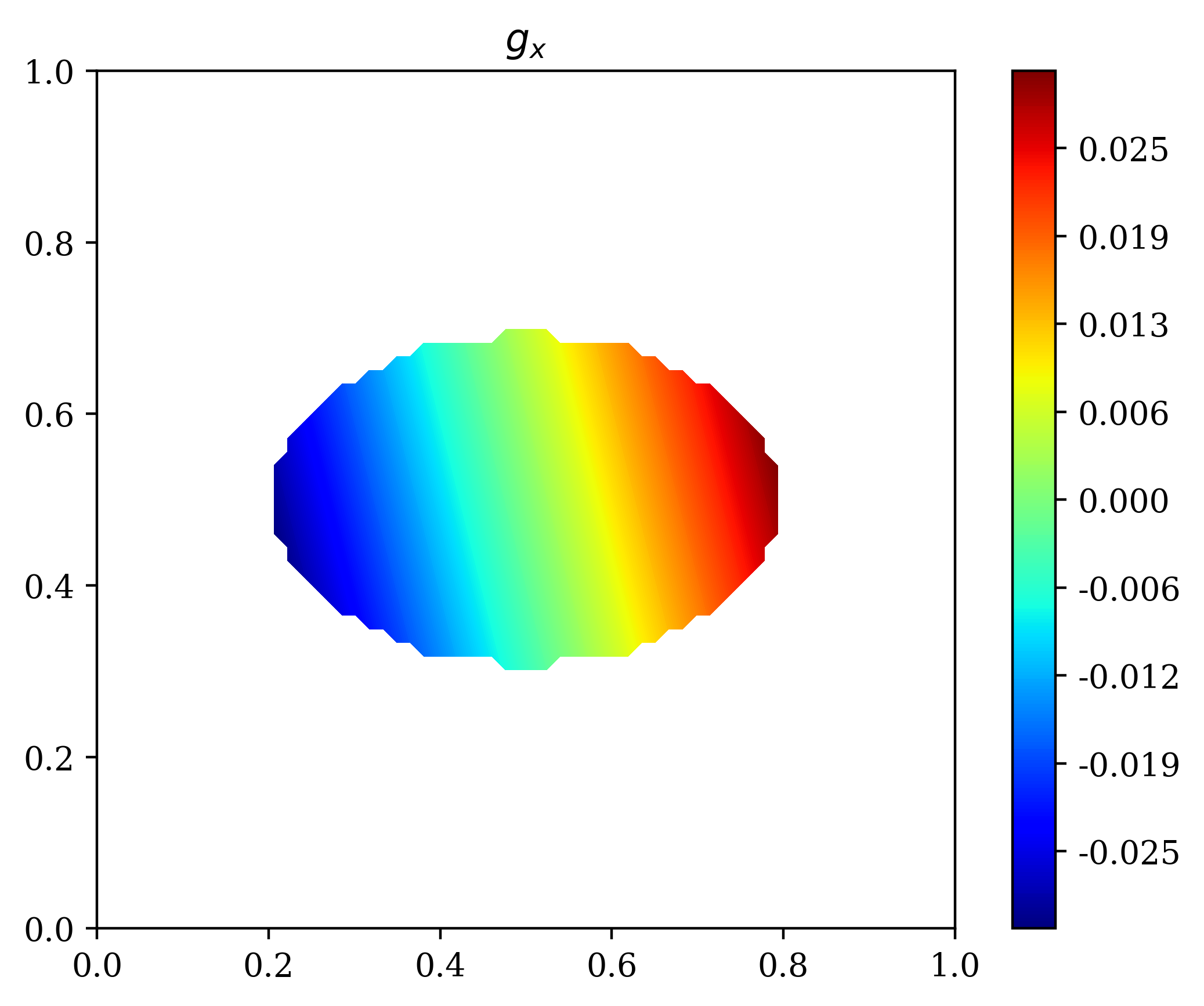}
		\includegraphics[width=0.18\textwidth]{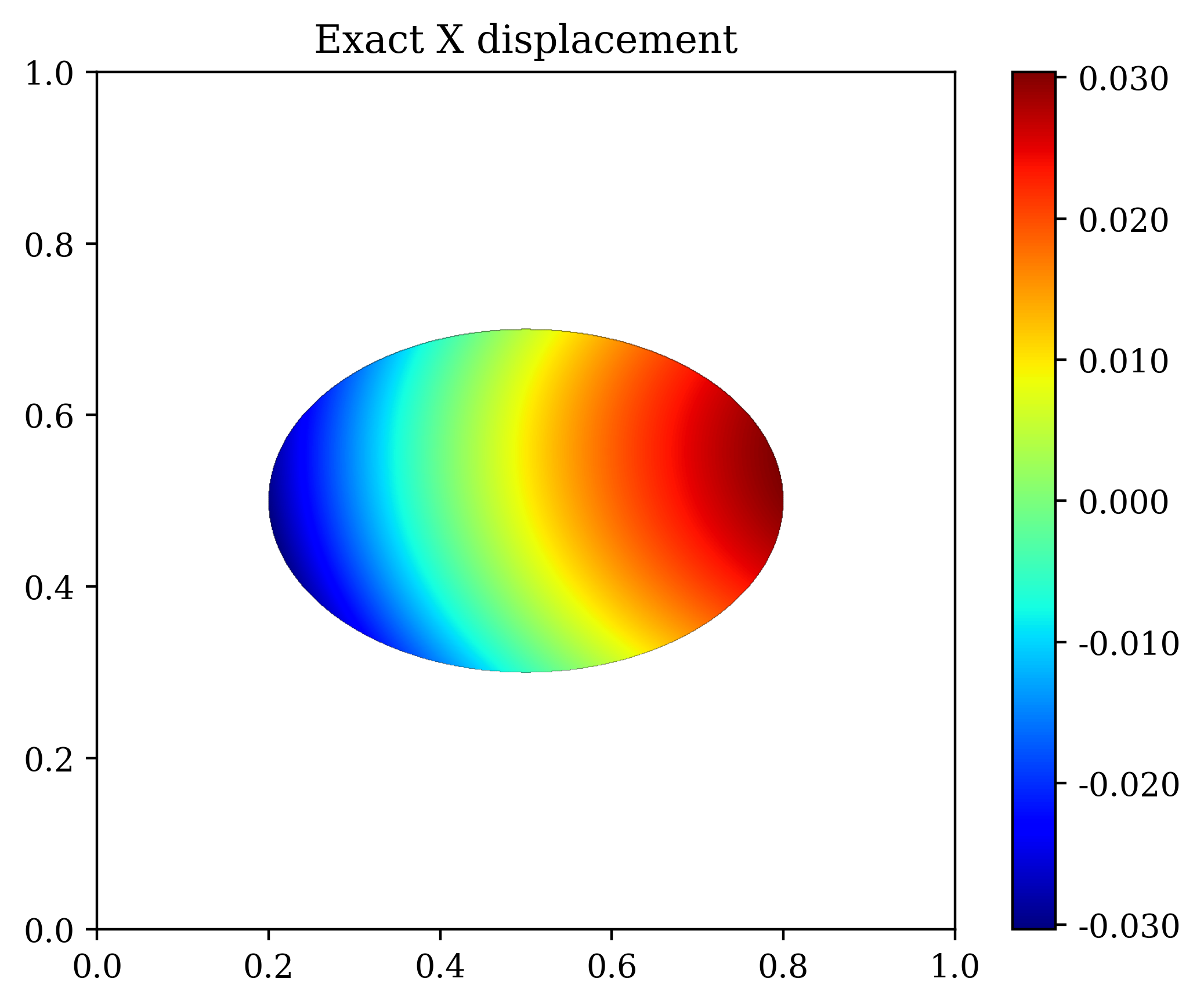}
		\includegraphics[width=0.18\textwidth]{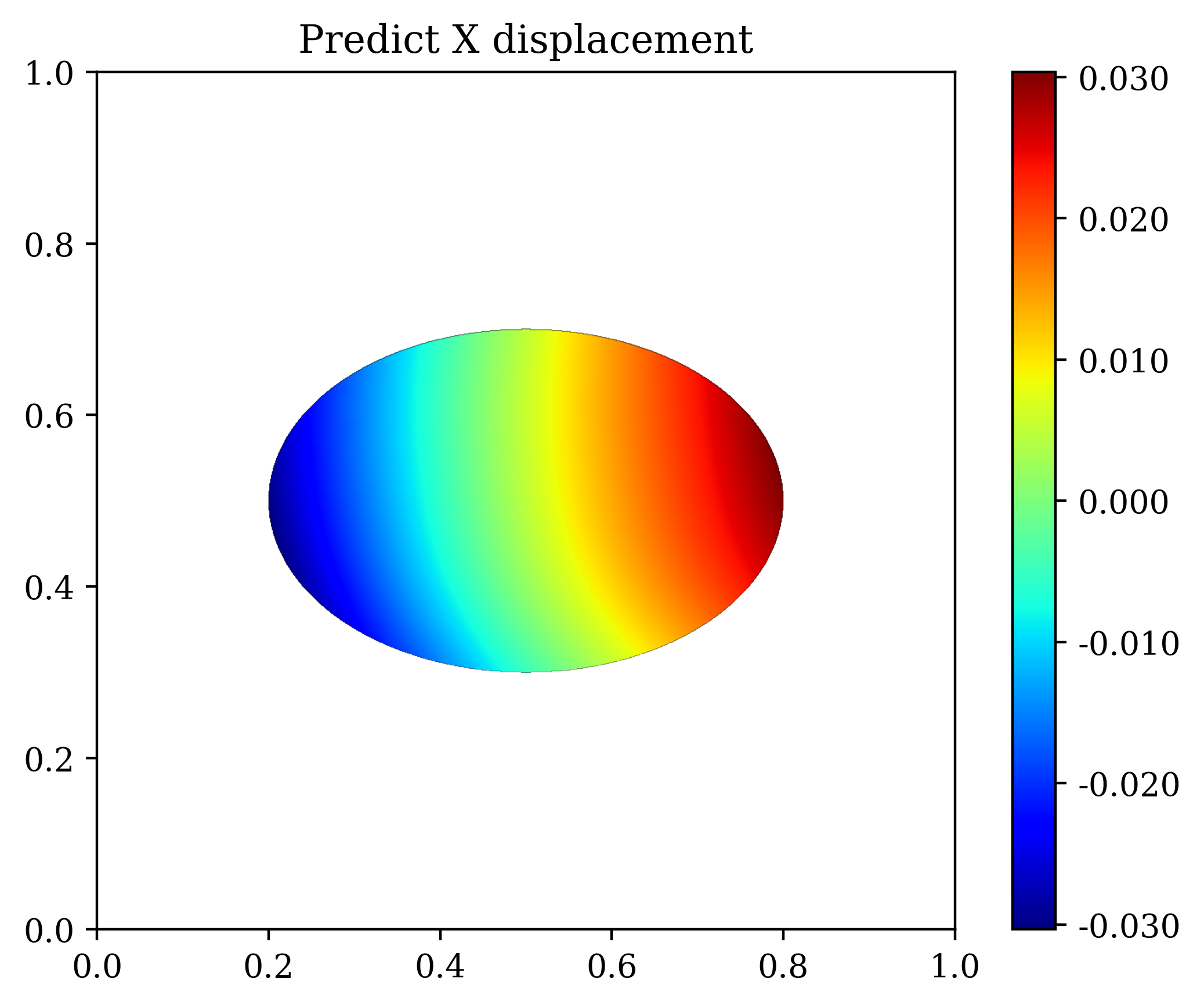}
		\includegraphics[width=0.185\textwidth]{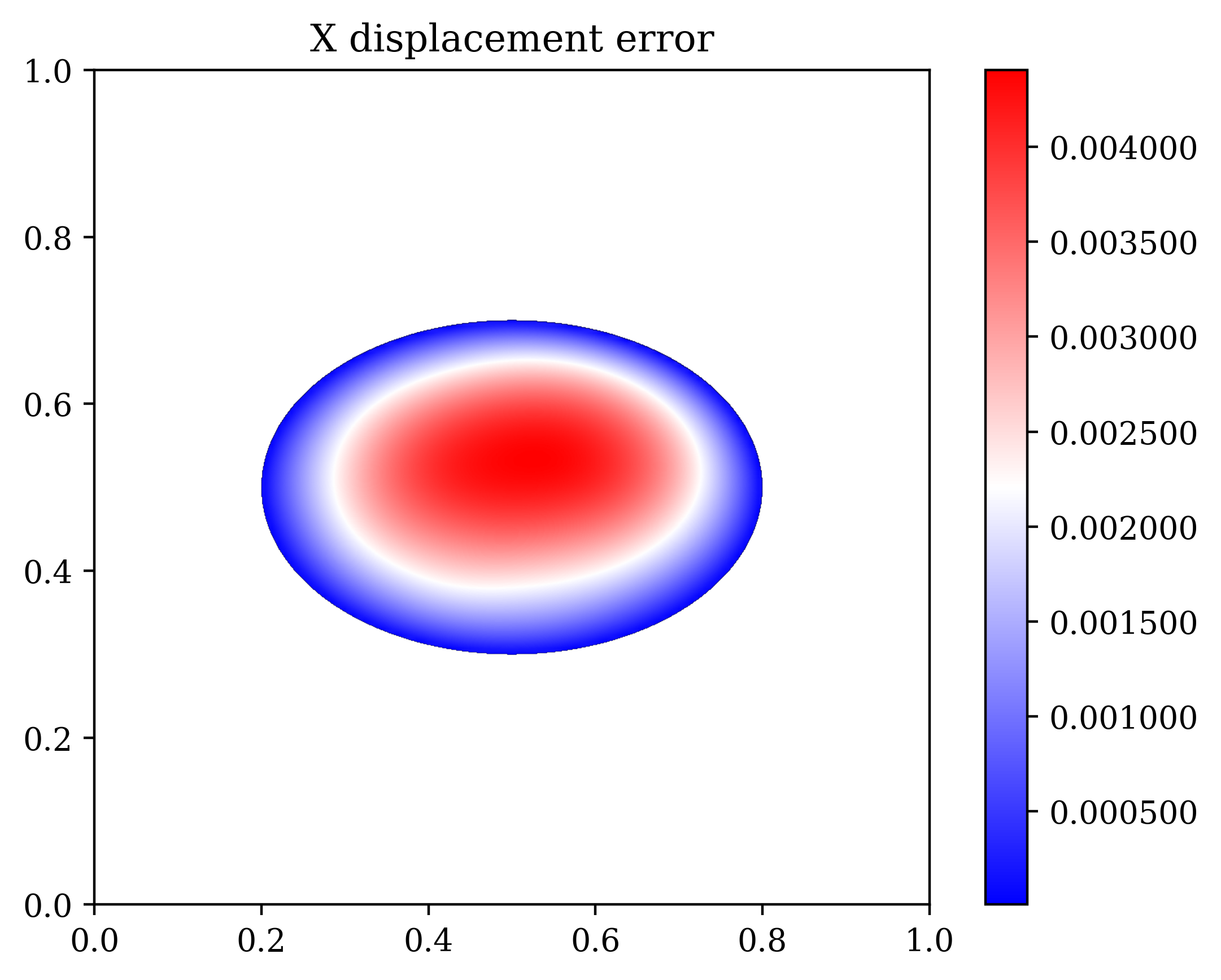}\\
		\vspace{0.1cm}
		\includegraphics[width=0.18\textwidth]{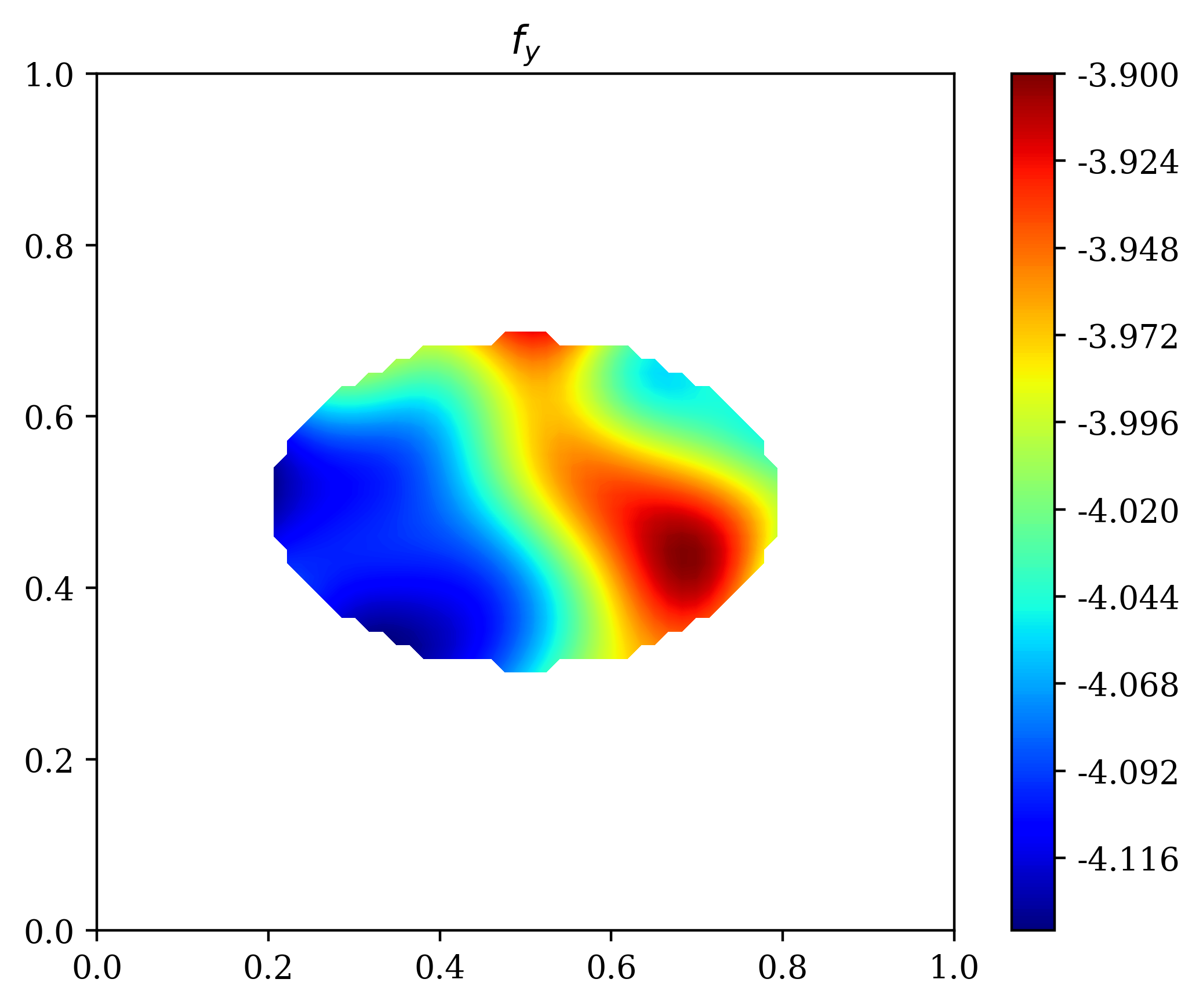}
		\includegraphics[width=0.18\textwidth]{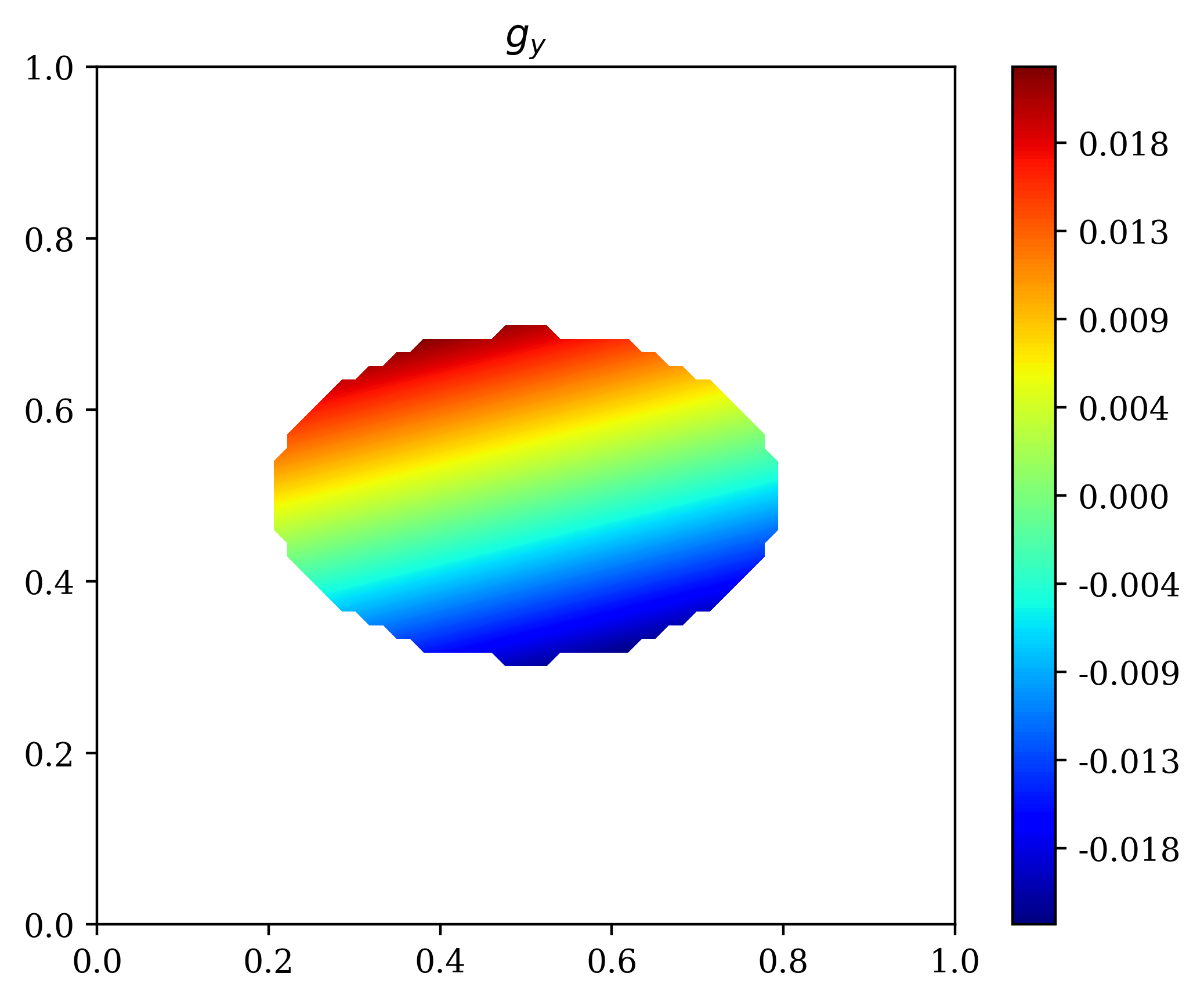}
		\includegraphics[width=0.18\textwidth]{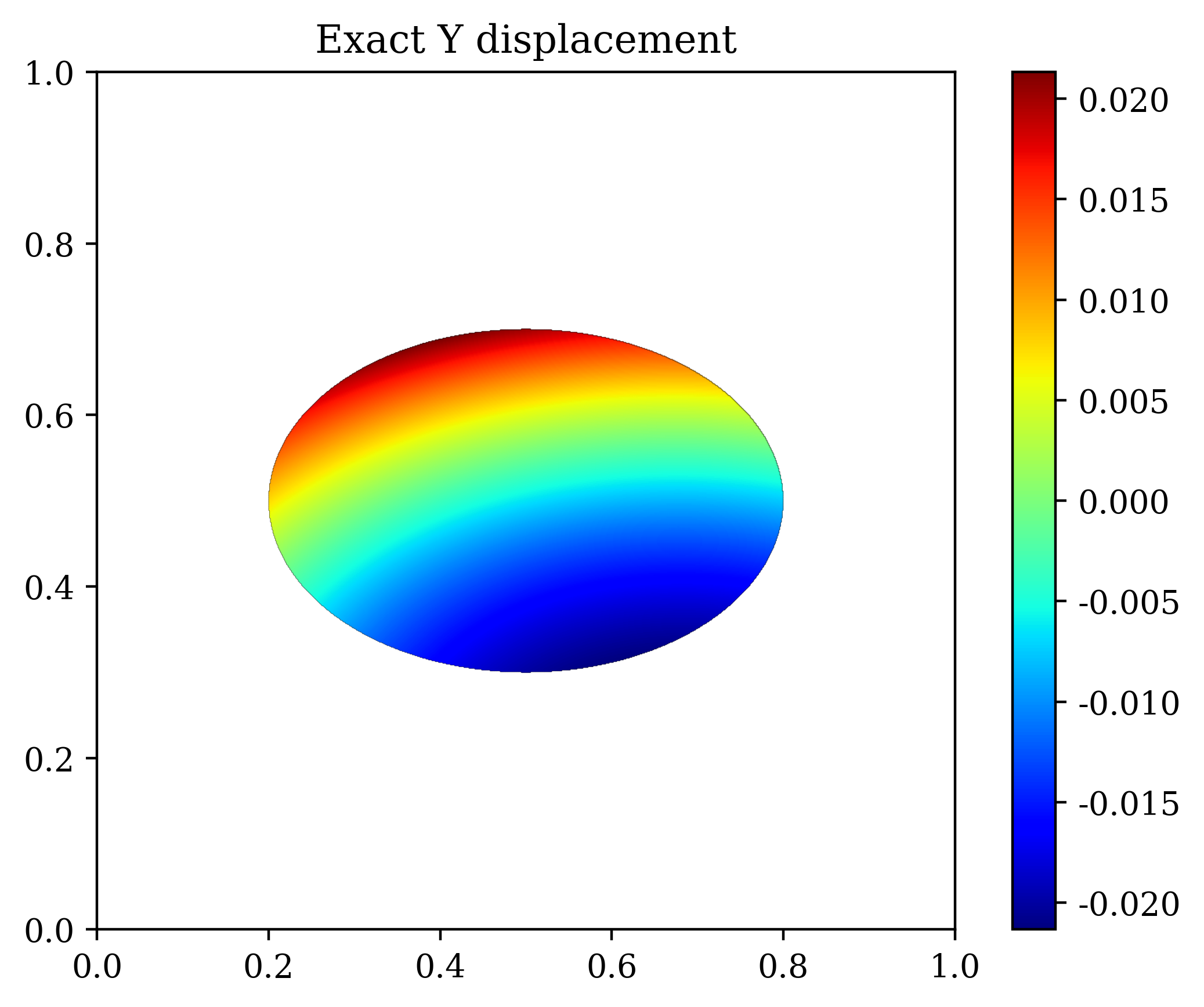}
		\includegraphics[width=0.18\textwidth]{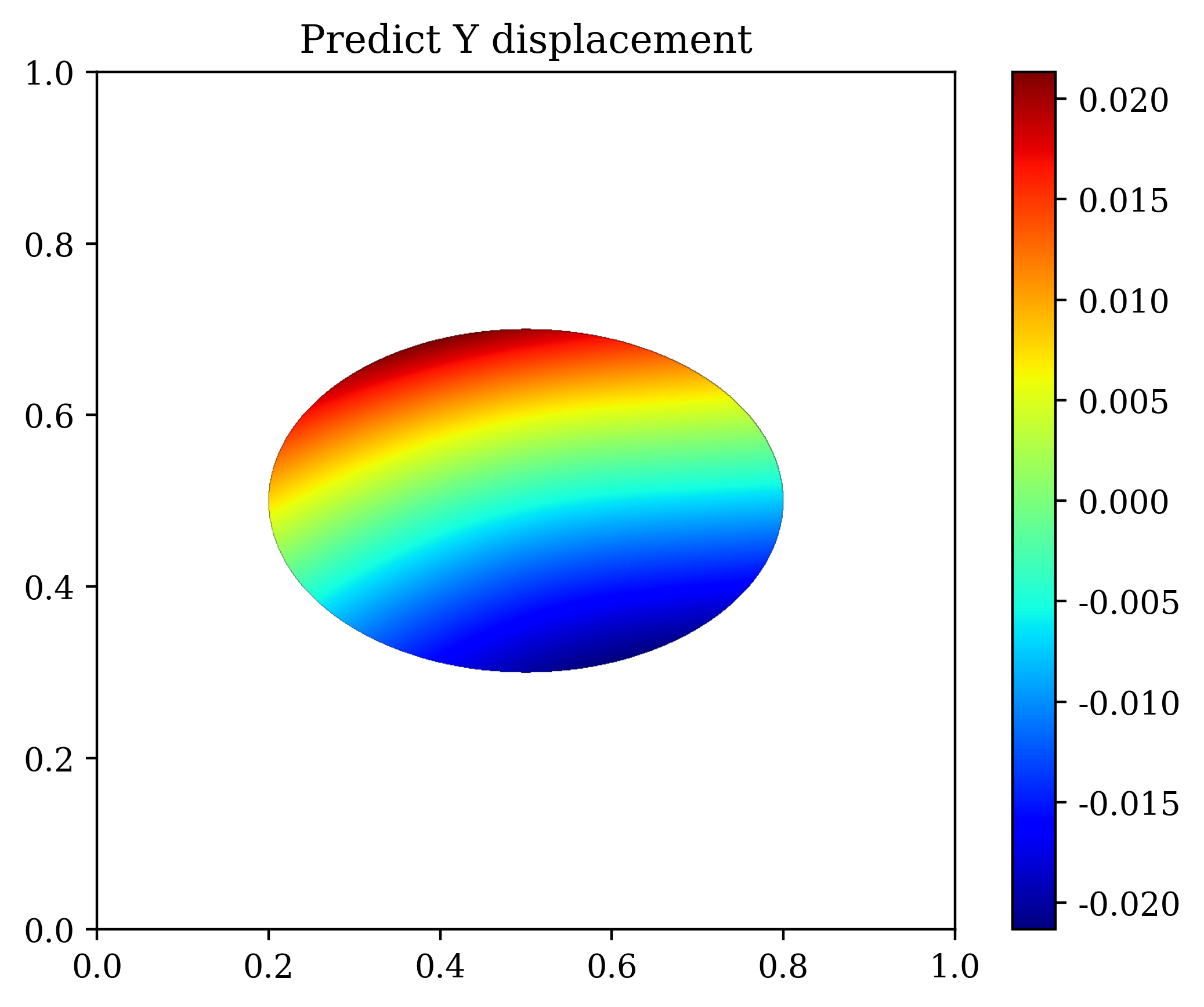}
		\includegraphics[width=0.185\textwidth]{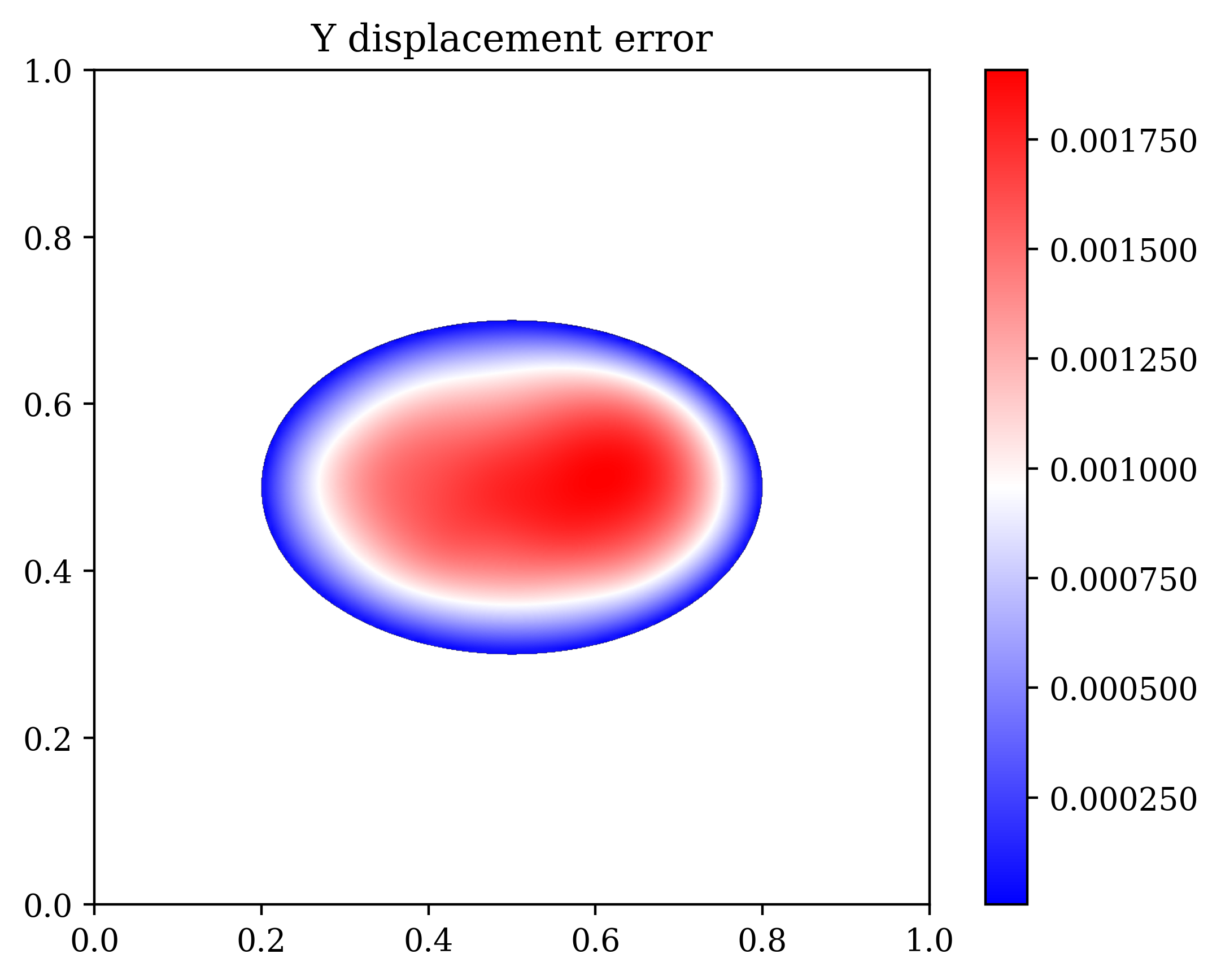}\par}
	\caption{WINO performance for the random shape case: (a) In-distribution test sample; (b) Out-of-distribution test sample. For rach panel: the top and bottom rows show the $x$- and $y$-components; the columns from left to right display $\mathbf{f}_h$, $\mathbf{g}_h$, exact displacement, WINO prediction, and absolute error contour, respectively.}
	\label{fig:case_arbit_performance}
\end{figure}

Fig.~\ref{fig:case_arbit_phi_error}b presents the epoch-wise evolution of the relative $L^2$ errors on the training and test sets, while Fig.~\ref{fig:case_arbit_phi_error}c reports box plots of the post-training WINO errors. WINO converges rapidly, and the relative $L^2$ errors on both the training and test sets decrease steadily without signs of overfitting.
A representative in-distribution test sample is shown in Fig.~\ref{fig:case_arbit_performance}a, including the input fields, exact displacement, and WINO prediction. The predicted field matches the reference solution closely, with pointwise absolute errors on the order of $\mathcal{O}(10^{-4})$. To further evaluate performance, Table~\ref{tab:performance_arbit} compares total computational cost (data generation + training) and accuracy for three approaches (WINO, purely data-driven $\varphi$-FEM-FNO, and WINO with labeled data) under the three metrics \eqref{eq:rel_err_L2}-\eqref{eq:rel_err_energy}. Because WINO does not require labeled reference-solution generation, its total cost is substantially lower than that of $\varphi$-FEM-FNO, and it also achieves smaller errors than $\varphi$-FEM-FNO across all three metrics. Incorporating labeled data into WINO further reduces these errors. To examine discretization invariance, we additionally evaluate WINO across mesh resolutions from $16\times16$ to $128\times128$ (Fig.~\ref{fig:case_arbit_phi_error}d). For the $16\times16$ and $32\times32$ settings, the retained Fourier modes in $\varphi$-FEM-FNO are set to 6 and 12, respectively. Even on the coarsest $16\times16$ mesh, most samples still stay within a relative error of $2\times10^{-2}$, indicating strong discretization invariance.

\begin{figure}[t]
	\centering
	\begin{minipage}[t]{0.58\textwidth}
		\centering
		\makebox[\linewidth][l]{\textbf{(a)}}\par\vspace{0.3em}
		\includegraphics[width=\linewidth]{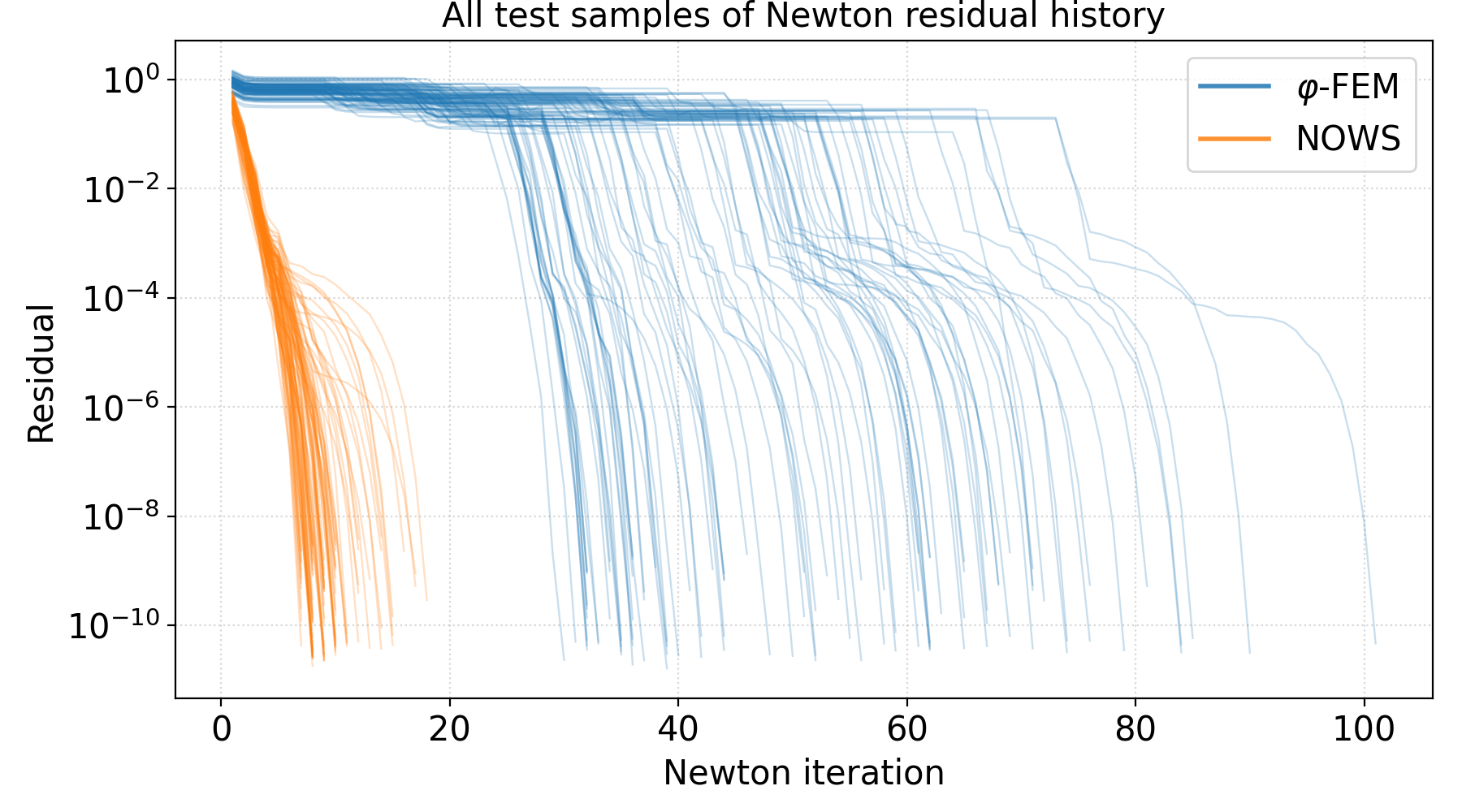}
	\end{minipage}
	\hfill
	\begin{minipage}[t]{0.39\textwidth}
		\centering
		\makebox[\linewidth][l]{\textbf{(b)}}\par\vspace{0.3em}
		\includegraphics[width=\linewidth]{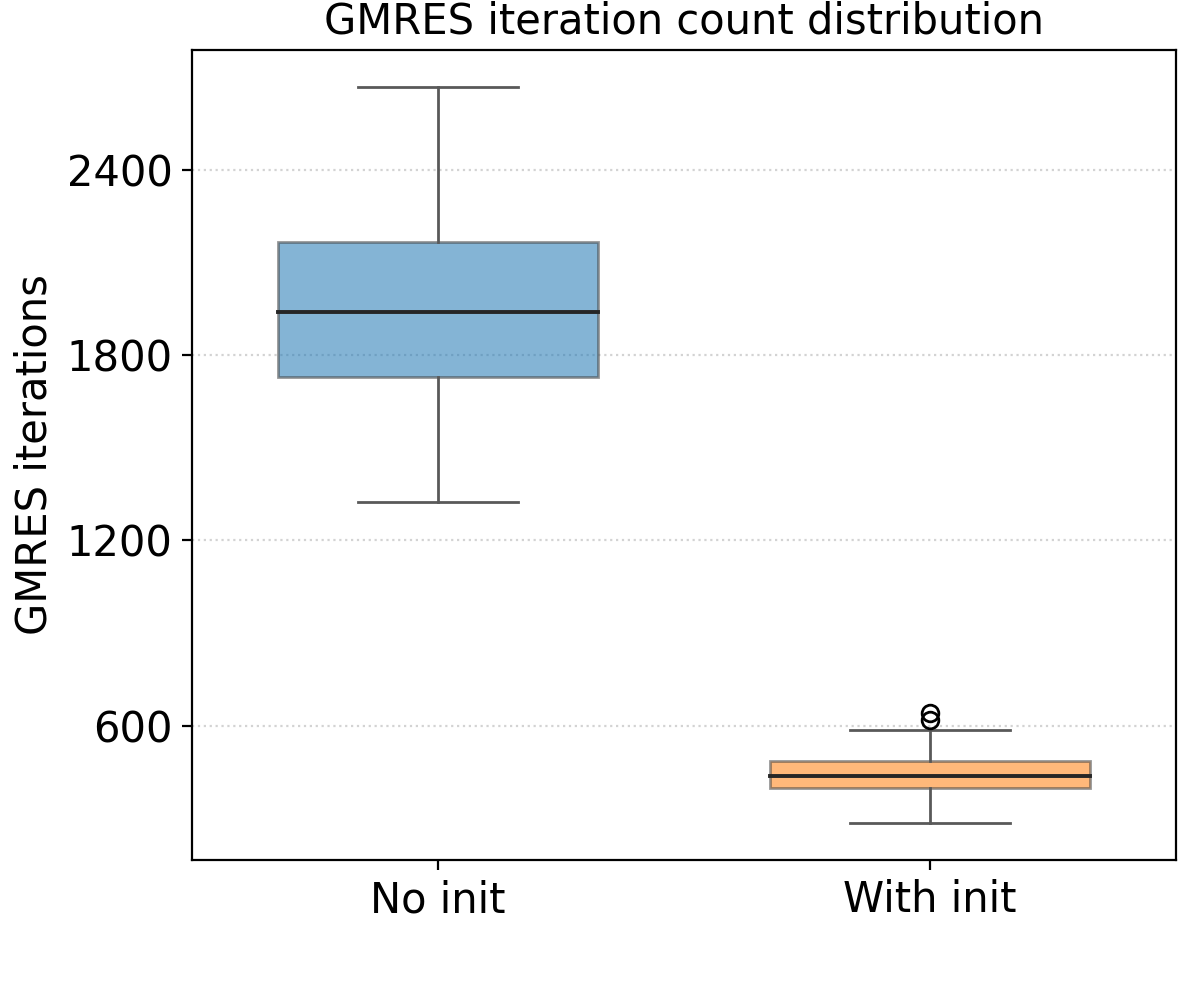}
	\end{minipage}
	\caption{NOWS comparison for the random shape case. (a) Outer Newton residual and total Newton iterations for classical $\varphi$-FEM (4 load increments) and NOWS. (b) Box plots of total GMRES iteration counts for all load increments.}
	\label{fig:case_arbit_nows}
\end{figure}

To further evaluate out-of-distribution (OOD) robustness, we test WINO on an unseen elliptical domain and unseen forcing realizations. The OOD level-set uses \eqref{eq:elliptical_phi} with $(c_x,c_y)=(0.5,0.5)$ and $(l_x,l_y)=(0.3,0.2)$. The displacement boundary function uses the same parametric form as in training \eqref{eq:case_ellip_g} with $(\alpha,\beta)=(1,1)$. For body-force $\mathbf{f}$ sampling, we use the Gaussian random field (GRF) \cite{rasmussen2003gaussian} process
\begin{equation}
	\eta(x,y)\sim\mathcal{GP}(\mu,k_\ell),
	\label{eq:grf_process}
\end{equation}
with squared-exponential covariance
\begin{equation}
	k_\ell(\mathbf{x},\mathbf{y})=s^2\exp\left(-\frac{\|\mathbf{x}-\mathbf{y}\|^2}{2\ell^2}\right).
	\label{eq:grf_cov}
\end{equation}
The two force components are sampled independently as $f_1\sim\mathcal{GP}(4.0,k_{l})$ and $f_2\sim\mathcal{GP}(-4.0,k_{l})$ with $s=0.1$ and $l=0.5$. The prediction and absolute error are shown in Fig.~\ref{fig:case_arbit_performance}b. Under this OOD setting, errors are moderately larger than in-distribution errors, while the predicted displacement fields remain physically consistent.

We next evaluate NOWS for this random-shape setting, using the WINO prediction as the initial iterate of the nonlinear $\varphi$-FEM solve. As a baseline, classical $\varphi$-FEM employs incremental load stepping (4 load increments) to improve Newton robustness, whereas NOWS directly solves the full-load problem without continuation. At each Newton step, the linearized system is solved by restarted GMRES with $m=120$ and relative tolerance $10^{-5}$. Fig.~\ref{fig:case_arbit_nows} compares Newton residual histories, total Newton iterations, and box plots of total GMRES iterations with and without neural initialization. Averaged over the test instances, the wall-clock time for the nonlinear $\varphi$-FEM solve without a neural warm start is $2.19\pm 0.67\,\mathrm{s}$, compared with $0.374\pm 0.08\,\mathrm{s}$ when WINO supplies the initial iterate. NOWS significantly reduces both Newton and GMRES iteration counts; a single WINO inference (about $5\,\mathrm{ms}$) is sufficient to provide a high-quality warm start and improve overall solver efficiency.

\subsection{Mixed Dirichlet and Neumann problems}

In this subsection, we consider hyperelastic problems with mixed Dirichlet and Neumann boundary conditions. The level-set function $\varphi$ is used only to represent the Neumann portion of the boundary, whereas Dirichlet segments are aligned with edges of the Cartesian background mesh $\mathcal{T}_h^{\mathcal{O}}$: no cut cells arise along those sides, and the prescribed displacements are imposed with a standard finite-element treatment. The $\varphi$-FEM weak formulation is \eqref{eq:phifem_weak_form}, and the WINO training objective built on that discretization is \eqref{eq:total_loss}. We present three numerical examples with different domains: a plate with an elliptic hole, Cook's membrane, and a pressure vessel. The penalty parameters selection details in the \eqref{eq:strong_loss}--\eqref{eq:total_loss} are described in the appendix \ref{app:penalty}.

\subsubsection{Plate with a hole}\label{sec:plate_with_hole}

\begin{figure}[t]
	\centering
	\includegraphics[width=0.385\textwidth]{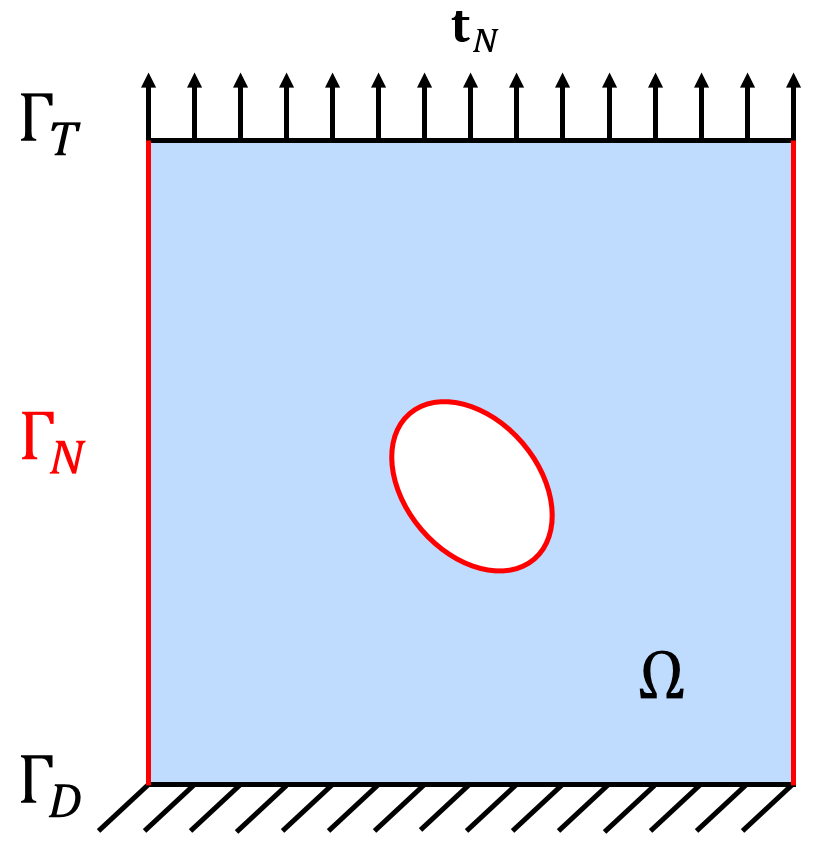}
	\hspace{0.5cm}
	\raisebox{4.1mm}{
		\includegraphics[width=0.31\textwidth]{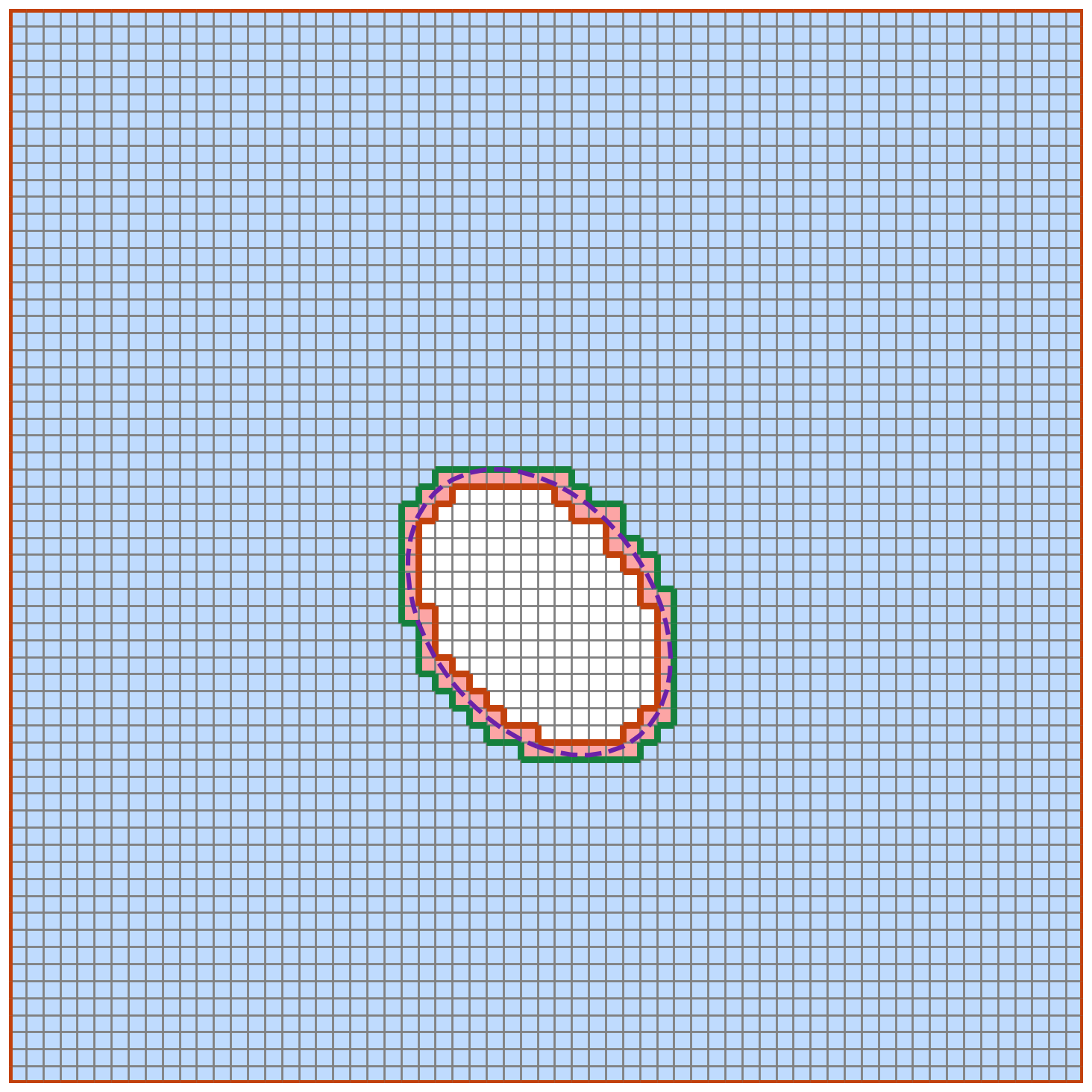}
	}
	\caption{The geometry and mesh tags for the plate with a hole case.}
	\label{fig:case_hole}
\end{figure}

In this subsection, we consider the hyperelastic problem on a square plate with an elliptic hole, which is formulated as
\begin{equation}
	\begin{split}
		-\nabla \cdot \mathbf{P}(\mathbf{F}(\mathbf{u})) &= \mathbf{0}, \quad \text{in } \Omega, \\
		\mathbf{u} &= \mathbf{0}, \quad \text{on } \Gamma_D, \\
		\mathbf{P}(\mathbf{F}(\mathbf{u})) \cdot \mathbf{n} &= \mathbf{t}_N, \quad \text{on } \Gamma_T, \\
		\mathbf{P}(\mathbf{F}(\mathbf{u})) \cdot \mathbf{n} &= \mathbf{0}, \quad \text{on } \Gamma_N,
	\end{split}
	\label{eq:case_hole}
\end{equation}
where the background domain is the square $\mathcal{O}=[0,1]^2$, and the physical domain $\Omega\subset\mathcal{O}$ is a square plate with an elliptic hole represented by the level-set function $\varphi$. Here, $\Gamma_D$ is the bottom side of $\Omega$, $\Gamma_T$ is the top boundary of $\Omega$, and $\Gamma_N$ is the remaining Neumann boundary, i.e., $\Gamma_N=\partial\Omega\setminus(\Gamma_D\cup\Gamma_T)$; see Fig.~\ref{fig:case_hole} (left). The background discretization uses a $64\times 64$ Cartesian grid on $[0,1]^2$, and the corresponding mesh and Neumann-boundary partition are shown in Fig.~\ref{fig:case_hole} (right). The elliptic hole is defined by the level set function
\begin{equation}
	\varphi(x, y) = 1 - \frac{((x-x_0)\cos(\theta) + (y-y_0)\sin(\theta))^2}{l_a^2} - \frac{((x-x_0)\sin(\theta) - (y-y_0)\cos(\theta))^2}{l_b^2},
\end{equation}
where $(x_0, y_0)$ is the center of the elliptic hole, $\theta$ is the rotation angle, and $l_a=r_{ab}l_b$ and $l_b$ are the lengths of the semi-major and semi-minor axes, respectively. In this test case, we randomly generate these parameters from
\[
x_0, y_0 \sim\mathcal{U}([0.4, 0.6]), \quad l_b \sim\mathcal{U}([0.1, 0.15]),\quad r_{ab} \sim\mathcal{U}([1.0, 1.5]), \quad \theta \sim\mathcal{U}([0, 2\pi]).
\]
Specifically, we sample $\mathbf{t}_N(x,y)\sim\mathcal{GP}(\mathbf{t}_c,k_\ell)$ at background-mesh nodes with $\mathbf{t}_c=[0,t_y]^T$, $t_y\sim\mathcal{U}([10,50])$, signal standard deviation $s=0.5$, and correlation length $\ell=0.2$. We set the Young's modulus and the Poisson ratio as $E=100$ and $\nu=0.3$.

\begin{figure}[t]
	\centering
	\begin{minipage}[t]{0.32\textwidth}
		\centering
		\makebox[\linewidth][l]{\textbf{(a)}}\par\vspace{0.3em}
		\includegraphics[width=\linewidth]{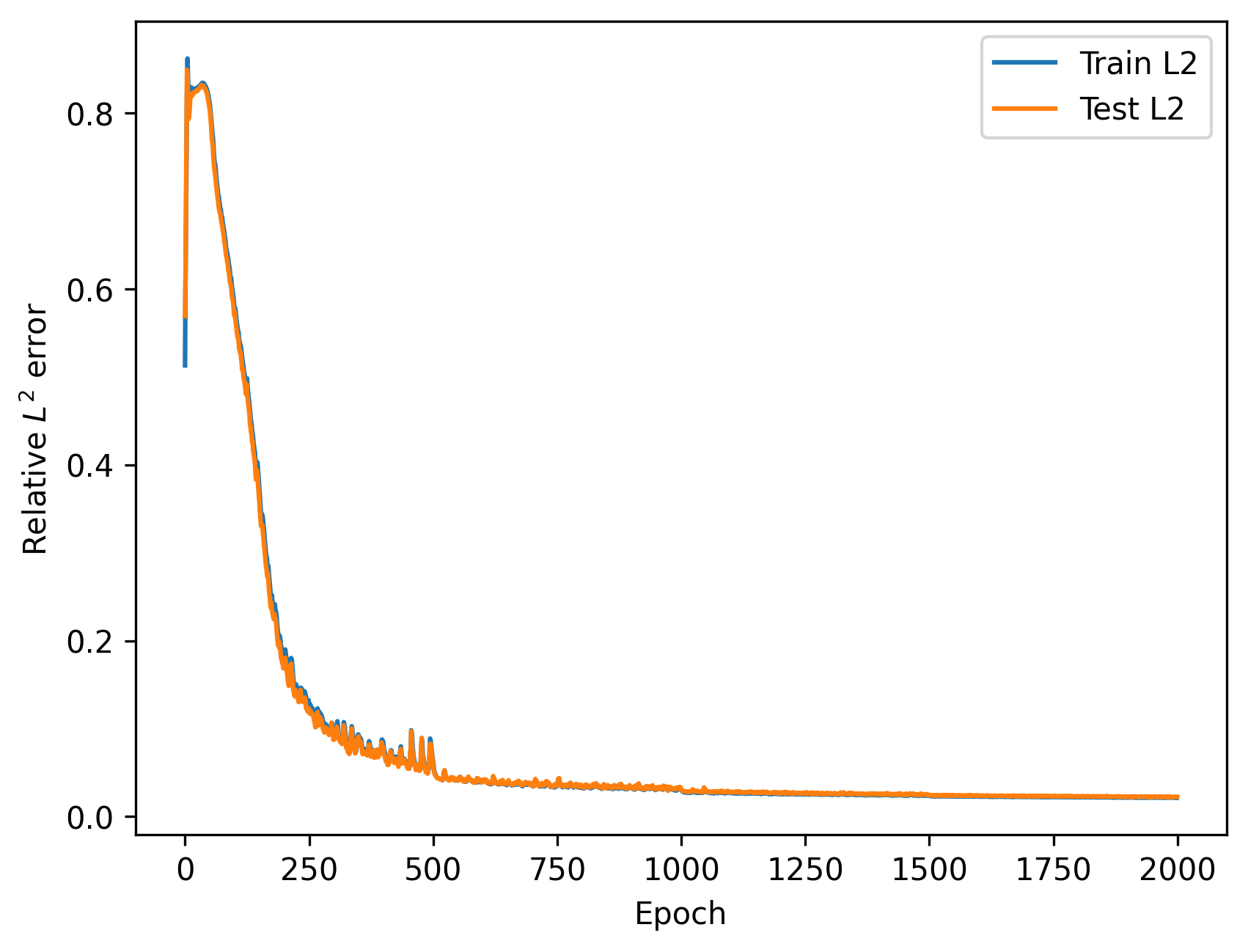}
	\end{minipage}
	\hfill
	\begin{minipage}[t]{0.32\textwidth}
		\centering
		\makebox[\linewidth][l]{\textbf{(b)}}\par\vspace{0.3em}
		\includegraphics[width=\linewidth]{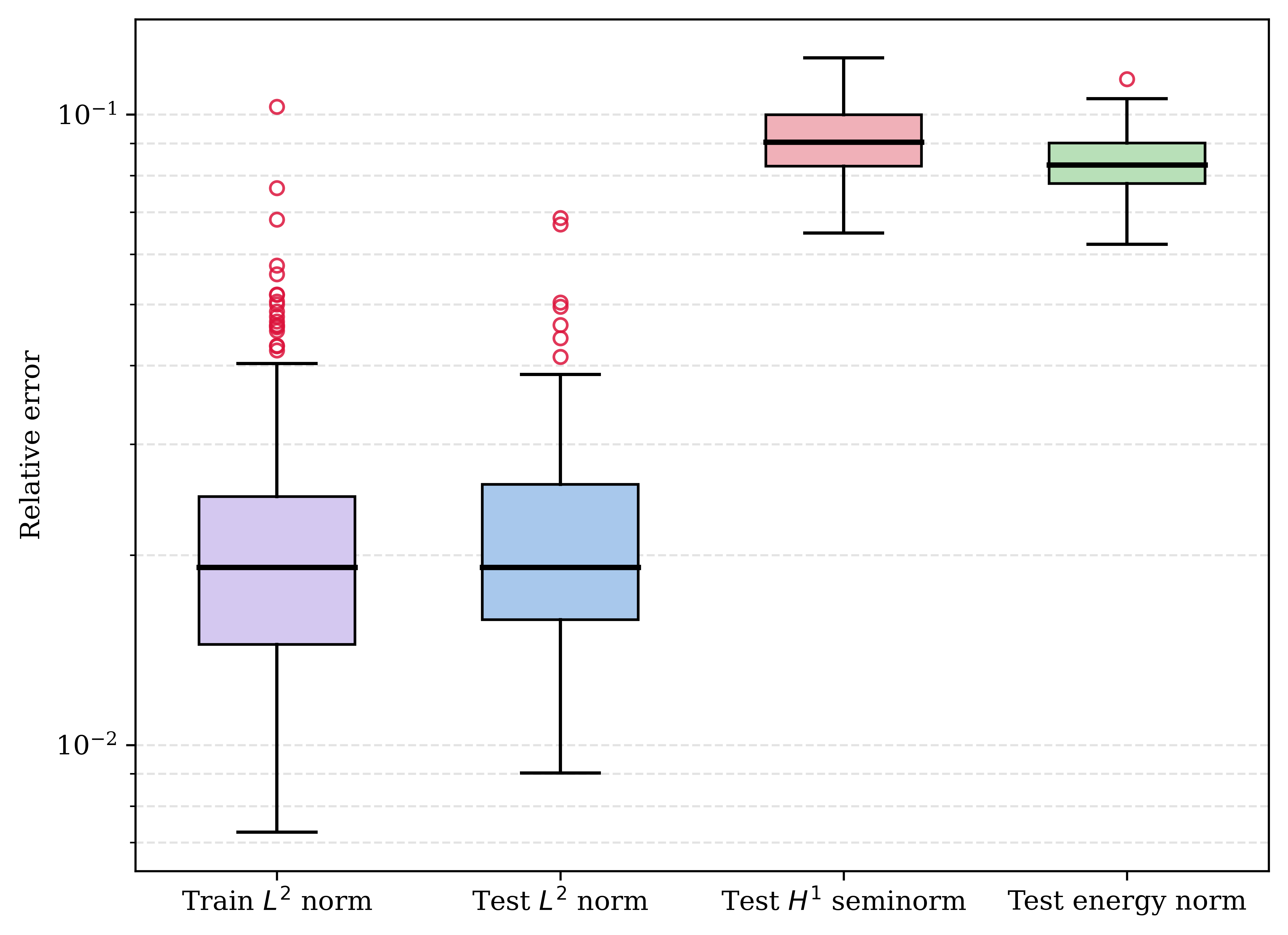}
	\end{minipage}
	\hfill
	\begin{minipage}[t]{0.328\textwidth}
		\centering
		\makebox[\linewidth][l]{\textbf{(c)}}\par\vspace{0.3em}
		\includegraphics[width=\linewidth]{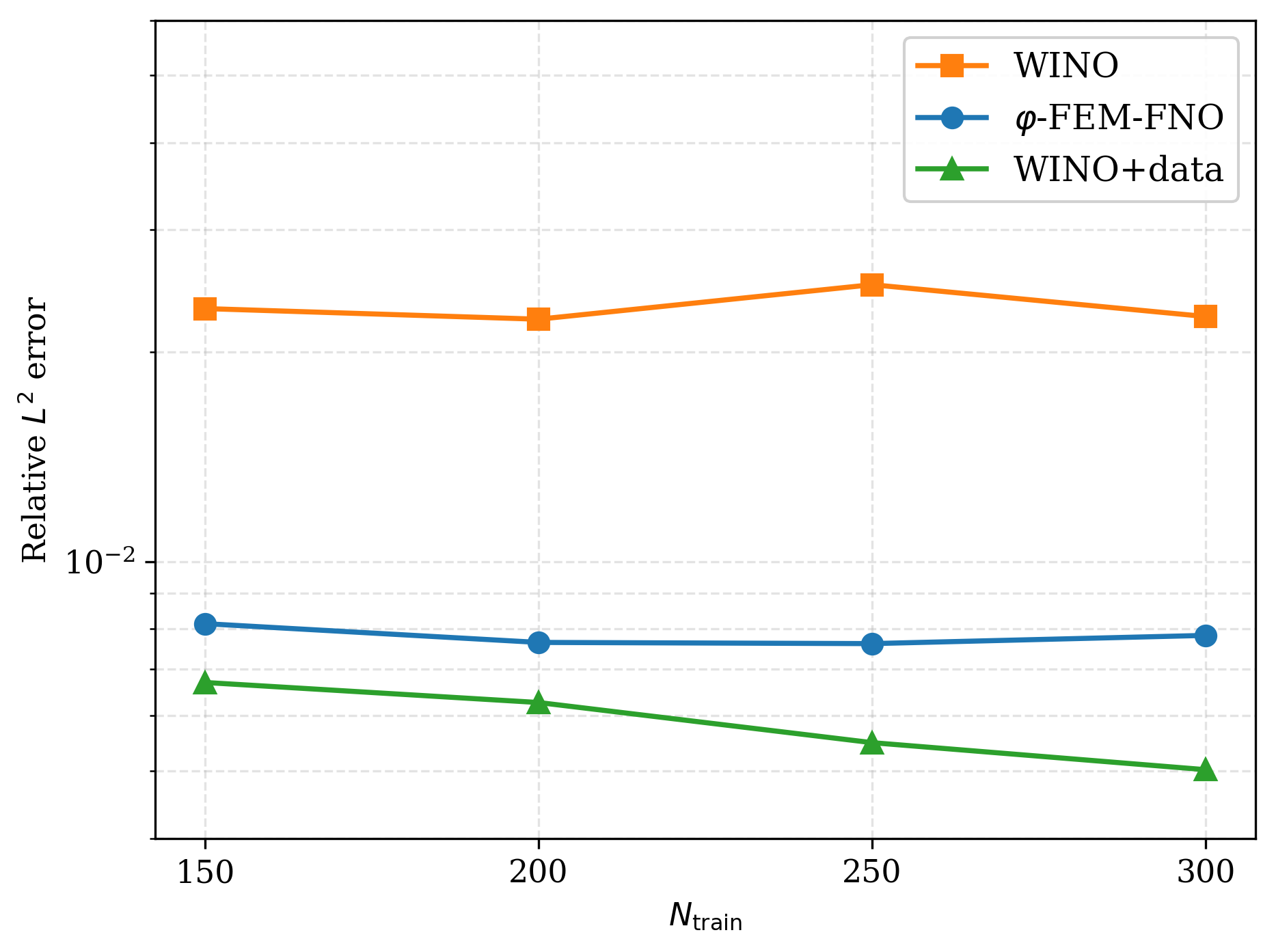}
	\end{minipage}
	\caption{Plate with a hole case. (a) Evolution of the mean relative $L^2$ errors during training; (b) Relative errors after training; (c) Influence of the number of training samples $N_{\mathrm{train}}$ on the relative $L^2$ error.}
	\label{fig:case_hole_training}
\end{figure}

\begin{table}[t]
	\centering
	\footnotesize
	\caption{Computational costs (data generation + training time) and relative errors of three methods for the plate with a hole case. Cost ratio is total time relative to WINO.}
	\label{tab:performance_hole}
	\begin{tabularx}{\linewidth}{@{}l >{\raggedright\arraybackslash}X >{\raggedright\arraybackslash}X >{\raggedright\arraybackslash}X >{\raggedright\arraybackslash}X >{\raggedright\arraybackslash}X@{}}
		\toprule
		Method & Total time & Cost ratio & $\|e\|_{L^2}$ & $\|e\|_{H^1}$ & $\|e\|_{E}$ \\
		\midrule
		WINO & $0.3 + 2912.8$s & $1.00\times$ & $2.26 \pm 1.10\%$ & $9.12 \pm 1.26\%$ & $8.39 \pm 0.99\%$ \\
		$\varphi$-FEM-FNO & $11635.8 + 1973.6$s & $4.67\times$ & $0.78 \pm 0.40\%$ & $15.3 \pm 6.27\%$ & $15.7 \pm 7.15\%$ \\
		WINO+data & $11635.8 + 2976.8$s & $5.02\times$ & $0.54 \pm 0.20\%$ & $8.12 \pm 1.18\%$ & $7.50 \pm 1.03\%$ \\
		\bottomrule
	\end{tabularx}
\end{table}

In this test case, the level-set function $\varphi$ is used only to represent the Neumann boundary along the elliptic hole; the outer boundary of the square plate lies on edges of the Cartesian background mesh and is therefore treated with standard conforming finite elements, without invoking the cut-cell $\varphi$-FEM boundary machinery on those edges. In WINO and $\varphi$-FEM-FNO, the homogeneous Dirichlet condition on $\Gamma_D$ is enforced by the multiplicative ansatz $\mathbf{u}_h=\mathbf{u}_{\theta}\,y$. The loss function of WINO is computed by \eqref{eq:total_loss} with the penalty parameters $\lambda_1=1\times10^{-3}, \lambda_2=1\times10^{-7}, \lambda_3=1\times10^{-6}$. In this case, we learn the mapping
\begin{equation}
	\mathcal{G}_{\theta}:(\varphi_h,\mathbf{t}_h)\mapsto(\mathbf{u}_h,\mathbf{y}_h,\mathbf{p}_h).
\end{equation}
We compute the loss function on a $64 \times 64$ mesh and train for 1{,}000 epochs using the SOAP optimizer with 300 training samples and 100 test samples. The ground-truth fields are generated by $\varphi$-FEM on the same mesh using second-order finite elements (biquadratic Q9 elements on quadrilateral meshes) by Newton iterations and GMRES methods in each Newton step.

\begin{figure}[t]
	\centering
	\begin{minipage}[t]{0.30\textwidth}
		\vspace{0pt}% anchor [t] to true top (avoids large gap above first graphic)
		\centering
		\makebox[\linewidth][l]{\textbf{(a)}}\par\vspace{0.3em}
		\includegraphics[width=0.9\linewidth]{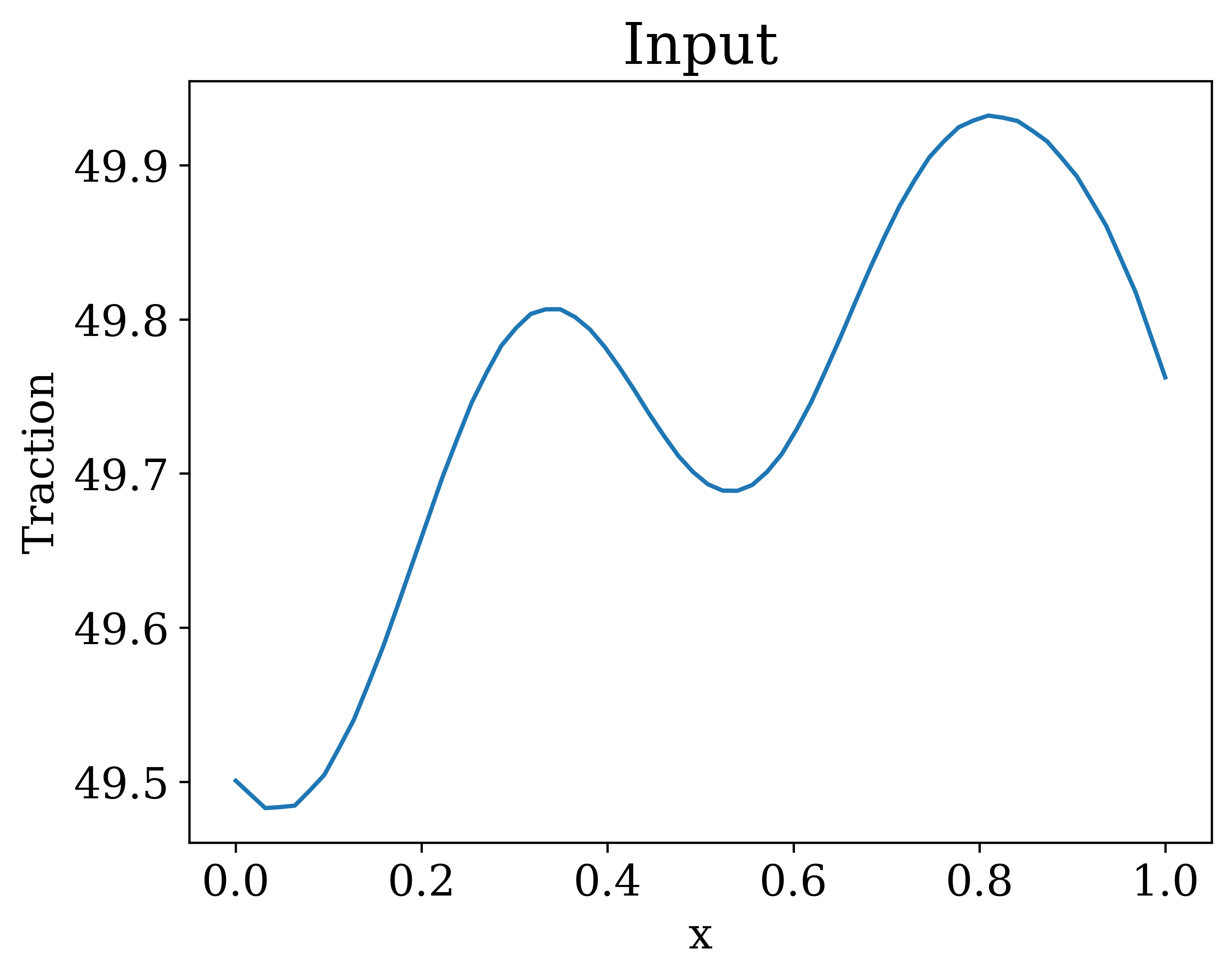}\\[0.6ex]
		\includegraphics[width=0.85\linewidth]{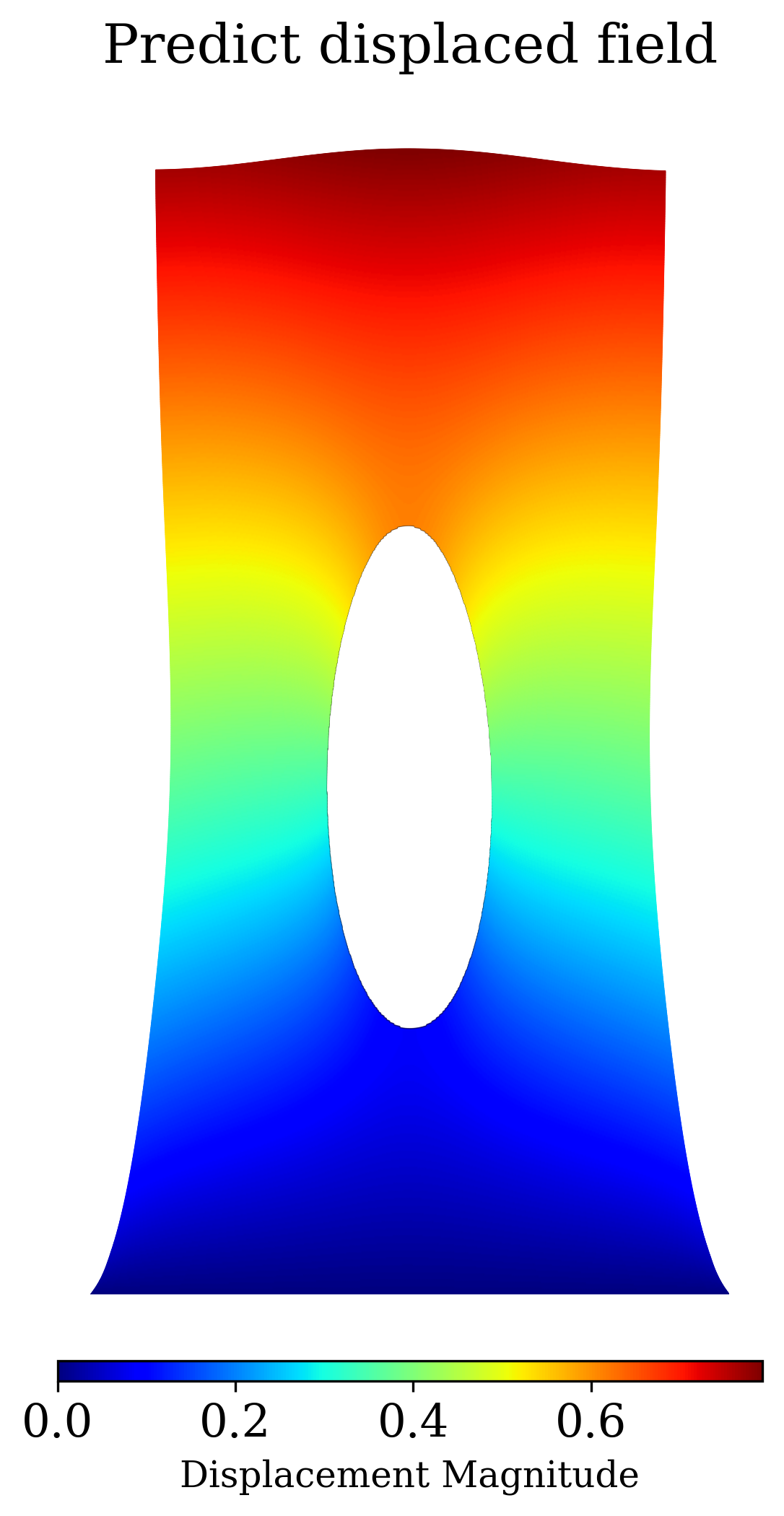}
	\end{minipage}%
	\hspace{0.08\textwidth}%
	\begin{minipage}[t]{0.6\textwidth}
		\vspace{0pt}% align top with left column
		\centering
		\makebox[\linewidth][l]{\textbf{(b)}}\par\vspace{0.3em}
		\begin{tabular}{@{}c@{\hspace{0.02\linewidth}}c@{}}
			\includegraphics[width=0.45\linewidth]{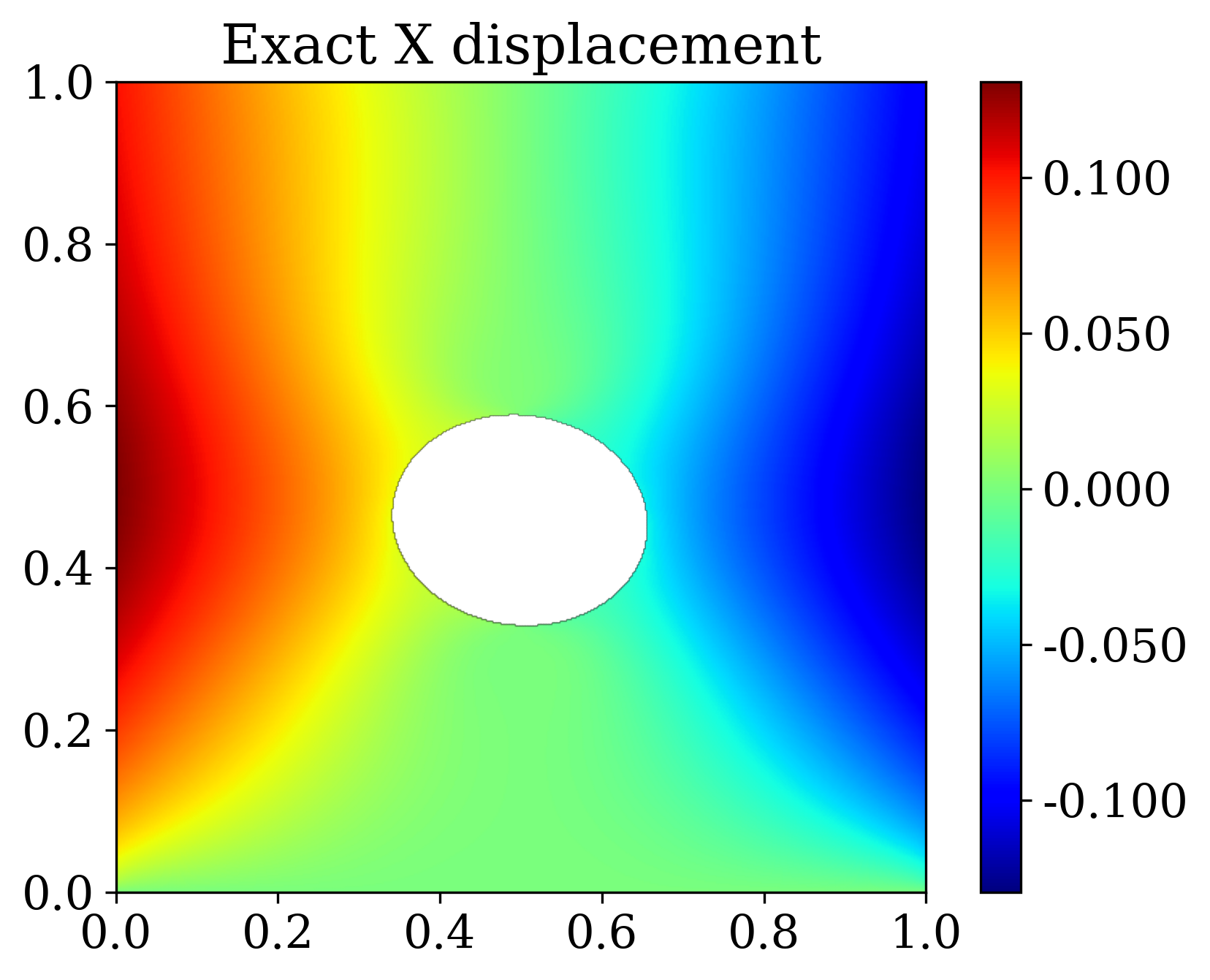} &
			\includegraphics[width=0.45\linewidth]{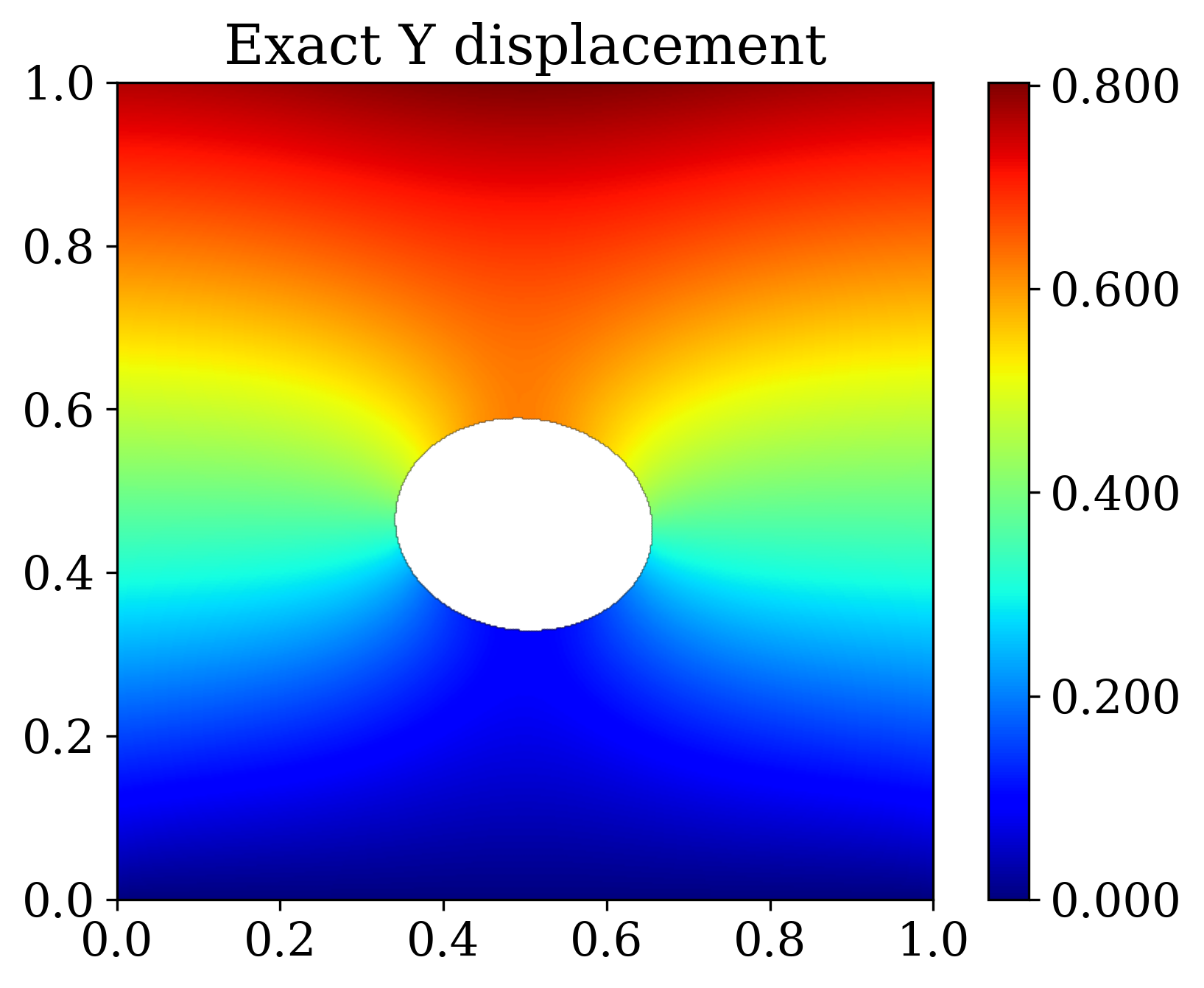} \\[0.5ex]
			\includegraphics[width=0.45\linewidth]{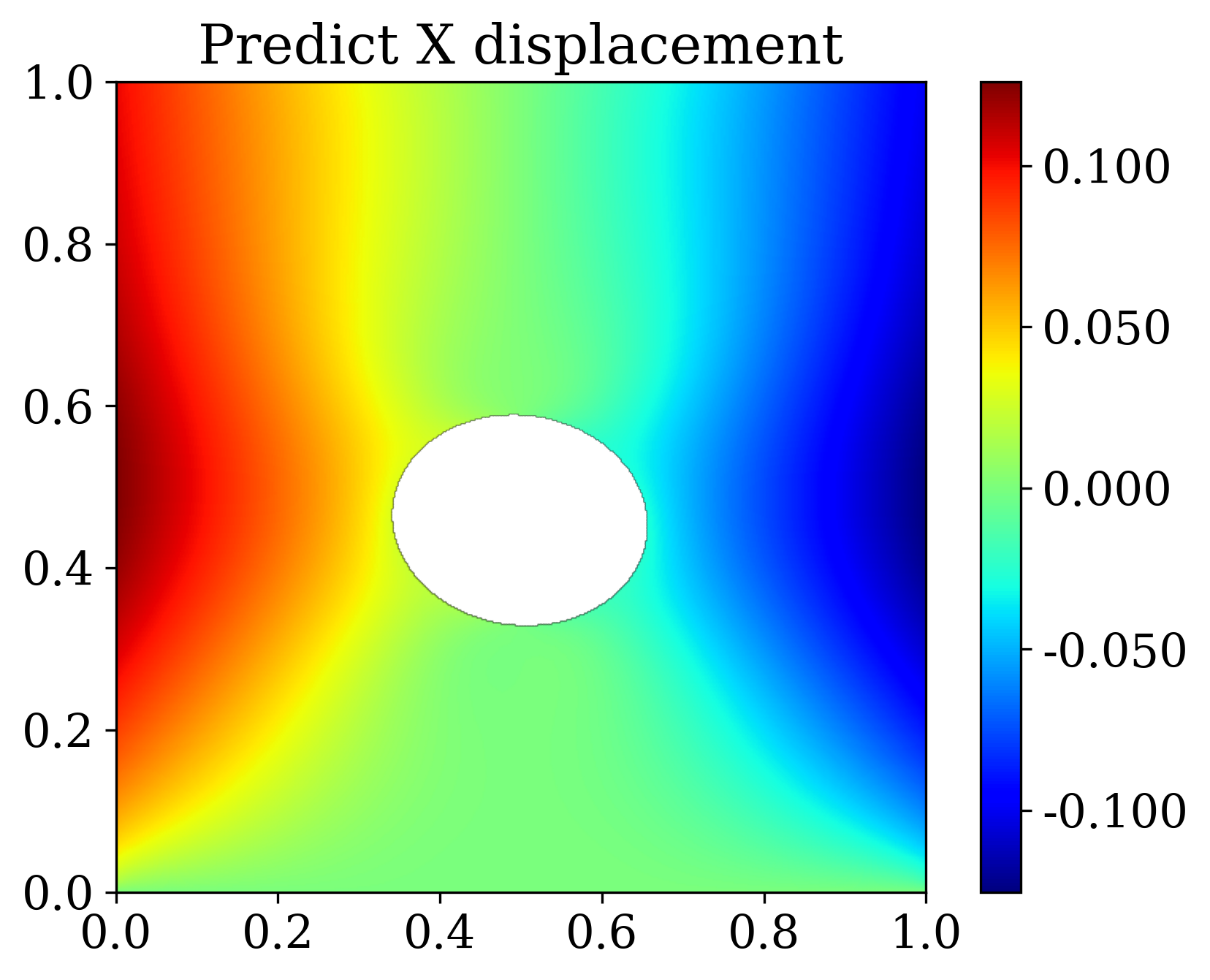} &
			\includegraphics[width=0.45\linewidth]{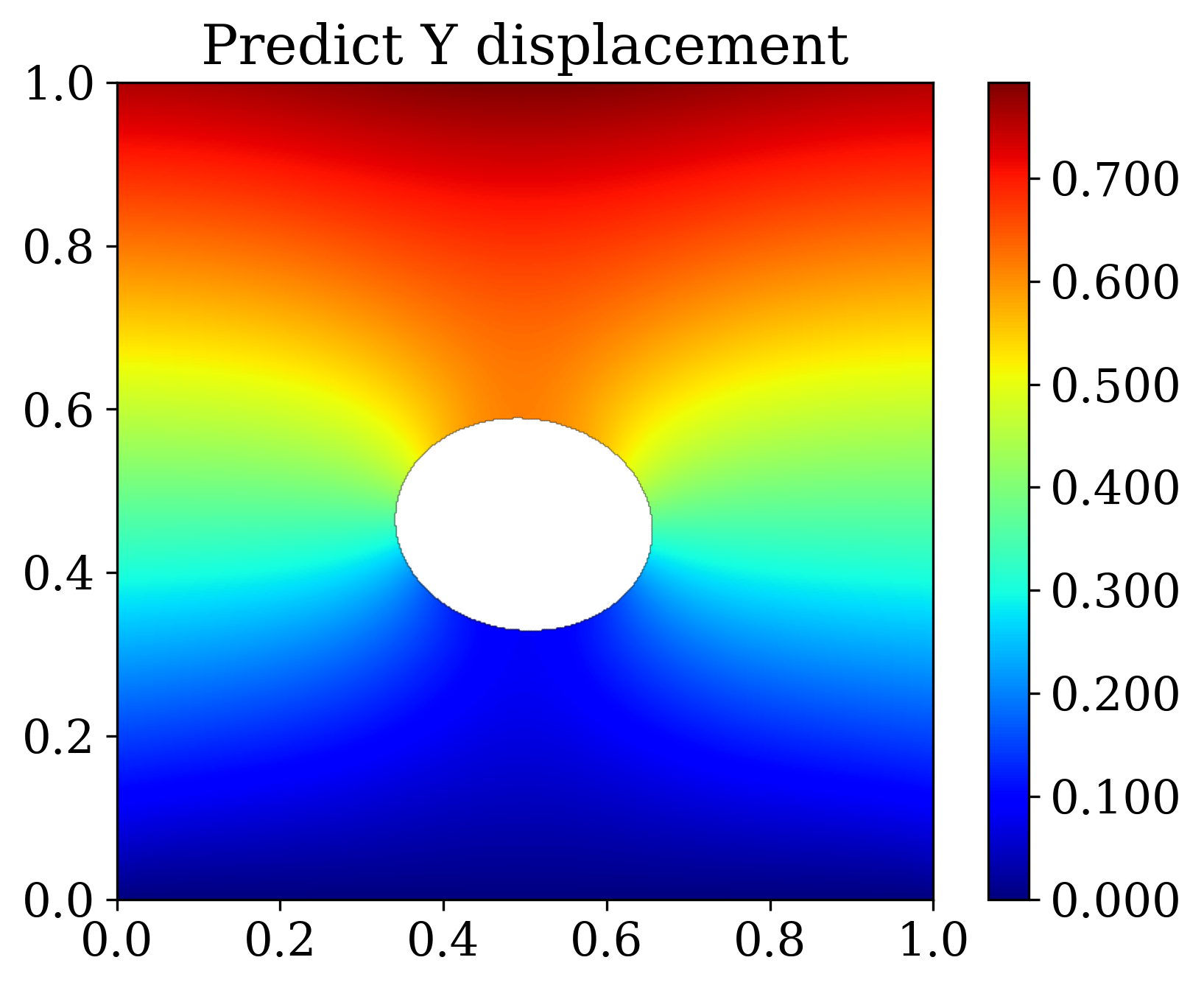} \\[0.5ex]
			\includegraphics[width=0.45\linewidth]{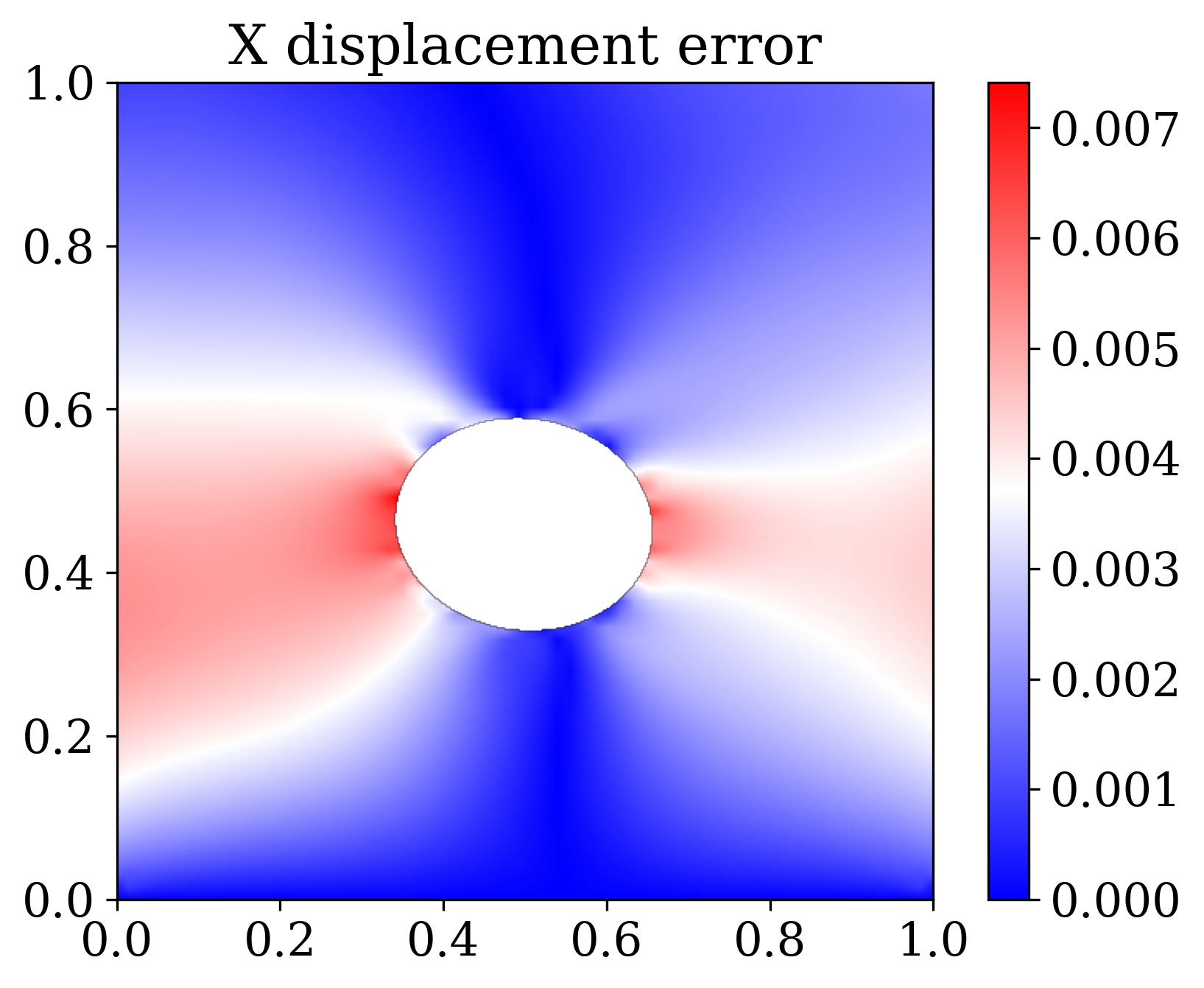} &
			\includegraphics[width=0.45\linewidth]{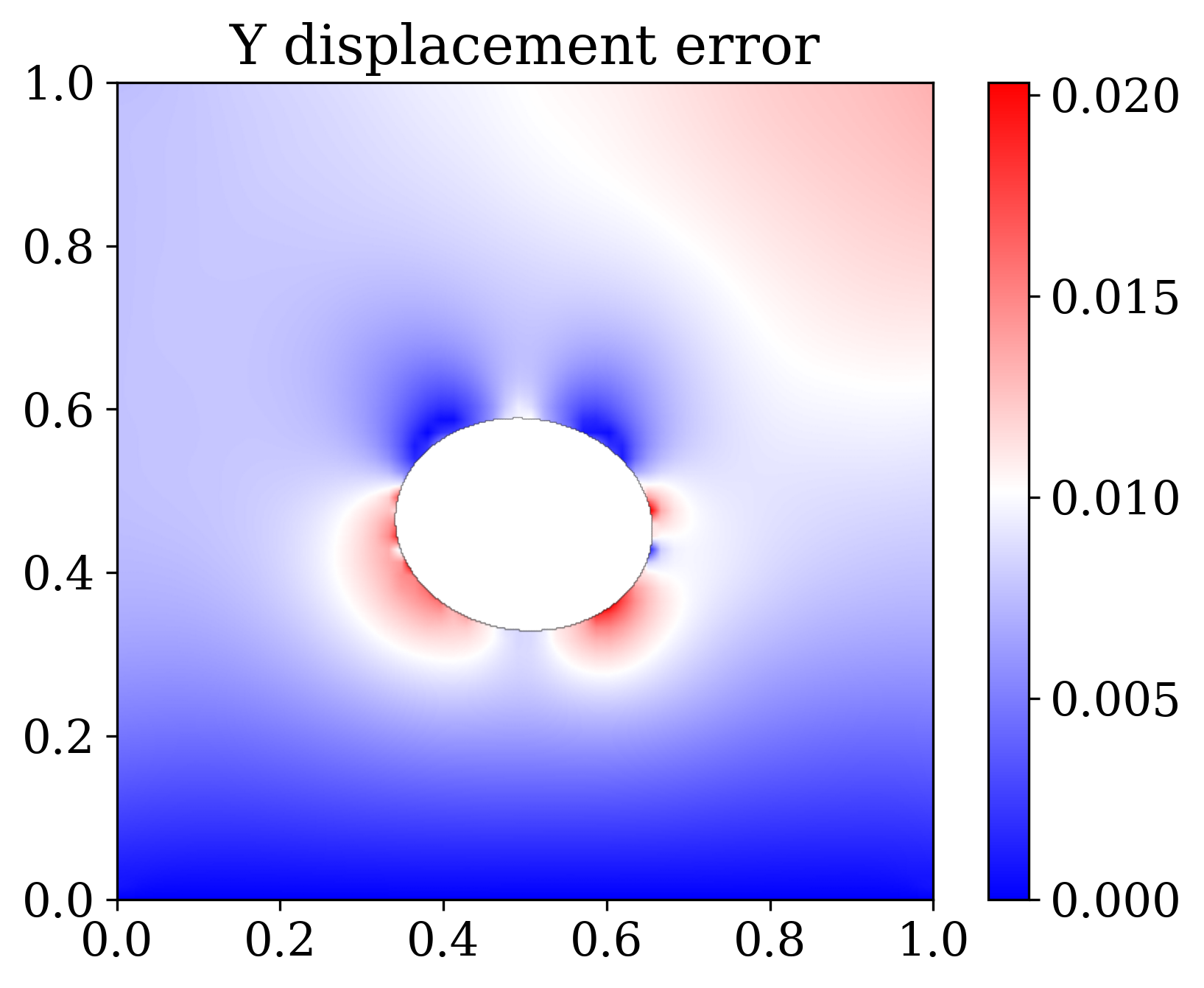}
		\end{tabular}
	\end{minipage}
	\caption{Representative test sample for the plate with a hole case. (a) Input boundary condition $\mathbf{t}_N$ and WINO-predicted deformed configuration; (b) $x$- and $y$-displacement components for the exact solutions (top row), WINO predictions (middle row), and pointwise absolute error contours (bottom row).}
	\label{fig:case_hole_displacement}
\end{figure}

\begin{figure}[t]
	\raggedright
	\makebox[\linewidth][l]{\textbf{(a)}}\par\vspace{0.3em}
	{\centering
	\includegraphics[width=0.22\textwidth]{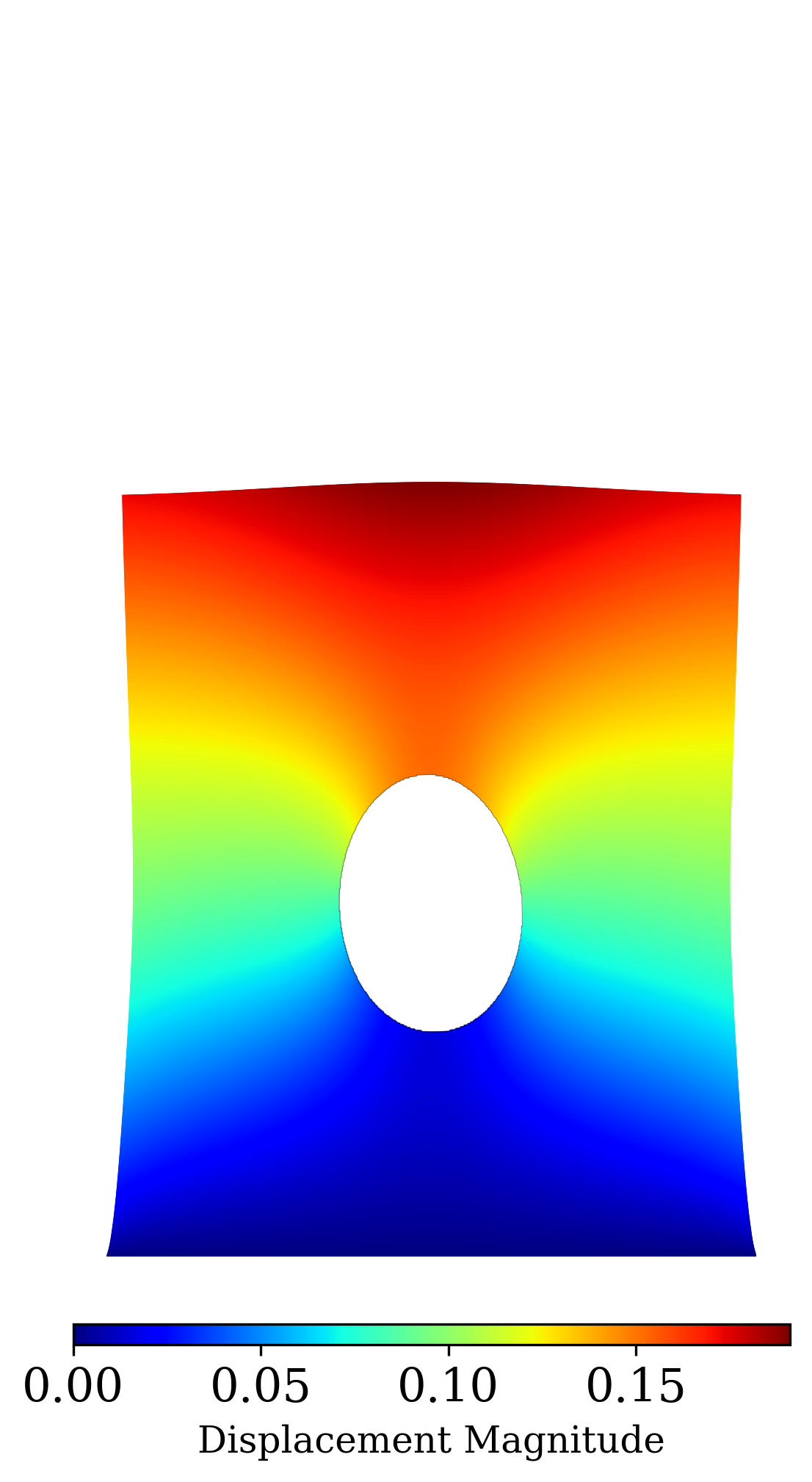}
	\hspace{0.02\textwidth}
	\includegraphics[width=0.22\textwidth]{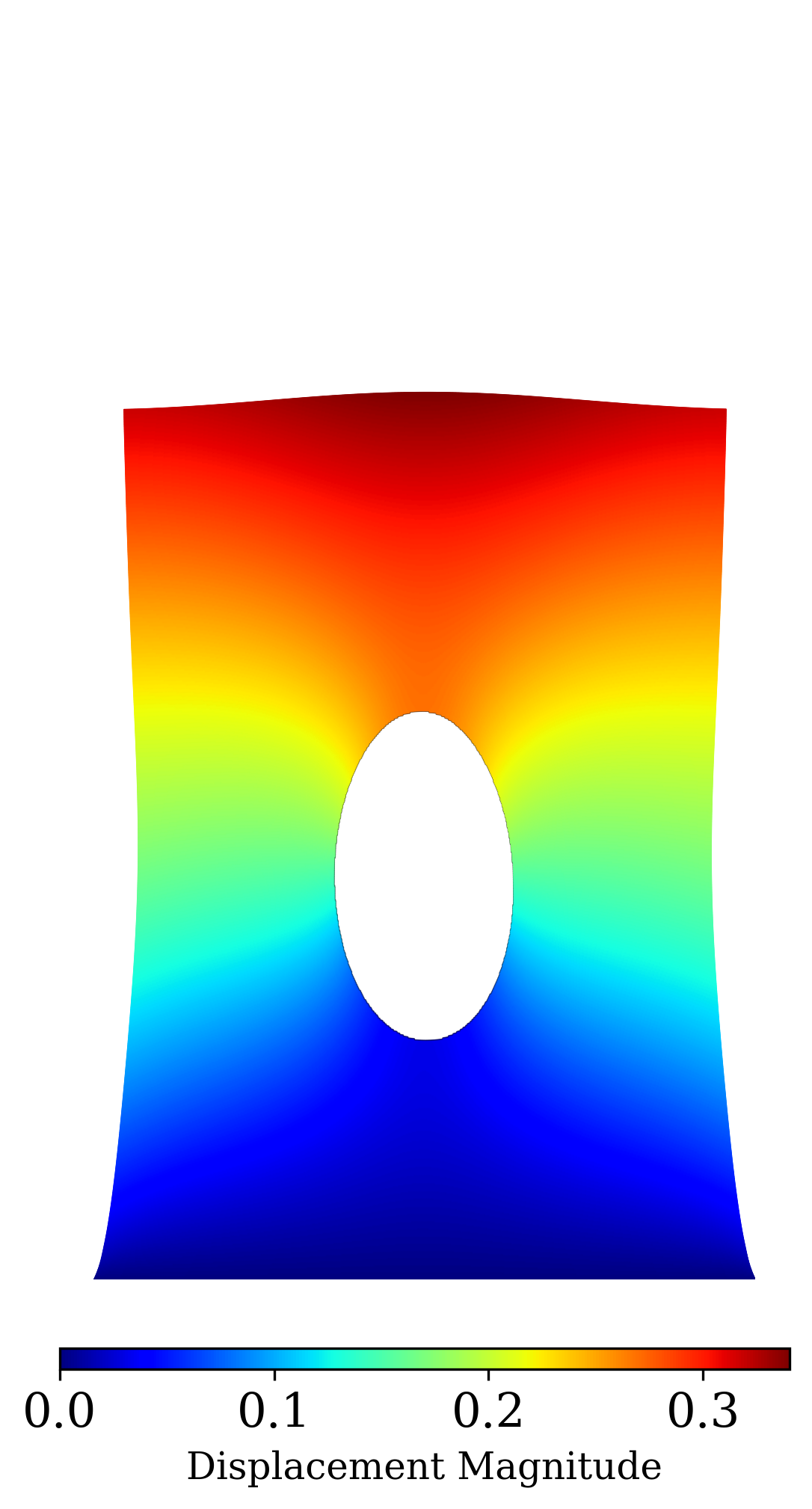}
	\hspace{0.02\textwidth}
	\includegraphics[width=0.22\textwidth]{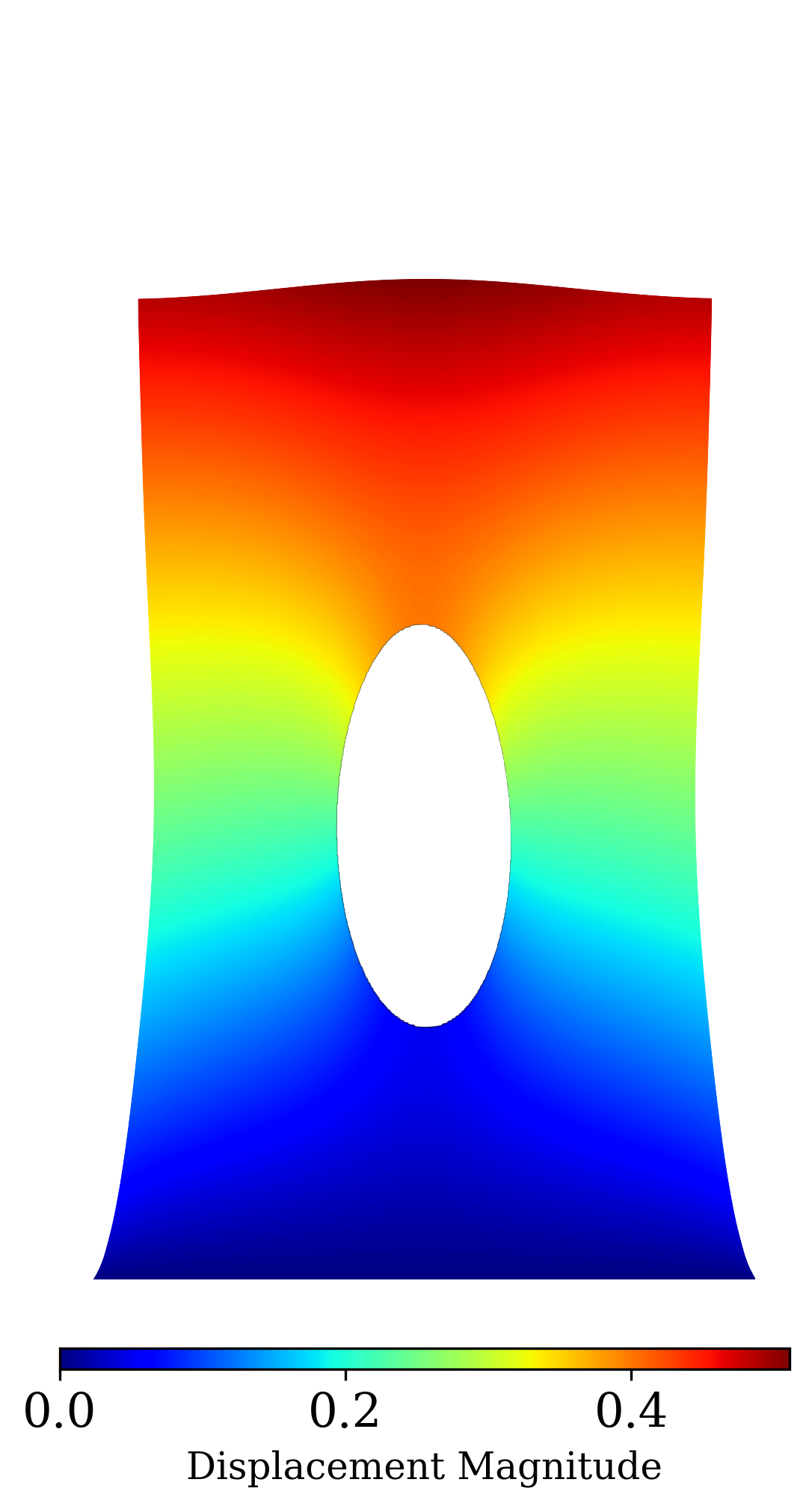}
	\hspace{0.02\textwidth}
	\includegraphics[width=0.22\textwidth]{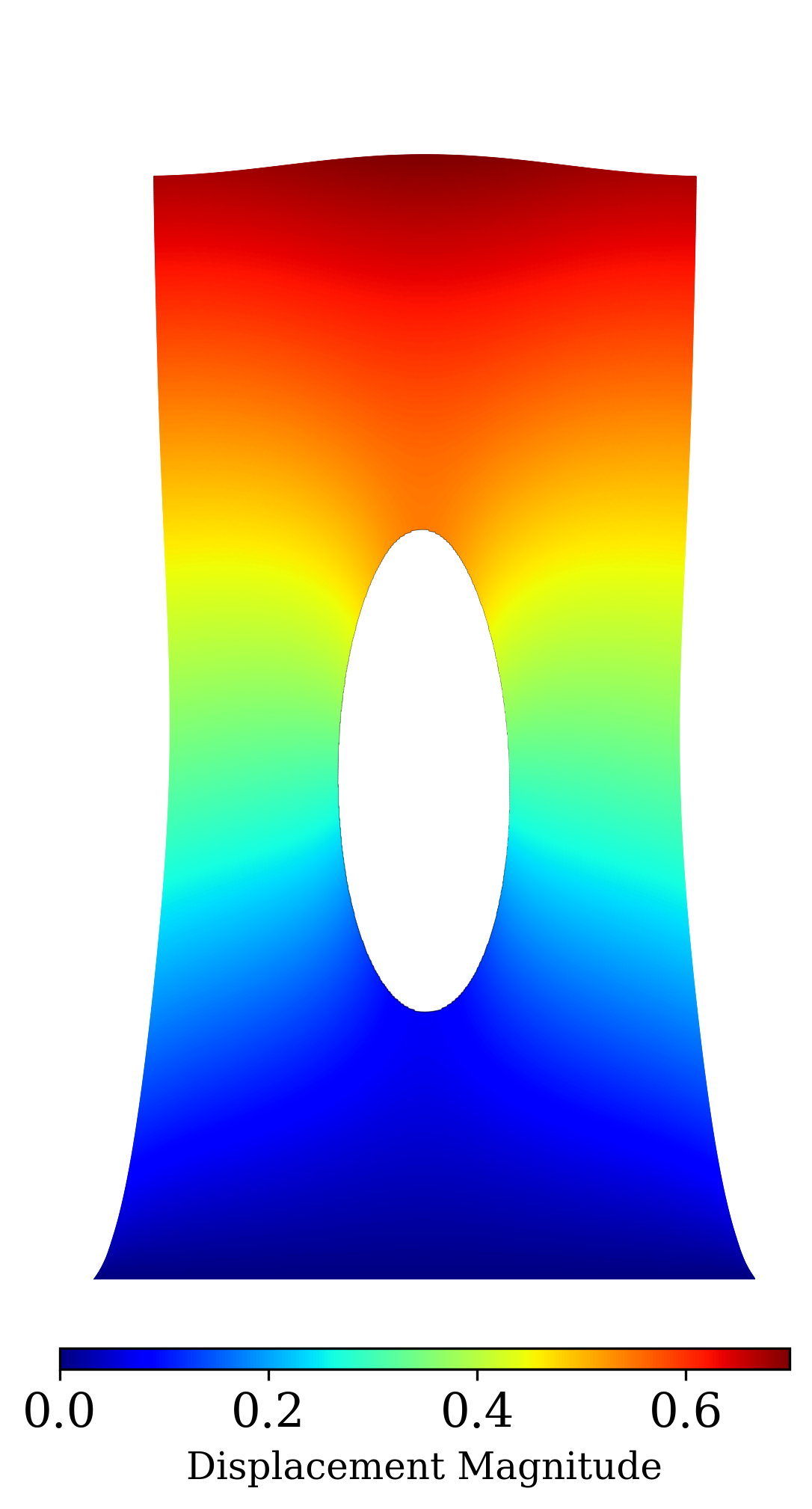}\par}
	\vspace{0.6em}

	\makebox[\linewidth][l]{\textbf{(b)}}\par\vspace{0.3em}
	{\centering
	\includegraphics[width=0.3\textwidth]{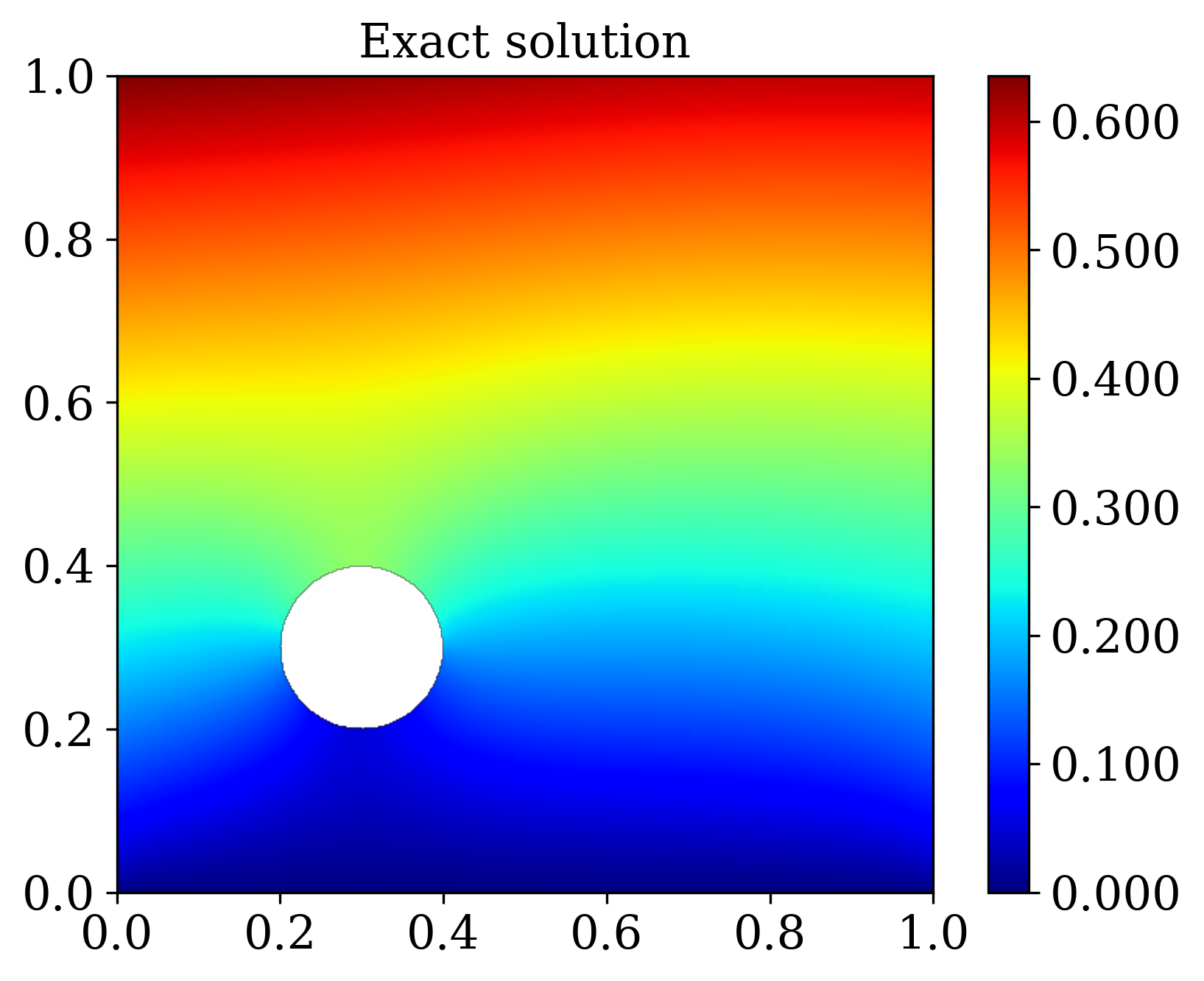}
	\hspace{0.02\textwidth}
	\includegraphics[width=0.3\textwidth]{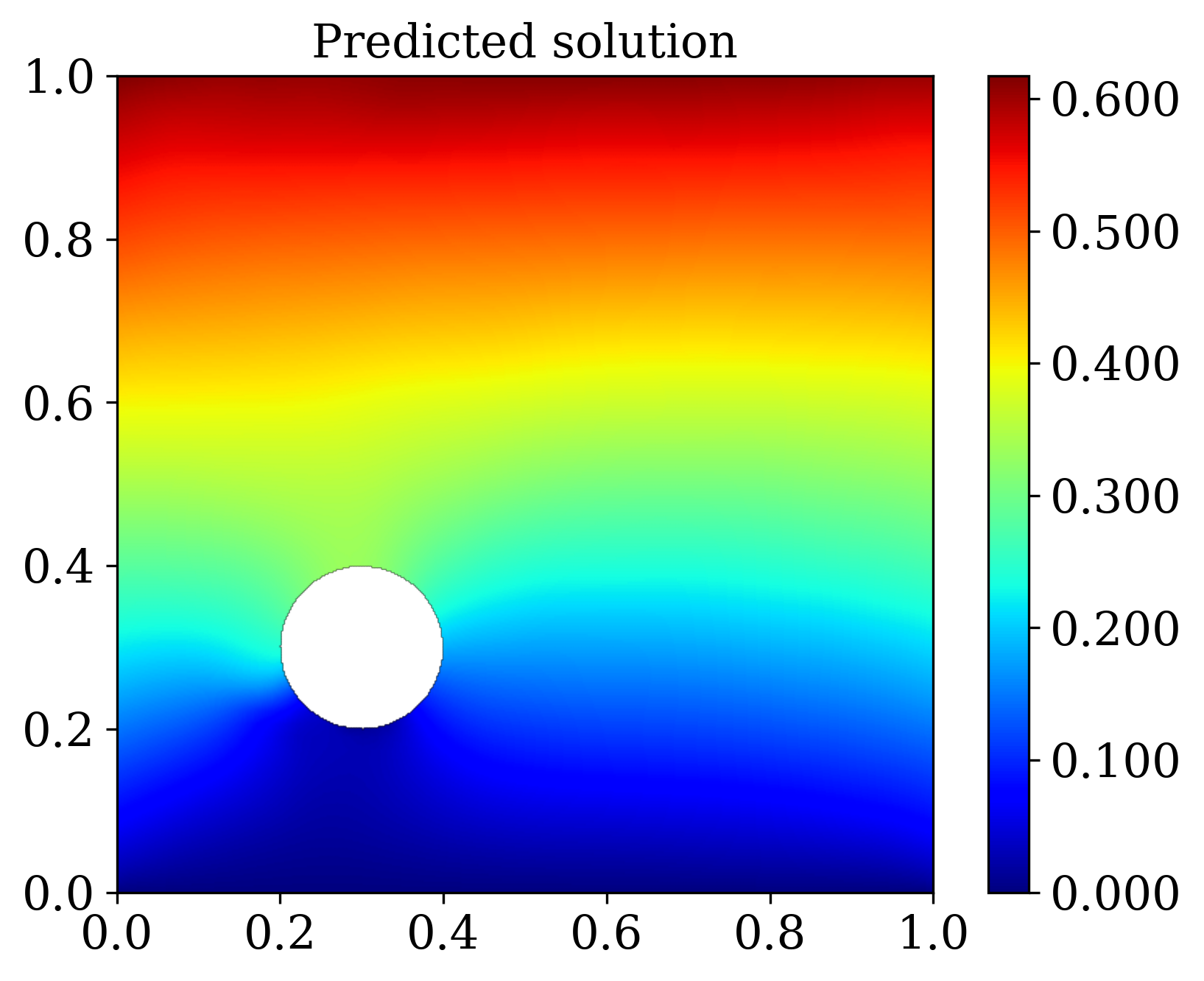}
	\hspace{0.02\textwidth}
	\includegraphics[width=0.3\textwidth]{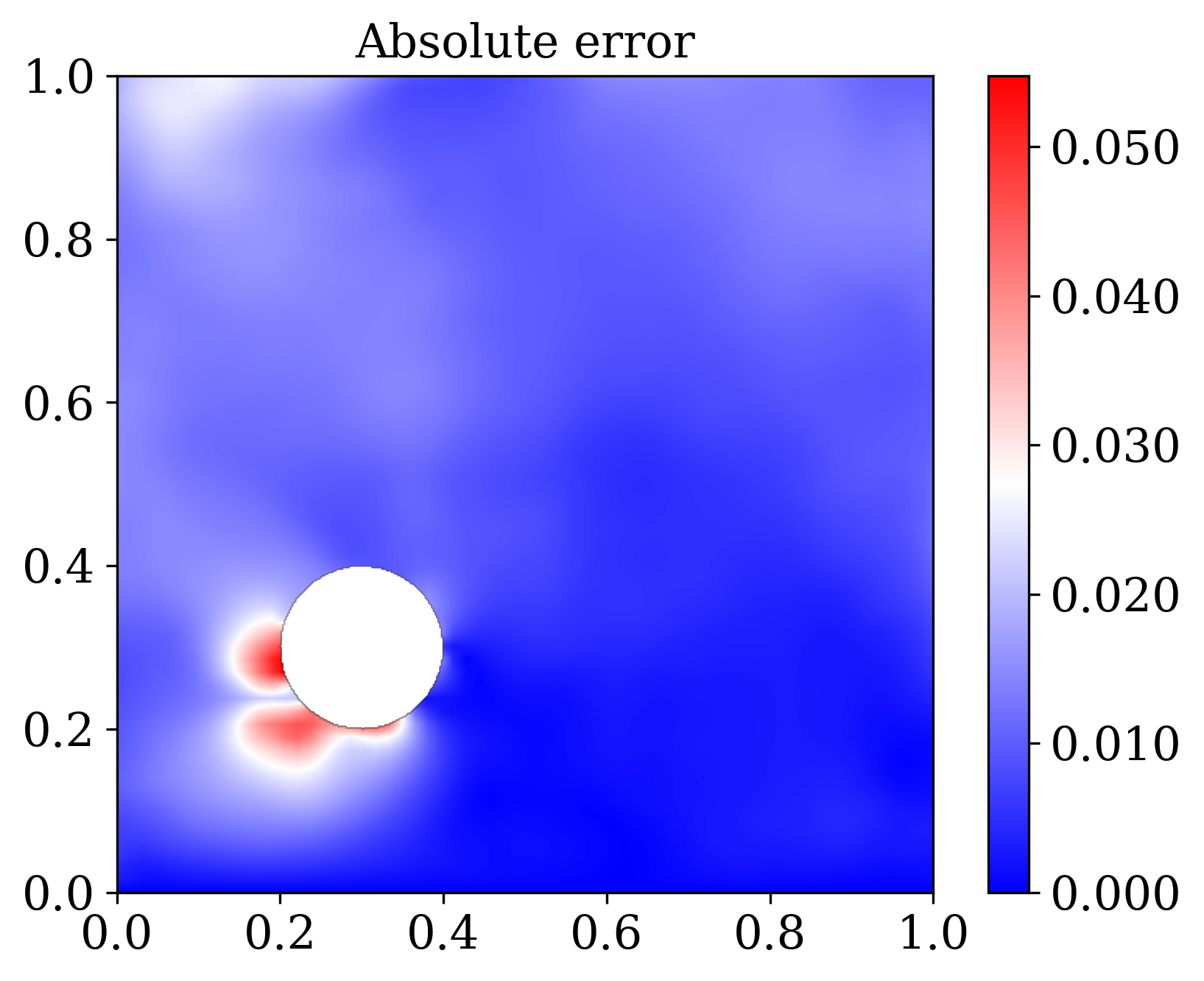}\par}
	\caption{WINO generalization for the plate with a hole case. (a) WINO-predicted displacement fields for increasing nominal traction levels $t_y\in\{15,25,35,45\}$; (b) Out-of-distribution test: exact solution (left), WINO prediction (middle), and absolute error contour (right).}
	\label{fig:case_hole_generalization}
\end{figure}

Fig.~\ref{fig:case_hole_training} reports the mean relative $L^2$ errors on the training and test sets: Fig.~\ref{fig:case_hole_training}a shows the evolution of the training and validation errors during optimization, and Fig.~\ref{fig:case_hole_training}b shows box plots of the errors after convergence. In each box, the horizontal line marks the median and the box spans the interquartile range. WINO attains satisfactory accuracy after a modest number of epochs, and the test-set error remains stable throughout training. Fig.~\ref{fig:case_hole_displacement} displays a representative WINO-predicted displacement field from the test set. The WINO solution agrees closely with the high-fidelity reference and accurately reproduces the displacement induced by the random traction loading. Since each training data uses a spatially GRF mean $\mathbf{t}_c=[0,t_y]^T$ with $t_y\sim\mathcal{U}([10,50])$, the learned map spans a wide range of nominal tractions. As a result, we interrogate the trained operator by evaluating WINO on a monotonically increasing sequence of traction levels and inspecting how the predicted deformation evolves. Fig.~\ref{fig:case_hole_generalization}a illustrates such a sweep for $t_y\in\{15,25,35,45\}$.

To further illustrate the computational efficiency of WINO, we compare it with the purely data-driven $\varphi$-FEM-FNO and with WINO augmented with labeled data. $\varphi$-FEM-FNO and the labeled-data variant of WINO use the same train/test splits, network architecture, and hyperparameters. For each method we record the total computational cost, including reference-data generation and training. Relative errors and these timings are reported in Table~\ref{tab:performance_hole}. Although $\varphi$-FEM-FNO attains a smaller relative $L^2$ error than label-free WINO, its reference-data generation cost is substantial (average 38.78s for one $\varphi$-FEM generation). In contrast, the data generation time of WINO is much lower since it doesn't need to compute the reference solution, and the weak-form physics regularization is reflected in substantially smaller relative $H^1$ and energy errors than for $\varphi$-FEM-FNO. Moreover, when labeled data are available, incorporating them into WINO yields further accuracy gains. Inference cost is similar for all three models (about $5\,\mathrm{ms}$ per forward pass in our setup), reflecting the fact that they share essentially the same network architecture. We also performed ten independent WINO training runs with different random seeds on the same training and test datasets. The seed-averaged mean relative $L^2$ error is $2.35\times10^{-2}$, and the $5$--$95\%$ percentile range is $[2.22\times10^{-2},\,2.45\times10^{-2}]$, which further supports the stability and reproducibility of the proposed training procedure. For data-driven methods, the data-generation time scales linearly with the number of training samples. We therefore test the influence of the training-set size on the three methods in Fig.~\ref{fig:case_hole_training}c. The relative $L^2$ errors remain essentially unaffected over the range considered (the batch size is set to $25$ for $N_{\mathrm{train}}=150,200$ and to $50$ for $N_{\mathrm{train}}=250,300$). Consequently, to ensure a fair cost comparison among the three methods, we use the same number of training samples for all of them.

We further examine the out-of-distribution (OOD) generalization of WINO. We construct auxiliary geometries featuring a circular hole near the lower-left corner of the plate (Fig.~\ref{fig:case_hole_generalization}b) and apply a spatially constant traction with vertical component $t_y=45$, i.e., a deterministic loading distinct from the random GRF-sampled tractions used in training. Together, the unseen hole placement/shape and this fixed loading define an OOD test regime relative to the training distribution. WINO still yields accurate predictions, with a relative $L^2$ error of $3.24\%$. The error is slightly larger than on in-distribution test cases, as expected under distribution shift, which suggests that the learned operator retains useful OOD robustness.

\begin{figure}[t]
	\centering
	\begin{minipage}[t]{0.59\textwidth}
		\centering
		\makebox[\linewidth][l]{\textbf{(a)}}\par\vspace{0.3em}
		\includegraphics[width=\linewidth]{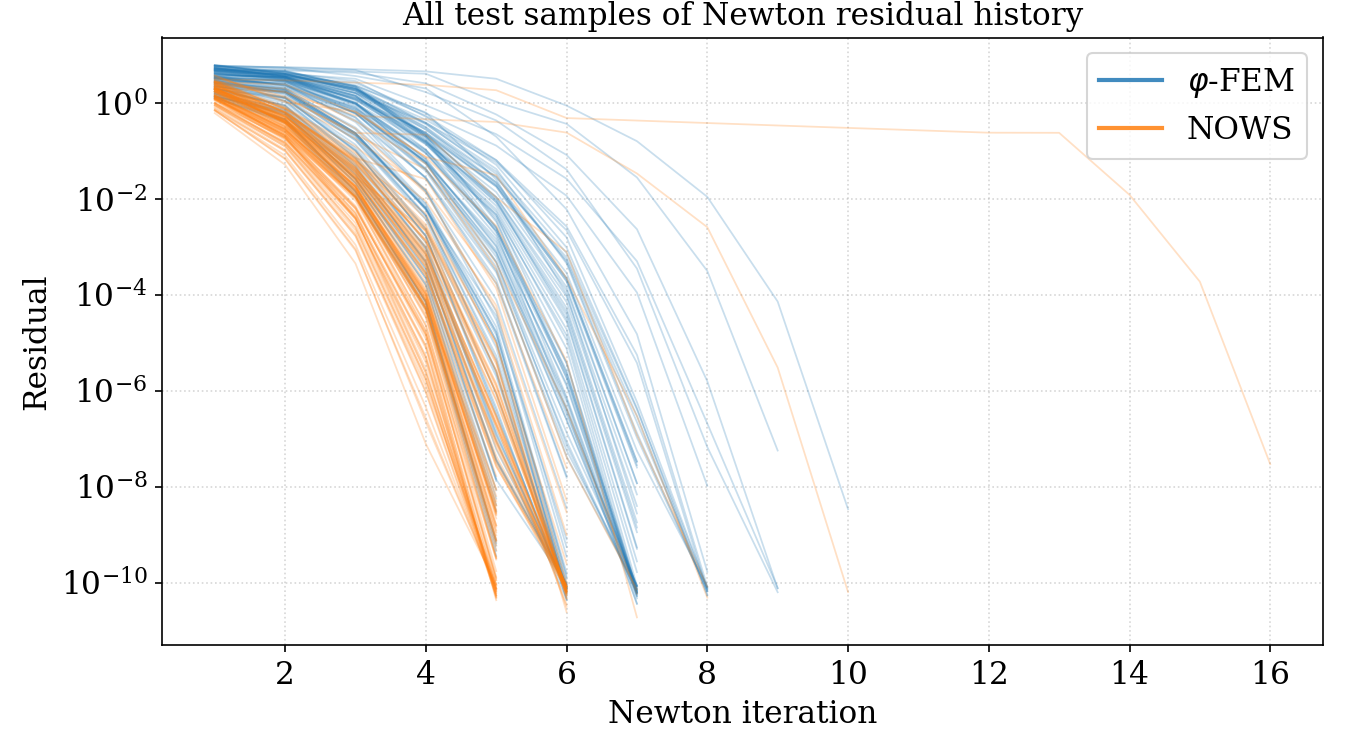}
	\end{minipage}
	\hfill
	\begin{minipage}[t]{0.4\textwidth}
		\centering
		\makebox[\linewidth][l]{\textbf{(b)}}\par\vspace{0.3em}
		\includegraphics[width=\linewidth]{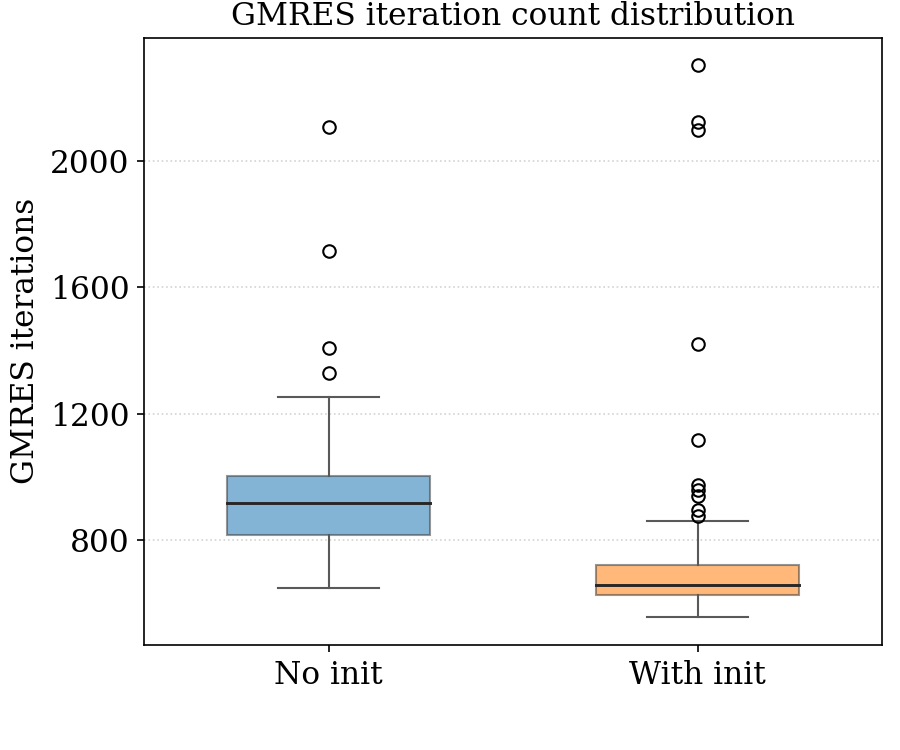}
	\end{minipage}
	\caption{NOWS comparison for the plate with a hole case. (a) Outer Newton residual history for $\varphi$-FEM and NOWS. (b) Box plots of the total GMRES iteration counts.}
	\label{fig:case_hole_nows}
\end{figure}

The WINO prediction can also be used as an initial guess to accelerate the nonlinear $\varphi$-FEM solve, i.e., the neural-operator warm-start (NOWS) strategy summarized in Section~\ref{sec:nows}. Hyperelasticity is governed by a nonlinear residual. Each Newton step linearizes that residual and yields a sparse tangent system, which is solved by restarted GMRES (restart $m=80$, relative tolerance $\mathrm{rtol}=10^{-5}$) in every Newton iteration in this case. In practice, the $64\times 64$ discrete fields returned by WINO, including the displacement $\mathbf{u}_h$, the auxiliary variables $\mathbf{y}_h$ and $\mathbf{p}_h$ from the $\varphi$-FEM formulation, are injected into the corresponding $\varphi$-FEM degrees of freedom and furnish the initial Newton iterate in the nonlinear solver. For consistency, $\varphi$-FEM is discretized in the same first-order finite-element space (bilinear $Q4$ elements on the Cartesian background mesh), so the neural outputs align nodewise with the solver unknowns without cross-mesh projection. Fig.~\ref{fig:case_hole_nows}a compares the outer Newton residual convergence of $\varphi$-FEM with a WINO-based warm start against a cold start, and Fig.~\ref{fig:case_hole_nows}b reports box plots of the total GMRES iteration counts over all Newton steps. NOWS consistently reduces both Newton iterations and total GMRES iterations, which directly lowers the overall computational cost. A single WINO forward pass costs only about $5\,\mathrm{ms}$ in our setup, which is negligible compared to the cost of subsequent finite element iterations.

\begin{remark}
	A natural question is whether a reverse-mapping strategy, commonly used in fluid-structure interaction, can be incorporated here. Although we explored this possibility, it was not adopted for two practical reasons: First, in the present setting the level-set field $\varphi$ serves only as a local descriptor of the Neumann boundary segment, rather than a global representation of the full geometry, making it difficult to define a consistent reverse map over the entire domain; Second, the relative $H^1$ error is significantly larger than the $L^2$ error, which leads to pronounced errors in the inverse deformation gradient $\mathbf{F}^{-1}$, and this sensitivity adversely affects the accuracy and robustness of the reverse-mapping procedure.
\end{remark}

\subsubsection{Cook's membrane}\label{sec:cooks_membrane}

\begin{figure}[t]
	\centering
	\begin{minipage}[t]{0.38\textwidth}
		\centering
		\makebox[\linewidth][l]{\textbf{(a)}}\par\vspace{0.3em}
		\includegraphics[width=0.9\linewidth]{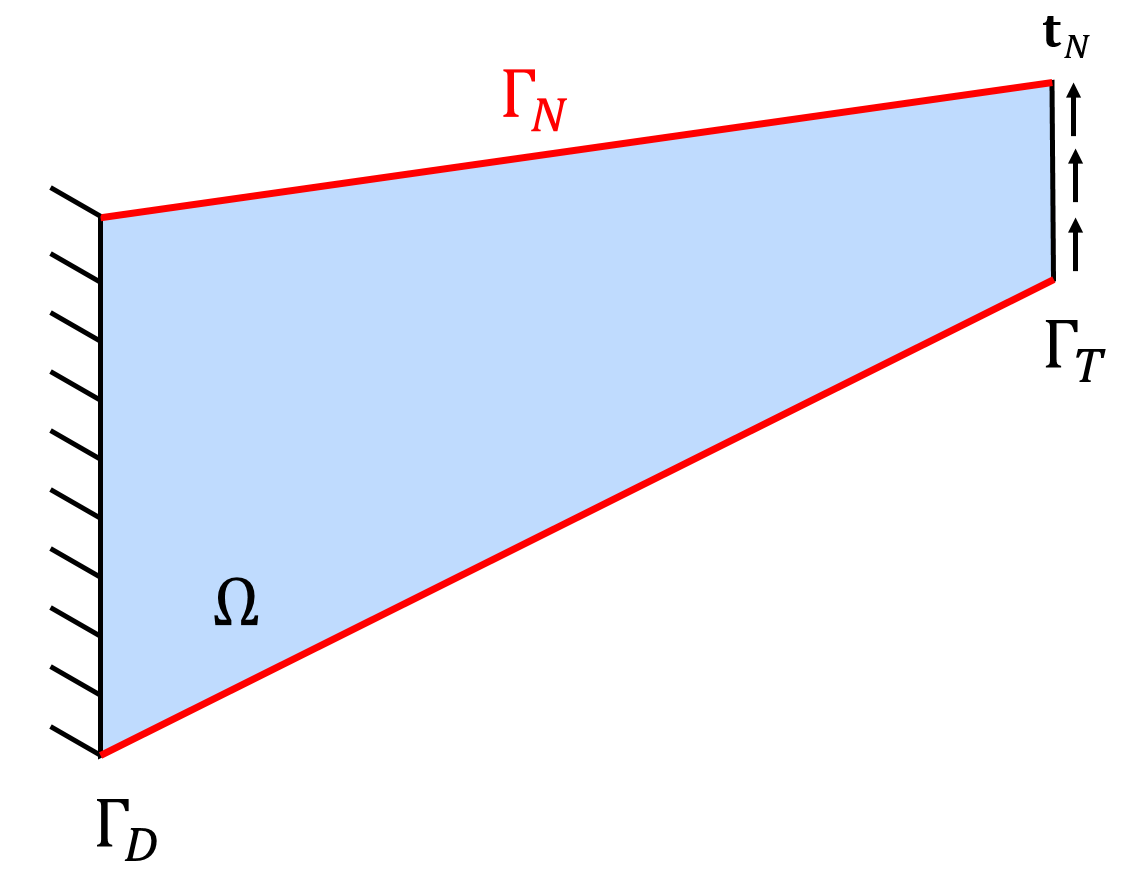}
	\end{minipage}
	\hspace{0.1cm}
	\begin{minipage}[t]{0.38\textwidth}
		\centering
		\makebox[\linewidth][l]{\textbf{(b)}}\par\vspace{0.3em}
		\includegraphics[width=\linewidth]{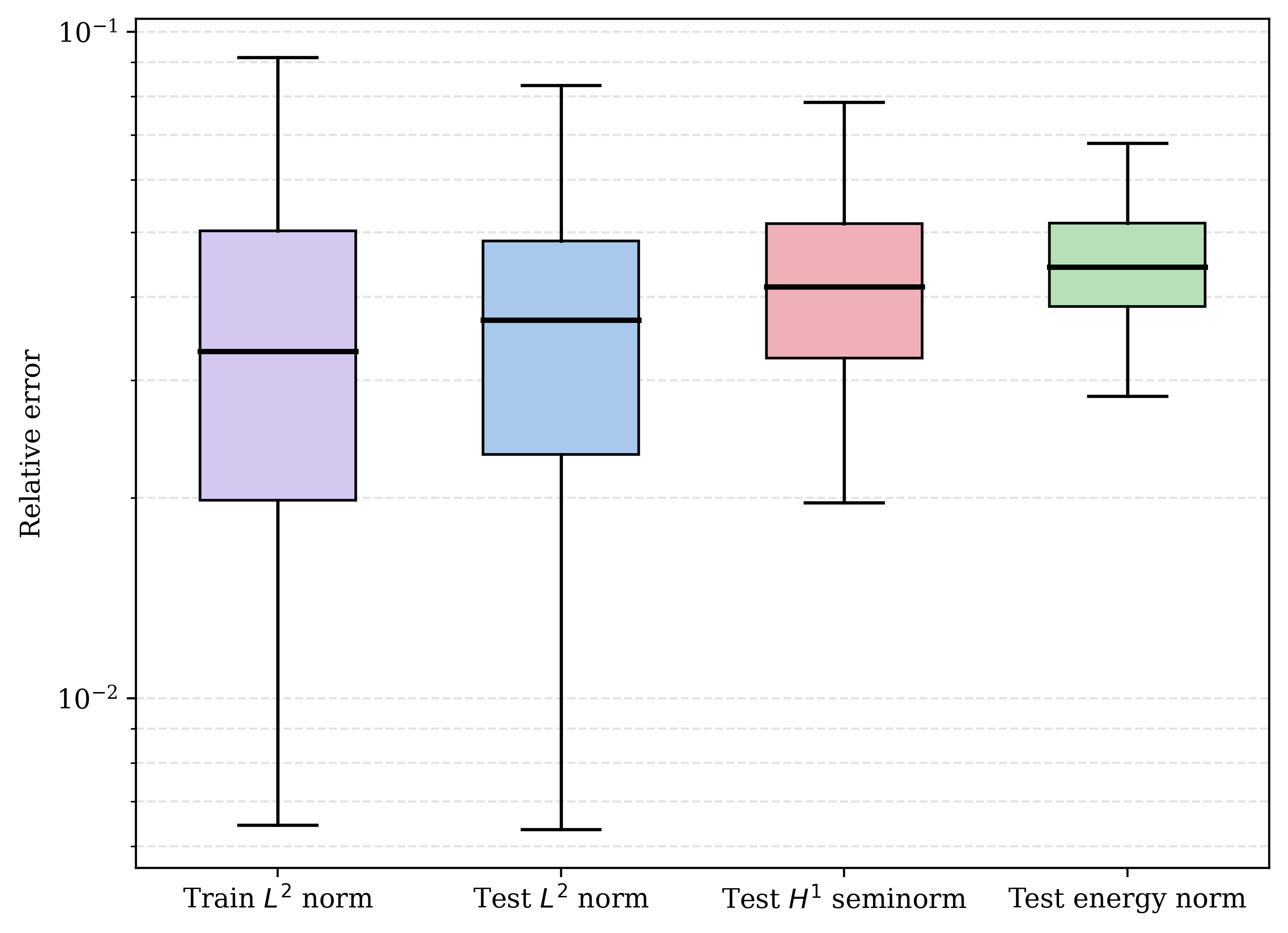}
	\end{minipage}\\
	\vspace{0.0cm}
	\begin{minipage}[t]{0.38\textwidth}
		\centering
		\makebox[\linewidth][l]{\textbf{(c)}}\par\vspace{0.3em}
		\includegraphics[width=\linewidth]{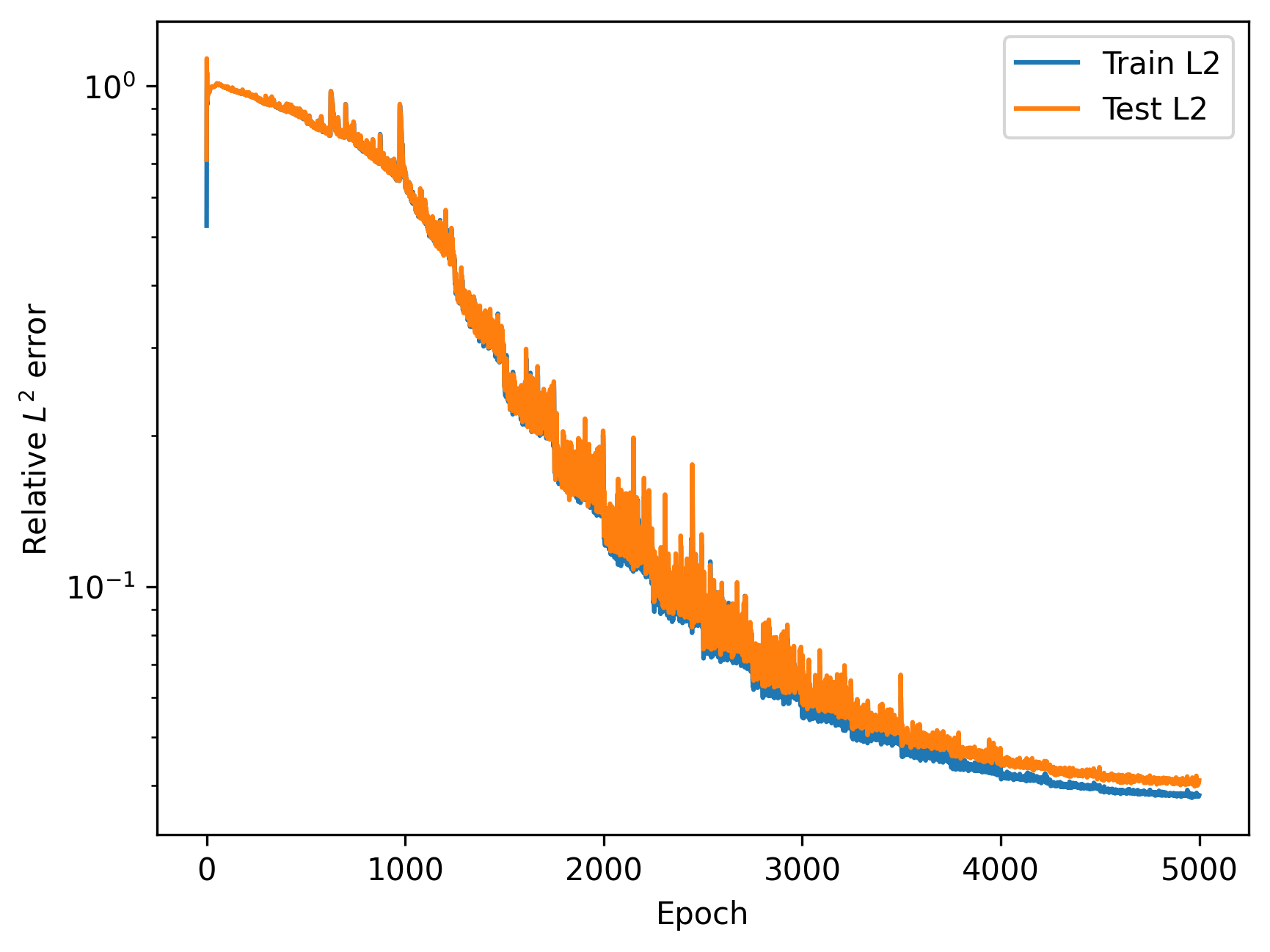}
	\end{minipage}
	\hspace{0.1cm}
	\begin{minipage}[t]{0.38\textwidth}
		\centering
		\makebox[\linewidth][l]{\textbf{(d)}}\par\vspace{0.3em}
		\includegraphics[width=\linewidth]{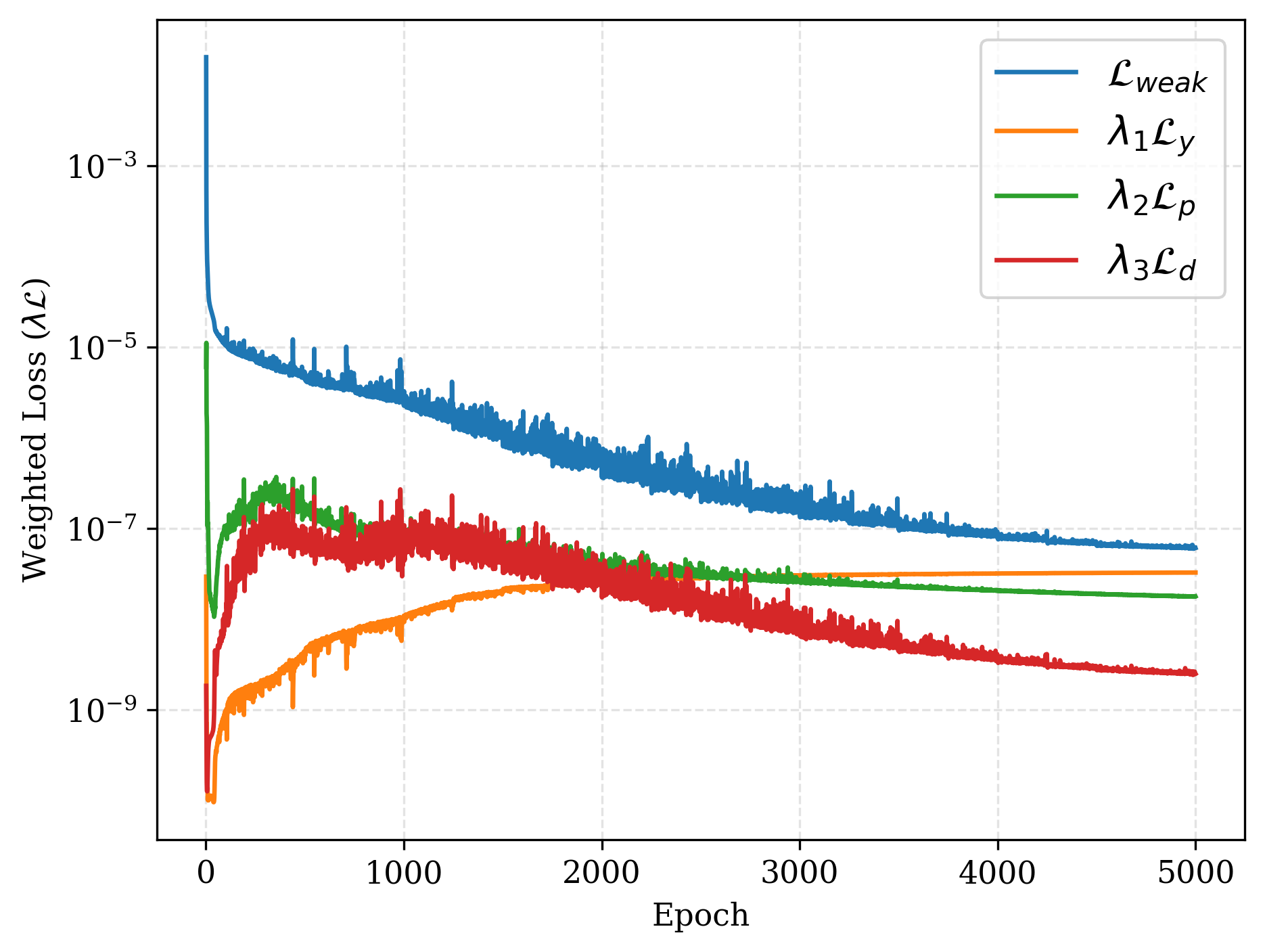}
	\end{minipage}
	\caption{Cook's membrane case. (a) Geometry of the Cook's membrane.  (b) Box plot of relative errors after training. (c) Evolution of the mean relative $L^2$ error during training. (d) The loss landscape of different residual terms.}
	\label{fig:case_cook}
\end{figure}

The next example is a shear-dominated Cook's membrane benchmark, which is more challenging than the plate with a hole case. We consider the same hyperelastic boundary value problem and loading setting as in Eq.~\eqref{eq:case_hole}, but embed a parametric Cook-type domain in the rectangular background $\mathcal{O}=[0,4]\times[0,2]$. All fields are represented on a fixed Cartesian mesh of size $51\times 51$ over $\mathcal{O}$. The right vertical boundary is fixed at $x=4$ with endpoints $(4,1.4)$ and $(4,2)$. The two slanted edges join these points to the left boundary $x=0$, where the lower and upper endpoint heights $y_b$ and $y_t$ are random,
\[
y_b\sim\mathcal{U}([0,0.2]),\qquad y_t\sim\mathcal{U}([1.8,2.0]).
\]
Since the sides $x=0$ and $x=4$ are mesh-aligned, $\varphi$ is used only for the two oblique boundaries. Define
\begin{equation}
\begin{aligned}
\psi_b(x,y) &= y-\left(y_b+\frac{1.4-y_b}{4}\,x\right), \\
\psi_t(x,y) &= \left(y_t+\frac{2-y_t}{4}\,x\right)-y,
\end{aligned}
\label{eq:cook_halfplanes}
\end{equation}
and set the level-set field as the product of these two edge factors,
\begin{equation}
\varphi(x,y)=\psi_b(x,y)\,\psi_t(x,y).
\label{eq:cook_phi}
\end{equation}
Then $\varphi=0$ on the two slanted boundaries, while the vertical sides are imposed directly on mesh-aligned edges, see Fig.~\ref{fig:case_cook}a. The Neumann traction is modeled by the same GRF construction as in Eqs.~\eqref{eq:grf_process} and \eqref{eq:grf_cov}. Specifically, we sample $\mathbf{t}_N(x,y)\sim\mathcal{GP}(\mathbf{t}_c,k_\ell)$ at background-mesh nodes with $\mathbf{t}_c=[0,t_y]^T$, $t_y\sim\mathcal{U}([5,15])$, signal variance $s^2=0.1$, and correlation length $\ell=0.2$.

In this benchmark, the homogeneous Dirichlet condition on $\Gamma_D$ is enforced by the multiplicative ansatz $\mathbf{u}_h=\mathbf{u}_{\theta}\,x$ for FNO models. The loss function of WINO is computed by \eqref{eq:total_loss} with the penalty parameters $\lambda_1=10^{-7},\lambda_2=1\times10^{-2}, \lambda_3=1$. We compute the loss function on a $51\times 51$ mesh and train for 5{,}000 epochs using the SOAP optimizer with 1,000 training samples and 100 test samples. The ground-truth fields are generated by $\varphi$-FEM on the $251\times 251$ mesh using second-order finite elements (biquadratic Q9 elements on quadrilateral meshes). Owing to strong nonlinearity in this membrane benchmark, the Newton linearizations are solved with sparse LU factorization rather than restarted GMRES: for the tangent systems arising here, GMRES did not converge within practical iteration budgets, whereas a direct solve based on LU decomposition remains robust. In this case, we learn the WINO mapping
\begin{equation}
	\mathcal{G}_{\theta}:(\varphi_h,\mathbf{t}_h)\mapsto(\mathbf{u}_h,\mathbf{y}_h,\mathbf{p}_h).
\end{equation}

\begin{figure}[t]
	\centering
	\begin{minipage}[t]{0.30\textwidth}
		\vspace{0pt}% anchor [t] to true top
		\centering
		\makebox[\linewidth][l]{\textbf{(a)}}\par\vspace{0.3em}
		\includegraphics[width=0.8\linewidth]{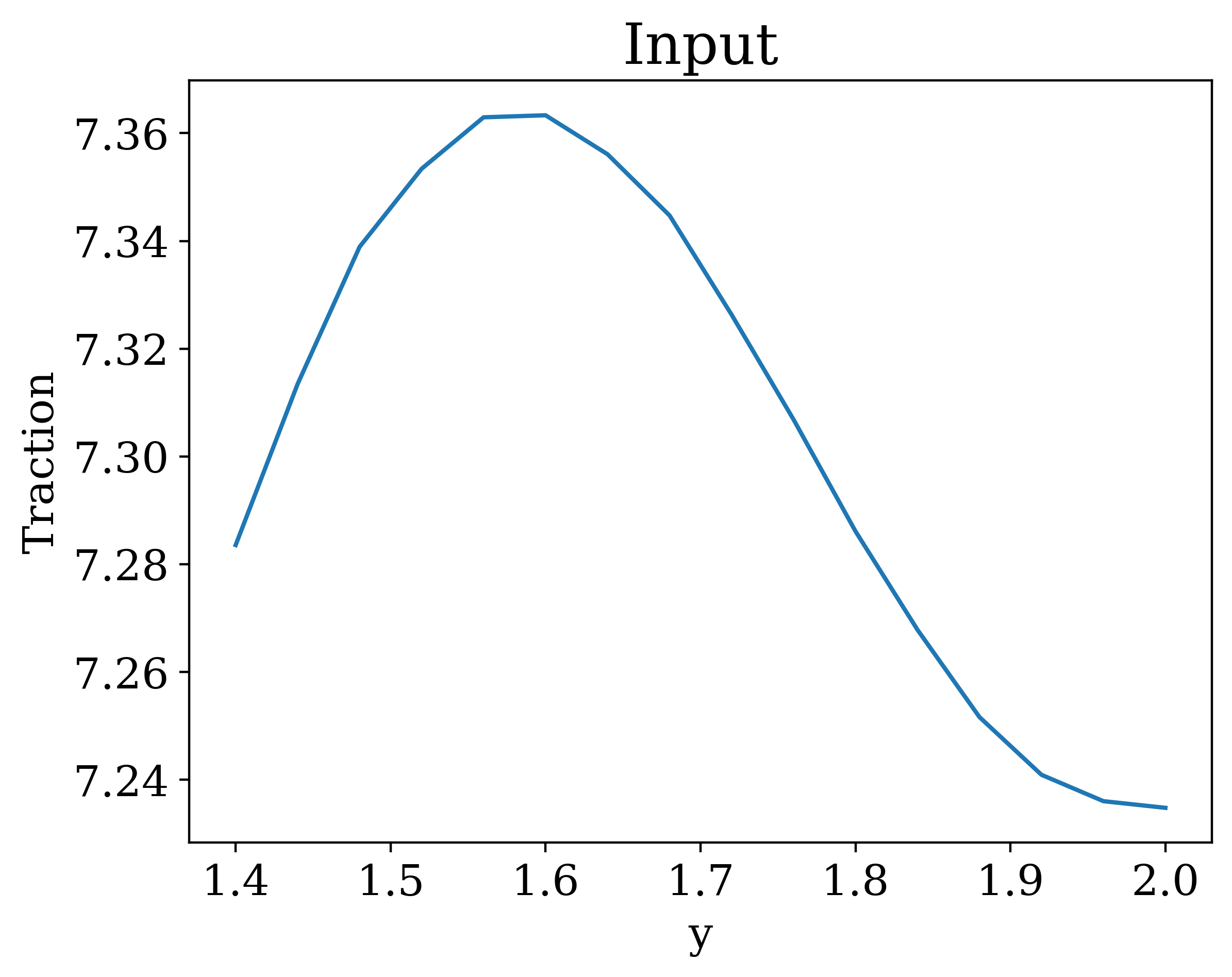}\\[0.6ex]
		\vspace{0.3cm}
		\includegraphics[width=\linewidth]{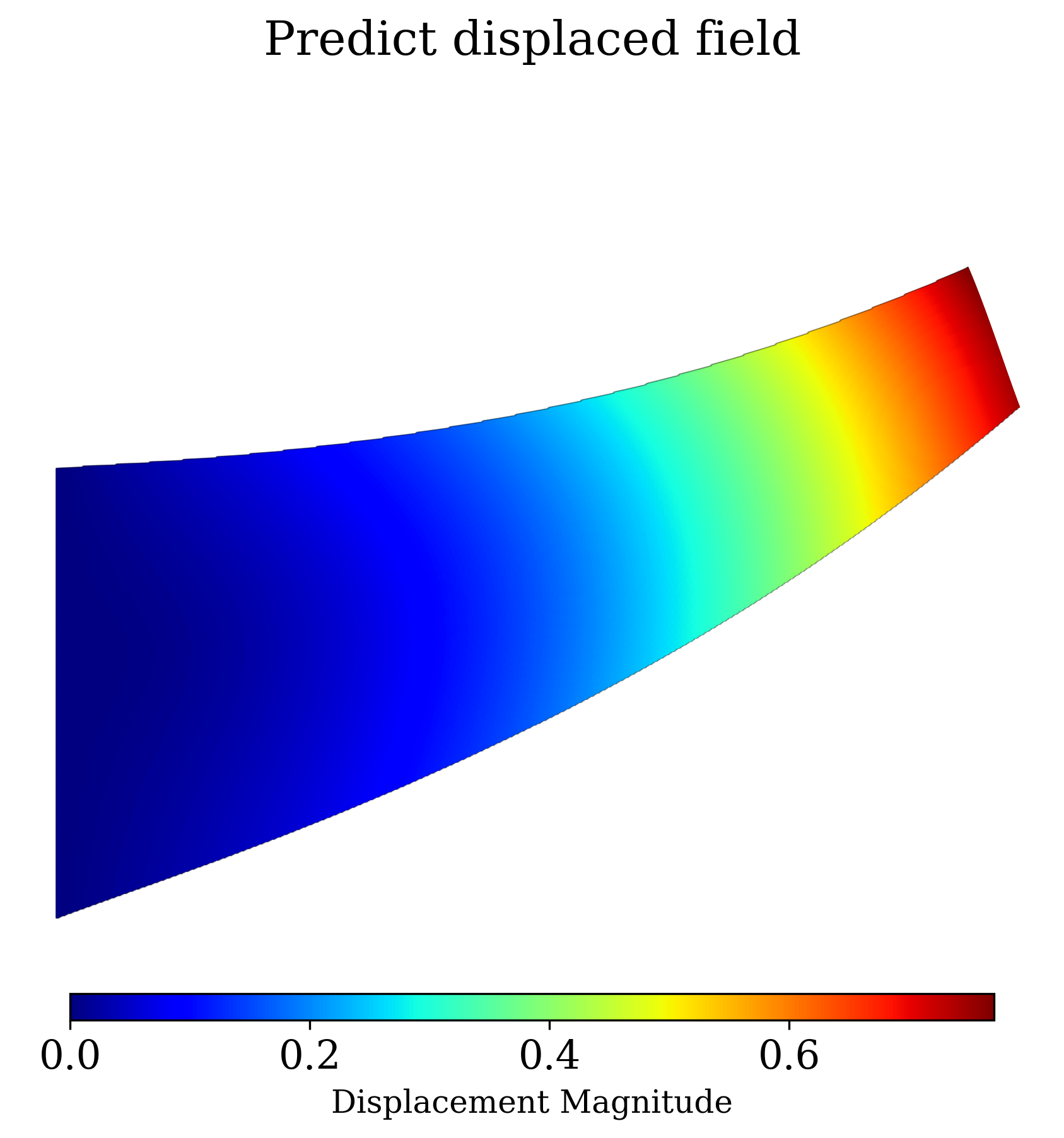}
	\end{minipage}%
	\hspace{0.03\textwidth}%
	\begin{minipage}[t]{0.57\textwidth}
		\vspace{0pt}% align top with left column
		\centering
		\makebox[\linewidth][l]{\textbf{(b)}}\par\vspace{0.3em}
		\begin{tabular}{@{}c@{\hspace{0.02\linewidth}}c@{}}
			\includegraphics[width=0.45\linewidth]{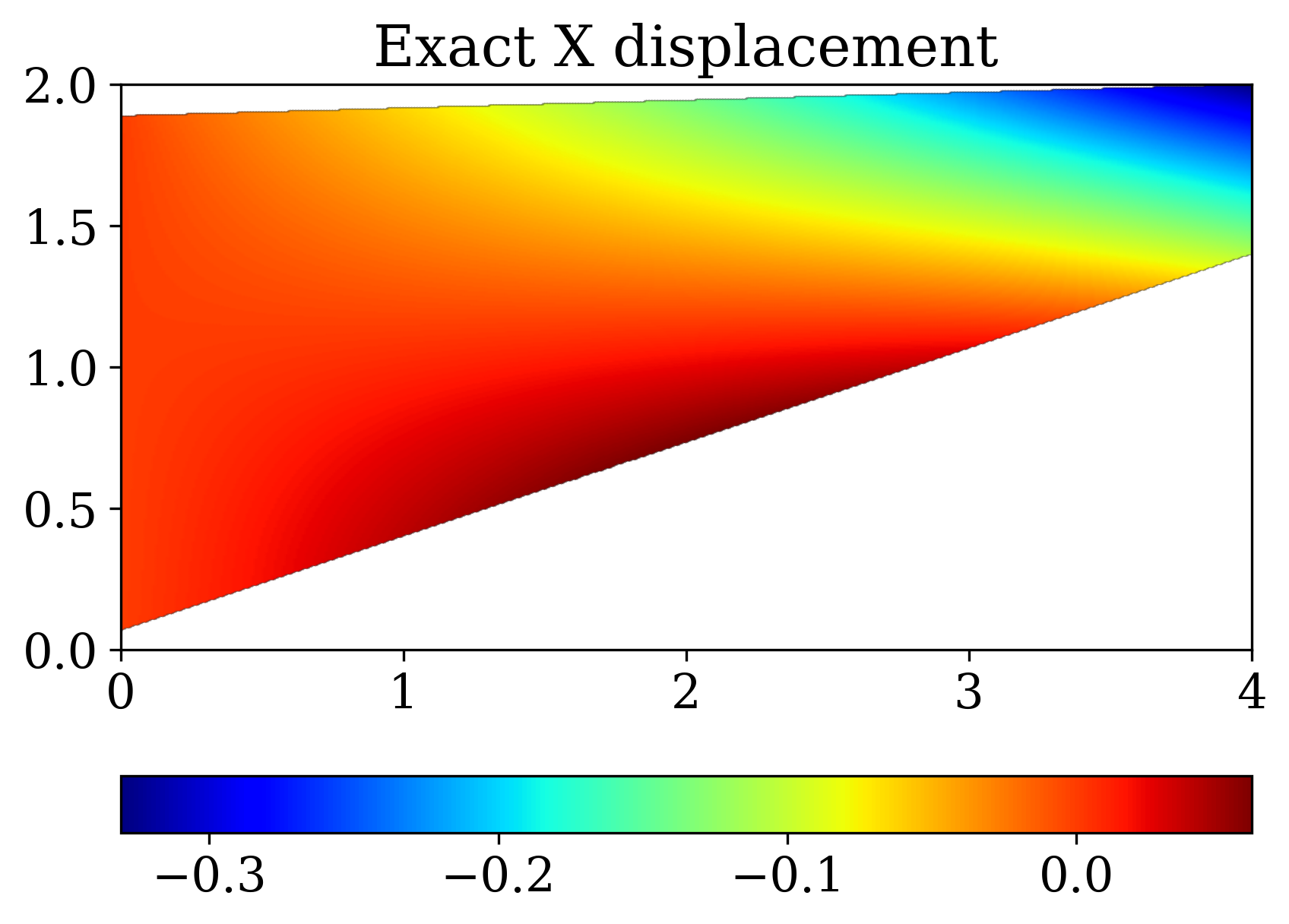} &
			\includegraphics[width=0.45\linewidth]{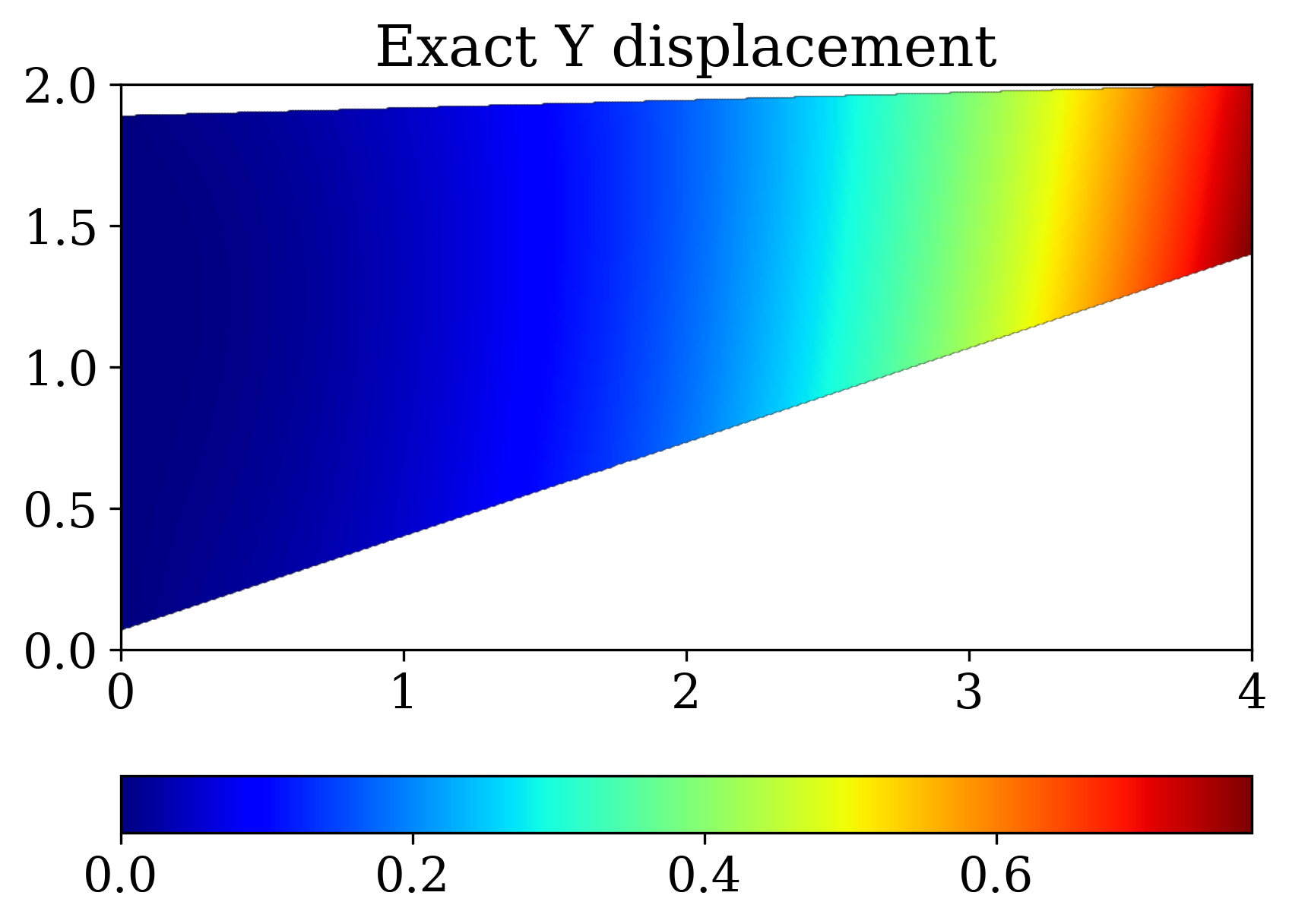} \\[0.5ex]
			\includegraphics[width=0.45\linewidth]{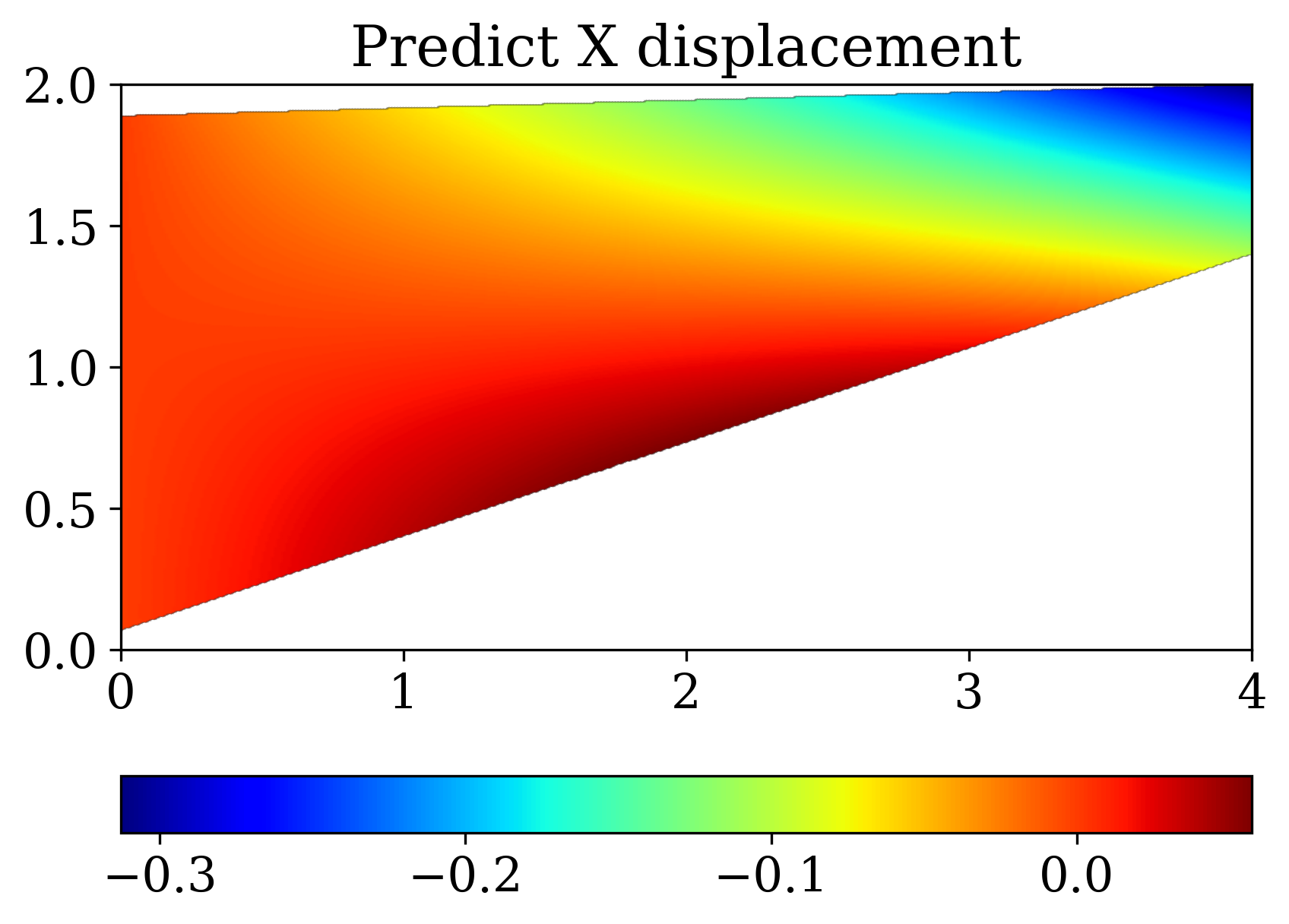} &
			\includegraphics[width=0.45\linewidth]{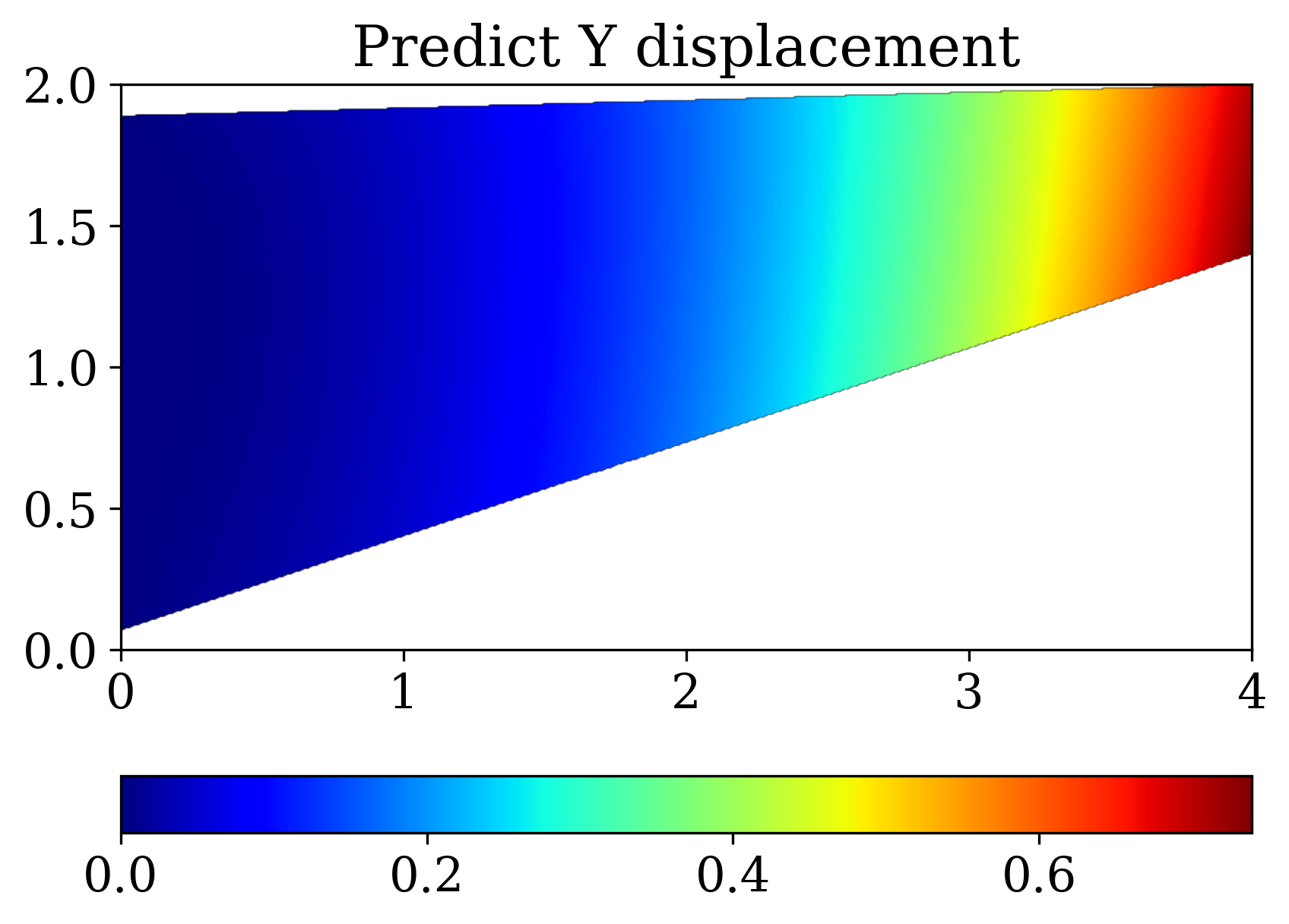} \\[0.5ex]
			\includegraphics[width=0.45\linewidth]{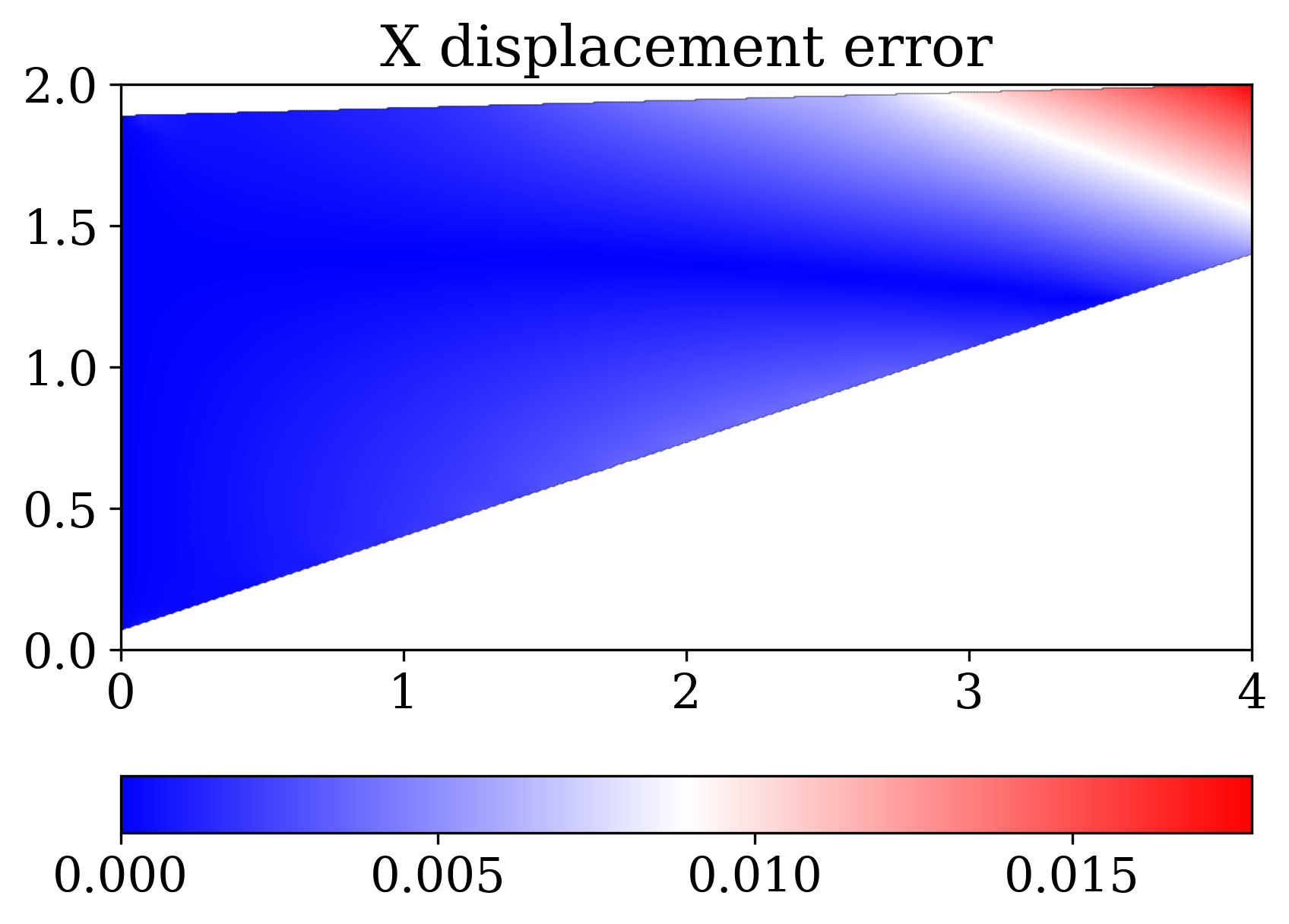} &
			\includegraphics[width=0.45\linewidth]{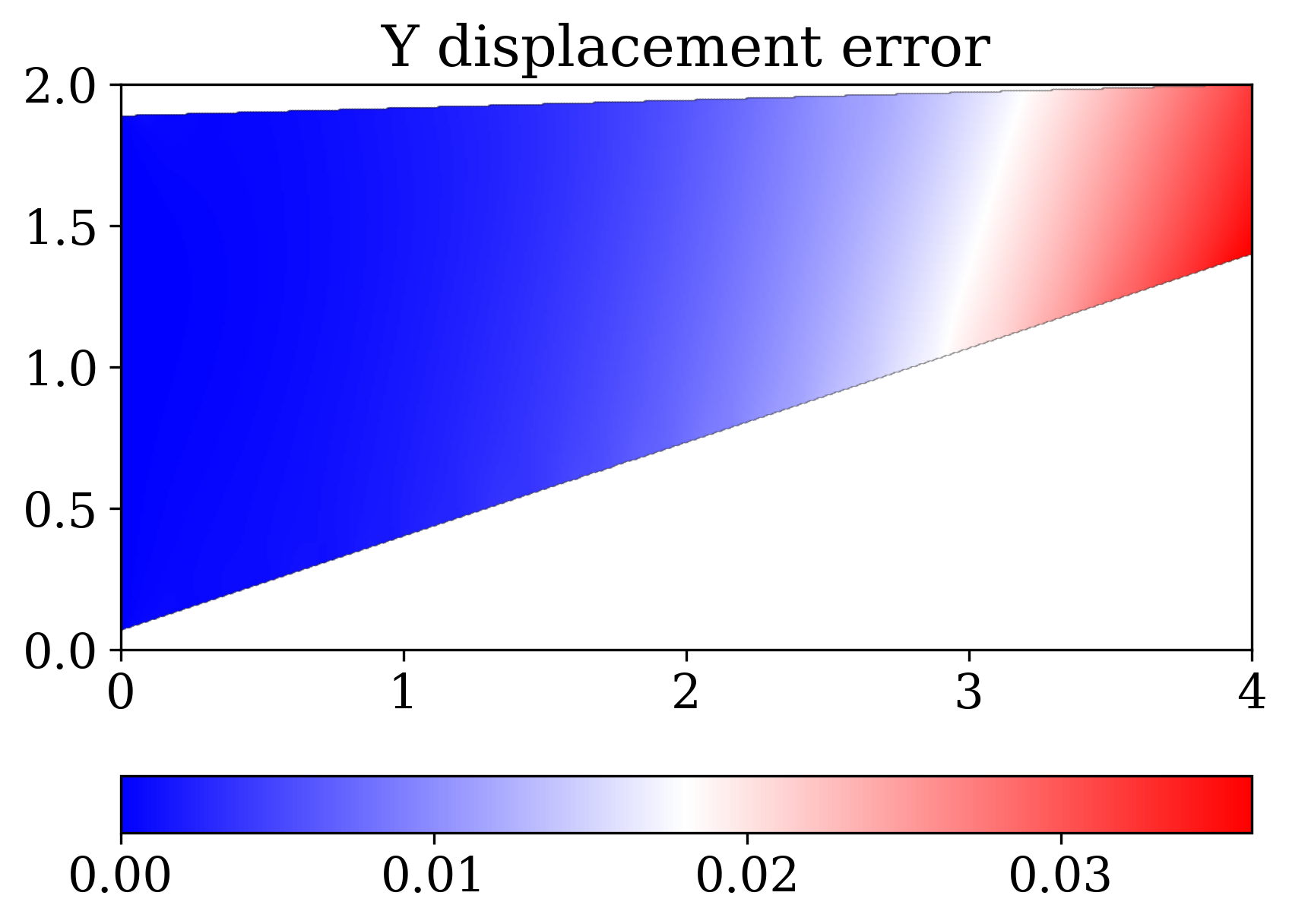}
		\end{tabular}
	\end{minipage}
	\caption{Representative test sample for the Cook's membrane case. (a) WINO-predicted deformed configuration; (b) $x$- and $y$-displacement components for the exact solutions (top row), WINO predictions (middle row), and pointwise absolute error contours (bottom row).}
	\label{fig:case_cook_displacement}
\end{figure}

\begin{table}[t]
	\centering
	\footnotesize
	\caption{Computational costs (data generation + training time) and relative errors of three methods for the Cook's membrane case. Cost ratio is total time relative to WINO.}
	\label{tab:performance_cook}
	\begin{tabularx}{\linewidth}{@{}l >{\raggedright\arraybackslash}X >{\raggedright\arraybackslash}X >{\raggedright\arraybackslash}X >{\raggedright\arraybackslash}X >{\raggedright\arraybackslash}X@{}}
		\toprule
		Method & Total time & Cost ratio & $\|e\|_{L^2}$ & $\|e\|_{H^1}$ & $\|e\|_{E}$ \\
		\midrule
		WINO & $0.6 + 19242.1$s & $1.00\times$ & $3.83 \pm 2.14\%$ & $4.32 \pm 1.68\%$ & $3.91 \pm 0.73\%$ \\
		$\varphi$-FEM-FNO & $73096.6 + 12757.2$s & $4.47\times$ & $0.05 \pm 0.02\%$ & $1.02 \pm 0.03\%$ & $2.40 \pm 0.54\%$ \\
		WINO+data & $73096.6 + 19679.1$s & $4.82\times$ & $0.07 \pm 0.01\%$ & $1.43 \pm 0.04\%$ & $3,18 \pm 0.65\%$ \\
		\bottomrule
	\end{tabularx}
\end{table}

\begin{figure}[t]
	\includegraphics[width=0.3\textwidth]{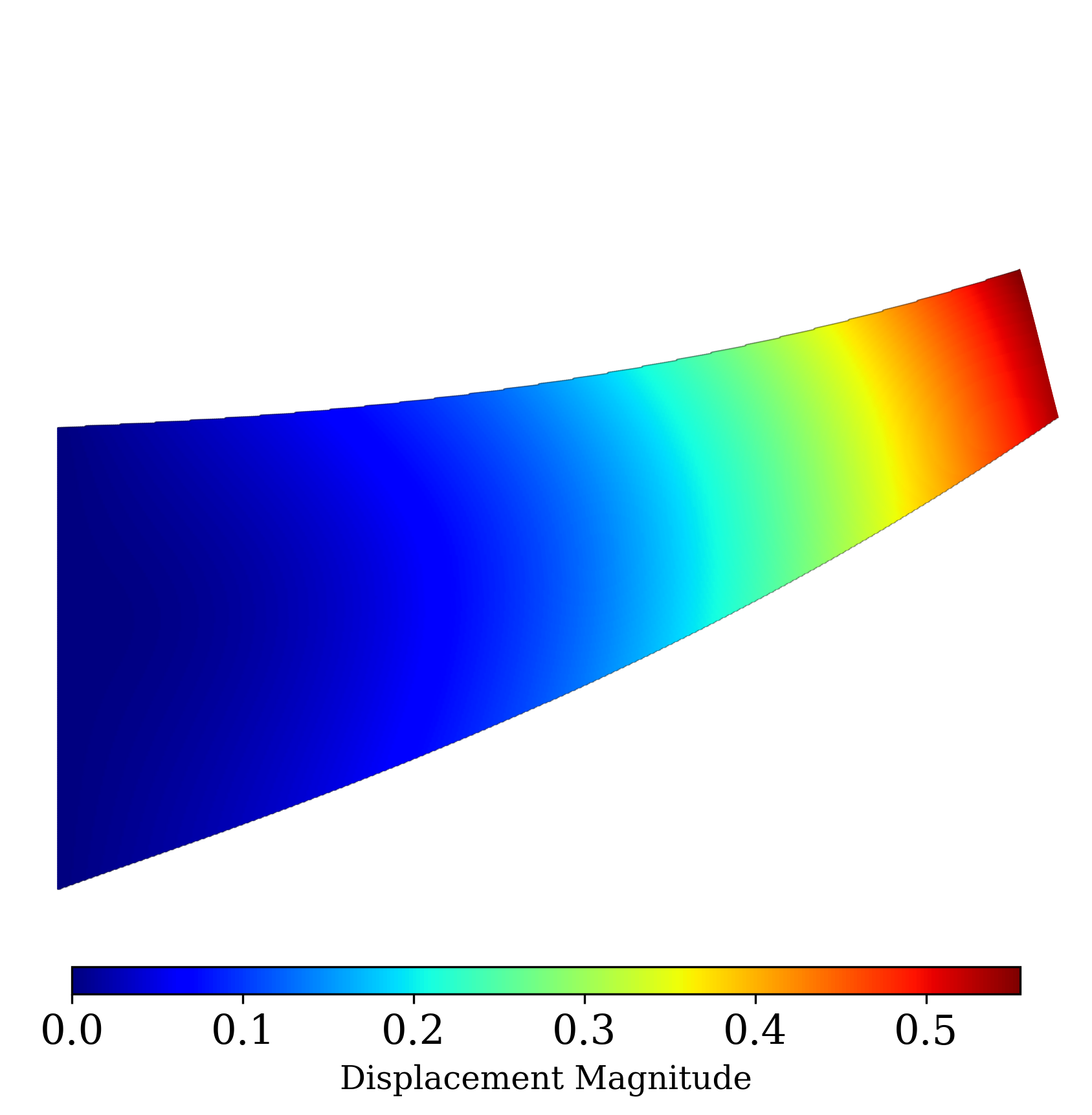}
	\hspace{0.02\textwidth}
	\includegraphics[width=0.3\textwidth]{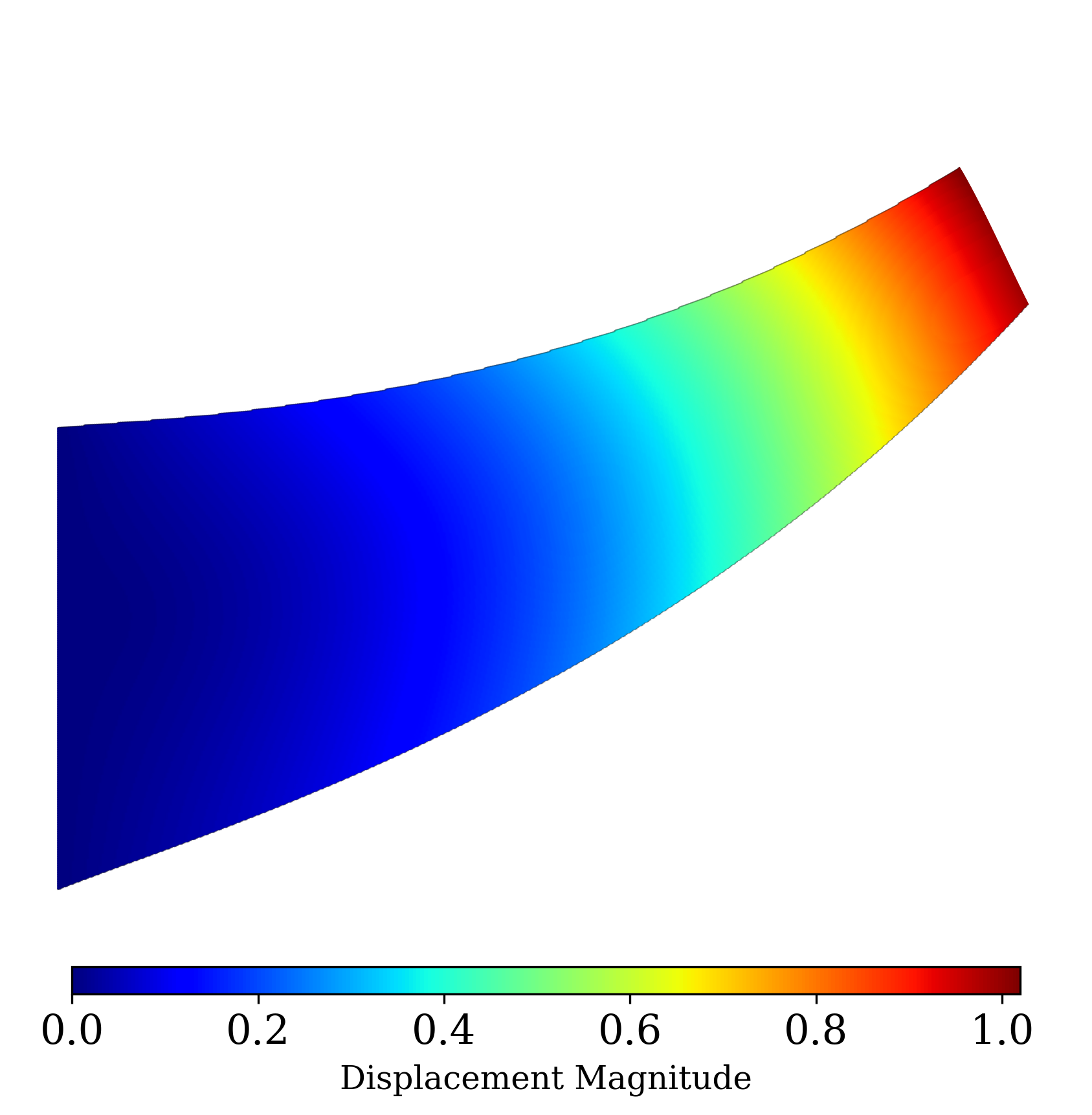}
	\hspace{0.02\textwidth}
	\includegraphics[width=0.3\textwidth]{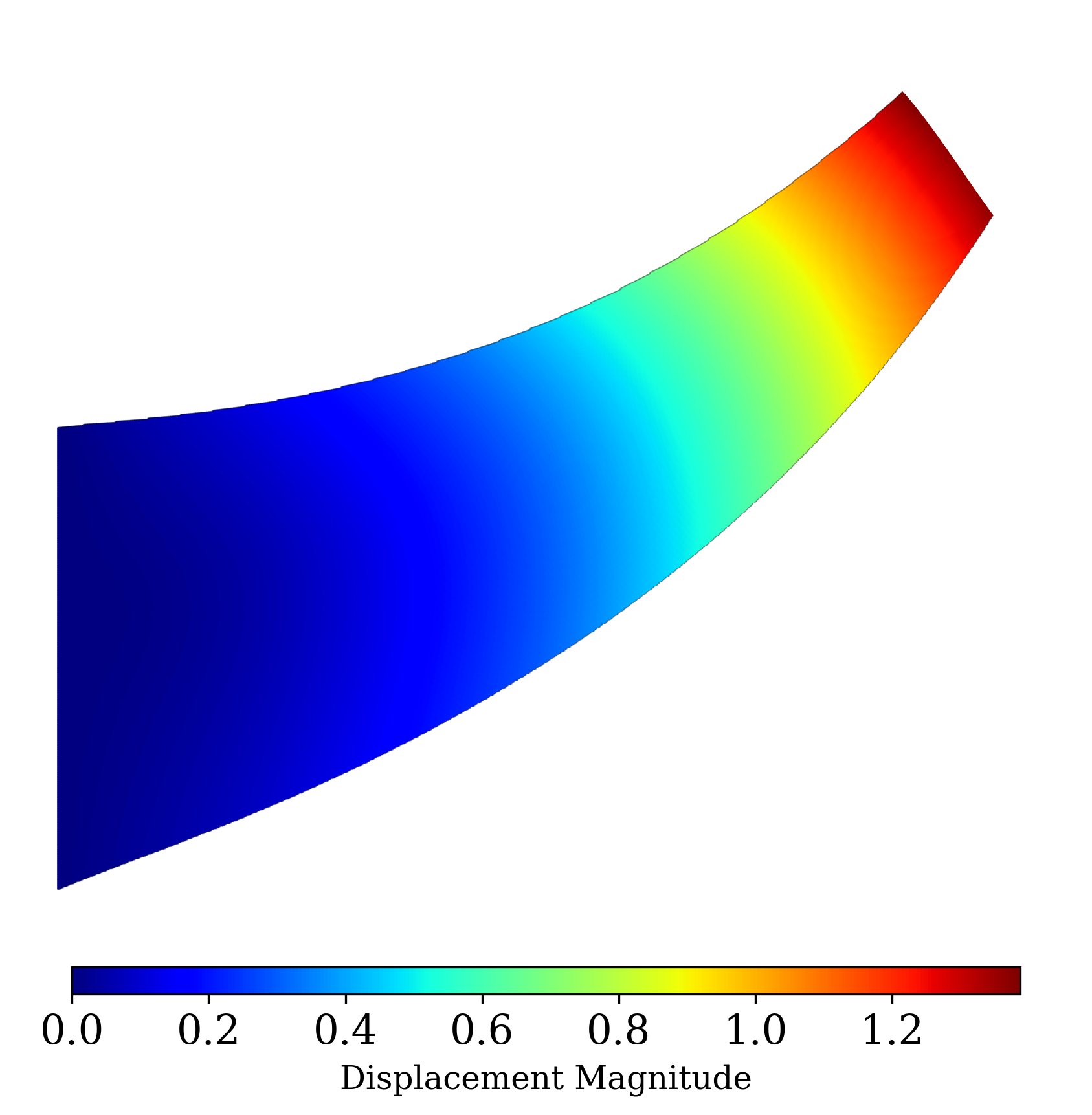}
	\vspace{0.6em}
	\caption{WINO predicted displacement fields for increasing traction levels $t_y=\{5,10,15\}$.}
	\label{fig:case_cook_predict5_15}
\end{figure}

Fig.~\ref{fig:case_cook}c reports the epoch-wise mean relative $L^2$ errors for the training and test sets, and both curves decrease steadily during optimization. Fig.~\ref{fig:case_cook}b summarizes the post-training error distributions using box plots, where the center line denotes the median and the box bounds indicate the interquartile range. A representative full-membrane test prediction is shown in Fig.~\ref{fig:case_cook_displacement}, and the WINO solution remains close to the high-fidelity reference while reproducing the displacement response under random pressure loading with good agreement. In addition, WINO can predict the displacement evolution under increasing traction loads, capturing how the membrane deformation grows with load magnitude; see Fig.~\ref{fig:case_cook_predict5_15}, where $t_y=5,10,15$ are applied. To assess efficiency, we further compare WINO, purely data-driven $\varphi$-FEM-FNO, and WINO with labeled data in terms of both accuracy and total computational cost (see Table~\ref{tab:performance_cook}). Although $\varphi$-FEM-FNO can reach very high accuracy in $L^2$ norm and requires only six training epochs ($\approx 15.31\,\mathrm{s}$) to attain approximately the same relative $L^2$ accuracy as label-free WINO, its reference-data generation is substantially more expensive (about $73.1\,\mathrm{s}$ per sample on average in this benchmark). In some benchmarks \cite{wang2026pretraining}, generating a single labeled sample can even exceed one hour, which motivates data-free training in such settings. Accordingly, the primary efficiency advantage of WINO in this comparison arises from its ability to circumvent large-scale upfront label generation, rather than from accelerating supervised operator training once such labels are already available. Fig.~\ref{fig:case_cook}d further illustrates the optimization landscape discussed in Remark~3.1 by tracking the evolution of the individual residual contributions $\mathcal{L}_{\mathrm{weak}}$, $\lambda_1\mathcal{L}_y$, $\lambda_2\mathcal{L}_p$, and $\lambda_3\mathcal{L}_d$ during training. While $\mathcal{L}_{\mathrm{weak}}$, $\mathcal{L}_p$, and $\mathcal{L}_d$ decrease overall, $\mathcal{L}_y$ can increase as optimization proceeds, indicating competing descent directions among the coupled auxiliary constraints. This multi-objective coupling makes the WINO loss more difficult to minimize than the single supervised regression used by $\varphi$-FEM-FNO and helps explain the slower convergence and higher residual accuracy plateau observed for label-free WINO on mixed boundary problems.

We next investigate mesh-transfer performance on the same 100 test realizations. Specifically, the operators previously trained on the $51\times51$ grid are applied, without retraining, to the corresponding inputs discretized on a substantially finer $251\times251$ mesh, and the resulting predictions are compared with $\varphi$-FEM reference solutions computed on the same $251\times251$ discretization (Table~\ref{tab:cook_resolution}). Each reference sample on the $251\times251$ mesh costs $730.23\,\mathrm{s}$ on average, which is substantially higher than the mean per-sample data cost on the $51\times51$ grid. The WINO results remain largely unchanged across resolutions, whereas the $\varphi$-FEM-FNO method exhibits more noticeable degradation on the refined mesh, particularly in the relative $H^1$ seminorm and energy norm. The reason is that $\varphi$-FEM-FNO tends to produce less accurate displacement gradients near domain boundaries, and such localized discrepancies are amplified in $H^1$-type and energy-based metrics. Although FNOs are theoretically grid-independent, their pointwise activation functions inevitably expand the function bandwidth beyond the Nyquist frequency, breaking continuous-discrete equivalence and inducing significant operator aliasing errors across different discrete grids \cite{bartolucci2023representation}. This supports the conclusion that WINO retains stable accuracy under mesh refinement beyond a mere change of supervised loss terms.

\begin{table}[t]
	\centering
	\footnotesize
	\caption{Performance of models on the Cook's membrane case with training on 51$\times$51 meshes and testing on 251$\times$251 meshes.}
	\label{tab:cook_resolution}
	\begin{tabularx}{\linewidth}{@{}l >{\centering\arraybackslash}X >{\raggedright\arraybackslash}X >{\raggedright\arraybackslash}X >{\raggedright\arraybackslash}X@{}}
		\toprule
		Method & Data cost (per sample) & $\|e\|_{L^2}$ & $\|e\|_{H^1}$ & $\|e\|_{E}$ \\
		\midrule
		WINO & \multirow{3}{*}{$730.23\,\mathrm{s}$} & $4.34 \pm 0.84\%$ & $4.79 \pm 1.17\%$ & $4.73 \pm 2.21\%$ \\
		$\varphi$-FEM-FNO &  & $0.55 \pm 0.14\%$ & $11.4 \pm 1.85\%$ & $35.5 \pm 1.51\%$ \\
		WINO+data &  & $0.36 \pm 0.04\%$ & $1.21 \pm 0.18\%$ & $0.35 \pm 0.04\%$ \\
		\bottomrule
	\end{tabularx}
\end{table}

\begin{remark}
	Since this benchmark involves mixed Dirichlet and Neumann conditions, the loaded right edge is taken as the vertical segment at $x=4$ between $(4,1.4)$ and $(4,2)$, while the two slanted sides are encoded by the level-set product \eqref{eq:cook_phi}. We use a uniform $51\times51$ Cartesian background on $\mathcal{O}=[0,4]\times[0,2]$ and fix the Neumann endpoints at $y=1.4$ and $y=2$ so that both points coincide with background-mesh nodes. Together with the mesh-aligned coordinate $x=4$, the traction boundary lies on a union of full element facets rather than cutting through element interiors; the Neumann boundary condition is therefore imposed by standard conforming facet assembly, without invoking the cut-cell $\varphi$-FEM auxiliary formulation on that edge.
\end{remark}

\subsubsection{Pressure vessel problem}

\begin{figure}[t]
	\centering
	\begin{minipage}[t]{0.23\textwidth}
		\centering
		\makebox[\linewidth][l]{\textbf{(a)}}\par\vspace{0.3em}
		\includegraphics[width=\linewidth]{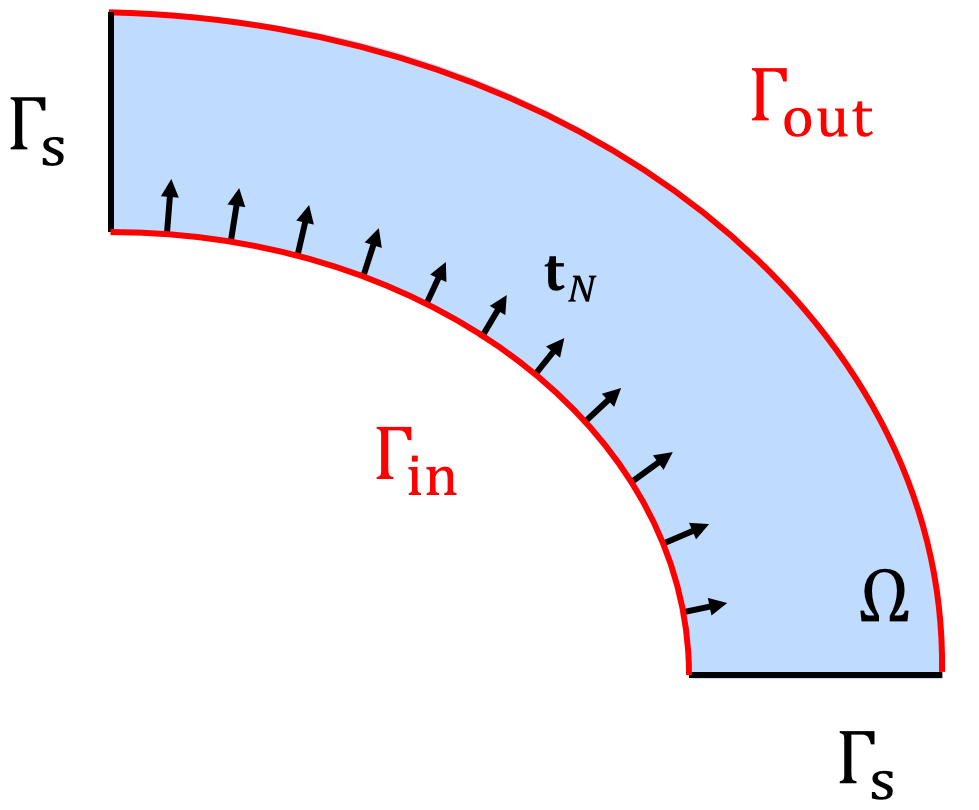}
	\end{minipage}
	\hfill
	\begin{minipage}[t]{0.23\textwidth}
		\centering
		\makebox[\linewidth][l]{\textbf{(b)}}\par\vspace{0.3em}
		\includegraphics[width=\linewidth]{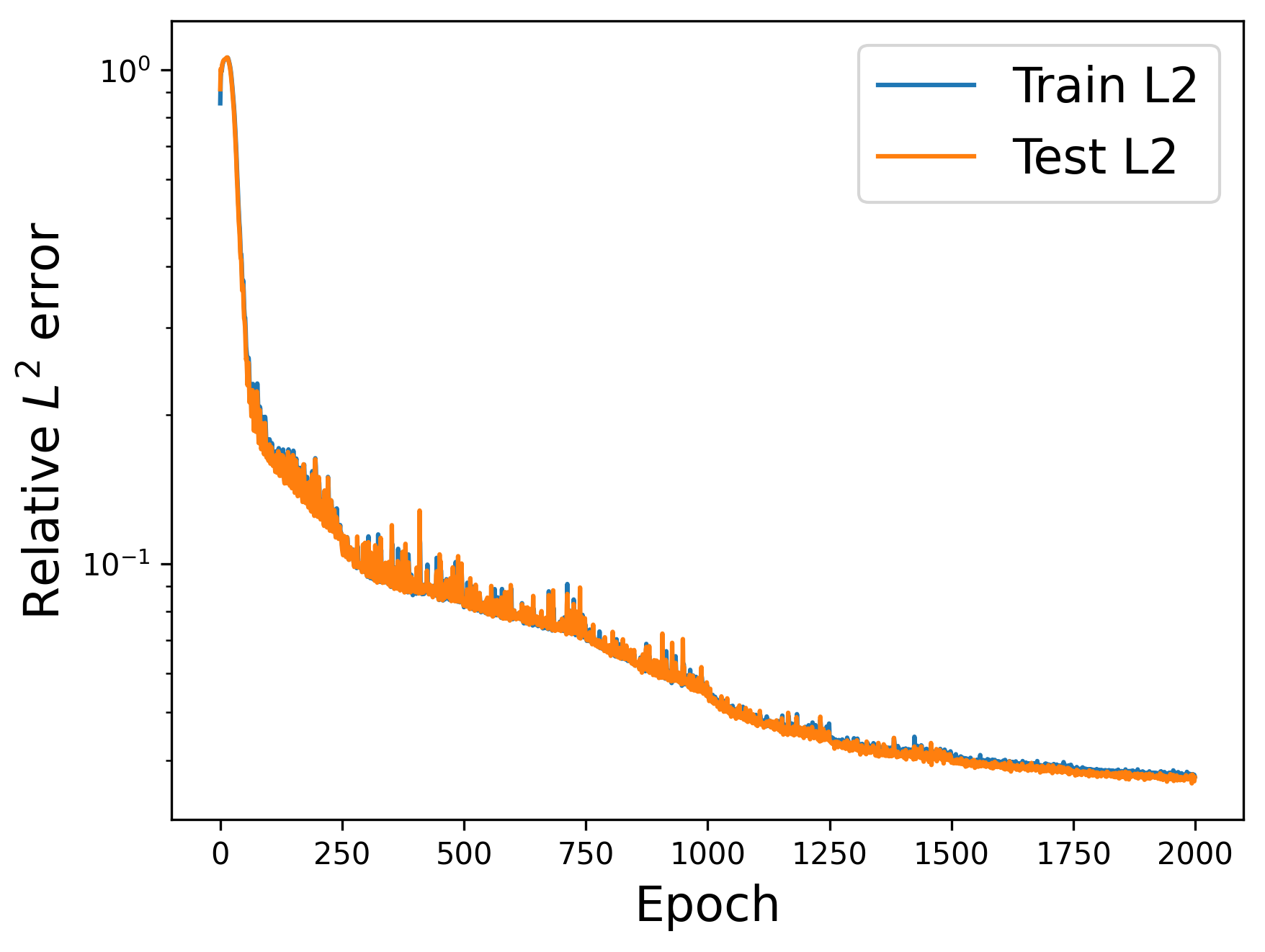}
	\end{minipage}
	\hfill
	\begin{minipage}[t]{0.23\textwidth}
		\centering
		\makebox[\linewidth][l]{\textbf{(c)}}\par\vspace{0.3em}
		\includegraphics[width=\linewidth]{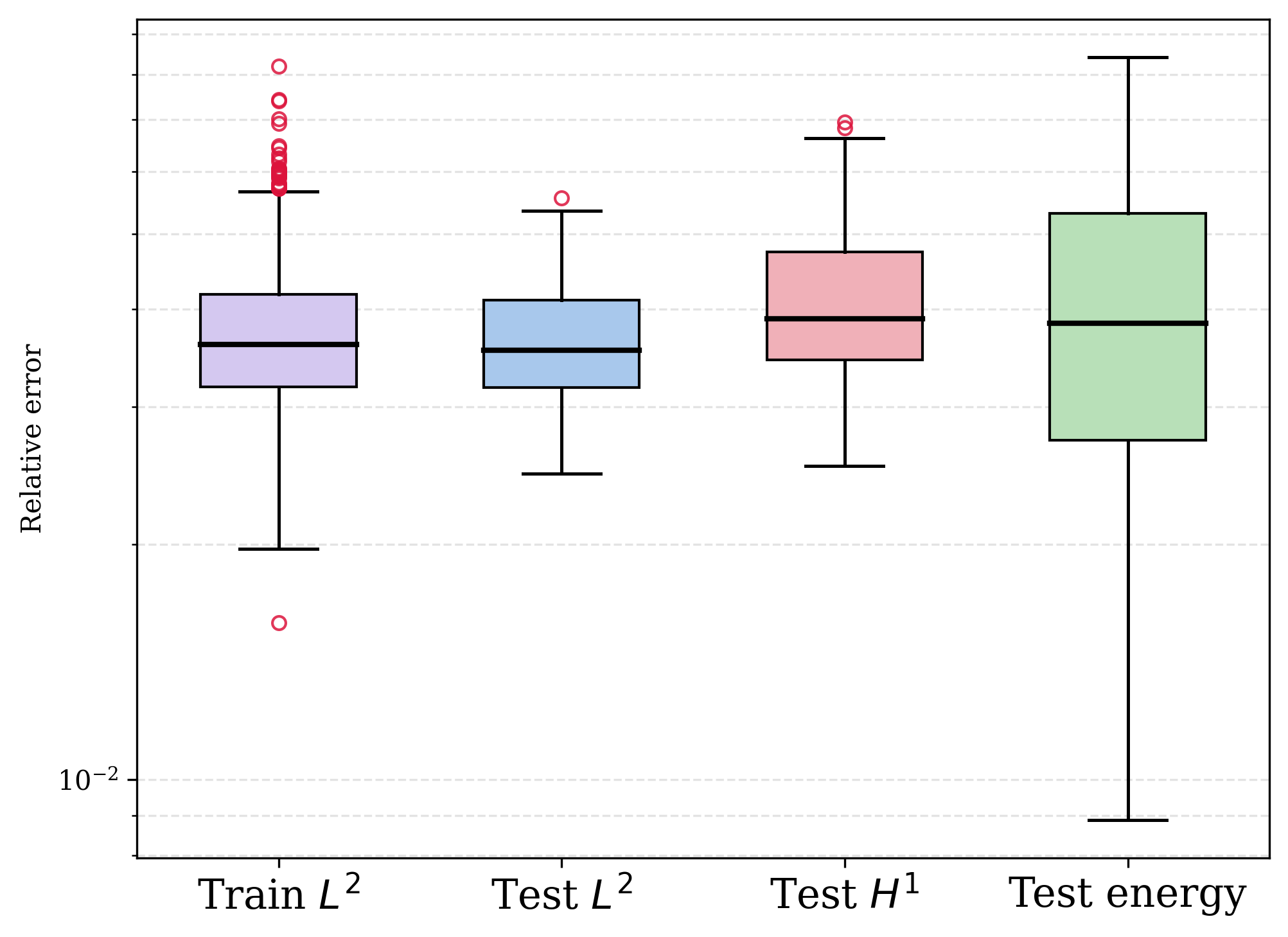}
	\end{minipage}
	\hfill
	\begin{minipage}[t]{0.23\textwidth}
		\centering
		\makebox[\linewidth][l]{\textbf{(d)}}\par\vspace{0.3em}
		\includegraphics[width=\linewidth]{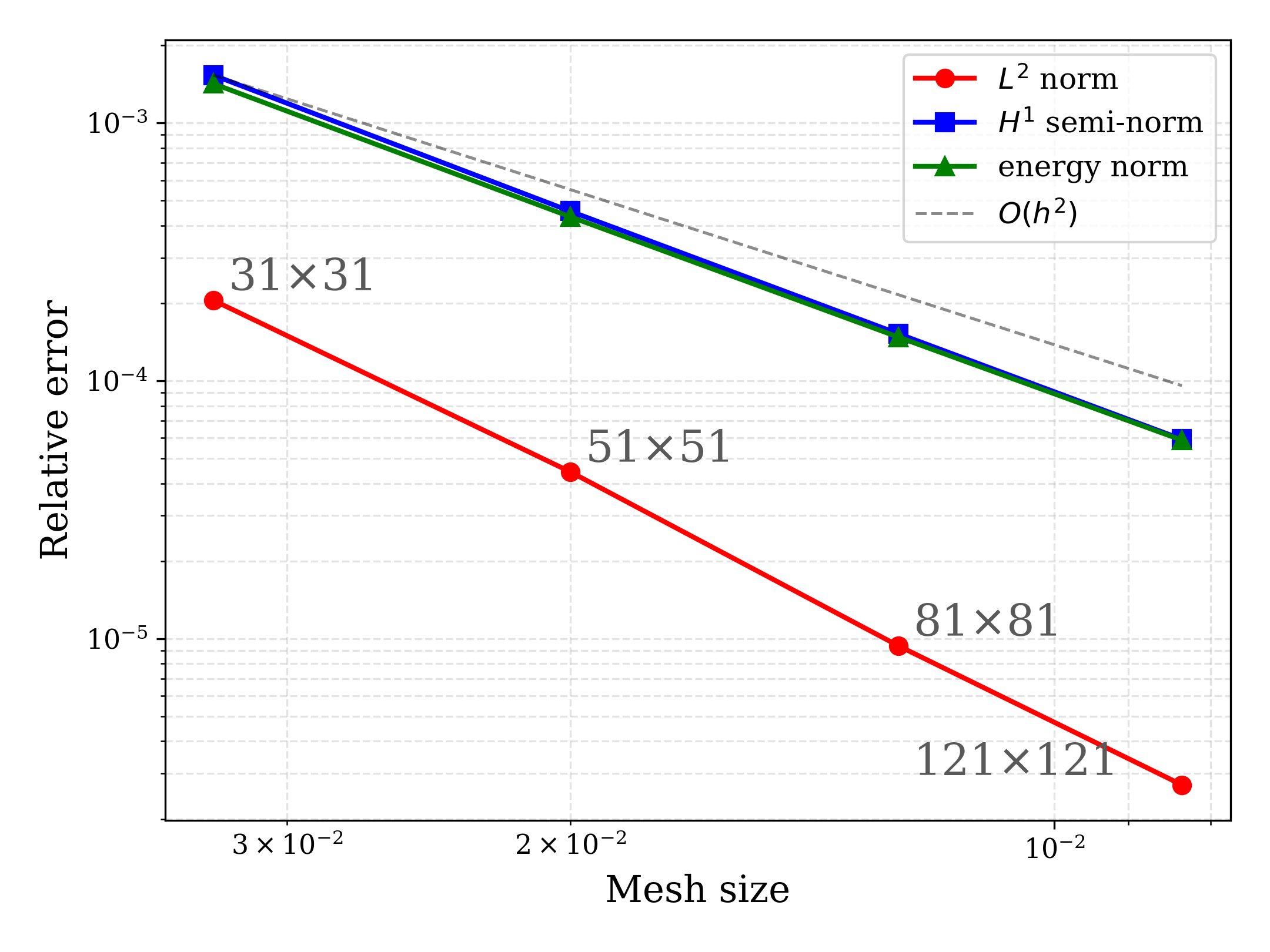}
	\end{minipage}
	\caption{Pressure vessel case. (a) Geometry of the pressure vessel. (b) Evolution of the mean relative $L^2$ error during training. (c) Box plot of relative errors after training. (d) convergence of $\varphi$-FEM: relative $L^2$, $H^1$ seminorm, and energy-norm errors versus mesh size.}
	\label{fig:case_vessel}
\end{figure}

We finally consider a quarter pressure-vessel benchmark in large-deformation hyperelasticity, which is formulated as follows:
\begin{equation}
	\begin{split}
		-\nabla\cdot\mathbf{P}(\mathbf{F}(\mathbf{u})) &= \mathbf{0}, \quad \text{in } \Omega,\\
		\mathbf{u}\cdot\mathbf{n} &= 0, \quad \text{on } \Gamma_{\mathrm{s}},\\
		\mathbf{P}(\mathbf{F}(\mathbf{u}))\cdot\mathbf{n} &= \mathbf{T}_N, \quad \text{on } \Gamma_{\mathrm{in}},\\
		\mathbf{P}(\mathbf{F}(\mathbf{u}))\cdot\mathbf{n} &= \mathbf{0}, \quad \text{on } \Gamma_{\mathrm{out}},
	\end{split}
	\label{eq:case_vessel}
\end{equation}
where $\Gamma_{\mathrm{s}}$ denotes the symmetry planes ($x=0$ and $y=0$), $\Gamma_{\mathrm{in}}$ is the inner pressurized wall, and $\Gamma_{\mathrm{out}}$ is the traction-free outer wall. The pressure load is a follower force: in the current configuration $\mathbf{t}_N=-p_0\mathbf{n}$, and by Nanson's law its pull-back to the reference configuration is 
\begin{equation}
\mathbf{T}_N=-p_0J\mathbf{F}^{-T}\mathbf{n},
\end{equation}
where $J=\det\mathbf{F}$ is the Jacobian determinant. The background mesh is defined on $[0,1]\times[0,1]$, and the symmetry boundary $\Gamma_{\mathrm{s}}$ is aligned with mesh-aligned edges. The inner and outer vessel walls are represented by the level-set function $\varphi(x,y)$:
\begin{equation}
	\varphi(x,y)=E_{\mathrm{in}}(x,y)\,E_{\mathrm{out}}(x,y),
\end{equation}
where $E_{\mathrm{in}}(x,y)$ and $E_{\mathrm{out}}(x,y)$ are the level-set functions for the inner and outer elliptical boundaries:
\begin{equation}
	E_{\mathrm{in}}(x,y)=\sqrt{\frac{x^2}{a_{\mathrm{in}}^2}+\frac{y^2}{b_{\mathrm{in}}^2}}-1,
	\label{eq:vessel_Ein}
\end{equation}
\begin{equation}
	E_{\mathrm{out}}(x,y)=\sqrt{\frac{x^2}{a_{\mathrm{out}}^2}+\frac{y^2}{b_{\mathrm{out}}^2}}-1.
	\label{eq:vessel_Eout}
\end{equation}
Here, $a_{\mathrm{out}}, b_{\mathrm{out}} \sim \mathcal{U}([0.9, 1.0])$, $a_{\mathrm{in}}=a_{\mathrm{out}}-r$, $b_{\mathrm{in}}=a_{\mathrm{out}}-r$, and $r \sim \mathcal{U}([0.2, 0.25])$. The geometry is shown in Fig.~\ref{fig:case_vessel}a.

\begin{table}[t]
	\centering
	\footnotesize
	\caption{Computational costs (data generation + training time) and relative errors of three methods for the pressure vessel case. Cost ratio is total time relative to WINO.}
	\label{tab:performance_vessel}
	\begin{tabularx}{\linewidth}{@{}l >{\raggedright\arraybackslash}X >{\raggedright\arraybackslash}X >{\raggedright\arraybackslash}X >{\raggedright\arraybackslash}X >{\raggedright\arraybackslash}X@{}}
		\toprule
		Method & Total time & Cost ratio & $\|e\|_{L^2}$ & $\|e\|_{H^1}$ & $\|e\|_{E}$ \\
		\midrule
		WINO & $0.613 + 17868.4$s & $1.00\times$ & $3.71 \pm 0.78\%$ & $4.21 \pm 1.14\%$ & $3.85 \pm 1.99\%$ \\
		$\varphi$-FEM-FNO & $13498.2 + 14308.4$s & $1.56\times$ & $0.03 \pm 0.01\%$ & $2.88 \pm 0.35\%$ & $0.31 \pm 0.09\%$ \\
		WINO+data & $13498.2 + 18755.4$s & $1.80\times$ & $0.01 \pm 0.01\%$ & $0.40 \pm 0.02\%$ & $0.06 \pm 0.01\%$ \\
		\bottomrule
	\end{tabularx}
\end{table}

The level-set function $\varphi$ is used only to represent the inner and outer elliptical boundaries. The symmetry planes are treated with standard conforming finite elements, without invoking the cut-cell $\varphi$-FEM boundary machinery on those edges. In WINO and $\varphi$-FEM-FNO, the homogeneous Dirichlet condition on $\Gamma_{\mathrm{s}}$ is enforced by the multiplicative ansatz $\mathbf{u}_x=\mathbf{u}_{\theta x}\,x, \mathbf{u}_y=\mathbf{u}_{\theta y}\,y$. The loss function of WINO is computed by \eqref{eq:total_loss} with the penalty parameters $\lambda_1=\lambda_3=1\times10^{-2}, \lambda_2=1\times10^{-3}$. We seek a neural operator $\mathcal{G}_\theta$ such that
\begin{equation}
\label{eq:vessel_mapping}
\mathcal{G}_\theta:(E_{\mathrm{in}},E_{\mathrm{out}},\partial_x E_{\mathrm{in}},\partial_y E_{\mathrm{in}},\partial_x E_{\mathrm{out}},\partial_y E_{\mathrm{out}},p_0)\to(\mathbf{u}_h,\mathbf{y}_h,\mathbf{p}_h).
\end{equation}
The channels $\partial_x E_{\mathrm{in}}$, $\partial_y E_{\mathrm{in}}$, $\partial_x E_{\mathrm{out}}$, and $\partial_y E_{\mathrm{out}}$ are obtained analytically by differentiating the closed-form expressions \eqref{eq:vessel_Ein} and \eqref{eq:vessel_Eout}. These quantities are evaluated node-wise on the background mesh together with $E_{\mathrm{in}}$ and $E_{\mathrm{out}}$. The pressure input $p_0$ is sampled as a Gaussian random field with mean $10$ and variance $0.1$. To better encode boundary information for WINO, on cut cells we approximate the level-set field using the Q9 interpolation (cf. \eqref{eq:element_interpolation}),
\[
\varphi \approx (\mathbf{N}^{Q9})^{T}\bar{\boldsymbol{\varphi}}^{e}.
\]
Here $\mathbf{N}^{Q9}$ is the Q9 shape-function matrix. Accordingly, the network inputs and outputs are represented on a $101\times 101$ grid, while the loss is computed on a $51\times 51$ mesh according to \eqref{eq:total_loss}. Q9-level sampling is applied only on intersected elements, and bilinear Q4 sampling is used on all uncut elements. The penalty parameters are $\lambda_1=1\times10^{-2}, \lambda_2=1\times10^{-3}, \lambda_3=1\times10^{-4}$. We train for 2{,}000 epochs using the SOAP optimizer with 1,000 training samples and 100 test samples. The ground-truth fields are generated by $\varphi$-FEM on the $256\times 256$ mesh using second-order finite elements (biquadratic Q9 elements on quadrilateral meshes). As in the Cook's membrane setup, reference solves use sparse LU factorization for each Newton linearization instead of GMRES, since strongly nonlinear tangent systems in this experiment prevented Krylov iterations from converging reliably \cite{golub2013matrix}.

\begin{figure}[t]
	\centering
	\begin{minipage}[t]{0.37\textwidth}
		\vspace{0pt}% anchor [t] to true top
		\centering
		\makebox[\linewidth][l]{\textbf{(a)}}\par\vspace{0.3em}
		\includegraphics[width=0.75\linewidth]{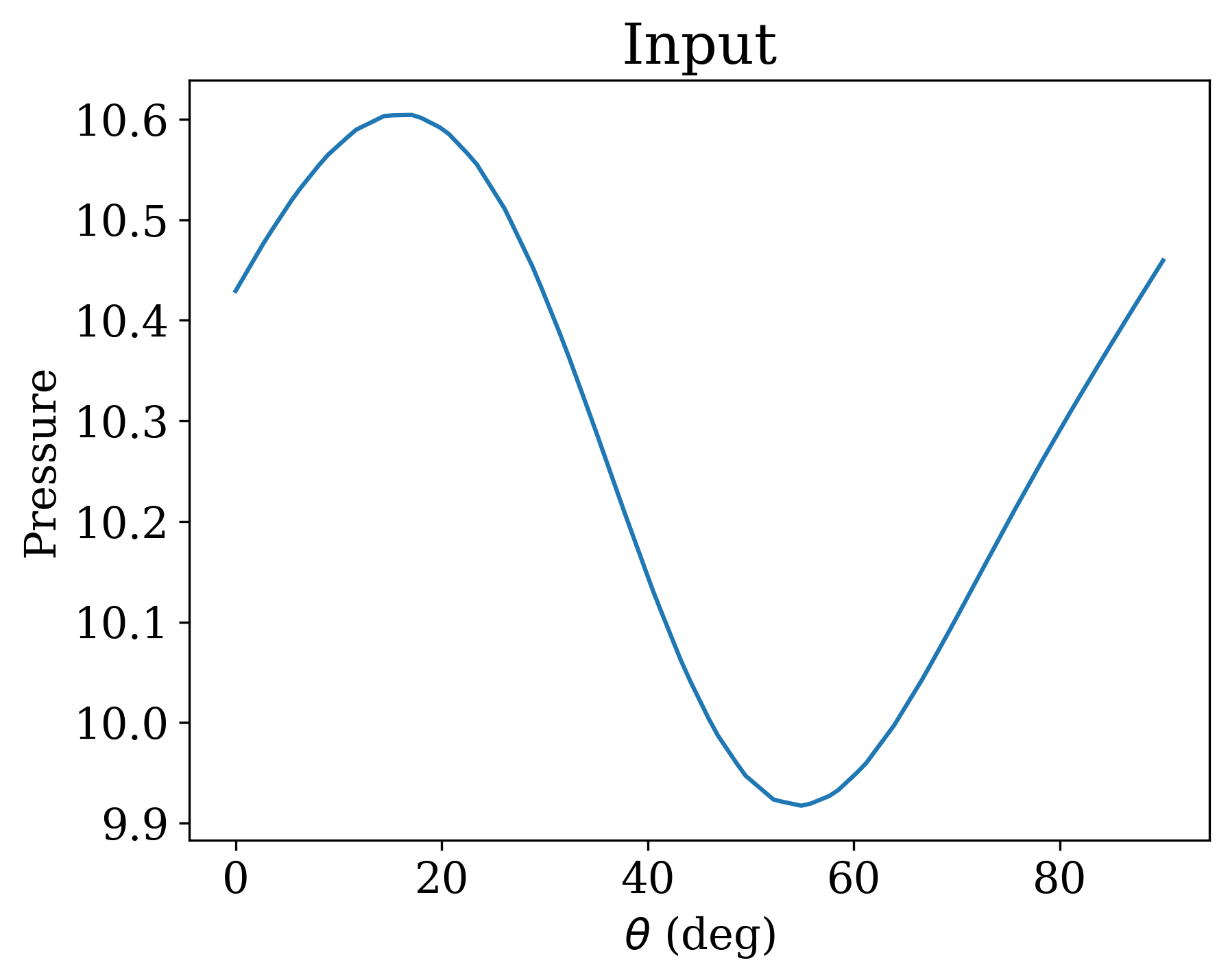}\\[0.6ex]
		\vspace{0.1cm}
		\includegraphics[width=0.9\linewidth]{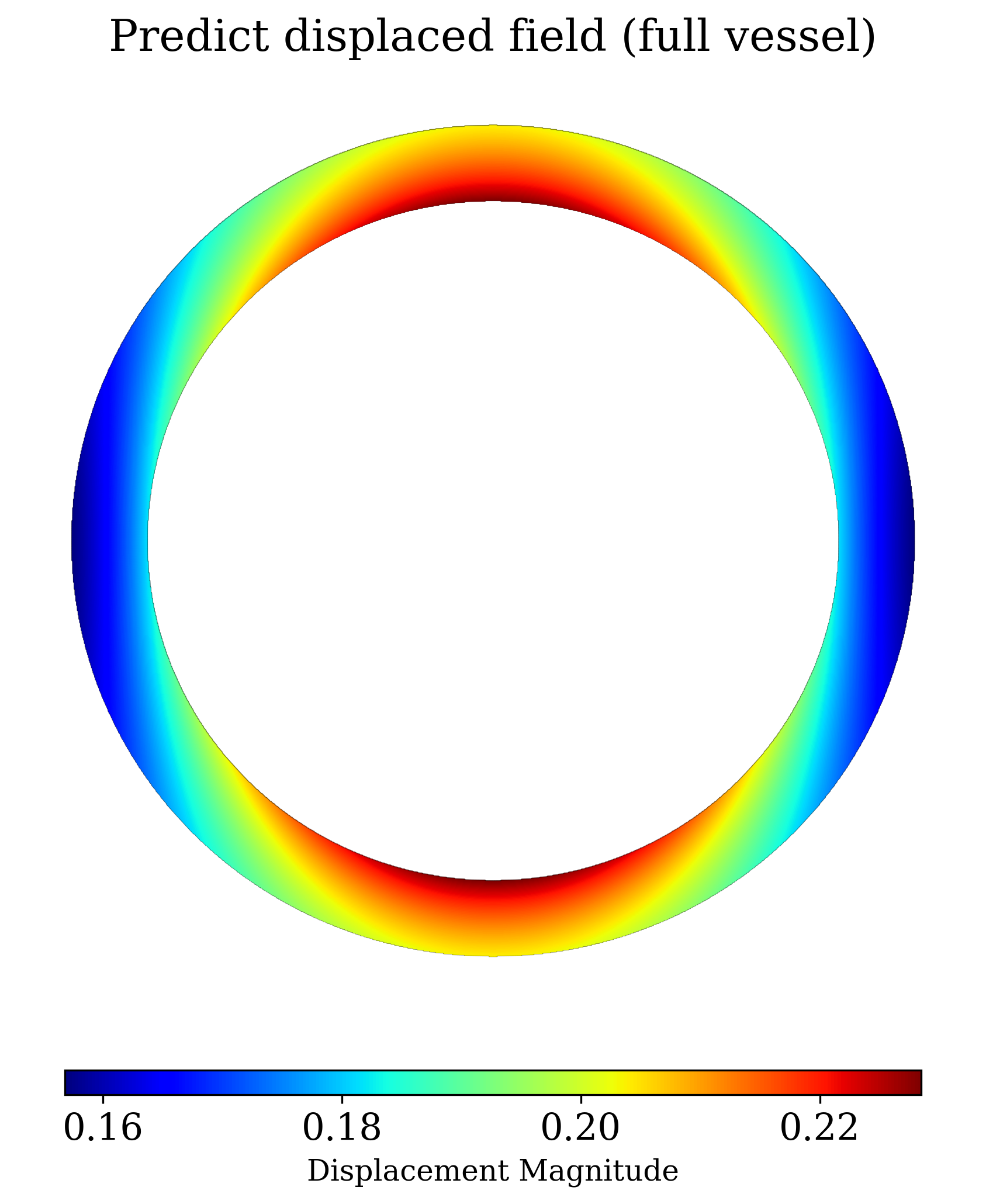}
	\end{minipage}%
	\hspace{0.03\textwidth}%
	\begin{minipage}[t]{0.57\textwidth}
		\vspace{0pt}% align top with left column
		\centering
		\makebox[\linewidth][l]{\textbf{(b)}}\par\vspace{0.3em}
		\begin{tabular}{@{}c@{\hspace{0.02\linewidth}}c@{}}
			\includegraphics[width=0.45\linewidth]{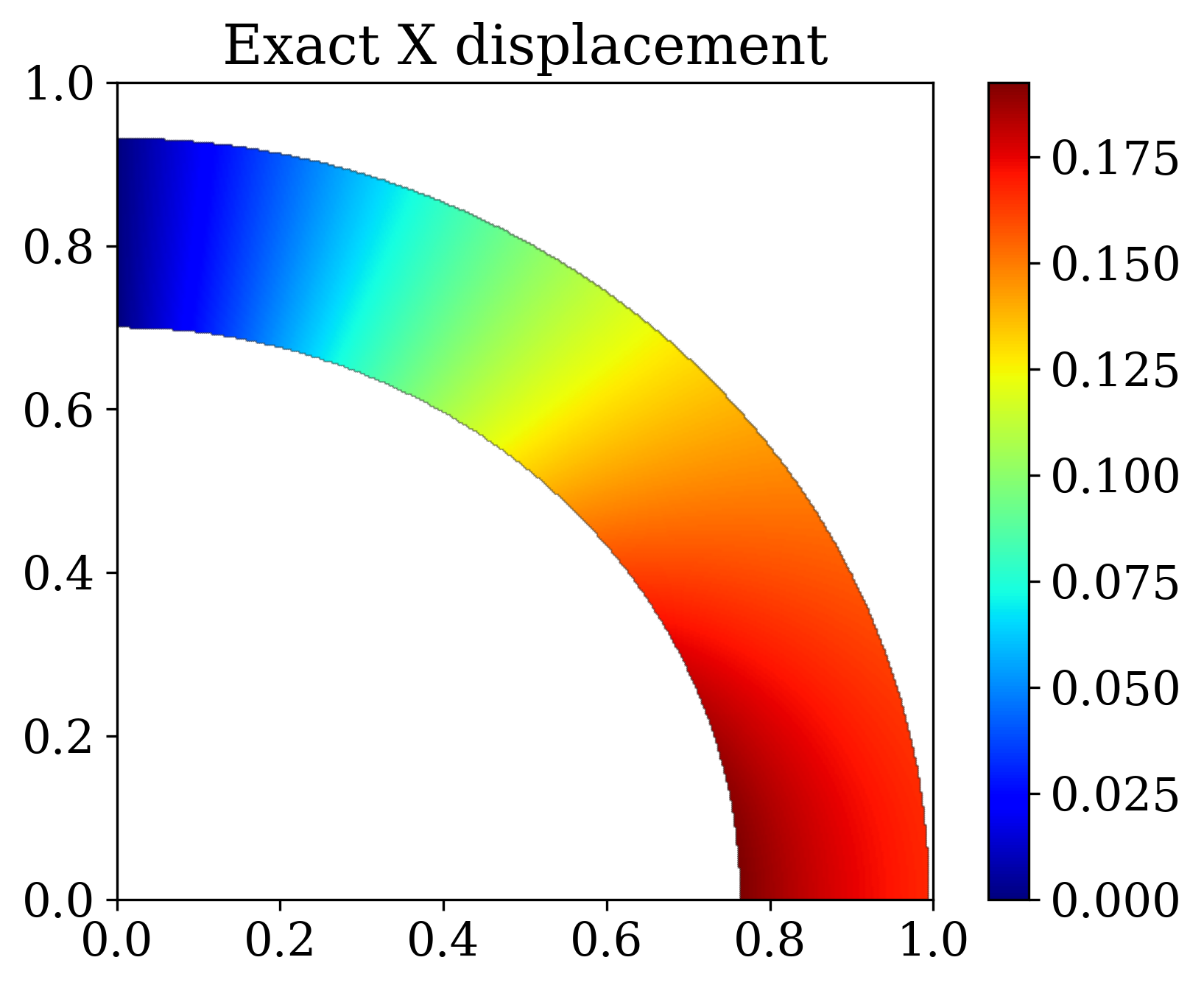} &
			\includegraphics[width=0.45\linewidth]{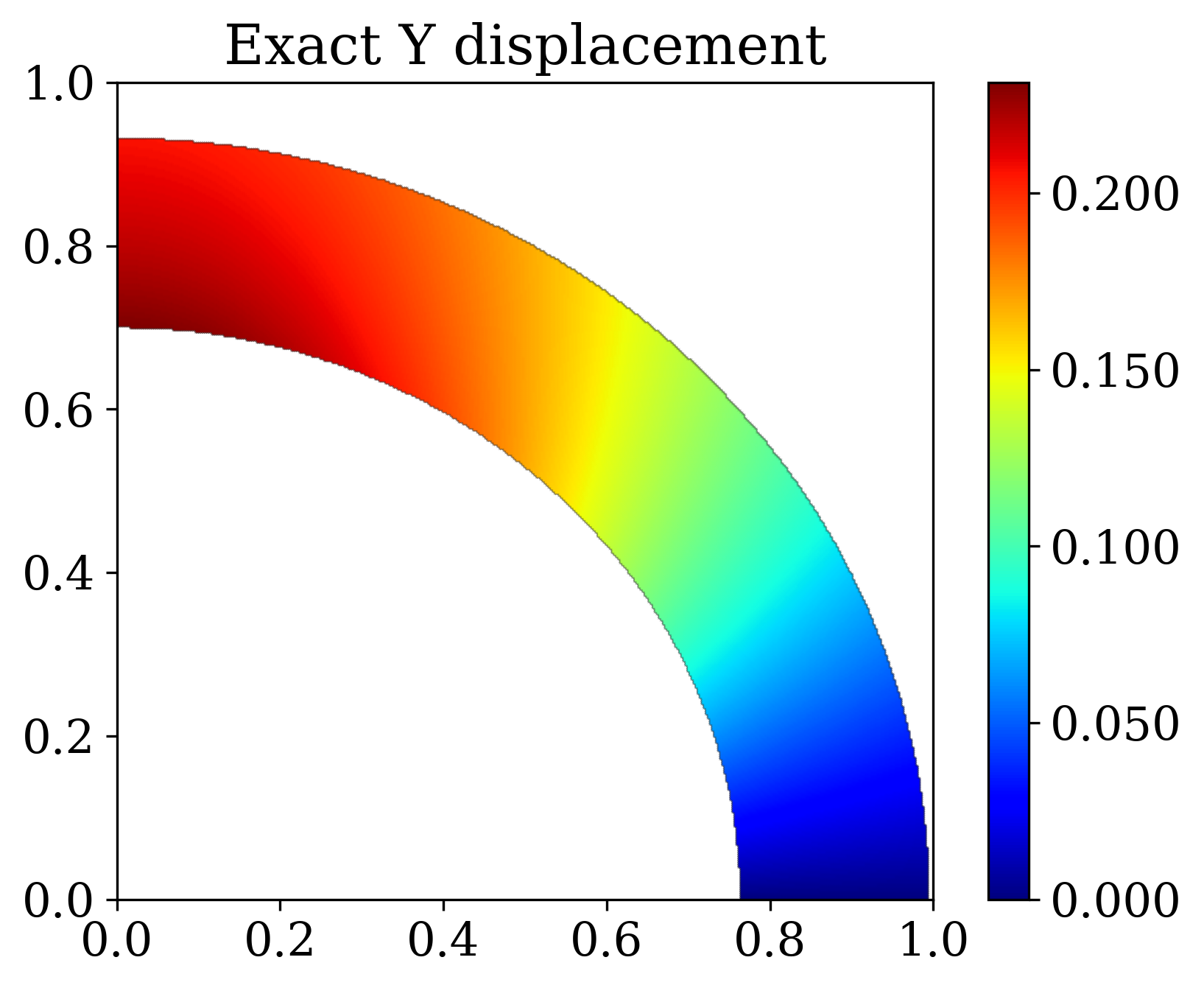} \\[0.5ex]
			\includegraphics[width=0.45\linewidth]{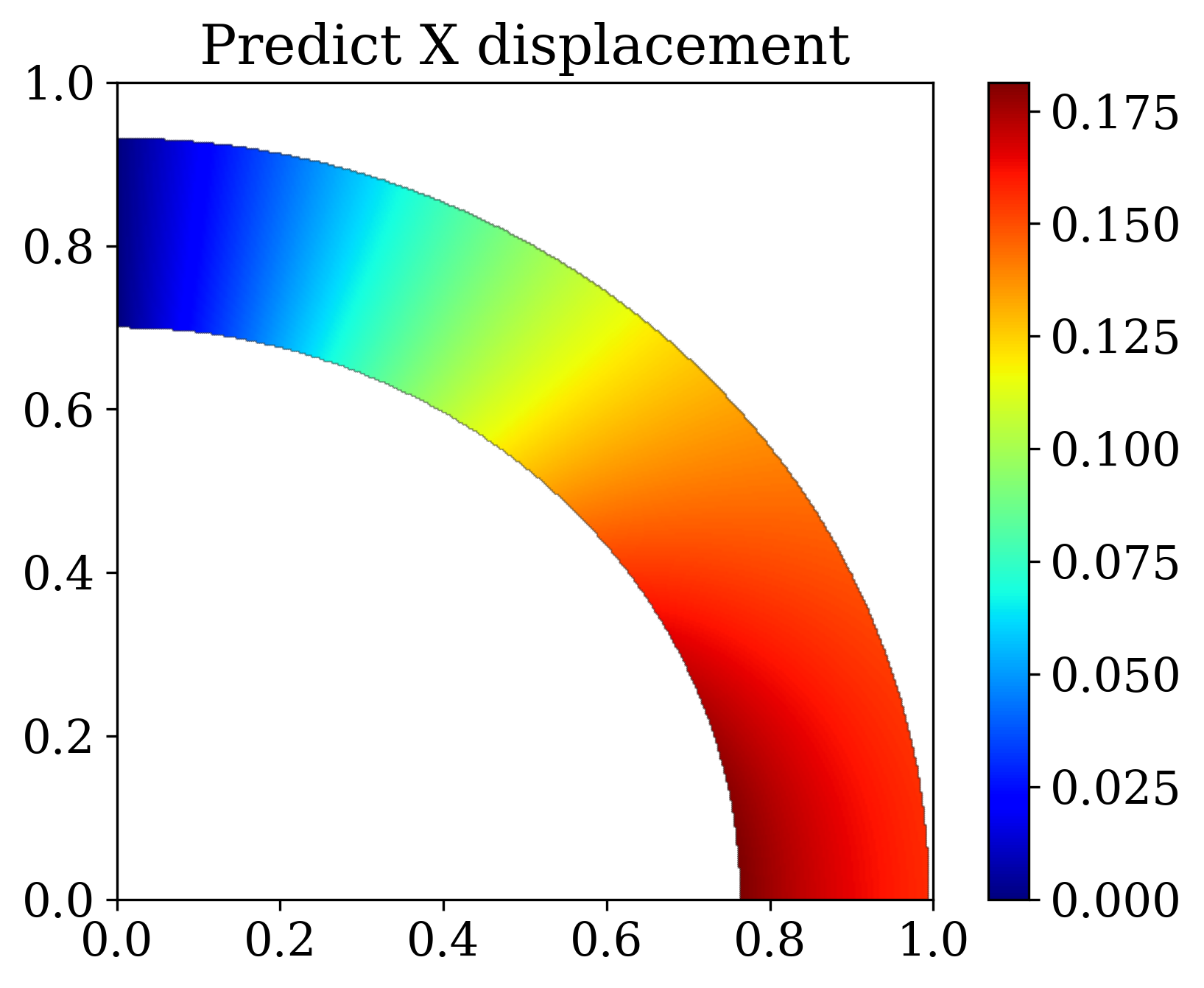} &
			\includegraphics[width=0.45\linewidth]{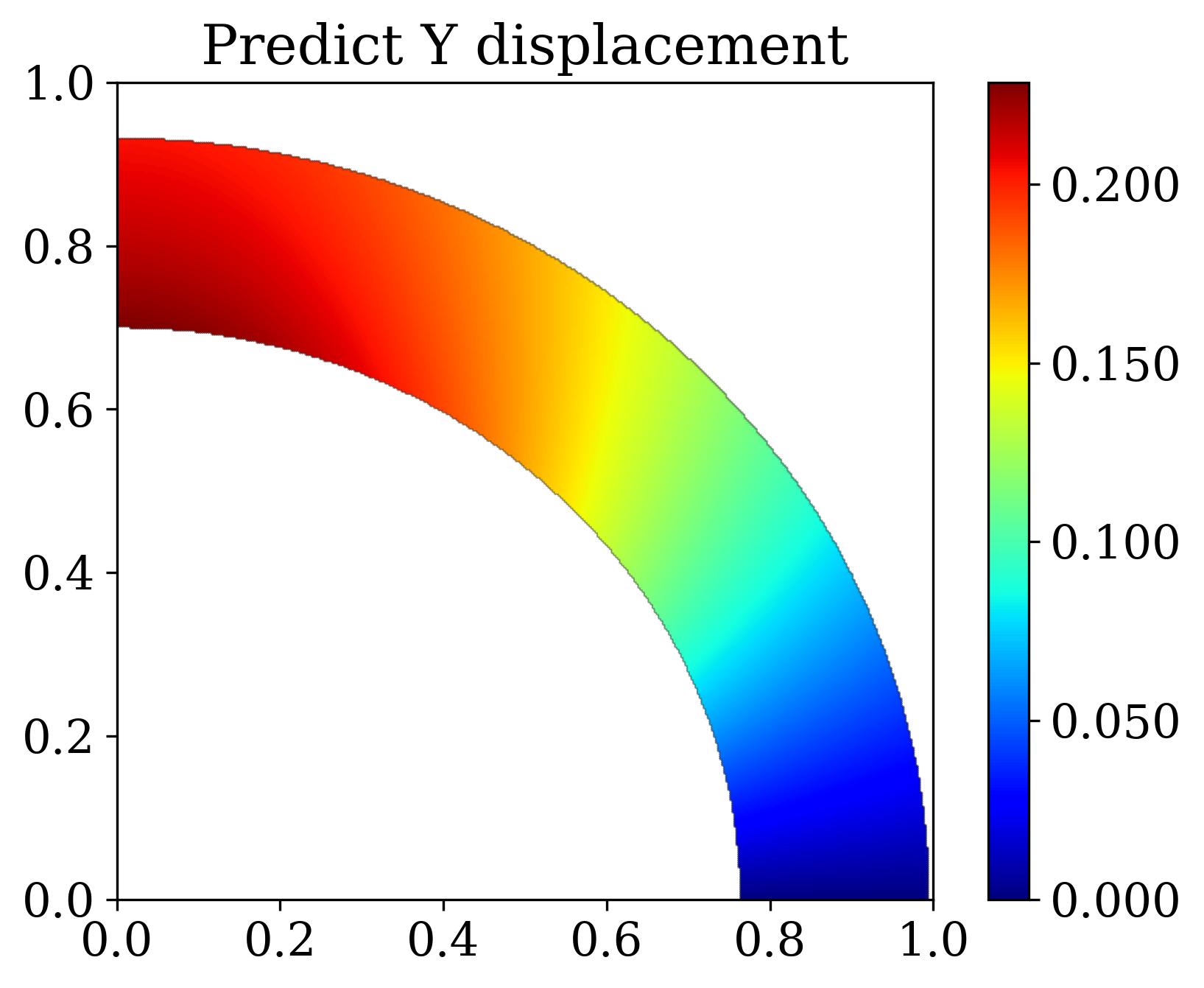} \\[0.5ex]
			\includegraphics[width=0.45\linewidth]{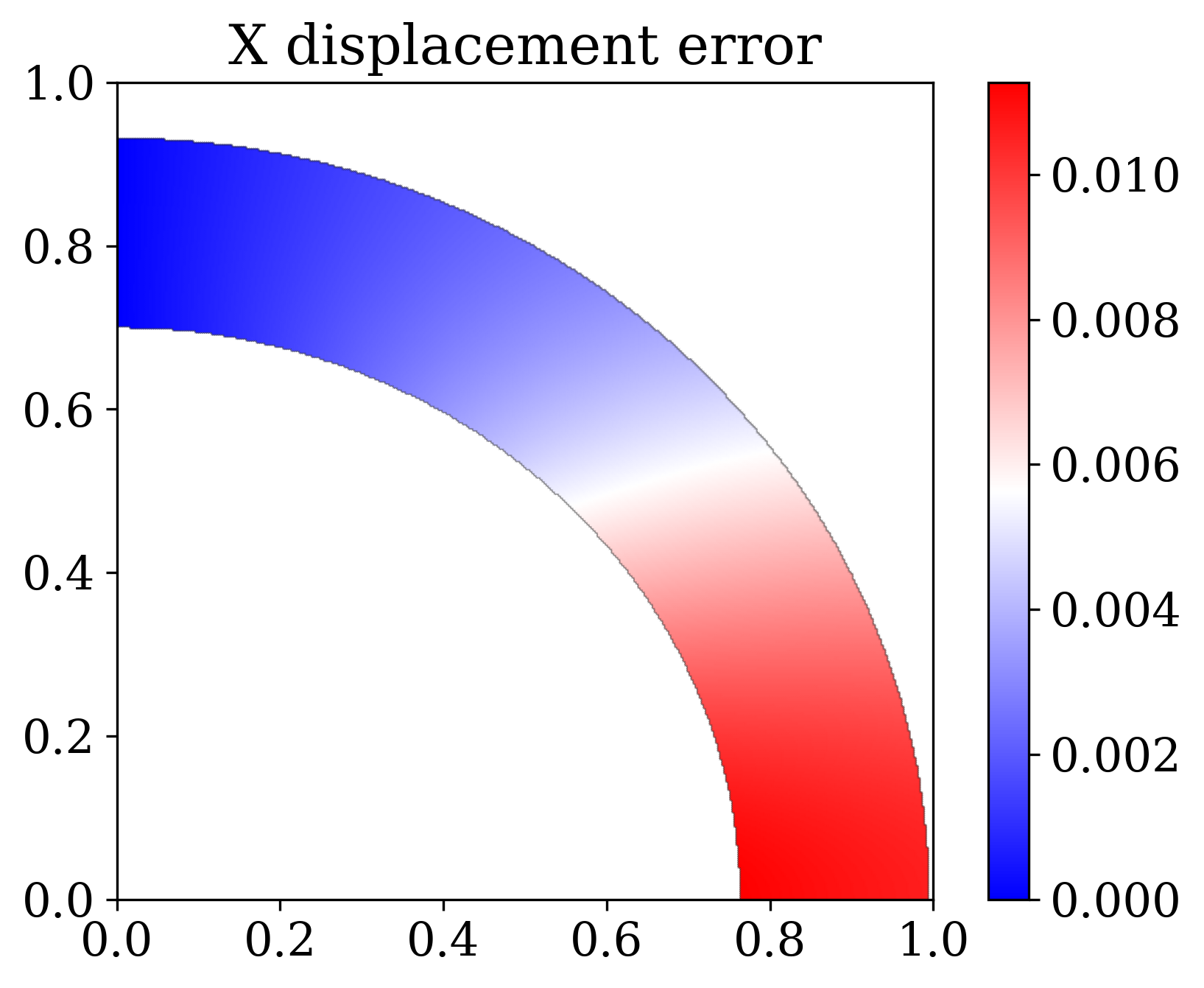} &
			\includegraphics[width=0.45\linewidth]{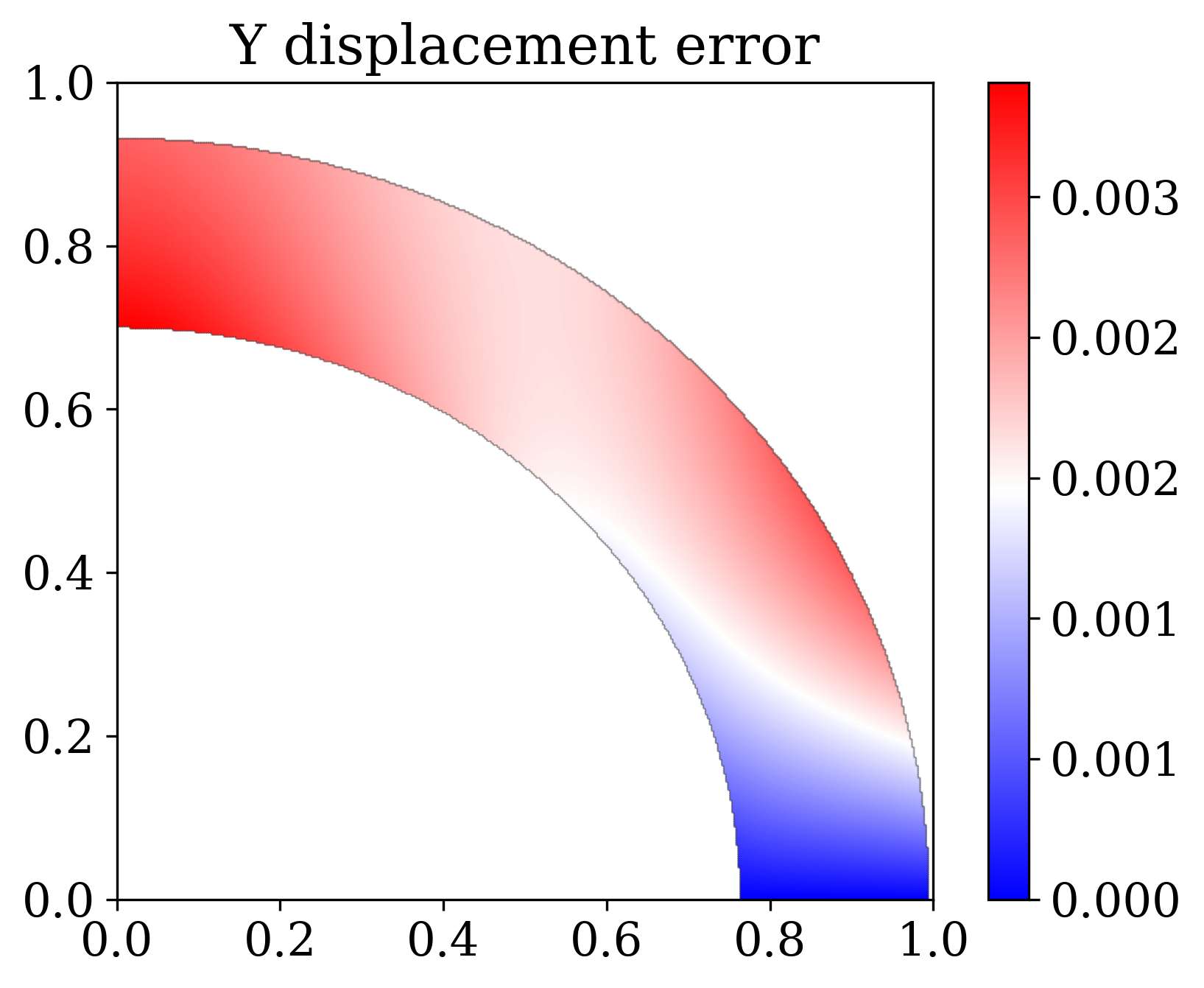}
		\end{tabular}
	\end{minipage}
	\caption{Representative test sample for the pressure vessel case. (a) WINO-predicted deformed configuration; (b) $x$- and $y$-displacement components for the exact solutions (top row), WINO predictions (middle row), and pointwise absolute error contours (bottom row).}
	\label{fig:case_vessel_displacement}
\end{figure}

To assess the convergence of the $\varphi$-FEM reference discretization, we plot the relative errors in the $L^2$ norm, the $H^1$ seminorm, and the energy norm as functions of the mesh size; see Fig.~\ref{fig:case_vessel}d. The reference solution is obtained by a conventional body-fitted finite element method on a finely resolved unstructured mesh of the quarter annular domain. The boundary geometry of the standard FEM method is discretized with 501 segments along each quarter-elliptic arc. The $\varphi$-FEM displacement is interpolated onto the reference finite element space, and relative $L^2$, $H^1$, and energy errors are evaluated on that body-fitted mesh. The results indicate second-order convergence, $\mathcal{O}(h^2)$, for $\varphi$-FEM on this benchmark. On the $51\times51$ mesh used in the WINO experiments, the corresponding relative errors are $4.57\times10^{-5}$, $4.55\times10^{-4}$, and $4.33\times10^{-4}$, respectively, confirming that the $\varphi$-FEM reference is sufficiently accurate for this case.

\begin{remark}
	The Neumann boundary condition of the pressure vessel case on the inner wall $\Gamma_{\mathrm{in}}$ is embedded in the cut cell region of the unfitted background mesh; its weak enforcement follows the $\varphi$-FEM Neumann branch and therefore requires the auxiliary fields $\mathbf{y}_h$ and $\mathbf{p}_h$ together with the displacement $\mathbf{u}_h$. To embed the Neumann boundary condition accurately on the inner wall, the level set is factorized as $\varphi=E_{\mathrm{in}}E_{\mathrm{out}}$, the analytic partial derivatives of $E_{\mathrm{in}}$ and $E_{\mathrm{out}}$ are concatenated as extra input channels (cf.~\eqref{eq:vessel_mapping}), and the discrete $\varphi$ on cut cells is interpolated with biquadratic Q9 shape functions.
\end{remark}

\begin{table}[t]
	\centering
	\footnotesize
	\caption{Performance of models on the pressure vessel case with training on 51$\times$51 meshes and testing on 256$\times$256 meshes.}
	\label{tab:vessel_resolution}
	\begin{tabularx}{\linewidth}{@{}l >{\centering\arraybackslash}X >{\raggedright\arraybackslash}X >{\raggedright\arraybackslash}X >{\raggedright\arraybackslash}X@{}}
		\toprule
		Method & Data cost (per sample) & $\|e\|_{L^2}$ & $\|e\|_{H^1}$ & $\|e\|_{E}$ \\
		\midrule
		WINO & \multirow{3}{*}{$42.56\,\mathrm{s}$} & $4.11 \pm 0.73\%$ & $4.55 \pm 1.07\%$ & $4.39 \pm 2.06\%$ \\
		$\varphi$-FEM-FNO &  & $0.44 \pm 0.03\%$ & $7.81 \pm 0.79\%$ & $0.74\pm 0.14\%$ \\
		WINO+data &  & $0.71 \pm 0.05\%$ & $0.98 \pm 0.03\%$ & $0.54 \pm 0.07\%$ \\
		\bottomrule
	\end{tabularx}
\end{table}

Fig.~\ref{fig:case_vessel}b shows the evolution of the mean relative $L^2$ errors on the training and test sets. The training error drops rapidly in the early epochs and then gradually stabilizes. The test error follows a similar trend, starting from a higher initial level and converging to a value close to that of the training error. Fig.~\ref{fig:case_vessel}c presents box plots of the post-training relative errors, where the horizontal line indicates the median and the box denotes the interquartile range. A representative full-vessel prediction from the test set is shown in Fig.~\ref{fig:case_vessel_displacement}, where the WINO solution closely matches the high-fidelity reference and accurately captures the displacement induced by random pressure loading. To further evaluate efficiency, we compare WINO with purely data-driven $\varphi$-FEM-FNO and WINO augmented with labeled data in terms of both accuracy and computational cost (see Table~\ref{tab:performance_vessel}). Although $\varphi$-FEM-FNO achieves very high accuracy in $L^2$ norm, its data-generation cost is substantially higher than that of WINO (about 13.4 seconds per sample on average), so WINO is more efficient overall. When labeled data are available, WINO augmented with labeled data also outperforms $\varphi$-FEM-FNO.

To assess the discretization invariance of WINO, we evaluate the operators previously trained on the coarse $51\times51$ mesh on the same 100 test samples. Specifically, these operators are applied, without retraining, to the corresponding inputs discretized on a significantly finer $256\times256$ mesh, and the resulting predictions are compared with $\varphi$-FEM reference solutions computed on the same $256\times256$ discretization (Table~\ref{tab:vessel_resolution}). WINO predictions and error levels remain stable under this cross-mesh evaluation, whereas $\varphi$-FEM-FNO does not exhibit comparable resolution invariance: its relative errors, especially in the $H^1$ seminorm, increase markedly on the fine reference mesh. This discrepancy is partly explained by larger $\varphi$-FEM-FNO errors in displacement gradients at various resolutions near the boundary, where mismatches are penalized strongly by the $H^1$ seminorm but contribute less to the discrete energy norm in \eqref{eq:rel_err_energy}, so the reported $H^1$ error can exceed the energy-norm error even when the relative $L^2$ error remains small. Moreover, WINO inference requires only $3.79\,\mathrm{ms}$ on the coarse $51\times51$ mesh and $40.36\,\mathrm{ms}$ on the fine $256\times256$ mesh, compared with about $42.56\,\mathrm{s}$ per $\varphi$-FEM reference solve.

\section{Conclusion}\label{section_conclusion}

This work introduces WINO, a weak-form physics-informed neural operator built on $\varphi$-FEM for large-deformation hyperelasticity on varying domains. By combining operator learning with element-wise weak-form residual evaluation, WINO avoids body-fitted remeshing and directly works on fixed Cartesian backgrounds while retaining physically meaningful supervision through finite-element test functions. Across the numerical benchmarks, we can find that WINO delivers stable and accurate predictions in $L^2$, $H^1$, and energy-type errors.

For boundary treatment, Dirichlet constraints are enforced either by lifting constructions in the network outputs, whereas Neumann conditions are naturally incorporated through the weak-form boundary terms and auxiliary $\varphi$-FEM variables on cut-cell regions. This hybrid design allows efficient handling of mixed and geometry-dependent boundary partitions without explicit boundary-fitted meshing. In addition, the framework preserves compatibility with standard finite-element assembly, making implementation practical for implicit boundaries and irregular cut configurations.

Compared with supervised $\varphi$-FEM-FNO, label-free WINO achieves substantial accuracy without requiring converged $\varphi$-FEM reference solutions for training. Across the reported benchmarks, the total computational cost (input sampling plus operator training) amounts to about $15\%$--$70\%$ of that of supervised $\varphi$-FEM-FNO. When labeled data are available, the WINO+data variant can further improve accuracy. These results position WINO as a practical label-free alternative for nonlinear mechanics tasks where solver-grade labels are costly, rather than as a universal replacement for supervised operator learning when high-fidelity data can be obtained at moderate cost.

Beyond the training objective, WINO can also accelerate the nonlinear $\varphi$-FEM solve through neural-operator warm starts (NOWS). Feeding WINO's displacement prediction into Newton's method as an initial iterate reduces reliance on incremental load stepping and mitigates slow or fragile convergence from cold starts. In our experiments, NOWS consistently lowers the total number of outer Newton steps and inner GMRES iterations relative to classical continuation, while the additional neural forward pass incurs only a small computational cost. Therefore, WINO serves as both a surrogate for offline solution generation and a practical initializer that couples learned operators with iterative high-fidelity solvers in strongly nonlinear hyperelastic settings.

The current study also has limitations. First, the underlying $\varphi$-FEM formulation is best suited to smooth boundaries; for non-smooth boundaries (e.g., corners or kinks), its approximation quality and robustness may degrade, which limits direct applicability. Second, the convergence speed of WINO training is sensitive to the penalty parameters in the loss function, so careful parameter tuning is often required in practice; to address this issue, we provide a penalty-parameter normalized scaling strategy in Appendix~\ref{app:penalty}. Third, the present implementation is restricted to 2D settings, structured Cartesian backgrounds, and relatively smooth constitutive responses, and the training cost remains nontrivial for high-resolution and strongly nonlinear regimes.

Future work will extend WINO to 3D and unstructured element types,  investigate adaptive/local refinement and multi-fidelity training, and develop reusable software and model libraries for broader engineering deployment. We also plan to incorporate more complex constitutive laws and extend the framework to additional problem classes, including phase-field-based variational damage model \cite{duan2025unified, duan2026unified}.

\section*{Declaration of competing interest}
The authors declare that they have no known competing financial interests or personal relationships that could have appeared to influence the work reported in this paper.

\appendix
\renewcommand{\thesection}{\Alph{section}}
\makeatletter
\renewcommand{\@seccntformat}[1]{\appendixname~\csname the#1\endcsname.~~}
\makeatother

\section{Influence of the ghost-penalty stabilization term}\label{app:WINO_Gh}

In this study, we focus on WINO augmented with ghost-penalty stabilization residuals for the plate-with-a-hole benchmark in Subsection~\ref{sec:plate_with_hole}. To incorporate these residuals into the loss function, we modify \eqref{eq:weak_residual_u} as
\begin{equation}
	\begin{aligned}
		\mathcal{R}_u^{(m)}(\mathbf{v}_{i,j}) &= \int_{\Omega_h^{(m)}} \mathbf{P}(\mathbf{F}(\mathbf{u}_h^{(m)})):\nabla \mathbf{v}_{i,j}\, d\mathbf{x}
		+\int_{\partial\Omega_h^{(m)}} (\mathbf{y}_h^{(m)} \mathbf{n})\cdot \mathbf{v}_{i,j}\, ds - \int_{\Omega_h^{(m)}} \mathbf{f}_{h}^{(m)}\cdot\mathbf{v}_{i,j}\, d\mathbf{x} \\
		&\quad + G_h^{(m)}(\mathbf{u}_h^{(m)}, \mathbf{v}_{i,j}),
	\end{aligned}
	\label{eq:weak_residual_u_Gh}
\end{equation}
where, following \eqref{Gh_definition},
\begin{equation}
	G_h^{(m)}(\mathbf{u}_h,\mathbf{v}_{i,j}) := \sigma_N h \int_{\Gamma_h^{G,(m)}} \bigl[\!\bigl[\mathbf{P}(\mathbf{F}(\mathbf{u}_h)) \mathbf{n}\bigr]\!\bigr] \cdot \bigl[\!\bigl[D_u(\mathbf{P} \circ \mathbf{F})(\mathbf{u}_h)[\mathbf{v}_{i,j}] \mathbf{n}\bigr]\!\bigr],
	\label{eq:Gh_sample}
\end{equation}
and $\Gamma_h^{G,(m)}$ denotes the ghost-penalty facet set of sample $m$, as defined in Section~\ref{subsec:phi_fem}. For the compressible Neo-Hookean model \eqref{eq:strain_energy}, the first Piola--Kirchhoff stress can be written as
\begin{equation}
	\mathbf{P}(\mathbf{F})=\mu\mathbf{F}+(\lambda\ln J-\mu)\mathbf{F}^{-\mathsf{T}}, \qquad J=\det\mathbf{F},
	\label{eq:PK_neo_hookean}
\end{equation}
and the directional derivative in the test direction $\mathbf{v}_{i,j}=\mathbf{e}_j v_i$ is evaluated with $\delta\mathbf{F}=\nabla\mathbf{v}_{i,j}$ by
\begin{equation}
		D_u(\mathbf{P}\circ\mathbf{F})(\mathbf{u}_h)[\mathbf{v}_{i,j}]
		=\mu\,\delta\mathbf{F}
		+\lambda\,(\mathbf{F}^{-\mathsf{T}}:\delta\mathbf{F})\,\mathbf{F}^{-\mathsf{T}} 
		+(\lambda\ln J-\mu)\,\bigl(-\mathbf{F}^{-\mathsf{T}}\,\delta\mathbf{F}^{\mathsf{T}}\,\mathbf{F}^{-\mathsf{T}}\bigr),
	\label{eq:Du_P_explicit}
\end{equation}
where $:$ denotes the Frobenius inner product and $\mathbf{F}$, $J$, and $\ln J$ are evaluated at $\mathbf{u}_h$. This is the constitutive linearization used in the ghost-penalty residual implementation, and the integral defining $G_h$ is discretized with three-point Gaussian quadrature. In this section, we set the parameter $\sigma_N=0.01$.

\begin{figure}[t]
	\centering
	\begin{minipage}[t]{0.42\textwidth}
		\vspace{0pt}
		\centering
		\makebox[\linewidth][l]{\textbf{(a)}}\par\vspace{0.05em}
		\includegraphics[width=\linewidth]{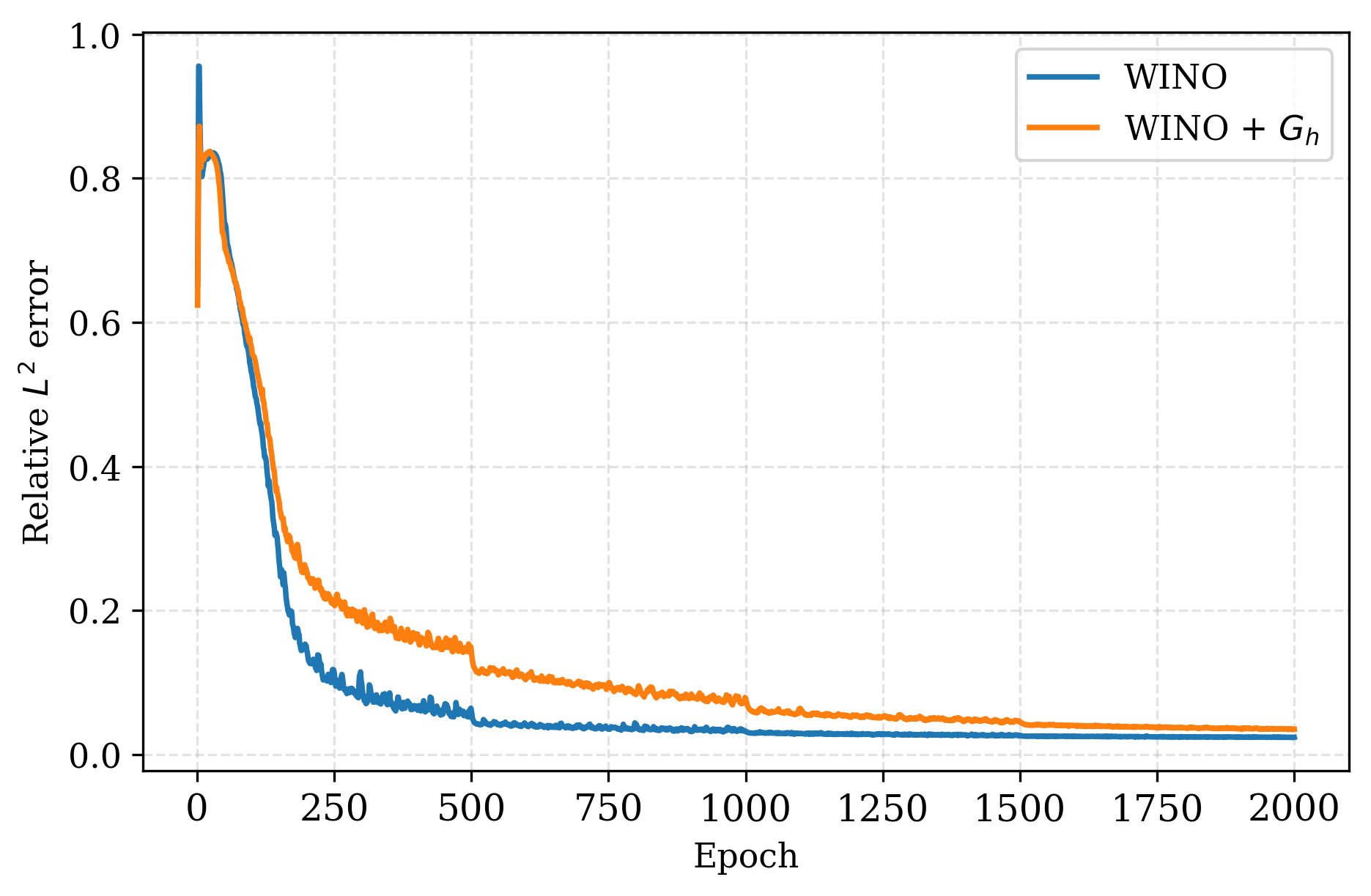}
	\end{minipage}
	\hfill
	\begin{minipage}[t]{0.38\textwidth}
		\vspace{0pt}
		\centering
		\makebox[\linewidth][l]{\textbf{(b)}}\par\vspace{0.3em}
		\includegraphics[width=\linewidth]{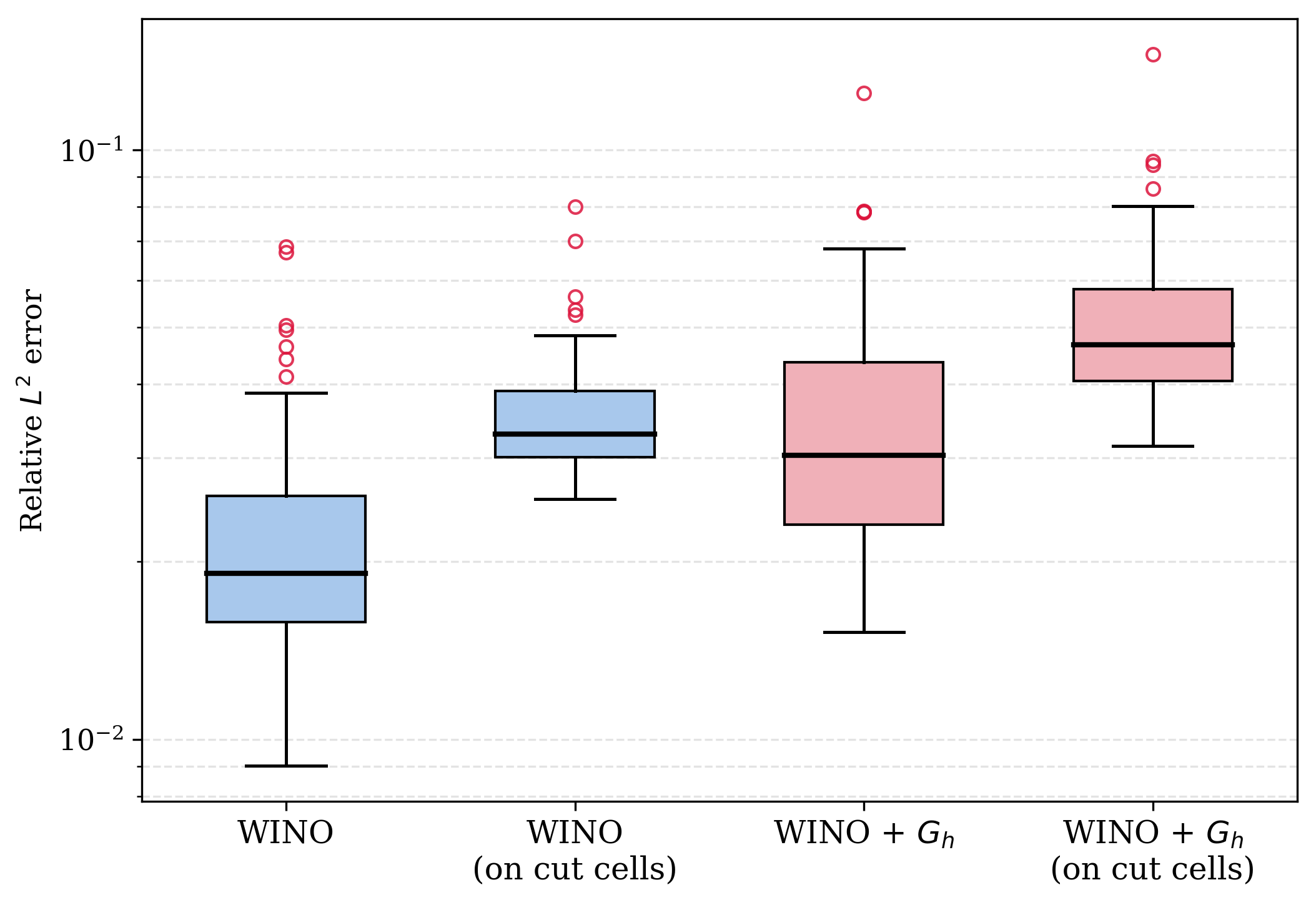}
	\end{minipage}
	\hfill
	\begin{minipage}[t]{0.163\textwidth}
		\vspace{0pt}
		\centering
		\makebox[\linewidth][l]{\textbf{(c)}}\par\vspace{0.25em}
		\includegraphics[width=\linewidth]{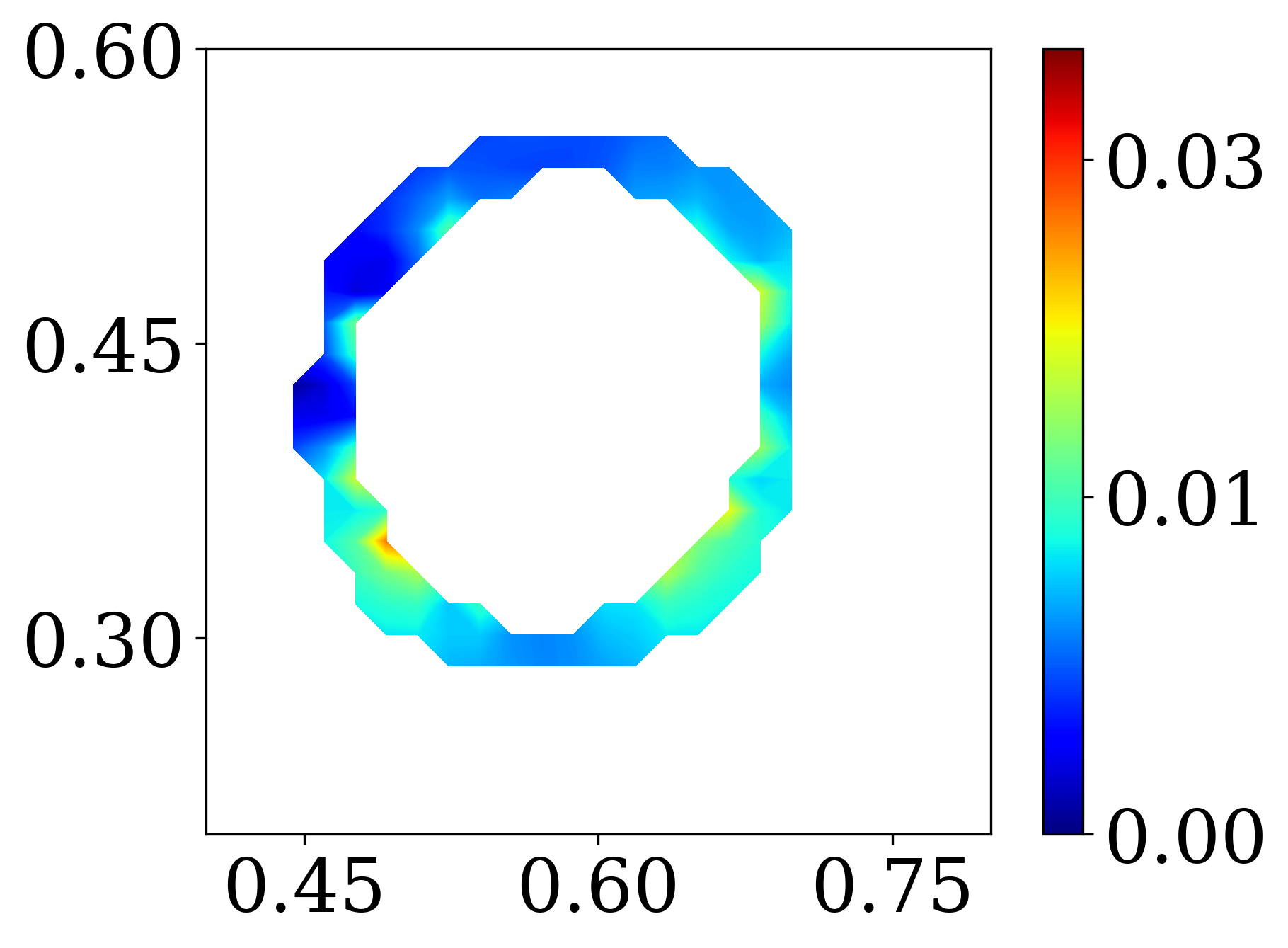}\\[0.1em]
		\includegraphics[width=\linewidth]{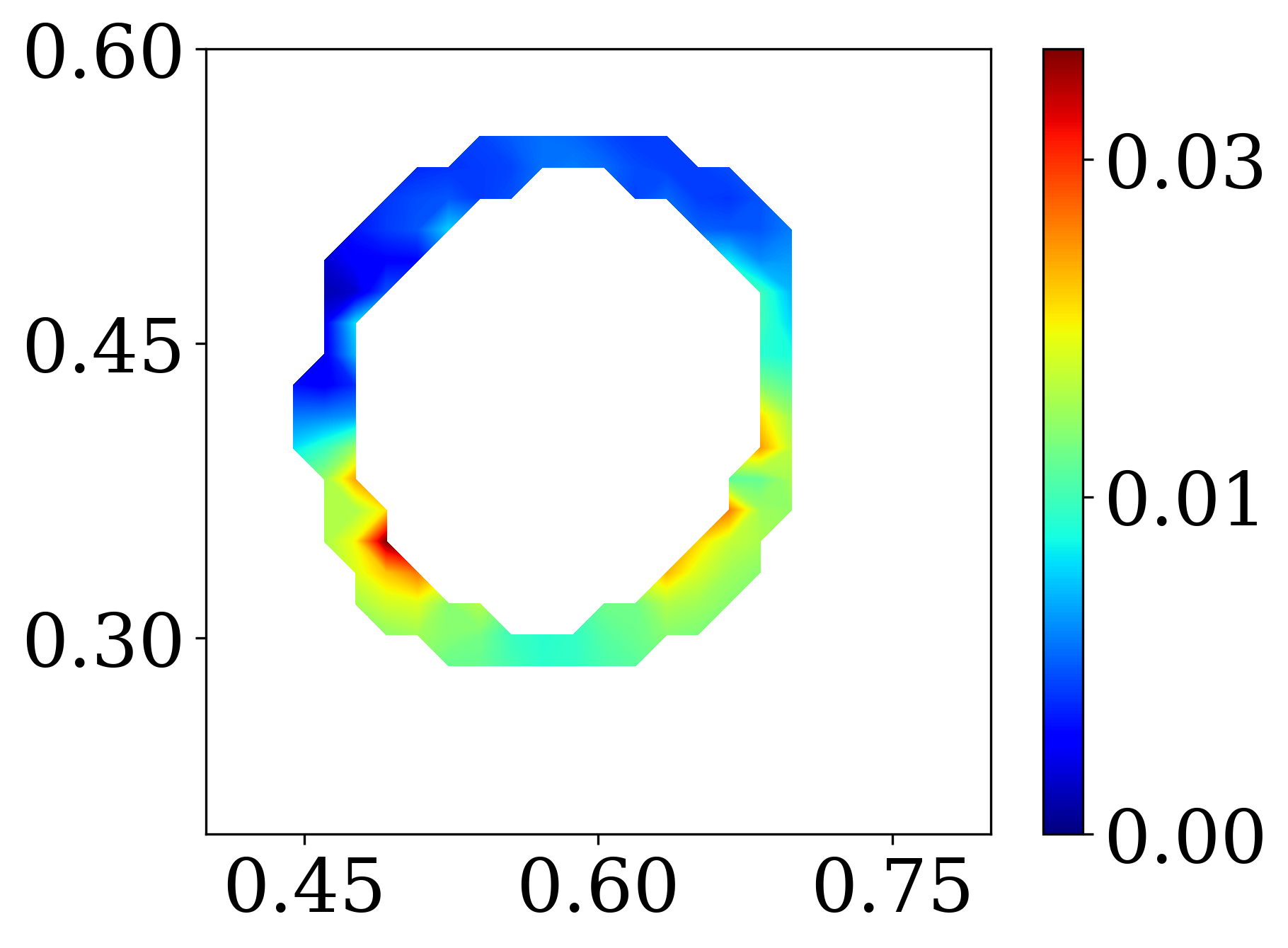}
	\end{minipage}
	\caption{Influence of ghost-penalty residuals for the plate with a hole case: (a)~Evolution of the mean relative $L^2$ error on the test set during training. (b)~Box plots of the relative $L^2$ errors after training (on the whole domain and on the cut-cell band $\mathcal{T}_h^\Gamma$). (c)~Pointwise absolute displacement errors restricted to $\mathcal{T}_h^\Gamma$ for a representative test sample (top: WINO; bottom: WINO+$G_h$).}
	\label{fig:appendix_Gh}
\end{figure}

\begin{table}[t]
	\centering
	\footnotesize
	\caption{Training time and relative errors of WINOs for the plate with a hole case. Cost ratio is total time relative to WINO.}
	\label{tab:appendix_Gh_performance}
	\begin{tabularx}{\linewidth}{@{}l >{\raggedright\arraybackslash}X >{\raggedright\arraybackslash}X >{\raggedright\arraybackslash}X >{\raggedright\arraybackslash}X >{\raggedright\arraybackslash}X@{}}
		\toprule
		Method & Training time & Cost ratio & $\|e\|_{L^2}$ & $\|e\|_{H^1}$ & $\|e\|_{E}$ \\
		\midrule
		WINO & $2912.8$s & $1.00\times$ & $2.26 \pm 1.10\%$ & $9.12 \pm 1.26\%$ & $8.39 \pm 0.99\%$ \\
		WINO+$G_h$ & $9139.1$s & $3.14\times$ & $3.49 \pm 1.71\%$ & $10.9 \pm 1.78\%$ & $9.69 \pm 1.25\%$ \\
		\bottomrule
	\end{tabularx}
\end{table}

The hyperparameters of WINO augmented with ghost-penalty residuals are the same as those of label-free WINO in Subsection~\ref{sec:plate_with_hole}. We train WINO+$G_h$ for 2{,}000 epochs; the error evolution is shown in Fig.~\ref{fig:appendix_Gh}a, and the training time together with the post-training relative errors are reported in Table~\ref{tab:appendix_Gh_performance}. Adding the $G_h$ residual neither accelerates convergence nor improves the final accuracy, most likely because the additional term makes the loss landscape more difficult to optimize. By contrast, the total training time of WINO+$G_h$ increases substantially, since evaluating \eqref{eq:Du_P_explicit} requires computing $\mathbf{F}^{-\mathsf{T}}$, which is computationally expensive. We also report the relative $L^2$ errors on the cut-cell band $\mathcal{T}_h^\Gamma$ in Fig.~\ref{fig:appendix_Gh}b and the pointwise absolute errors in Fig.~\ref{fig:appendix_Gh}c. The cut-cell errors are consistent with those on the full domain, indicating that omitting $G_h$ from the loss does not disproportionately degrade accuracy near cut cells; moreover, WINO+$G_h$ is generally less accurate than label-free WINO. A similar increase in computational cost would also arise for pure Dirichlet problems if $G_h$ were included. Therefore, the ghost-penalty residuals are not used in the numerical experiments reported in the main text.

\section{Newton iteration and Krylov methods}
\label{app:nows}

Consider finding $\mathbf{U}\in\mathbb{R}^n$ such that the nonlinear balance
\begin{equation}
	\mathbf{K}(\mathbf{U})=\mathbf{F}
\end{equation}
holds, where $\mathbf{F}\in\mathbb{R}^n$ is a prescribed load vector and $\mathbf{K}:\mathbb{R}^n\to\mathbb{R}^n$ is a nonlinear operator. Defining the residual map $\mathbf{R}(\mathbf{U}):=\mathbf{F}-\mathbf{K}(\mathbf{U})$, the balance condition is equivalent to the root-finding problem
\begin{equation}
	\mathbf{R}(\mathbf{U})=\mathbf{0},
	\label{eq:newton_residual}
\end{equation}
which is the standard form solved by Newton iteration.

Starting from an initial iterate $\mathbf{U}^{(0)}$, Newton's method generates the sequence
\begin{equation}
	\mathbf{U}^{(k+1)}=\mathbf{U}^{(k)}+\Delta\mathbf{U}^{(k)}, \qquad k=0,1,2,\ldots,
	\label{eq:newton_update}
\end{equation}
where the increment $\Delta\mathbf{U}^{(k)}$ solves the linearized system
\begin{equation}
	\mathbf{J}(\mathbf{U}^{(k)})\,\Delta\mathbf{U}^{(k)}=-\,\mathbf{R}(\mathbf{U}^{(k)}).
	\label{eq:newton_linear}
\end{equation}
Here $\mathbf{J}(\mathbf{U})=\partial \mathbf{R}/\partial \mathbf{U}$ is the Jacobian matrix. Each Newton step thus reduces to a sparse linear system, and the quality of the initial guess $\mathbf{U}^{(0)}$ strongly influences the number of outer Newton iterations required.

When the tangent matrix $\mathbf{J}(\mathbf{U}^{(k)})$ in \eqref{eq:newton_linear} is symmetric positive definite (SPD), the correction equation can be written in the standard form
\begin{equation}
	\mathbf{A}\mathbf{x}=\mathbf{b},
	\label{eq:krylov_linear}
\end{equation}
with $\mathbf{A}=\mathbf{J}(\mathbf{U}^{(k)})$, $\mathbf{b}=-\mathbf{R}(\mathbf{U}^{(k)})$, and $\mathbf{x}=\Delta\mathbf{U}^{(k)}$. Krylov methods approximate $\mathbf{x}$ using only matrix--vector products with $\mathbf{A}$, and no explicit inverse is required. The conjugate gradient (CG) is a common Krylov method for \eqref{eq:krylov_linear} when $\mathbf{A}$ is SPD. The algorithm can be interpreted as a projection method that minimizes the following quadratic functional
\begin{equation}
	\mathcal{Q}(\mathbf{x})=\tfrac{1}{2}\mathbf{x}^T\mathbf{A}\mathbf{x}-\mathbf{b}^T\mathbf{x},
	\label{eq:cg_quadratic}
\end{equation}
over the affine Krylov space $\mathbf{x}_0+\mathcal{K}_k(\mathbf{A},\mathbf{r}_0)$, where $\mathbf{x}_0$ is an initial guess for \eqref{eq:krylov_linear}, and its residual is $\mathbf{r}_0=\mathbf{b}-\mathbf{A}\mathbf{x}_0$,
\begin{equation*}
	\mathcal{K}_k(\mathbf{A},\mathbf{r}_0)=\mathrm{span}\bigl\{\mathbf{r}_0,\mathbf{A}\mathbf{r}_0,\ldots,\mathbf{A}^{k-1}\mathbf{r}_0\bigr\}.
\end{equation*}
At iteration $k$, the iterate $\mathbf{x}_k$ lies in $\mathbf{x}_0+\mathcal{K}_k(\mathbf{A},\mathbf{r}_0)$ and minimizes $\mathcal{Q}$ on that space. The method generates mutually $\mathbf{A}$-conjugate search directions $\{\mathbf{p}_k\}$ and updates the approximation and the residual as
\begin{equation}
	\begin{aligned}
		\mathbf{x}_{k+1} &= \mathbf{x}_k+\alpha_k \mathbf{p}_k,\\
		\mathbf{r}_{k+1} &= \mathbf{r}_k-\alpha_k \mathbf{A}\mathbf{p}_k.
	\end{aligned}
	\label{eq:cg_xr_updates}
\end{equation}
The scalar $\alpha_k$ is chosen so that $\mathcal{Q}$ is minimized along $\mathbf{p}_k$. The full recursion also defines $\mathbf{p}_{k+1}$ from $\mathbf{r}_{k+1}$ and the previous direction so that conjugacy is maintained. In exact arithmetic, CG terminates in at most $n$ steps. The error $\mathbf{e}_k=\mathbf{x}_k-\mathbf{x}^\ast$, with $\mathbf{x}^\ast=\mathbf{A}^{-1}\mathbf{b}$, satisfies the classical bound in the energy norm $\|\mathbf{v}\|_{\mathbf{A}}=\sqrt{\mathbf{v}^T\mathbf{A}\mathbf{v}}$:
\begin{equation}
	\|\mathbf{e}_k\|_{\mathbf{A}}\le 2\left(\frac{\sqrt{\kappa(\mathbf{A})}-1}{\sqrt{\kappa(\mathbf{A})}+1}\right)^{k}\|\mathbf{e}_0\|_{\mathbf{A}},
	\label{eq:cg_convergence}
\end{equation}
where $\kappa(\mathbf{A})$ is the condition number of $\mathbf{A}$ \cite{saad1981krylov}. Thus convergence accelerates when $\mathbf{A}$ is well conditioned or effectively preconditioned.

When $\mathbf{A}$ is nonsymmetric or indefinite, CG is generally not applicable. The generalized minimal residual (GMRES) method is a standard Krylov solver for nonsingular systems \eqref{eq:krylov_linear}. GMRES builds the same Krylov space $\mathcal{K}_k(\mathbf{A},\mathbf{r}_0)$ as above but chooses $\mathbf{x}_k\in\mathbf{x}_0+\mathcal{K}_k(\mathbf{A},\mathbf{r}_0)$ to minimize the Euclidean residual norm,
\begin{equation}
	\|\mathbf{b}-\mathbf{A}\mathbf{x}_k\|_2 = \min_{\mathbf{x}\in\mathbf{x}_0+\mathcal{K}_k(\mathbf{A},\mathbf{r}_0)} \|\mathbf{b}-\mathbf{A}\mathbf{x}\|_2.
	\label{eq:gmres_minres}
\end{equation}
The Arnoldi process produces $\mathbf{V}_k\in\mathbb{R}^{n\times k}$ with orthonormal columns spanning $\mathcal{K}_k(\mathbf{A},\mathbf{r}_0)$ and an extended upper Hessenberg matrix $\widetilde{\mathbf{H}}_k\in\mathbb{R}^{(k+1)\times k}$ such that
\begin{equation}
	\mathbf{A}\mathbf{V}_k = \mathbf{V}_{k+1}\widetilde{\mathbf{H}}_k,
	\label{eq:gmres_arnoldi}
\end{equation}
where $\mathbf{V}_{k+1}\in\mathbb{R}^{n\times(k+1)}$ has orthonormal columns. Writing $\mathbf{x}_k=\mathbf{x}_0+\mathbf{V}_k\mathbf{y}_k$ with $\mathbf{y}_k\in\mathbb{R}^k$, the minimization \eqref{eq:gmres_minres} reduces to the small least-squares problem
\begin{equation}
	\mathbf{y}_k = \arg\min_{\mathbf{y}\in\mathbb{R}^k}\bigl\|\,\beta\mathbf{e}_1-\widetilde{\mathbf{H}}_k\mathbf{y}\,\bigr\|_2,\qquad \beta=\|\mathbf{r}_0\|_2,
	\label{eq:gmres_lsq}
\end{equation}
where $\mathbf{e}_1$ is the first canonical basis vector in $\mathbb{R}^{k+1}$. In practice, restarted GMRES($m$) \cite{morgan2002gmres} limits memory and orthogonalization work by truncating the Krylov basis after $m$ inner iterations and restarting from the current iterate. 

As for CG, performance depends strongly on the initial guess $\mathbf{x}_0$ and on \emph{preconditioning}. Rather than iterating on $\mathbf{A}$ directly, one introduces a nonsingular matrix $\mathbf{M}$ whose inverse or triangular solve approximates $\mathbf{A}^{-1}$ at moderate cost, and applies GMRES to the left-preconditioned system $\mathbf{M}^{-1}\mathbf{A}\,\mathbf{x}=\mathbf{M}^{-1}\mathbf{b}$ (or to an equivalent right-preconditioned formulation). A good preconditioner reduces the effective condition number and clusters eigenvalues of the iteration operator, which typically lowers the Krylov iteration count for a given tolerance. For sparse Jacobian systems, algebraic multigrid is a common choice. GAMG (geometric--algebraic multigrid), the smoothed-aggregation / graph-based multigrid preconditioner available in PETSc \cite{balay2019petsc}, builds a hierarchy of coarse operators directly from the sparse matrix: neighboring degrees of freedom are grouped through graph coarsening, and restriction and prolongation maps combine smoothed interpolation with relaxation sweeps so that a multigrid cycle approximates the action of $\mathbf{A}^{-1}$ at cost that scales nearly linearly with the number of unknowns for problems whose principal part is elliptic. Here ``geometric--algebraic'' refers to combining algebraic coarsening from matrix connectivity with intergrid operators informed by geometric-multigrid practice. In practice, each Newton linearization may be solved with restarted GMRES preconditioned by GAMG applied to the Jacobian matrix.

\section{Normalized scaling strategy for penalty parameters}
\label{app:penalty}

For the mixed Dirichlet--Neumann benchmarks, the WINO training objective \eqref{eq:total_loss} combines the weak-form momentum residual with squared auxiliary residuals on cut cells. For notational convenience in this appendix, we introduce a global weight $\lambda_0$ and collect the four scalar weights in the vector
\begin{equation}
	\boldsymbol{\lambda}=\lambda_0\bigl[1,\;\lambda_1,\;\lambda_2,\;\lambda_3\bigr]=\bigl[\lambda_0,\;\lambda_1\lambda_0,\;\lambda_2\lambda_0,\;\lambda_3\lambda_0\bigr],
	\label{eq:penalty_vector}
\end{equation}
which is equivalent to multiplying the original weight pattern in \eqref{eq:total_loss}--\eqref{eq:strong_loss} by $\lambda_0$. Here $\lambda_0$ weights the weak loss $\mathcal{L}_{\mathrm{weak}}$ in \eqref{eq:weak_loss}, and $\lambda_1$, $\lambda_2$, and $\lambda_3$ are the dimensionless relative multipliers of the auxiliary loss contributions $\mathcal{L}_y$, $\mathcal{L}_p$, and $\mathcal{L}_d$ in \eqref{eq:strong_loss}. Once the scaled vector $\boldsymbol{\lambda}$ is fixed, the penalty parameters used in \eqref{eq:strong_loss} are recovered by
\begin{equation}
	\lambda_i=\frac{(\boldsymbol{\lambda})_{i+1}}{\lambda_0}, \qquad i=1,2,3.
	\label{eq:penalty_recovery}
\end{equation}

Our normalized scaling strategy is to set
\begin{equation}
	\boldsymbol{\lambda}=\boldsymbol{\lambda}_r \oslash \mathbf{L}_{\mathrm{scale}},
	\label{eq:penalty_scaling}
\end{equation}
where $\oslash$ denotes element-wise division and $\mathbf{L}_{\mathrm{scale}}\in\mathbb{R}^4$ is determined from the unweighted loss contributions evaluated before optimization; representative values are listed in Table~\ref{tab:penalty_scale}. The goal is to keep the individual loss terms on comparable magnitudes during training, typically within $0.1$--$10$.
For the plate-with-a-hole and pressure-vessel benchmarks, we use the reference vector
\begin{equation}
	\boldsymbol{\lambda}_r=\bigl[10^{2},\;1,\;1,\;10^{-7}\bigr].
	\label{eq:lambda_r_general}
\end{equation}
The first component weights the weak-form residual assembled over all test functions in $\Omega_h$ and is therefore assigned a larger reference value, whereas the fourth component associated with $\mathcal{L}_d$ plays only an auxiliary role during training and is assigned a smaller value. For the more strongly nonlinear Cook's membrane benchmark, we instead adopt
\begin{equation}
	\boldsymbol{\lambda}_r=\bigl[10^{2},\;10^{-4},\;10^{-3},\;10^{-7}\bigr],
	\label{eq:lambda_r_cook}
\end{equation}
obtained by adjusting the second and third components relative to the reference above; this choice was found to accelerate convergence of the mean relative $L^2$ error on the test set in that benchmark. The penalty parameters $\lambda_1$, $\lambda_2$, and $\lambda_3$ quoted in Section~\ref{section_Results} are obtained from the scaled vector $\boldsymbol{\lambda}$ in \eqref{eq:penalty_scaling} via \eqref{eq:penalty_recovery}.

\begin{table}[t]
	\centering
	\footnotesize
	\caption{Pre-training unweighted loss contributions, scaling vectors $\mathbf{L}_{\mathrm{scale}}$, reference vectors $\boldsymbol{\lambda}_r$, and scaled penalty vectors $\boldsymbol{\lambda}$ for the mixed boundary benchmarks. Each entry lists the four components associated with $\mathcal{L}_{\mathrm{weak}}$, $\mathcal{L}_y$, $\mathcal{L}_p$, and $\mathcal{L}_d$, respectively.}
	\label{tab:penalty_scale}
	\begin{tabularx}{\linewidth}{@{}l >{\raggedright\arraybackslash}X >{\raggedright\arraybackslash}X >{\raggedright\arraybackslash}X@{}}
		\toprule
		 & Plate with a hole & Cook's membrane & Pressure vessel \\
		\midrule
		Initial loss
		& $[1\mathrm{e}\text{-}2,\,3\mathrm{e}\text{-}2,\,2\mathrm{e}3,\,2\mathrm{e}\text{-}6]$
		& $[1\mathrm{e}\text{-}2,\,9\mathrm{e}\text{-}1,\,3\mathrm{e}\text{-}4,\,5\mathrm{e}\text{-}10]$
		& $[1\mathrm{e}\text{-}1,\,4\mathrm{e}\text{-}1,\,9\mathrm{e}0,\,2\mathrm{e}\text{-}6]$ \\
		$\mathbf{L}_{\mathrm{scale}}$
		& $[1\mathrm{e}\text{-}2,\,1\mathrm{e}\text{-}1,\,1\mathrm{e}3,\,1\mathrm{e}\text{-}6]$
		& $[1\mathrm{e}\text{-}1,\,1\mathrm{e}0,\,1\mathrm{e}\text{-}4,\,1\mathrm{e}\text{-}10]$
		& $[1\mathrm{e}\text{-}1,\,1\mathrm{e}\text{-}1,\,1\mathrm{e}0,\,1\mathrm{e}\text{-}6]$ \\
		$\boldsymbol{\lambda}_r$
		& $[1\mathrm{e}2,\,1\mathrm{e}0,\,1\mathrm{e}0,\,1\mathrm{e}\text{-}7]$
		& $[1\mathrm{e}2,\,1\mathrm{e}\text{-}4,\,1\mathrm{e}\text{-}3,\,1\mathrm{e}\text{-}7]$
		& $[1\mathrm{e}2,\,1\mathrm{e}0,\,1\mathrm{e}0,\,1\mathrm{e}\text{-}7]$ \\
		$\boldsymbol{\lambda}$
		& $[1\mathrm{e}4,\,1\mathrm{e}1,\,1\mathrm{e}\text{-}3,\,1\mathrm{e}\text{-}1]$
		& $[1\mathrm{e}3,\,1\mathrm{e}\text{-}4,\,1\mathrm{e}1,\,1\mathrm{e}3]$
		& $[1\mathrm{e}3,\,1\mathrm{e}1,\,1\mathrm{e}0,\,1\mathrm{e}\text{-}1]$ \\
		$[1,\lambda_1,\lambda_2,\lambda_3]$
		& [1, 1e-3, 1e-7, 1e-5] & [1, 1e-7, 1e-2, 1e0] & [1, 1e-2, 1e-3, 1e-4] \\
		\bottomrule
	\end{tabularx}
\end{table}

\section*{Data availability}
The paper's datasets were created using the accompanying code. The dataset is available for free access at https://github.com/bokai-zhu/WINO.

\section*{Code availability}
The code for the numerical experiments is available for free access at https://github.com/bokai-zhu/WINO.

\section*{Declaration of generative AI and AI-assisted technologies in the manuscript preparation process}

During the preparation of this work, the authors used ChatGPT in order to improve language clarity and assist with grammar refinement. After using this tool, the authors reviewed and edited the content as needed and take full responsibility for the content of the published article.

\bibliographystyle{IEEEtran}
\bibliography{WINO_ref.bib}

\end{document}